\documentclass[11pt]{amsbook} 
\usepackage{color} 
\usepackage{amssymb, amsmath, amsthm, amsfonts} 
\usepackage{mathtools}
\usepackage{caption}
\usepackage[ansinew]{inputenc}
\usepackage{graphicx}

\usepackage[margin=1.25in]{geometry}

\def \sgn {\hbox {sgn}}
\newcommand{\R}{\mathbb{R}}
\newcommand{\C}{\mathbb{C}}

\newcommand{\N}{\mathbb{N}}
\newcommand{\Z}{\mathbb{Z}}

\renewcommand{\P}{{P}} 
\newcommand{\I} {{I}}\newcommand{\HH}{\mathbb{H}}
\newcommand{\F}{\mathcal{F}} \newcommand{\E}{{E}}
\newcommand{\eps}{\epsilon}

\newcommand{\U}{\mathbb D}
\newcommand{\1}[1]{\mathbf{1}_{\left\{ #1\right\}}}

\newcommand{\oldcomments}{}

\numberwithin{figure}{chapter}
\numberwithin{section}{chapter}

\newtheorem{theorem}{Theorem}[chapter]
\newtheorem{lemma}[theorem]{Lemma}
\newtheorem{proposition}[theorem]{Proposition}
\newtheorem{corollary}[theorem]{Corollary}

\newtheorem{remark}[theorem]{Remark}
\newtheorem{definition}[theorem]{Definition}
\newtheorem{conjecture}[theorem]{Conjecture}
\newtheorem{reminder}[theorem]{Reminder}

\newtheorem{exercise}[theorem]{Exercise}

\addtocounter{chapter}{-1}

\frontmatter

\title {Lecture notes on the Gaussian Free Field}
\author{Wendelin Werner \and Ellen Powell} 
\address{ETH Z\"urich \and Durham University}
\email{wendelin.werner@math.ethz.ch, ellen.g.powell@durham.ac.uk}

\date {}
\begin {document}

\begin {abstract}
The Gaussian Free Field (GFF) in the continuum appears to be the natural generalisation of Brownian motion, when one replaces time by a multidimensional continuous 
parameter. While Brownian motion can be viewed as the most natural random real-valued function defined on $\R_+$ with $B(0)=0$, the GFF in a domain $D$ of 
$\R^d$ for $d \ge 2$ is a natural random real-valued  generalised function defined on $D$ with zero boundary conditions on $\partial D$. In particular, it is not a random continuous function. 

The goal of these lecture notes is to describe some aspects of the continuum GFF and of its discrete counterpart defined 
on lattices, with the aim of providing a gentle self-contained introduction to some recent developments on this topic, such as the relation between the continuum GFF, Brownian loop-soups and the Conformal Loop Ensembles CLE$_4$.

This is an updated and expanded version of the notes written by the first author (WW) for graduate courses at ETH Z\"urich in 2014 and 2018. It has benefited from the comments and corrections of students, as well as of a referee; we thank them all very much. The exercises that are interspersed in the first half of these notes mostly originate from the exercise sheets prepared by the second author (EP) for this course in 2018. 
\end {abstract}

\maketitle

\vskip 10cm 

  \subsubsection*{Acknowledgements} {\small The authors acknoweldge the support of the grant 175505 of the Swiss National Science Foundation. EP also thanks the FIM of ETH Z\"urich, for support during visits in 2019-2020.}
 
\tableofcontents

\mainmatter

\chapter*{Overview} 

Let us start with a very very sketchy and necessarily incomplete historical overview in order to try to explain the scope of these lecture notes. 

One simple way to think of the Gaussian Free Field (GFF) is that it is the most natural and tractable model for a random function defined on either a discrete graph (each vertex of the graph 
is assigned a random real-valued height, and the distribution favours configurations where neighbouring vertices have similar heights) or on a subdomain of $\R^d$. We will refer to these 
two cases as the discrete GFF and the continuum GFF respectively. 

The Gaussian free field, in both its discrete and continuum versions, has been one of the main building blocks in mathematical physics at least since the early 1970s. Many of its important features were pointed out and used in a number of seminal works by Symanzik, Nelson, Brydges, Fr\"ohlich, Spencer, Simon and many others. Often these works were connected with questions originating from Quantum Field Theory -- see for instance \cite {Nelson, Simon, GRS, BFS} or \cite {FFS} and the references therein. 
In the theoretical physics community, a number of later developments (such as Conformal Field Theory -- CFT, Liouville Quantum Gravity -- LQG) in the 1980s, used the continuum GFF as an essential building block, together with a number of other new fundamental ideas, in order to describe aspects of random systems in two dimensions. 

While the discrete GFF is indeed a random function defined on the vertices of a graph, the continuum GFF is a somewhat more complicated object when $d \ge 2$. Indeed, it is not a random continuous function -- it is only a random generalised function. The height of the GFF at a given point is not well-defined, but the ``mean height'' of a realisation of the GFF on some given bounded open set \emph{is} a well-defined Gaussian random variable. The fact that the continuum GFF is not 
a proper function is not a problem in CFT or LQG, as in these theories the focus is put on correlation functions (leading to results on critical exponents for example) and these turn out to be well-defined. On the other hand, it makes it seem almost impossible to detect random geometric structures (i.e. random fractal objects) in a sample of the GFF. 

Just before the turn of the century, Oded Schramm \cite {Schramm} constructed  Schramm-Loewner Evolutions (SLE): a 
family of random curves in the plane providing a direct mathematical approach to the random geometric objects (random interfaces, random domains) that appear in these two-dimensional statistical physics questions. This was quite a novel perspective.  
In fact, in order to connect SLE with random fields, it is natural to consider the ``entire'' collection of interfaces that are present in the system (not only the particular interface described by one SLE). This gives rise to the Conformal Loop Ensembles (CLE) introduced and studied in \cite {Sheffield,SheffieldWerner}, that can be defined using appropriate generalisations of SLE.

Another important SLE-related development initiated in \cite {SchrammSheffield} -- see also \cite {Dubedat,IMG1} -- is that one particular SLE (the SLE$_4$) and one particular CLE (the CLE$_4$) can be directly related to the continuum GFF, and interpreted as ``level-lines'' of this random generalised function. 
This led many authors to revisit some of the basic features of the GFF, such as its Markov property, leading to a novel and alternative understanding of the continuum GFF in two dimensions. 

A central role in some developments around SLE and CLE is played by the so-called Brownian loop-soup introduced in \cite {LawlerWerner}. This is a random gas of non-interacting Brownian loops defined in a domain $D$. The law of the 
Brownian loop-soup is described by its positive intensity $c$; a loop-soup with intensity $2c$ is then the union of two independent loop-soups with intensity $c$. It is also possible to define a discrete analogue of these gases of Brownian loops: a so called random walk loop-soup. Both in the discrete and the continuum setting (when $d=2$) the loop-soup with intensity $c=1$ turns out to be 
directly connected to the GFF, while the loop-soup with intensity $c=2$ is directly related to uniform spanning trees (for instance via Wilson's algorithm).  
It also turns out -- \cite {SheffieldWerner} -- that letting $c$ vary between $0$ and $1$ one can construct  many CLEs directly as the collection of outer boundaries of clusters of Brownian loops in a loop-soup. 

Some of the striking properties of the random walk loop-soup with $c=1$ (the one that is related to CLE$_4$ and the GFF) correspond to combinatorial type identities that were, for instance, instrumental in the pioneering works of Brydges, Fr\"ohlich and Spencer \cite {BFS}. 
Again, one main difference in more recent developments is to use these gases of loops to construct random geometric objects such as clusters of loops, and not just to compute relevant interesting quantities. 
\medbreak 

The goal of these lecture notes is not to go through all the aforementioned items. It is rather
to provide a self-contained introduction to the Gaussian Free Field and its main properties, 
with an emphasis on more geometrical aspects (i.e., on some random geometric sets that can be coupled to the GFF): 

\begin {itemize}
 \item 
We will start with a gentle introduction to the discrete GFF; we will discuss its various resampling properties and decompositions. We will then study its spatial Markov property
and the closely related concept of local sets. We will also discuss features of its partition function, with a special role played by the determinant of the Laplacian, and its direct relation to random walk loop-soups. There will be one little detour via the GFF and loop-soups on cable-graphs, as recently worked out by Titus Lupu, and another via Wilson's algorithm to construct a uniform spanning tree. 

\item 
We will then move on to the continuum GFF. We will start by explaining what sort of random object (i.e, generalised function) it actually is, and how to make sense of various properties that generalise those of the discrete GFF. 
This can be somewhat tricky due to the fact that the continuum GFF is not defined pointwise.
In Chapter \ref {Ch4}, we will spend some time describing the Markov property and the important concept of local sets for the continuum GFF. 

\item 
In the subsequent chapter, we will focus on the continuum GFF in two dimensions, and describe some of its special features, such as its relation to SLE$_4$ and CLE$_4$. 
In particular, we will describe the main ideas that lead to the construction of the GFF via a family of nested CLE$_4$ loops: providing a topographic description of the field. In this chapter, we will 
try to provide most main ideas for the proofs, but will not go through all of the technical details (and this chapter should not be viewed as an introduction to SLE). 

\item 
In the final chapter, we very superficially browse without proofs through some further related topics, such as the Liouville Quantum Gravity area measure and its relation to SLE, the GFF with Neumann boundary conditions or the scaling limit of the uniform spanning tree in two dimensions. 
\end {itemize}

We stress that this is not a comprehensive survey of {\em all} the questions related to the GFF -- many important GFF-related questions (such as the question of which discrete models -- other than the discrete GFF -- have been shown to give rise to the continuum GFF in the scaling limit, or the recent developments related to constructive Conformal Field theory) will not be discussed or addressed. 

Some pointers to papers that discuss the results that we do present in the notes are given at the end of each chapter, but our bibliography is not meant to be a full list of all the relevant material present in the literature either.

\chapter {Warm-up} 
\label {warmup}

\section {Conditioned random walks} \label{section::one_dimensional_case}
Let us first recall some features of random walks and Brownian motions (more specifically, Brownian bridges) that will guide us as we try to construct the Gaussian Free Field.

\begin {reminder}
Recall that when $(B(t))_{t \in [0,1]}$ is a one-dimensional Brownian motion, then the process $( \beta_t := B_t - tB_1 )_{t \in [0,1]}$ is called a 
Brownian bridge. Basic considerations on covariance functions and Gaussian processes show that the process $\beta$ is a centred Gaussian process that is independent of the random variable $B_1$, so that its law can be interpreted as the law of Brownian motion ``conditioned to be equal to $0$ at time $1$''. The covariance structure of $\beta$ is $E ( \beta_t \beta_s ) = t (1-s)$ when $0 \le t \le s \le 1$. 
\end {reminder}

One-dimensional Brownian motion is known to be 
the scaling limit of a rather large class of random walks with independent and identically distributed increments (as soon as the laws of the individual steps of the walks have expectation $0$ and variance $1$). 
Similarly, the Brownian bridge is known to be the scaling limit of a rather large class of random walks, when they are {\em conditioned} to be back at $0$ after a large number of steps. For instance: 
\begin {enumerate}
\item Choose a path $(S(0), \ldots, S(N))$ with $N$ steps, when $N$ is even, with values in $\Z$,  uniformly from the set ${\mathcal S}_N$ of walks such that 
$$ S(0)= S(N)=0 \hbox { and } | S(j)- S(j-1)| = 1 \hbox { for all } 1 \le j \le N.$$  Then the law of $(S_{[Nt]} / \sqrt {N} )_{t \in [0,1]}$ 
is known to converge weakly (for the topology of the sup-norm on the space of real-valued right-continuous functions on $[0,1]$) to the law of the Brownian bridge  (here and in the sequel $[u]$ denotes the integer part of the real number $u$).  
Note that ${\mathcal S}_N$ has $N! / ((N/2)!)^2$ elements, as one only needs to choose the times of the $N/2$ upwards steps.
\begin{figure}[ht]
	\centering
	\includegraphics[scale=0.4]{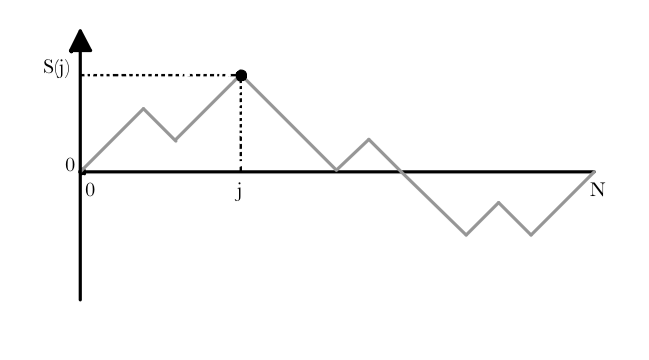}
	\caption{Linear interpolation of $(S(0),S(1),\cdots, S(N))$, with $N=16$.}
\end{figure} 

\item Take a symmetric density function $h(x)$  on $\R$ such that $\int x h(x) dx = 0$ and $\int x^2 h(x) = 1$, and consider the random vector $(S(1), \ldots, S(N-1))$ with density \oldcomments{(with respect to Lebesgue measure on $\R^{N-1}$)}
proportional to 
$$ \prod_{j=1}^{N} h( \gamma_j- \gamma_{j-1}) $$ 
at $(\gamma_1, \ldots, \gamma_{N-1})$ (with the convention $\gamma_0 = \gamma_N = 0 $). Then again, one can show that 
the law of $(S_{[Nt]} / \sqrt {N} )_{t \in [0,1]}$ converges  (in the same sense as above, which implies in particular the weak convergence of the finite-dimensional distributions) to the law of the Brownian bridge. 
\end {enumerate}
The proofs of these facts are not very difficult, but they do not fall into the scope of the present lectures. The results do illustrate however that Brownian bridges (and
constant multiples of the Brownian bridge) are indeed natural universal objects describing the fluctuations of a random function $f$ on $[0,1]$, constrained to satisfy $f(0)=f(1)=1$. 

\begin{figure}[ht]
	\centering
	\includegraphics[scale=0.5]{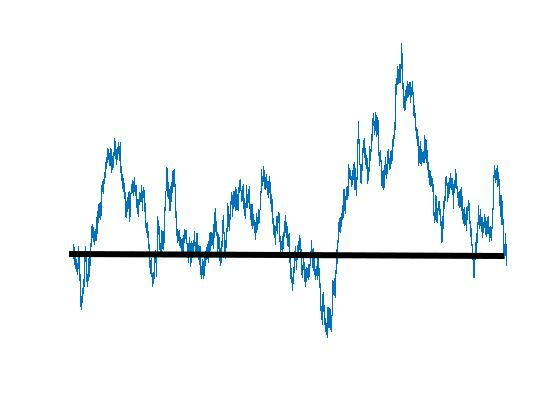}
	\caption{A Brownian bridge from zero to zero.}
\end{figure} 

\begin {remark}
It is worth noticing that for each given $N$, the laws of conditioned random walks of the type (1) or (2) can be viewed as the unique stationary measures of simple Markov chains on the space of ``admissible'' paths. For instance, in case (1) and when $N\ge 4$,
the natural dynamics on the space \oldcomments{${\mathcal S}_N$} can be described as follows. When we are given a path $\gamma$ in \oldcomments{${\mathcal S}_N$},  the Markovian \oldcomments{algorithm} to produce a new path $\gamma'$ is:   

(a) Choose a point $J$ uniformly at random in $\{1, \ldots , N-1 \}$. The new path $\gamma'$ will then be equal to $\gamma$ except possibly at time $J$. 

(b) If $\gamma (J-1) = \gamma (J+1)$, define $\gamma'$ to be equal to $\gamma$ except at time $J$, and set 
$$\gamma' (J) = \gamma  (J-1) - (\gamma (J)- \gamma (J-1)).$$  If $\gamma (J-1) \not= \gamma (J+1)$ (which means that 
$| \gamma (J+1) - \gamma (J-1) | = 2$), then keep $\gamma$ unchanged, i.e., set $\gamma' = \gamma$. 

It is then a simple exercise to check that this Markov chain is irreducible, aperiodic and that the uniform measure on ${\mathcal S}_N$ is reversible (indeed, if the probability to jump from $\gamma$ to $\gamma'$ when $\gamma' \not= \gamma$ in one step 
is positive, then it is equal to $1/(N-1)$, and equal to the probability to jump from $\gamma'$ to $\gamma$). \oldcomments{Hence the law of the conditioned random walk in case (1) is equal to the unique stationary law of this Markov chain.} 

In fact, an even more natural alternative to (b) is to toss a fair coin in the case where  $\gamma (J-1) = \gamma (J+1)$ in order to decide whether $\gamma' (J) - \gamma'(J-1)$ is equal to $+1$ or $-1$. Again, the 
uniform measure on ${\mathcal S}_N$ is the unique stationary measure for this dynamic. 
\end {remark}

\begin{figure}[ht]
	\centering
	\includegraphics[scale=0.45]{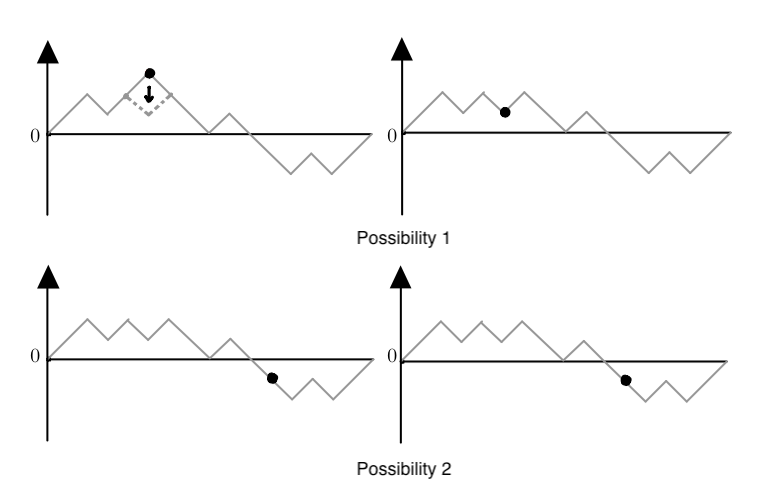}
	\caption{Two example steps in the Markov chain. The top figures illustrate the first possibility described in (b), and the bottom figures the second. The vertex selected in step (a) is marked with a dot. }
\end{figure} 

\begin {exercise} 
Describe a similar natural irreducible Markov chain on the state \oldcomments{space} of functions from $\{ 1, \ldots, N-1 \}$ into $\R$, such that the law described in (2) is an invariant stationary measure for this chain. 
\end {exercise}

There is one special case of type (2) conditioned walks that is worth highlighting. 
This is when one takes $h$ to be the Gaussian distribution function with variance $1$ i.e., $h(x)= \exp (- x^2/2) / \sqrt {2\pi}$. 
Then  $(S(1), \ldots, S(N-1))$ is a centred Gaussian vector, and its covariance function is easily shown to be given by 
$$ E [ S(j) S(j') ] =  j (N-j') / N$$ 
when $1 \le j \le j' < N$. In particular, 
$$E [ (S(j) / \sqrt {N}) \times (S(j') / \sqrt {N} ) ] = (j/N) \times (1- (j'/N)).$$ 
Note that in fact, if $\beta = ( \beta_t, t \in [0,1])$ is itself a Brownian bridge, then the vector
$(\sqrt {N} \beta (1/N), \ldots, \sqrt {N} \beta ((N-1)/N))$ is distributed exactly like $(S(1), \ldots, S(N-1))$. In this case, the convergence in distribution of
the conditioned walk to the Brownian bridge is then a direct consequence of this observation and of the almost sure continuity of the Brownian bridge.

\section {Concrete examples in the discrete square.} 
\label{section::examples_discrete_square}

What is the corresponding object \oldcomments{describing fluctuations}, when instead of considering a one-dimensional string, one looks at some tambourine skin? 
In other words, what happens in the previous cases when one replaces the one-dimensional 
time-segment $[0,1]$ by a two-dimensional set $D$ (that plays the role of the shape of the tambourine),
and tries to look at random functions from $D$ into $\R$? 

Let us start with discrete models, defined on grid approximations of $D$. 
To be specific, let us consider $N \ge 2$ and define $\overline \Lambda_N := \{0, \ldots, N \}^2$ to be the closed $N \times N$ discrete square. We let $\Lambda_N := \{1, \ldots, N -1 \}^2$ be the {\em inside} of the square and 
$\partial_N := \overline \Lambda_N \setminus {\Lambda}_N$ be its {\em boundary}. We denote by $E_N$ the set of (unoriented) edges that join two neighbouring points (i.e., at distance $1$) in $\overline\Lambda_N$. 
Let us consider the family of functions $f$ from the discrete square $\overline \Lambda_N$ into $\R$, with the constraint that $f$ is equal to zero on $\partial_N$.
Here are some concrete ways to choose such a function $f$ at random: 
\begin {enumerate} 
\item Choose $f$ uniformly among the finite set of all integer-valued functions $f$ such that (a) $f=0$ on the boundary of the square, and (b) for any $x$ in $\{1, \ldots , N-1 \}^2$ and any $y$ neighbouring $x$ (i.e. in $\overline\Lambda_N$ and at distance 
$1$ from $x$), $f(x)- f(y)  \in \{ -1, 0, 1 \}$. This is somehow the  analogue of the discrete random walk (1) from Section \ref{section::one_dimensional_case}, when it is also allowed to stay constant (it is useful to use this variant here in order to avoid parity constraints due to the boundary conditions). 
\item One can also consider the following continuous analogue: choose a function uniformly (i.e., with respect to the Lebesgue measure on $\R^{ \Lambda_N}$)  
in the set of all {\em real-valued} functions $f$ such that for any $x$ in $\{1, \ldots , N-1 \}^2$ and any $y$ neighbouring $x$, $|f(x)- f(y) |  \le 1$
(with the convention that $f=0$ on the boundary of the square). This is the analogue of a discrete random walk bridge, where the steps of the walk are 
chosen uniformly in $[-1, 1]$. 
\item More generally: when $h$ is the density function of a symmetric $L^2$ random variable with zero mean, one can choose $f$ in such a way that the random vector $(f(x))_{x \in \Lambda_N}$ has density \oldcomments{(with respect to Lebesgue measure on $\R^{(N-1)^2}$)}
proportional to 
$$\prod_{e \in E_N} h (| \nabla \gamma(e)  | )$$ 
at $(\gamma_x)_{x \in \Lambda_N}$, where here and in the sequel, $| \nabla  \gamma (e)  |$ denotes the absolute value of the difference between the two values of $\gamma$ at the two extremities of the edge $e$. We use this for vectors $(\gamma_x)_x$ and functions $(f(x))_x$ interchangeably (with the obvious interpretation).
\end {enumerate}
One way to think about it is that each edge $e \in E_N$ consists of a little spring (so that the tambourine skin is actually made of a little trampoline web of springs). Each point $x$ in $\Lambda_N$ (in the horizontal plane) is allowed to move vertically (in some third direction perpendicular to $\overline \Lambda_N$) to the position 
$(x, \gamma(x))$ in three-dimensional space, while the boundary points $x \in \partial_N$ are stuck to height $0$. 
The spring on the edge $e$ puts some constraints on the height-difference between the two extremities of $e$, and in particular tends to prevent this difference from being very large. 

\begin{figure}[h]
	\centering
	\includegraphics[scale=0.35]{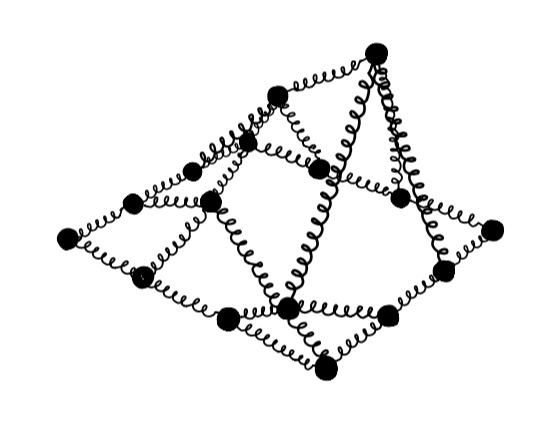}
	\caption{An illustration when $N=3$ and $d=2$.}
\end{figure} 
Just as in the previous one-dimensional case, each of these measures can be viewed as the stationary measure of some rather simple Markov chain on the state \oldcomments{space} of functions from $\Lambda_N$ into $\R$, where at 
each step of the chain, one resamples the value \oldcomments{(height)} of the function at at most one site, according to the conditional distribution of that height given those of its neighbours. 

Then, by analogy with the previous one-dimensional case, one would like to argue that all of these models, when $N \to \infty$ and when appropriately rescaled, do converge to the same random object, that is some sort 
of random function from $[0,1]^2$ into $\R$. 
For instance, one can first transform any of these discrete random functions $f_N$ on $\overline \Lambda_N$ into a function defined on $[0,1]^2$ simply by rescaling (and making the function constant on each square): 
$$ \hat f_N (x_1, x_2) := f_N (  [Nx_1], [Nx_2 ]).$$ 
Then, the hope is that for some good choice of sequence $\eps_N$, the law of $\eps_N \hat f_N$ will converge to that of some ``universal random function'' $f$ from $[0,1]^2$ to $\R$. 

As we will see very soon, the story turns out to be a little more subtle due to the actual nature of this universal random function $f$, but the conjecture is roughly that this should be correct. Loosely speaking:  

\begin {conjecture} 
For each of the aforementioned models \oldcomments{(1)-(3) of Section \ref{section::examples_discrete_square}}, 
one can find a sequence $\eps_N$ (actually we will see that in this two-dimensional case $\eps_N$ should be constant) such that in some appropriate sense, $\eps_N \hat f_N$ converges in distribution 
to a universal non-trivial random generalised function. 
\end {conjecture}

This is actually still a conjecture for most of the examples mentioned above! 
There exist a couple of cases where this is known to be true (for instance when $h$ is the exponential of a uniformly concave function), but for case (1), this is (to our knowledge) an open problem.  
In these lectures, we will actually not discuss these universality questions at all. Rather, we will first focus on the special Gaussian subcase of example (3), for which one can: 
\begin {itemize}
\item  say a lot in the discrete case, 
\oldcomments{which already gives rise to}  combinatorially very rich mathematical objects; 

\item show very easily that (when suitably rescaled), the \oldcomments{discrete models} converge in distribution \oldcomments{as $N\to \infty$} to their counterparts in the continuum. 
\end {itemize}

This particular example is that of the {\em discrete Gaussian Free Field} (we will use the acronym GFF \oldcomments{for Gaussian Free Field} throughout these notes). 
This is the case where the function $h(u)$ is the distribution function of a Gaussian random variable i.e.,   $\exp (- u^2 / 2 \sigma^2)$ for some choice of $\sigma^2$. 
So, the discrete GFF is the probability measure on $\R^{\Lambda_N}$ with density at $(\gamma_x)_{x \in { \Lambda_N}}$ a constant multiple of 
$$ \exp ( - \sum_{e \in E_N} | \nabla  \gamma (e) |^2 / (2 \sigma^2) )$$ 
with the convention that $\gamma =0$ on $\partial_N$.

In this case,  the obtained random  function $f_N$ is a centred Gaussian vector. Hence, its law is fully described via its covariance function, and if one controls this covariance function well in 
the limit when $N \to \infty$, one will obtain convergence to some Gaussian object in the continuum space (with covariances given by limit of the covariances).
\oldcomments{Hence, we can determine what the continuous object that we are looking for should be.} 

The structure of the lectures will be the following. In the next two chapters, we will define and study some features of this discrete GFF, focusing especially on those that will have 
a natural analogue in the continuum. 
We will then discuss the definition of the continuum GFF in an arbitrary number of dimensions and describe some of its properties. Finally, we will restrict to the case of two dimensions, and survey some of the special results that hold in this setting.

\chapter {The discrete GFF} 
\label{Ch1}
\label{chapter:discrete GFF}

\section {Definition}

\subsection {Notation} 

Before defining the discrete GFF, let us first introduce some notation that we will use throughout these notes. We suppose that $d\ge 1$. 

When $f$ is a function from $\Z^d$ into $\R$, we define $\overline f(x)$ to be the average value of $f$ at the $(2d)$ neighbours of $x$. In other words,  
$$ \overline f (x) =  \frac 1 {2d}  \sum_{y  :  y\sim x} f(y), $$
where here and in the sequel, $ \sum_{y  :  y\sim x}$ means that we sum over the $2d$ neighbours of $x$ in $\Z^d$.

\begin {definition}[Discrete Laplacian -- \textbf{careful, this is not the standard definition}] \label{def:dlap}

We define the discrete Laplacian 
$\Delta f$ of $f$ to be the function 
$$\Delta f (x) := \overline f(x) - f(x).$$  
\end {definition} 

\begin {remark} 
We would like to emphasise that throughout these lecture notes, we are going to use Definition \ref{def:dlap} of $\Delta$ to be our discrete Laplacian. This is \emph{not} the standard 
definition that one finds in most textbooks, where the discrete Laplacian is often defined as $\sum_{y : y \sim x } ( f(y) - f(x))$ (so it differs by the multiplicative factor $2d$). 
\end {remark}

When $D$ is a subset of $\Z^d$, we define its (discrete) boundary 
$$ \partial D :=\{x\in \Z^d: d(x,D)=1 \} \hbox { and } \overline D := D\cup \partial D.$$ 

We will denote by $ \F_{(D)} $ the set of functions from $\Z^d$ into $\R$ that are equal to $0$ outside of $D$. 
When $D$ is finite and has $n$ elements, then $\F_{(D)}$ is of course a real vector space of dimension $n$. 

When $F$ is a function from $\overline{D}$ into $\R$ (which is not defined outside of $D \cup \partial D$) then we can 
still define $\overline F (x)$ and $\Delta F (x)$ for all $x \in D$ just as before.

  We define the set {$E_{\overline{D}}$} to be the set of edges of $\Z^d$ such that at least one end-point of the edge is in $D$. 
  For each $F \in \F_{(D)}$ and each unoriented edge $e \in E_{\overline{D}}$, we define 
$| \nabla F (e) | := | F(x)- F(y)|$ as before, where $x$ and $y$ are the two endpoints of $e$. Note that to decide about the sign of $ \nabla F $, we would need to consider oriented edges, but that
$|\nabla F (e)|$ and its square do not depend on the orientation of $e$. Similarly, when $F_1$ and $F_2$ are in $\F_{(D)}$, we can define unambiguously the 
product $\nabla F_1 (e) \times \nabla F_2 (e)$. 
Finally, when $D$ is finite we define 
$$ {\mathcal E}_D (F) := \sum_{e \in E_{\overline{D}}} |\nabla F (e)|^2.$$ 
This quantity (or half of this quantity) is often referred to as the Dirichlet energy of the function $F$. 
\begin{figure}[h]
	\centering
	\includegraphics[scale=0.3]{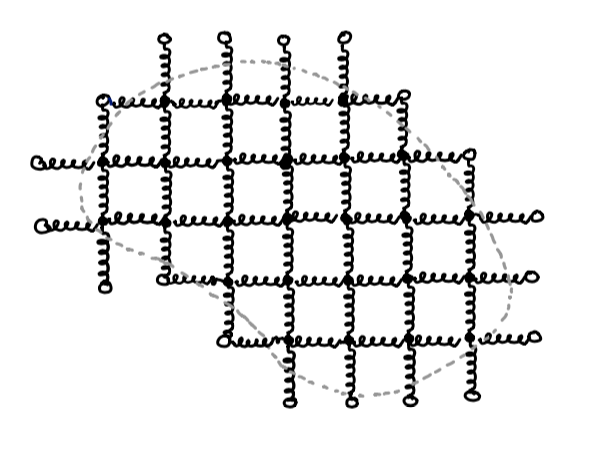}
	\caption{A domain $D\subset \Z^2$, formed by taking all $z\in \Z^2$ that lie inside a domain $\Omega\subset \R^2$ (the boundary of $\Omega$ is represented by the dotted line). Solid discs represent points of $D$, and open discs points of $\partial D$. Each edge in $E_{\overline{D}}$ is depicted as a ``spring''.}
\end{figure} 

\subsection {Definition via the density function} 

\begin {definition}[Discrete GFF via its density function]
The discrete GFF in $D$ with Dirichlet boundary conditions (also sometimes referred to as zero boundary conditions) on $\partial D$ is the centred Gaussian vector  $(\Gamma(x))_{x \in D}$ whose density function on $\R^D$ at $(\gamma_x)_{x \in D}$  
is a constant multiple of 
$$ \exp ( - \frac {1}{2} \times \frac {  {\mathcal E}_D (\gamma)} {2d}  ) =  \exp ( - \frac {1}{2 } \times \frac {1}{2d}  \sum_{e \in E_{\overline{D}}} | \nabla \gamma (e) |^2  )  $$ 
with the convention that $\gamma = 0$ on $\partial D$. 
\end {definition} 

 \begin {remark}
 {We use the notation $(\gamma_x)_{x\in D}$ rather than $(\gamma(x))_{x\in D}$ to distinguish it as a fixed vector. The quantity $|\nabla \gamma(e)|$ when $e$ has endpoints $\{x,y\}$ is equal to $|\gamma_x-\gamma_y|$.}
\end {remark} 

Note that by definition
$(\gamma_x)_{x \in D} \mapsto {\mathcal E}_D (\gamma) $ is a bilinear form, and it is also  positive definite (indeed if ${\mathcal E}_D (\gamma)$ is $0$, it means that $|\nabla \gamma (e) |= 0 $ on all edges, so that $\gamma$ is identically $0$). Thus, 
the exponential above is indeed a multiple of the density function of some Gaussian vector on $\R^D$, and this definition makes sense.  

\begin {remark} 
We could also introduce a positive parameter $\sigma$ to the model, in order to heuristically describe the ``stiffness'' of springs associated with the edges in $E_{\overline{D}}$. This would lead us to consider the random field with density function instead given by a multiple 
of 
$$  \exp ( - \frac {1}{2} \times \frac {  {\mathcal E}_D (\gamma)} {2d \sigma^2}  ) .$$ 
The random process $(\Gamma_{(\sigma)} (x))_{x \in D}$ obtained in this way is clearly equal in distribution to $( \sigma \Gamma(x))_{x \in D}$. 
\end {remark}

Recall that the law of a centred Gaussian vector is completely determined by its covariance function. It will turn out that the covariance function of the Gaussian Free Field is very nice, and we will come back to this later.

\subsection {Resampling procedure and consequences} 

Suppose that $x$ is a given point in $D$. What is the conditional distribution of $\Gamma(x)$ given $(\Gamma(y))_{y \in D \setminus \{ x \}}$? 
An inspection of the density function of $\Gamma$ shows that the conditional distribution of $\Gamma(x)$ given $(\Gamma(y))_{y \in D \setminus \{ x \}} = (h(y))_{y \in D \setminus \{ x \}}$
has a density at $(\gamma_x)_{x\in D}$ that is proportional to 
$$ \exp \big( -  \frac {1}{2 \times (2d)} \sum_{y: y \sim x} | \gamma_x - h(y) |^2 \big) .$$ 
Expanding this sum over $y$, we get that this   
is equal to 
$$ \exp(-\frac{1}{2}(\gamma_x-\overline h(x))^2)$$
times some normalising function that depends only on $h$.
In other words, this conditional law is that of the Gaussian distribution ${\mathcal N} ( \overline h (x), 1 )$. 

A first  feature worth stressing (which is due to the interaction via nearest-neighbours only) is that this conditional distribution depends only on the values $h(y)$ at the neighbours $y$ of $x$. 
A second feature is that in fact, the conditional law of $\Gamma (x) - \overline h(x)$ is a standard normal Gaussian (for all choices of $\overline h (x)$). This means that, for all $x$,
$ \Gamma (x) - \overline \Gamma (x)$ is a standard Gaussian random variable that is independent of $(\Gamma (y))_{y \in D \setminus \{ x \}}$. This fact has a number of important consequences.

A first consequence is that it indicates what the natural Markov chain (on the space of functions) is, for which the law of the GFF is stationary. \oldcomments{For this chain, the Markovian step can be described as follows:} 
if we are given a function $h$ in ${\mathcal F}_{(D)}$, then we choose a point $x \in D$ uniformly at random, and replace the value of $h(x)$ by $\overline h(x) + N$ where 
$N$ is a standard Gaussian random variable.  

A second consequence is that it allows us to derive immediately some interesting properties of the covariance function of $\Gamma$. 
For all $x$ and $y$ in $D$, let us denote this covariance function by 
$$ \oldcomments{\Sigma} (x,y) = \oldcomments{\Sigma}_x (y)  := E [ \Gamma (x) \Gamma(y) ].$$ 
In this way, one can view for each given $x$, $y \mapsto \Sigma_x (y)$ as a function in ${\mathcal F}_{(D)}$. 

Then, when $x \not= y$ are both in $D$, 
\begin {eqnarray*}
\lefteqn {\Sigma_x (y) = E [ \Gamma (x) \overline \Gamma(y) ] + E [ \Gamma (x) (\Gamma (y) - \overline \Gamma (y) )]}\\
&&= E [ \Gamma (x) \overline \Gamma (y) ]
=\frac 1 {2d} \sum_{z: z \sim y} E[ \Gamma (x) \Gamma(z) ]
= \overline {\Sigma_x}  (y) .\end {eqnarray*}
Similarly, 
\begin {eqnarray*}
\lefteqn { 
	\Sigma_x (x) = E [ \Gamma (x) \Gamma (x) ] = E [ \Gamma (x) \overline \Gamma (x) ] +  E [ (\Gamma (x) - \overline \Gamma (x)) \Gamma (x)] }\\
&&= (2d)^{-1} \sum_{z: z \sim x} E[ \Gamma(z) \Gamma (x)] +  E [ (\Gamma (x) - \overline \Gamma (x))^2 ] + 
E [   (\Gamma (x) - \overline \Gamma (x)) \overline \Gamma (x) ]\\
&& 
=  \overline{  \Sigma_x} (x)  + 1 + 0.\end {eqnarray*}
In other words, the function $\Sigma_x$ satisfies 
$$ \Delta \Sigma_x (y) = \oldcomments{-} \1{y=x} $$ for all $y$ in $D$. 
Note that (for each given $x$) this provides as many linear equations as there are entries for $\Sigma_x ( \cdot)$ (both sets have the cardinality of $D$). 
As we will see in a moment, these equations are clearly \oldcomments{linearly} independent, so that these relations fully determine $\Sigma_x$.

\subsection {The discrete Green's function}

The previous analysis leads us naturally to quickly review and  browse through some basic definitions and properties 
related to the discrete Laplacian and the discrete Green's function. 

\medbreak 
\subsubsection*{The discrete Laplacian}
Recall that for all  $F \in \F_{(D)}$, we defined for $x \in D$,
$$
\Delta F(x):= \frac{1}{2d} \sum_{y  :  y\sim x} (F(y)-F(x))=\overline F(x)-F(x)
$$ 
By convention, we will denote by $\Delta_D F$ the function that is equal to $\Delta F$ in $D$ and is equal to $0$ outside of $D$ 
(mind that here we do not care about the value of $\Delta F$ outside of $D$, in particular on $\partial D$). Again, we stress that this is not the most standard 
definition of the discrete Laplacian (our definition is $1/ (2d)$ times the usual one). 

Clearly, we can then view $\Delta_D$ as a linear operator from $\F_{(D)}$ into itself. 
It is easy to check that $\Delta_D$ is injective using the maximum principle. [If $\Delta_D F =0$, then choose $x_0 \in D$ so that $|F(x_0)| = \max_{x\in D} |F(x)|$, and because $\Delta_D F (x_0) = 0$,
this implies readily that the value of $F$ on all the neighbours of $x_0$ are all equal to $F( x_0)$ (as otherwise, their mean value could not be equal to $F(x_0)$). 
But then, this also holds for all neighbours of neighbours of $x_0$ as well. 
Eventually, since $D$ is finite, this means that we will find a boundary point $y$ for which $F (y) = F(x_0)$. Since $F = 0$ on the boundary, it follows that $\max_{x \in D} |F(x)| = |F(x_0)| = 0$]. 

Hence, $\Delta_D$ is a bijective linear map from the vector space $\F_{(D)}$ into itself. 
One can therefore define its linear inverse map: for any choice of function $u: D \to \R$, there exists exactly one function $F \in \F_{(D)}$ such that $\Delta_D F(x) = u(x)$ for all $x \in D$.  

If we apply this to the previous analysis, it shows that indeed, 
$y \mapsto \Sigma_x (y)$ is the unique function in $\F_{(D)}$ such that its Laplacian $\Delta_D$ in $D$ is the function $y \mapsto \oldcomments{-} \1{y=x}$. As we will see in a moment, this function has a name...

\medbreak

\subsubsection* {The Green's function}
Let $(X_n)_{n\geq 0}$ be a simple random walk in $\Z^d$, with law denoted  by $\P_x$ when it is started at $x$. Let $\tau= \tau_D := \inf\{ n \geq 0: X_n\notin D\}$ be its first exit time from $D$.

\begin {definition}[Green's function]
We define the Green's function 
$G_D$ in $D$ to be the function defined on $D \times D$  by
$$
G_D(x,y):= \E_{x}\Bigl[\sum_{k=0}^{\tau-1}\1{X_k=y} \Bigr].
$$ 
\end {definition}
By convention, we will set $G_D (x,y) = 0$ as soon as one of the two points $x$ or $y$ is not in $D$.
It is sometimes convenient to reformulate this definition in a more symmetric way that highlights that $G_D (x,y) = G_D( y,x)$:  
\begin{align*}
&G_D(x,y) = \E_{x}\Bigl[\sum_{k\geq 0} \1{X_k=y,k <\tau} \Bigr] = \sum_{k\geq 0} \P_x(X_k=y,k<\tau)\\
&= \sum_{k\geq 0} \# \{\text{paths $x\to y$ in $k$ steps within $D$} \} \times \Bigl[ \frac{1}{2d} \Bigr]^k
\end{align*}
and this last expression is clearly symmetric in $x$ and $y$ (the number of paths from $x$ to $y$ with $k$ steps in $D$ is equal to the number of paths from $y$ to $x$ with $k$ steps in $D$).

Let us now explain why the following result holds.
\begin {proposition}
The Green's function $G_D$ is the inverse of $-\Delta_D$, and it is equal to $\Sigma$.   
\end {proposition}

\begin {proof} 
We will use a slightly convoluted, but hopefully instructive, strategy to prove this (see the remark below for a more direct approach). 
The idea is that the Markov property of the simple random walk immediately enables us 
to determine the Laplacian of the function  $g_{D,x} (\cdot)= G_D(\cdot, x)$ in $D$ (note that $g_{D,x} \in \F_{(D)}$, as $g_{D,x}$  is equal to zero outside of $D$). Indeed, we have that for all $y \neq x$ in $D$, $\Delta_D g_{D,x}(y)=0$, simply because
$$
G_D(y,x)= \E_{y}\Bigl[\sum_{k\geq 1} \1{X_k=x,k <\tau} \Bigr]
= \sum_{z : z\sim y} \frac{1}{2d} G_D (z, x) ,
$$
where we have used the Markov property at time $1$ in the first identity.
Also, the very same observation (but noting that at time $0$, the random walk starting at $x$ is at $x$) shows that $\Delta_D g_{D,x}(x)=-1$. 
Hence, $g_{D,x}$ is a function in $\F_{(D)}$ satisfying 
$$\Delta_D g_{D,x}(y)=-\1{x=y}$$ for all $y$ in $D$. 
Since $\Delta_D$ is a bijection of $\F_{(D)}$ onto itself, the function $g_{D,x}$ is in fact the unique function in $\F_D$ with this property. We therefore conclude that for all $x$ and $y$ in $D$, 
$$ \oldcomments{\Sigma} (x,y) = G_D (x, y).$$
\end {proof}

\begin {remark}
\label {Pn}
For the record, let us also mention that there is a two-line proof of the fact that $G_D$ is the inverse of $-\Delta_D$, that does not rely on any of our previous considerations. Note that 
the matrix $P_D := I + \Delta_D$ is the transition matrix of the simple random walk on  $\Z^d$, when restricted to $D$, since $P_D (x,y)$ corresponds to the probability to jump from $x$ to $y$. Let us label the $n$ points of $D$ by $\{ x_1, \ldots, x_n\}$, so that we can view (and we will use this type of notation on numerous occasions in these notes) the functions $G_D$, $-\Delta_D$ and $\Sigma$ defined on $D \times D$ as $n \times n$ matrices. Then, it is clear that for all $k \ge 0$, 
$$ \P_x [ X_k = y , k < \tau ] = (P_D)^k (x, y),$$
where $(P_D)^k$ is the $k$-th power of the matrix $P^D$. Hence, 
$$ G_D (x, y) = \sum_{k \ge 0} (P_D)^k (x,y) $$ 
from which it follows that $G_D$ is equal to the inverse of $(I - P_D)$, that is, $- \Delta_D$. 
\end {remark}

This provides the following equivalent definition of the discrete Gaussian Free Field: 
\begin {definition}[Discrete GFF via the covariance function]
The discrete Gaussian Free Field in $D$ with Dirichlet boundary conditions on $\partial D$ is the centred Gaussian process $(\Gamma(x))_{x \in D}$ with covariance function $G_D (x,y)$ on $D \times D$. 
\end {definition}

\begin {remark} 
We see that, as opposed to the first definition, this second equivalent definition actually also works when $D$ is infinite, so long as $G_D$ is well-defined. That is, as long as the random walk in $D$, killed when it reaches $\partial D$, is transient. 
In other words, the definition can also be used for any infinite subset of $\Z^d$ when $d \ge 3$ (because the simple random walk on $\Z^d$ is transient), or for any infinite subset $D \not= \Z^d$ when $d=1,2$.
\end {remark}

\begin {remark} 
The two definitions are equivalent. It is a matter of taste whether one prefers to use the more hands-on (and maybe more intuitive) approach via density functions or the slightly more general setting of Gaussian processes, when one wants to derive properties of the GFF.   
\end {remark}

\section {Informal comments about the possible scaling limit} 
\label{sec::informal_sl}

In this section, 
we use the above definition of the discrete Gaussian free field to formulate some heuristics about how a ``continuum Gaussian free field'' on a subset of $\R^d$ could be defined. 
This section is non-rigorous, and can be viewed as an appendix to the warm-up chapter. It serves only as an appetiser to the actual study of the continuum GFF later on. 

Suppose that $\Omega$ is some open subset of $\R^d$ for $d \ge 1$. The idea is to approximate the continuum process $(\Gamma (x))_{x \in \Omega}$ that we would want to define, using the 
GFF on a fine grid approximation of $\Omega$. For each positive $\delta$, one can for instance define $\tilde D_\delta = \delta \Z^d \cap  \Omega$ and $D_\delta = \delta^{-1} \tilde D_\delta = \Z^d \cap (\delta^{-1} \Omega)$ so that $\tilde{D}_\delta$ is a subset of the fine grid $\delta \Z^d$, which is a good approximation to $\Omega$, and $D_\delta$ is its $(1/\delta)$ blow-up: a subset of $\Z^d$.
We can therefore define the GFF $\Gamma_\delta$ on $D_\delta$ as in the previous section, and a GFF $\tilde \Gamma_\delta$ on $\tilde D_\delta$ by setting $\tilde \Gamma_\delta (x) = \Gamma_\delta (x\delta^{-1})$. In other words, $\tilde \Gamma_{\delta}$ is a GFF on the grid approximation $\tilde{D}_\delta$
of $\Omega$ in $\delta \Z^d$, normalised in such a way that the variance of the difference between $\tilde \Gamma_\delta (x)$ and the mean value of its $2d$ neighbours in $\tilde{D}_\delta$ is equal to $1$ for all $x \in \tilde D_\delta$. 

We can extend this random function $\tilde \Gamma_\delta$ to all of $\R^d$ by (for instance) choosing $\tilde \Gamma (y) = \tilde \Gamma (x)$ for all $y = (y_1, \ldots, y_d) \in [x_1, x_1 +\delta) \times \ldots \times [x_d, x_d +\delta)$ 
when $x \in \delta \Z^d$. 

Now the philosophy is the following: when a centred Gaussian process converges in law (which is exactly when all its finite-dimensional distributions converge), then the limiting law is bound to be a centred Gaussian process as well, and the covariances 
of the limit are the limits of the covariances. 

\begin{exercise} Let $V$ be a finite set and let $(\Gamma_n(x))_{x\in V}$ be a centred Gaussian process for every $n\in \N$ with $E[\Gamma_n(x)\Gamma_n(y)]=:\Sigma_n(x,y)$. Suppose that for every $x,y\in V$, $\Sigma_n(x,y)\to \Sigma(x,y)$ for some positive definite bilinear form $\Sigma: V\times V\to \R$. Show that $\Gamma_n$ converges in distribution to $\Gamma$: the centred Gaussian process $(\Gamma(x))_{x\in V}$ with covariance matrix $\Sigma$ \end{exercise}

So, it is natural to look at what happens to the covariance function of $\tilde \Gamma_\delta$ as $\delta \to 0$. Let us collect here some observations and facts, leaving out any detailed proof: 
\begin {enumerate} 
\item When $x \not= y$ in $\Omega$, then it turns out that as $\delta \to 0$, 
$$ G_{D_\delta} ( x \delta^{-1} , y \delta^{-1} ) \sim \delta^{d-2} G_\Omega (x, y),$$
where $G_\Omega(x, y)$ is some positive function of $x$ and $y$ (called the continuum Green's function, but we will not discuss this here). 
The main point to note is that this quantity converges when $d = 2$, but tends to $0$ when $d > 2$.  A simple way to understand the formula above is to note that 
the mean number of steps spent by the random walk before exiting a compact portion of $\Omega$ is of the order of $\delta^{-2}$ (this $2$ comes from the central limit theorem renormalisation). 
On the other hand, in expectation, this time is spread rather regularly among all points $y$ (when $y$ is not too close to $x$), and the number of such points $y$ is of the order of $\delta^{-d}$. 
Hence, we should not be surprised by the coefficient $\delta^{d-2}$. 
\item 
As a consequence, when $d = 2$, we see that the covariances converge to something non-trivial {\em without any rescaling}. In other words, one would like to simply take the limit of $(\tilde \Gamma_\delta (x))_{x \in \Omega}$ to define the continuum GFF in $\Omega$. 
We already see that such a limit is unlikely to be a continuous function (which will be why we refer to it as the ``continuum Gaussian free field'' -- this name coming from the fact that it is defined in the \emph{continuum}), because the variance of the difference between the values of $ \tilde \Gamma_\delta$ at 
two points that are $\delta$ apart in $\Omega$ will be of order $1$, and in particular will not go to $0$. In fact, $\mathbb{E}[(\tilde{\Gamma}_\delta(x))^2]$ will grow like $\log(1/\delta)$ as $\delta\to 0$: see Exercise \ref{ex:disc_evalues} for an example.
\item 
When $d\ge 3$, things are even worse! In order to get a limit for the covariance function, point (1) implies that we need to rescale $\tilde \Gamma_\delta$ and to look instead at $\delta^{1- d/2} \tilde \Gamma_\delta $. 
This time, it means that the variance between the value of $\delta^{1- d/2} \tilde \Gamma_\delta $ at a point $x$ and its mean-value among the $2d$ neighbours of $x$ in $\tilde{D}_\delta$ is not only going to stay positive as $\delta\to 0$, but will actually blow up. Hence, the stiffness of the springs in our intuitive picture is going to vanish quickly as $\delta \to 0$. It therefore seems that in the limit, any obtained process must be unbounded everywhere, and equal to $\pm \infty$ simultaneously at each point of $\Omega$!
\item 
We finally observe that for $x \in \Omega$ 
the variance of $\delta^{1- d/2} \tilde \Gamma_\delta (x)$ tends to infinity as $\delta \to 0$ (when $d=2$, this follows from recurrence of random walk in $\Z^2$). 
So, any limiting process cannot possibly be defined as a random function, as it would then be a centred Gaussian with infinite variance. 
We could try to fix this by renormalising $\tilde \Gamma_\delta$ by some constant $\eps (\delta)$, so that the variance of $\eps(\delta) \tilde \Gamma_\delta (x)$ converges to something finite, and the process has a proper Gaussian limit. However, the covariance function of the limit would then be $0$ on $\{ (x, y) \in D \times D , \ x \not= y \}$, so that the limiting process would consist 
of a collection of independent Gaussian random variables (one for each point in the domain $D$). This  is clearly not the interesting process that we are looking for! 
\end {enumerate}

As we shall see, in a later chapter, the proper way to define the Gaussian free field in the continuum will be to view it 
as a random \emph{generalised function} rather than as a normal (point-wise defined) function.

In the remainder of this chapter and in the next chapter, we will actually continue to focus on aspects of the discrete GFF. These will turn out to have natural counterparts for the continuum GFF later on.

\begin{figure}[h]
	\centering
	\includegraphics[scale=1]{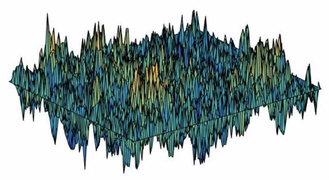}
	\caption{A simulation of $\tilde{\Gamma}_\delta$ on a square.}
\end{figure}

\section {Variations on the Markov property} 

Now we would like to ask: is there an analogue of the Markov property for the simple random walk that extends to the setting of the discrete GFF? In this section we will use the more hands-on definition of the GFF via density functions, as it provides a little more insight. However the Gaussian process setting is also very well suited to elegantly derive some of the Markovian properties that we discuss here.

We remark at this point that in the previous sections we did define the discrete GFF  in any finite subset of $\Z^d$ (i.e.,  we did not assume this set to be connected).

\subsection {The GFF with non-zero boundary conditions} 

In view of our intuitive description of the GFF, it is natural to generalise our definition to the case of non-zero boundary conditions. 
More precisely, suppose that $f$ is some given real-valued function defined on $\partial D$. Then, the definition of the GFF via its density function can be extended as follows: 
\begin {definition}[Discrete GFF with non-zero boundary conditions, via its density function]
The discrete GFF in $D$ \emph{with boundary condition $f$} on $\partial D$ is the Gaussian vector  $(\Gamma(x))_{x \in D}$ whose density function on $\R^D$ at $(\gamma_x)_{x \in D}$ 
is a constant multiple of 
$$ \exp ( - \frac {1}{2} 
\times  \frac {{\mathcal E}_D (\gamma)}{2d} ), $$
with the convention that $\gamma = f$ on $\partial D$. Note that the values of $f$ on $\partial D$ are implicitly used in the expression of ${\mathcal E}_D (\gamma)$ via the terms $|\nabla \gamma (e)|$ for those edges $e \in E_{\overline{D}}$ having one endpoint in $\partial D$. 
\end {definition} 
In other words, instead of fixing the height of $\Gamma$ on $\partial D$ to be $0$, we now fix it to be $f$. Then $\Gamma$ is still a Gaussian process, but it is not necessarily centred.

\medbreak 

Let us now make a few simple comments. A first, obvious, observation is that when $f$ is constant and equal to $c$ on $\partial D$, then if $(\Gamma(x))_{x \in D}$ is a GFF with boundary condition $f$, $(\Gamma(x) -c)_{x \in D}$ is a GFF with Dirichlet boundary conditions. A second immediate observation, that can be deduced directly from the expression of the density function for $\Gamma$ is the following: suppose that $(\Gamma(x))_{x \in D}$ is a GFF in $D$ with boundary condition $f$ on $\partial D$ and that $O$ is some given subset of $D$. Then, the conditional law of 
$(\Gamma (x))_{x \in O}$ given $(\Gamma(x))_{x \notin O}$ will be a GFF in $O$ with boundary conditions given by the (random) function $f_O$ on $\partial O$ that is equal to the observed values of $\Gamma$ on $\partial O$. 
We can rephrase this in a form that will be reminiscent of the simple Markov property of random walks, except that one replaces the time-set $[0, t]$ by the subset $O$ of $D$: 
\begin {proposition}[Markov property, version 1]
The conditional law of $(\Gamma (x))_{x \in O}$ given that $(\Gamma(x))_{x \notin O}$ is equal to $ (f(x))_{x \notin O}$ is that of a  GFF in $O$ with boundary condition $f|_{\partial O}$.  
\end {proposition}
From this we see why it is so natural to consider the GFF with non-zero boundary conditions. 

\begin{figure}[h]
	\centering
	\includegraphics[scale=0.5]{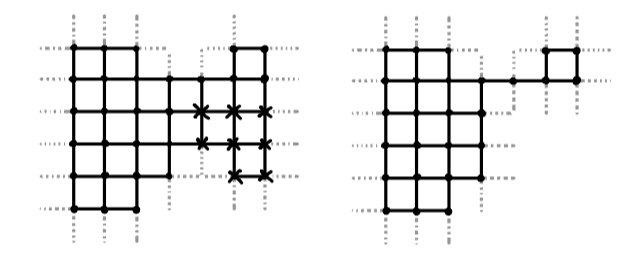}
	\caption{The left-hand side is an example of $D\subset \Z^2$ and $O\subset D$, where the vertices of $D\setminus O$ are marked with a cross, and the vertices of $O$ are marked with a disc. The edges of $\Z^2$ joining two points in $D$ are represented by solid lines, and the edges with one endpoint in $D$ and one endpoint in $\partial D$ are represented by dotted lines. The right-hand side illustrates $O$, where here solid lines are edges joining two vertices in $O$ and dotted lines are edges with one endpoint in $O$ and one endpoint in $\partial O$. The Markov property says that if $\Gamma$ is a GFF on the left graph, and we are given the values of $\Gamma$ ``on the crosses'', then $\Gamma$ restricted to the right graph has the law of a GFF in that graph with non-zero boundary conditions.   }
\end{figure} 
\begin {reminder}\label{reminder_ip}
Let us also  recall the following very elementary fact: when $F_1$ and $F_2$ are two real-valued functions defined on $\Z^d$ and with finite support, then if we define 
$$
(F_1,F_2)= \frac{1}{2} \times \frac{1}{2d}\times \sum_{x \in \Z^d} \sum_{y \in \Z^d, y\sim x} (F_1(y)-F_1(x))(F_2(y)-F_2(x) )
$$ 
we have
\begin {eqnarray*}
 \lefteqn {
 (F_1, F_2) = \frac{1}{2d}\times \sum_{x\in \Z^d} \sum_{y : y\sim x}  \bigl[ - F_1(x) ( F_2 (y) - F_2(x)) \bigr] }\\
 &&
= -\sum_{x\in  \Z^d} F_1(x)\Delta F_2(x)=-\sum_{x\in  \Z^d} F_2(x)\Delta F_1(x), 
\end {eqnarray*}
where we have deduced the last equality by symmetry.  

In particular if for some $B\subset D$, $F_1$ is equal to $0$ outside of $B$ and $F_2$ is harmonic in $B$ (meaning that $\Delta F_2 (x)=0$ for all $x\in B$), then the product $F_1 (x) \Delta F_2 (x)$ is zero everywhere, so that 
$(F_1, F_2)= 0$ and  
\begin {equation} 
(F_1, F_1) + (F_2, F_2) = ( F_1 + F_2, F_1+ F_2).
\label {remark}
\end {equation}
\end {reminder}

We will also use the following definition: 
when $f$ is a real-valued function defined on $\partial D$, we define {\em the harmonic extension $F$ of $f$ to $D$} to be the unique function defined in  $D \cup \partial D$ such that 
$F=f$ on $\partial D$ and $\Delta F  = 0$ in $D$. 

\begin {proposition} 
If $(\Gamma (x))_{x \in D}$ is a GFF with Dirichlet boundary conditions in $D$, and if $F$ is the harmonic extension to $D$ of some given function $f$ on $\partial D$, then 
$(\Gamma (x) + F(x))_{x \in D}$ is a GFF in $D$ with boundary condition $f$ on $\partial D$.
\end {proposition}

Equivalently, one can of course restate this as: 

\begin {proposition}[Markov property, version 2]\label{prop:MPv2}
If $(\Gamma (x))_{x \in D}$ is a GFF in $D$ with boundary conditions $f$ on  $\partial D$, and if $F$ is the harmonic extension to $D$ of $f$, then 
$(\Gamma (x) - F(x))_{x \in D}$ is a GFF in $D$ with Dirichlet boundary conditions. 
\end {proposition}
Hence, the Gaussian vector $(\Gamma (x))_{x \in D}$ is characterised by its 
expectation $(F(x))_{x \in D}$ and its covariance function $\Sigma (x,y) = G_D (x, y)$. 
The effect of the non-zero boundary conditions is only to tilt the expectation of the GFF, but it does not change its covariance structure. 

\begin {proof}[Proof of Proposition \ref{prop:MPv2}] 
The proof is an immediate consequence of the equation (\ref {remark}). Let us consider a GFF $ \Gamma$ in $D$ with Dirichlet boundary conditions, and let 
$F$ be the harmonic extension of $f$ to $D$. Then if we define $\tilde \Gamma = F +  \Gamma$, by a simple change of variables, $\tilde{\Gamma}$ will have a density at $(\gamma_x)_{x \in D}$ which 
is a multiple of 
$$ \exp (  - (  \gamma - F  ,  \gamma - F ) ) ,   $$ 
with the convention that $\gamma=f$ on $\partial D$. This (given that $F$ is deterministic, and using (\ref {remark})) 
is a multiple of 
$$ \exp (  - (\gamma, \gamma ) ) = \exp ( - \frac {1}{2} 
\times  \frac {{\mathcal E}_D (\gamma)}{2d} )   $$ 
(using the same convention on $\gamma$), so that $\Gamma$ is indeed a GFF in $D$ with boundary conditions $f$ on $\partial D$.  
\end {proof}

Let us now introduce some notation that we will be using quite a lot. 
Suppose that $\Gamma$ is a GFF in a finite subset $D$ of $\Z^d$ with boundary conditions given by some real-valued function $f$ on $\partial D$. 
Suppose that $B$ is some finite subset of $D$. We define $O= O(B) := D \setminus B$ and then define the following two new processes: 

\begin {definition} 
\emph{(The processes $\Gamma_B$ and $\Gamma^B$)}
\begin {itemize}
\item
$(\Gamma_B(x))_{x \in D}$ is the process that is equal to $\Gamma$ in $B$ and in $O(B)$, it is defined to be the harmonic extension to $O$ of the values of $\Gamma$ on $\partial O$. 
So the process $\Gamma_B$ can be constructed in a deterministic way given $f$ and the values of $\Gamma$ on $B$.
\item 
The process $(\Gamma^B(x))_{x\in D}$ is then defined to be equal to $\Gamma - \Gamma_B$. Clearly, $\Gamma^B (x) = 0$ as soon as $x \notin O$, and $\Gamma_B + \Gamma^B = \Gamma$. 
\end {itemize}
\end {definition} 

Combining our previous observations readily implies the following alternative statement of the Markov property:
\begin {proposition}[Markov property, version 3] 
The processes $\Gamma_B$ and $\Gamma^B$ are independent, and $\Gamma^B$ is a GFF in $O=D \setminus B$ with Dirichlet boundary conditions. 
\end {proposition}
One main feature in the statement above is the independence of $\Gamma^B$ from $\Gamma_B$, i.e., that fact that $\Gamma^B$ does not depend on the values of $\Gamma$ in $B$. 
Another equivalent way to reformulate this result is therefore that
conditionally on $(\Gamma (x))_{x \in B}$, the conditional law of $(\Gamma (x))_{x \in D \setminus B}$ is that of a GFF in $D \setminus B$ with boundary conditions given by the values of $\Gamma$ on $\partial (D \setminus B)$.  

Note that the special case where $D \setminus B$ is a singleton point $\{x\}$ is exactly the resampling property of the GFF that we mentioned earlier: the conditional law of the GFF at $x$ given its values at all 
other points is equal to a Gaussian random variable with variance $1$ and mean given by the mean value of the GFF at the neighbours of $x$. 

\begin {remark}
\label {GD-GO} 
Since $\Gamma_B$ and $\Gamma^B$ are independent, and since we know that the covariance functions of $\Gamma$ and $\Gamma^B$ are $G_D$ and $G_O$ respectively, we get that 
$$ G_D (x, y) = E [ \Gamma (x) \Gamma (y) ] = E [ \Gamma_B (x) \Gamma_B (y) ] + E [ \Gamma^B (x) \Gamma^B (y) ] = E [ \Gamma_B (x) \Gamma_B(y) ] + G_O (x, y),$$
so that the covariance function of $\Gamma_B$ is 
$$ E[ \Gamma_B (x) \Gamma_B (y) ] = G_D (x,y) - G_O (x,y) $$
for all $x, y $ in $D$.
\end {remark}

\subsection {Deterministic and algorithmic discoveries of the GFF} 
\label {sec:alg_discoveries}
Suppose that $\Gamma$ is a GFF in $D$ with Dirichlet boundary conditions.
We are going to iteratively apply the Markov property described in the previous section in order to discover the values of the GFF in $D$ one by one. More precisely,  
suppose that $D= \{ x_1, \ldots, x_n \}$, and for each $j$, define $B_j = \{ x_1, \ldots, x_j \}$ and 
$O_j = \{ x_{j+1}, \ldots, x_n \}$. The discovery then proceeds as follows:
\begin {itemize} 
\item 
We first discover $\Gamma (x_1)$. This is a centred Gaussian random variable with variance $G_D (x_1, x_1)$. We can therefore write it as 
$N_1 \times \sqrt { G_D (x_1, x_1) }$ where $N_1$ is a centred Gaussian variable with variance $1$. Note that $\Gamma^{B_1}$ is a GFF 
in $O_1$ that is independent of $\Gamma (x_1)$. 
\item 
We then discover $\Gamma^{B_1} (x_2)$. Given that we already know $\Gamma (x_1)$ and therefore the function $\Gamma_{B_1}$, we can then recover $\Gamma (x_2) = \Gamma^{B_1} (x_2) + \Gamma_{B_1} (x_2)$. 
Since $\Gamma^{B_1} (x_2)$ 
is a centred Gaussian random variable with variance $G_{O_1} (x_2, x_2)$, we can write it as $N_2 \times \sqrt {G_{O_1} (x_2, x_2)}$. The Markov property ensures 
that $N_2$ and $N_1$ are independent. Note that at this point we know $\Gamma (x_1)$ and $\Gamma (x_2)$, and can therefore determine the whole function $\Gamma_{B_2}$. 
\item 
We then discover $\Gamma^{B_2} (x_3)$, which allows us to recover $\Gamma (x_3) = \Gamma^{B_2} (x_3) + \Gamma_{B_2} (x_3)$, and continue iteratively. 
\end {itemize}
In this way, we discover $n$ independent identically distributed centred Gaussian random variables $N_1, \ldots, N_n$, and these 
$n$ variables fully describe the GFF $\Gamma$.
 
\begin{exercise} Conclude that we can write 
$$ \Gamma ( \cdot) = \sum_{j=1}^n N_j  \times \sqrt { G_{O_{j-1}} ( x_j, x_j) } \times v_j ( \cdot)$$
for some functions $(v_j)_{1\le j \le n}$. Describe explicitly the form of these functions.
\end{exercise}

In this way, we have constructed the 
$n$-dimensional Gaussian vector $\Gamma$ as a linear combination of $n$ independent Gaussian variables 
(which we can of course always do for Gaussian vectors -- there is nothing special happening here, see Exercise \ref{ex:stand_gauss} below). Notice that if we had chosen another exploration order for $D$, then 
we would have obtained a different decomposition of $\Gamma$ (in fact, corresponding to a different choice of orthonormal basis for the bilinear form $(\cdot, \cdot)$ from Reminder \ref{reminder_ip}) . 
So, in a way, the iterative discovery of the GFF that we just described corresponds 
to the usual way to find an orthogonal basis for a positive definite bilinear form. 

In fact, there is an interesting probabilistic variant that is worth highlighting here. It is actually possible to use some other kind of algorithm in the above exploration, that will make us discover 
the points of $D$ in a random order. We will not give an abstract definition here of what such algorithmic discoveries are, but we will rather 
illustrate it with concrete examples.
For instance, suppose that as before $x_1, \ldots, x_n$ is some deterministic labelling of the $n$ points of $D$. We could instead discover the GFF at these $n$ points 
in an order $\tilde x_1, \ldots, \tilde x_n$ described as follows.  
After having discovered $\Gamma (x_1)$, we know that the conditional law of 
$\Gamma^{B_1}$ is that of a GFF in $O_1$ (and that this process is in fact independent of $\Gamma (x_1)$). So, if we would then like to discover the GFF $\Gamma^{B_1}$, we could actually 
use information that was revealed when we discovered $\Gamma (x_1)$ to decide on an ordering of the points in $O_1$. For example, we could choose the point $\tilde x_2$, depending on the sign of $\Gamma (x_1)$: for instance, by 
deciding that $\tilde x_2$ is $x_2$ if $\Gamma (x_1)$ is positive, and that $\tilde x_2 = x_3$ otherwise. We could then choose $\tilde x_3$ to be $x_4$ if $\Gamma (x_1) + \Gamma (\tilde x_2)  \in [0, 1]$ and 
$\tilde x_3= x_5$ otherwise, and so on. 
Moreover, we are clearly allowed to use additional randomness (that is not generated by $\Gamma$) in our exploration mechanism. For instance, we could have chosen $\tilde x_1$ uniformly at random in $D$. 
In all such explorations, a simple iteration argument shows that for all $j < n$, if we define the random sets 
$$\tilde B_j := \{  \tilde x_1, \ldots, \tilde x_j \} 
\hbox { and } \tilde O_j = D \setminus \tilde B_j,$$ then the conditional law of $\Gamma$ restricted to $\tilde O_j$, given $\tilde B_j$ and 
the values of $\Gamma$ on $ \tilde B_j$, is the law of a GFF in $\tilde O_j$ with boundary conditions given by the values of $\Gamma$ on $\partial \tilde O_j$.

Finally we observe that for each $j$, the set $\tilde B_j$ can take only finitely many (or countably many if $D$ is infinite) values. Hence the previous statement can be rephrased as follows: 
for any given finite $B$ with $j$ elements, the GFF $\Gamma^B$ is independent of the filtration generated by 
the event $\{ B= \{ \tilde x_1, \ldots, \tilde x_j \} \}$ and $\Gamma_B$. 

\begin{exercise}\label{ex:stand_gauss}
	Suppose that $V$ is a finite dimensional real vector space equipped with a positive definite inner product $(\cdot, \cdot)$. Let $\mu$ be the law of a random variable, whose density with respect to Lebesgue measure $dv$ on $V$ is proportional to $e^{-(v,v)/2}$. Show that for any deterministic orthonormal basis $(f_1,\cdots, f_n)$ of $V$ with respect to $(\cdot, \cdot)$, if $(\alpha_1,\cdots, \alpha_n)$ are i.i.d $\mathcal{N}(0,1)$ random variables, then 
	\begin{equation}\label{eq::fd} \sum_{i=1}^n \alpha_i f_i \end{equation} 
	has law $\mu$. Show that $\mu$ is the unique law such that if $X\sim \mu$ then $(X,v)\sim \mathcal{N}(0,(v,v))$ for any fixed $v\in V$. \end{exercise}

\begin{exercise} \label{ex:disc_evalues} Consider the subset $\Lambda_N=[1,N-1]\times [1,N-1]$ of $\Z^2$ for $N\in \N$. Show that for suitable $(m_1, m_2) \in \N^2$ \[ \psi_{m_1,m_2}(x_1,x_2)=\sin(\frac{\pi}{N}x_1 m_1)\sin(\frac{\pi}{N}x_2 m_2)\] is an eigenvector of $\Delta_{\Lambda_N}$, and determine its eigenvalue. Use this to write an expression for a Gaussian free field in $\Lambda_N$ with Dirichlet boundary conditions, as a sum of the form \eqref{eq::fd}, where the $f_i$'s are multiples of an appropriate collection of the $\psi_{m_1, m_2}$'s. For a challenge: use this to show that $G_{D_N}((N/2,N/2),(N/2,N/2)) \asymp \log N$ as $N\to \infty$. 
\end{exercise} 

\begin{exercise}[The classical infinite dimensional example: Brownian motion.] Consider the space $L^2[0,1]$ of square integrable functions from $[0,1]$ to $\R$ equipped with the usual inner product $(f,g)=\int f(x) g(x) \, dx$. Suppose that $(f_i;\; i\ge 1)$ are an orthonormal basis of $L^2([0,1])$ and that we have an infinite sequence $(\alpha_i;\; i\ge 1)$ of independent $\mathcal{N}(0,1)$ random variables defined on some probability space $(\Omega, \mathcal{F}, P)$. Show that 
	\[ W^{(n)}(\cdot)= \sum_{i=1}^n \alpha_i(\I_{[0,\cdot]},f_i)  \] converges (as an element of $L^2[0,1]$) in $\mathcal{L}^2(P)$ to a random variable $W(\cdot)$ and that $W$ is a centred Gaussian process with $E[W(s)W(t)]=s\wedge t$ for every $s,t\in [0,1]$.
	In other words, $W$ is a Brownian motion on $[0,1]$ and we have the decomposition 
	$W(\cdot)= \sum_{i=1}^\infty \alpha_i(\I_{[0,\cdot]},f_i)$.
\end{exercise}

\subsection{Local sets of the GFF} 

Inspired by the previous examples of algorithmic discoveries of a GFF, we are now going to introduce a more abstract class of random subsets of $D$ that are coupled 
with the GFF, and for which one can generalise the simple Markov property.  In other words, we will
define a class of random sets that are the GFF analogue of stopping times for random walks.

Suppose that $D$ is a finite fixed subset of $\Z^d$ and that $\Gamma$ is a GFF in $D$. 
We will use the notation $B$ for {\em deterministic} subsets of $D$, and continue to write
$\Gamma_B$ and $\Gamma^B$ as before. Recall that the simple Markov property of the GFF states that for any deterministic $B$, $\Gamma^B$ is a GFF in $D\setminus B$ that is independent of $\Gamma_B$.

\begin {definition}[Local sets]\label{def:local_set}
When a random set $A \subset D$ is defined on the same probability space as a discrete GFF $\Gamma$ on $D$, we say that 
{\em the coupling $(A, \Gamma)$ is local} if for all fixed $B \subset D$, the GFF $\Gamma^{B}$ in $D \setminus B$
is independent of the $\sigma$-field generated by $(\Gamma_{B}, \{ A= B \})$. 
\end {definition}

Note that this is a property of the joint distribution of $(A, \Gamma)$. 
Sometimes, this property is referred to by saying that ``$A$ is a local set of the free field $\Gamma$'' but we would like to stress that 
this definition does not imply that $A$ is a deterministic function of $\Gamma$; the $\sigma$-algebra on which the coupling is defined can be larger than $\sigma(\Gamma)$. For instance, if $A$ is a random set that is independent of $\Gamma$, 
then the coupling $(A, \Gamma)$ is clearly local. 

A simple criteria implying that a random set is local is the following. 
\begin {proposition} 
If for all fixed $B \subset D$, the event $\{ A =B \}$ is measurable with respect to the field generated by $\Gamma_B$, then $A$ is local.  
\end {proposition}
\begin {proof}
Indeed, if the criteria is satisfied, then the $\sigma$-field generated by  $(\Gamma_{B}, \{ A= B \})$ is just the $\sigma$-field generated by $\Gamma_B$, which is 
independent of that generated by $\Gamma^B$ by the simple Markov property. 
\end {proof}

\medbreak
Here are two instructive examples that we can keep in mind for later on:
\begin{enumerate}
	\item Let $D= \{1, 2, \ldots, n \} \subset \Z$, so that the GFF with Dirichlet boundary conditions on $D$ can be viewed as a Gaussian random walk conditioned to be back at the origin at time $n+1$. Suppose that  $x \in \{ 1, \ldots, n\}$ is chosen uniformly at random and independently of $\Gamma$. 
	Then, let $$y_+ := \max \{ y \ge  x\ : \  \Gamma (z) \times \Gamma (x)  > 0 \;\; \forall x\le z \le y\};$$
	$$y_- := \min \{ y \le  x \ : \ \Gamma (z)  \times  \Gamma (x) > 0 \;\; \forall y\le z \le x \}.$$
	Roughly speaking, the set $A=[y_- -1, y_+ +1 ]$ can be interpreted as an excursion of $\Gamma$ above or below $0$. 
	It is then a simple exercise to check that $A$ is a local set (using a mild variation of the criteria above).  
	
	\item We can do exactly the same when $D \subset \Z^d$ for $d >1$: First choose $x$ at random independently of $\Gamma$, and let $E$ be the connected component containing $x$ of the set of points $y$ in $D$ such that $ \Gamma (y)  \Gamma (x) > 0$. Then $A:= \overline E \cap D$ will be a local set of $\Gamma$ (note that, on the other hand, $E\cap D$ is not a local set, unless the connected components of $D$ are singletons). 
\end{enumerate}

Let us also remark that if $(A, \Gamma)$ is a local coupling, then for any $B \subset B'$ (since one can decompose $\Gamma^{B}$ further into 
$\Gamma^B = ( \Gamma^B)^{B' \setminus B} + (\Gamma^B)_{B' \setminus B}$ so that $(\Gamma^B)^{B' \setminus B} = \Gamma^{B'}$), the GFF $\Gamma^{B'}$ is independent of $(\Gamma_{B'}, 1_{A = B })$. 
In particular, the GFF $\Gamma^{B'}$ is independent of the $\sigma$-field generated by 
$\Gamma_{B'}$ and the event  $ \{ A \subset B' \}$. 
We will use this fact in the proof of the following lemma:

\begin {lemma}\label{disc_unionls}
Suppose that $(A_1, \Gamma)$ and $(A_2, \Gamma)$ are two local couplings (with the same GFF and on the same probability space) such that conditionally on $\Gamma$, the sets $A_1$ and $A_2$ are independent. 
Then, $(A_1 \cup A_2, \Gamma)$ is a local coupling. 
\end {lemma}

It is worthwhile stressing the fact that the conditional independence assumption cannot be dispensed with. Consider for instance the case where $d=1$, $D= \{ -1, 0, 1 \}$ and where $\xi$ is a random variable independent of $\Gamma$ with $\P (\xi = 1) = \P (\xi = -1) = 1/2$. 
Then we define $A_1 = \{ \xi \}$ and $A_2 = \{ \xi \times \sgn (\Gamma (0) ) \}$. Clearly, $A_1$ is independent of $\Gamma$, and $A_2$ is independent of $\Gamma$, so that $(A_1, \Gamma)$ and $(A_2, \Gamma)$ are both local couplings. Yet, $(A_1 \cup A_2, \Gamma )$ is not a local coupling (because $\Gamma (0)$ is positive as soon as $A_1 \cup A_2$ has only one element). 

\begin{proof}
	Let $U$ and $V$ denote measurable sets of $\R^D$. Then, writing $B=B_1\cup B_2$ for any $B_1$ and $B_2$ (again omitting reference to $D$ in the following to simplify notation),
	\begin{align*}
	\P &\Bigl[{\Gamma^{B} \in U}, \ {\Gamma_{B} \in V}, \ A_1 = B_1, \ A_2 = B_2 \Bigr]\\
	&=  \E \Bigl[ \P ( {\Gamma^{B} \in U}, \ {\Gamma_{B} \in V}, \ A_1 = B_1, \ A_2 = B_2 \mid \Gamma )  \Bigr] \\
	&=   \E \Bigl[  \1{{\Gamma^{B} \in U},  {\Gamma_{B} \in V}} \P(  A_1 = B_1, \ A_2 = B_2 \mid \Gamma )  \Bigr] \\
	&=  \E \Bigl[  \1{{\Gamma^{B} \in U},  {\Gamma_{B} \in V}} \P(  A_1 = B_1 \mid \Gamma) P( A_2 = B_2 \mid \Gamma )  \Bigr]
	\end{align*}
	where the last line follows from the assumption of conditional independence. However we know that $\Gamma^B$ is independent of $( \Gamma_B, 1_{ A_1 = B_1 })$ (since $B_1 \subset B$), from which it follows that 
	$$ \P(  A_1 = B_1 \mid \Gamma) = \P ( A_1 = B_1 \mid \Gamma_B )$$
	is a measurable function of $\Gamma_B$, and that the same is true for $\P ( A_2 = B_2 \mid \Gamma)$. 
	Hence, since $\Gamma_B$ and $\Gamma^B$ are independent, we have 
	\begin {eqnarray*}
	\lefteqn { \P \Bigl[ {\Gamma^{B} \in U}, \ {\Gamma_{B} \in V}, \ A_1 = B_1, \ A_2 = B_2 \Bigr] } \\
	& =& \P ( \Gamma^{B} \in U) \times  \P \Bigl[ {\Gamma_{B} \in V}, \ A_1 = B_1, \ A_2 = B_2 \Bigr]  
	\end {eqnarray*}
	If we now fix $B$ and sum over all $B_1$ and $B_2$ such that $B_1  \cup B_2 =B$, we conclude that 
	$$ { \P \Bigl[ {\Gamma^{B} \in U}, \ {\Gamma_{B} \in V}, \ A_1 \cup A_2 = B \Bigr] } 
	= \P ( \Gamma^{B} \in U) \times  \P \Bigl[ {\Gamma_{B} \in V}, \ A_1 \cup A_2 = B \Bigr].
	$$
	This is sufficient to deduce that $\Gamma^B$ is independent of the $\sigma$-algebra generated by $\Gamma_B$ and by the event $\{ A_1 \cup A_2  = B \}$ (because 
	this $\sigma$-algebra is generated by the family of events of the type $\{ \Gamma_B \in V , A_1 \cup A_2  = B \}$ which is a family that is stable under finite 
	intersections). Hence,  $(A_1\cup A_2, \Gamma)$ is a local coupling.  
\end{proof}

\begin {remark} 
The following simple example shows that not all local sets can be discovered in an algorithmic way. Consider $D = \{ 1, 3, 5 \} \subset \Z$. The GFF in $D$ therefore consists of 
three independent centred Gaussian random variables $\Gamma (1), \Gamma( 3)$ and $ \Gamma (5)$ with variance $1$. We denote their respective signs by $\sigma(1)$, $\sigma(3)$  and $\sigma(5)$. 
We will use some extra randomness to choose our random set $A$: 
\begin {itemize} 
\item When $\sigma(1)= \sigma (3) = \sigma (5)$, we choose $A = \{ 1, 3, 5 \} $. 
\item When $\sigma (i_1) = \sigma (i_2) \not= \sigma (i_3)$ for $\{ i_1, i_2, i_3 \} = \{1, 3, 5 \}$, 
we choose $A = \{ i_1, i_3 \}$ with probability $1/2$ and  $A = \{ i_2, i_3 \}$ with probability $1/2$.
\end {itemize}
It is easy to see that $A$ is indeed a local set: the only case to check in Definition \ref{def:local_set} is when $B$ is a two-point set, and then given that $A=B$ and given $\Gamma_B$, we see that the conditional distribution of the sign of the third point must be symmetric, so that the conditional distribution of the GFF at this point ($=\Gamma^B$) is still a centred Gaussian with variance $1$ . 
It is also clear that $A$ must have at least two elements, and that with probability $3/4$, it consists of two elements at which the GFF has opposite signs. 
On the other hand, for any set obtained by an algorithmic exploration as in Section \ref{sec:alg_discoveries} (with at least two elements), the probability that the second revealed value of the GFF has the same sign as the first one is always $1/2$. Thus $A$ cannot possibly be obtained in such a way. 

We remark, however, that this example  of a ``non-algorithmic'' local set is not really something inherently related to the GFF (since it is actually based on a percolation type model with i.i.d. inputs). 
\end {remark}

\section {Determinant of the Laplacian}

We are now going to give various equivalent definitions of an important quantity: the determinant of the Laplacian. 
Recall that when $D \subset \Z^d$ is finite with $n$ elements, we can view the Laplacian as a bijective 
 linear operator from ${\mathcal F}_{(D)}$ into itself. We will denote this operator by $\Delta_D$, as before.
If we write $D = \{ x_1, \ldots, x_n \} $, one can represent $- \Delta_D$ as an $n \times n$ symmetric matrix $(-\Delta_D (x_i, x_j))_{i,j \le n}$, with only $1$'s on the diagonal,
and  off-diagonal terms equal to $0$ or $-1/(2d)$. One can therefore define its determinant, which is a non-zero real number.  Note the sum of the values of $-\Delta_D$ on a line (corresponding to the vertex $x$) can be either $0$ (if all the neighbours of $x$ are in $D$) or 
positive (if at least one neighbour of $x$ is in $\partial D$). We can also note that the matrix $-(2d) \Delta_D$ is integer-valued, so that $(2d)^n \det(- \Delta_D)$ is necessarily an integer (we will see in the next chapter that this integer is actually the number of spanning trees that one can draw in $D$ with wired boundary conditions on $\partial D$). 

The Green's function $G_D$ is a symmetric function defined on $D \times D$, so that it can be also written as a square symmetric matrix $(G_D (x_i, x_j))_{i,j \le n}$. This matrix 
is the inverse matrix of $- \Delta_D$, 
because for all $x$ and $y$ in $D$, 
$$ \sum_{z \in D} \Delta_D (x,z) G_D (z, y) = (\Delta_D \Sigma_y) (x)  =  -\1{x=y}.$$ 
Hence, we have in particular that the determinants of $-\Delta_D$ and $G_D$ are not equal to $0$ and satisfy 
$$ \det G_D  = 1 / \det ( - \Delta_D).$$ 
The matrix $G_D$ is that of a positive definite bilinear form because for all $\lambda_1, \ldots, \lambda_n$, 
$$ \sum_{i,j} \lambda_i \lambda_j G_D (x_i, x_j) = E [ ( \sum_i \lambda_i \Gamma (x_i))^2 ] \ge 0,$$ (i.e., because $G_D$ is a covariance function). This means that its determinant is necessarily positive, and so the determinant of $- \Delta_D$ is therefore positive as well. Of course, one could have seen this from properties of the matrix $\Delta_D$ directly. 

Now let us recall some simple facts about Gaussian vectors. 
\begin {reminder} 
The classical relationship between the density and the covariance function of a centred Gaussian vector is as follows. 
\begin {itemize}
\item 
When $X$ is a centred Gaussian vector $(X_1, \ldots, X_n)$ with non-degenerate covariance 
matrix $\Sigma = (\Sigma_{i,j})_{i,j \le n}$, then its density on $\R^n$ can be written as 
$$ \frac { 1} {(2 \pi )^{n/2} \sqrt { \det \Sigma } } \exp \left\{ - \frac 1 2 \times \sum_{i,j} \gamma_i \gamma_j \Sigma^{-1}_{i,j} \right\} d\gamma_1 \ldots d\gamma_n, $$
where $\Sigma^{-1}$ is the inverse matrix of $\Sigma$. 
\item
Conversely, when $X$ is a centred Gaussian vector $(X_1, \ldots, X_n)$ with density of the form  
$$  C \exp ( - \frac 1 2 \times \sum_{i,j} \gamma_i \gamma_j (-\Delta_{i,j}) ) d\gamma_1 \ldots d\gamma_n, $$ 
where $(\gamma_j) \mapsto - \sum_{i,j \le n} \Delta_{i,j} \gamma_i \gamma_j$ is a positive definite bilinear form, 
then the covariance matrix of $X$ is $\Sigma := - \Delta^{-1}$, and the coefficient $C$ satisfies
$$  C = \frac 1 { (2 \pi )^{n/2} \sqrt { \det \Sigma  }} =   \frac {\sqrt {\det (-\Delta)}} { (2 \pi )^{n/2}} .$$  
\end {itemize}
\end {reminder}

Applying this to our GFF set-up could have provided us a more direct (but maybe less instructive than the resampling route we chose) 
way to see that the covariance function of the GFF is given by the Green's function. 

Now, we see that when $D$ has $n$ elements, the density of the GFF in $D$ is exactly 
$$ \frac { \sqrt { \det (- \Delta_D) } }  {(2 \pi )^{n/2}  }\exp \left\{ - \frac{{\mathcal E}_D (\gamma)}{2 \times (2d)} \right\} d\gamma_1 \ldots d\gamma_n,$$ 
which provides the first following intuitive interpretation for the quantity $\det ( - \Delta_D )$: it somehow measures how ``constrained'' the springs are 
by the condition that they are chained together, compared to if they were independent and identically distributed. 
Another interpretation is the following:  
\begin {proposition} 
The quantity $\sqrt {\det ( - \Delta_D)} / ({2 \pi })^{n/2}$ is the density of the GFF distribution at the point $(0, \ldots, 0)$. 
\end {proposition} 
In other words, the quantity $\sqrt {\det G_D}$ describes how costly it is to ask the GFF to be very small everywhere: 
$$ \lim_{\eps \to 0}  \eps^{-n} P \Bigl[ \forall i \le n, \ | \Gamma(x_i) | \le \eps \sqrt { \pi / 2 } \Bigr] =  1 / \sqrt { \det G_D}. $$  

Let us now combine this with the explicit decomposition of the GFF in $D=\{x_1,\ldots , x_n \}$, where one first discovers $\Gamma (x_1)$ and is then left to discover the 
GFF in $D \setminus \{ x_1 \}$ etc. On the one hand, 
since $\Gamma (x_1)$ is a centred Gaussian random variable with variance 
$G_D (x_1, x_1)$, we know that as $\eps \to 0$,
$$ P \Bigl[ | \Gamma(x_1) | \le \eps \sqrt { \pi / 2 } \Bigr] \sim \eps / \sqrt { G_D (x_1, x_1 )}.$$
On the other hand, since 
$\Gamma^{\{ x_1 \}}$ is independent of $\Gamma (x_1)$ (together with the fact that $\Gamma = \Gamma_{\{x_1\}} + \Gamma^{\{ x_1 \} }$ and that 
the density of $\Gamma(x_1)$ is smooth), we 
readily see that as $\eps \to 0$,
$$ P\Bigl[ \forall i \in \{2, \ldots, n \}, \ | \Gamma(x_i) | \le \eps \sqrt { \pi / 2 } \quad   \big|  \,  \ | \Gamma(x_1) | \le \eps \sqrt { \pi / 2 } \Bigr] 
\sim \eps^{n-1} /  \sqrt {\det G_{D \setminus \{ x_1 \}} }.$$ 
Hence, we can conclude that 
$$  { \det G_{D} } = {G_D (x_1, x_1)}  \times { \det G_{D \setminus \{ x_1 \}} },$$
and it then follows by induction that: 
\begin {proposition}
\label {G_D}
\[\det G_D =  \prod_{j=1}^n G_{D \setminus \{ x_1, \ldots , x_{j-1} \} } (x_j, x_j).\]
\end {proposition}
In particular, we observe that the product on the right-hand side does not depend on the ordering $\{x_1,x_2, \ldots, x_n\}$ that we gave to the points of $D$.
This fact will be useful in our description of Wilson's algorithm in the next chapter. 

\begin {remark}
It is easy to check by other simple means 
that this product does not depend on the order of the $\{x_j\}$. For instance,  by proving  the simple identity 
$$
G_{B} ( x, x) G_{B \setminus \{ x \} } (x', x') 
= G_{B} ( x', x') G_{B \setminus \{ x' \} } (x, x) 
$$ 
for all finite sets $B$, and all $x$ and $x'$ in $B$ (this can be viewed as a general property of a Markov chain on a state-space with three elements). 
\end {remark} 

Let us now explain how the previous considerations allow us to provide an expression for the Laplace transform of (the square of) a GFF in terms of determinants.
Suppose that $\Gamma$ is a GFF with Dirichlet boundary conditions in $D = \{ x_1, \ldots, x_n \} \subset \Z^d$ as before, and for all $k := (k(x_1)), \ldots, k(x_n) ) \in (\R_+)^n$, let $I_k$ be the  
diagonal matrix with $I_{i,i} = k(x_i)$ for each $i$. 

\begin {proposition}[Laplace transform of the square of the GFF]
\label {LaplaceGFF}
Suppose that $\Gamma$ is a GFF in $D$ with Dirichlet boundary conditions. Then, for all $k \in (\R_+)^n$, 
$$ E [ \exp ( - \frac 1 2 \sum_{j=1}^n k(x_j) \Gamma(x_j)^2 ) ] =  \sqrt { \frac { \det ( - \Delta_D  )}{\det (- \Delta_D + I_k) }}. $$
\end {proposition}

\begin {proof}
This is a straightforward consequence of the previous considerations: the matrix $- \Delta_D$ is positive definite so that $- U := - \Delta_D + I_k$ is 
positive definite as well (recall that the diagonal terms of $I_k$ are all non-negative). Moreover, we have 
\begin {eqnarray*}
\lefteqn { E [ \exp ( - \frac 1 2  \sum_{j=1}^n {k(x_j) \Gamma(x_j)^2} ) ] }\\
&=& 
\frac {\sqrt {\det (-\Delta_D)}} { (2 \pi )^{n/2}} \times  \int_{\R^n}
\exp ( -  \sum_{j=1}^n \frac {k(x_j) \gamma_j^2} 2 ) \times \exp ( -  ( \sum_{i,j} \frac {\gamma_i \gamma_j}{2} (- \Delta_D) (x_i, x_j))) d\gamma_1 \ldots d\gamma_n \\
&=& 
\frac {\sqrt {\det (-\Delta_D)}} { (2 \pi )^{n/2}} \times  \int_{\R^n}
\exp ( -   ( \sum_{i,j} \frac { \gamma_i \gamma_j}{2} (- U_{i,j}) )) d\gamma_1 \ldots d\gamma_n \\
&=& 
\sqrt { \frac { \det ( - \Delta_D  )}{\det (- \Delta_D + I_k) }}. 
\end {eqnarray*}
\end {proof}

\section {GFF on other graphs} 
\label{GFFonothergraphs}

\subsection {The massive Gaussian Free Field} 
\label {S51}

Let us first describe one particular generalisation of the GFF in $D \subset \Z^d$ that will be useful in the next chapter. 
A more general set-up (including this particular case) will be presented in Section \ref{elnet}.

Just as before, we are given a finite subset $D$ of $\Z^d$, and we define the energy  $ {\mathcal E}_D$ of a function in ${\mathcal F}_{(D)}$ in the same way. 
We are now also given a non-negative function $k=(k(x))_{x\in D}$ on $D$. Given $k$ and $D$, we define the following:

\begin {definition}[Massive GFF] 
The massive GFF in $D$ (with Dirichlet boundary condition and mass function $k$) is the centred Gaussian random vector $(\Gamma (x))_{x \in D}$ with density 
at the point $(\gamma_x)_{x \in D}$ that is proportional to 
$$ \exp \Bigl[ - \frac {1}{2} \times \Big( \frac {  {\mathcal E}_D (\gamma)} {2d}  + \sum_{x \in D} k(x) \gamma_x^2  \Big)  \Bigr] $$ 
with the convention that $\gamma = 0$ on $\partial D$. 
\end {definition}

Heuristically at each site $x$, one adds a little ``vertical'' spring with ``intensity'' $k(x)$ that tries to pull the height of the GFF back to $0$. Note that the proportionality constant in front of 
this density must be equal to $\sqrt {  \det ( - \Delta_D + I_k )} / { (2 \pi )^{n/2}}$ (see Proposition \ref{LaplaceGFF}). Of course if $k\equiv 0$, then the massive GFF is the same as the standard, massless version that we have discussed so far.

\medbreak

In this set up it is natural to consider instead of the Laplacian $\Delta_D$, the operator $U = U_{D,k}$ on ${\mathcal F}_{(D)}$ defined by 
$$ [UF]  (x) = \Delta_D F (x) - k(x) F (x).$$ It then follows from Proposition \ref{LaplaceGFF} that the covariance function 
$\Sigma$ of the massive GFF is given by the inverse matrix of $-U = - \Delta_D + I_k$. 

An alternative way to see this is through the following resampling property (that may be checked by simply inspecting the density function of the massive GFF $\Gamma$): for any $x\in D$, the conditional law of $\Gamma (x)$ given $(\Gamma (y))_{ y \ne x}$ 
is a Gaussian with mean $\overline \Gamma (x) / (1+ k (x))$ and variance $1/(1 +  k (x))$. 
Just as in the case where $k=0$, one can then use this to characterise the covariance function $\Sigma$ of this massive field by the fact that for all $x,y$ in $D$,
\begin {equation}
\label {EqU} 
U\Sigma_x (y) = -\mathbf{1}_{y=x}, 
\end {equation}
where $\Sigma_x ( \cdot ) = \Sigma (x, \cdot )$. This shows that $(-U) \times \Sigma$ is the identity matrix.  

We now explain how, just as for the (non-massive) Green's function, the function $\Sigma (x,y)$ can be interpreted in terms of certain random walks. These are the discrete-time 
and continuous-time random walks $(X_n)_{n \ge 0}$ and $(Y_t)_{t \ge 0}$ with ``killing rate given by $k$''.  As we will see, the relation will be neater for the  continuous-time walk: a feature that will also show up when we will discuss the relation with the GFF itself.

To define the walks with killing, we create an additional ``cemetery'' state $\partial$, and then $X$ and $Y$ are the discrete (resp. continuous) time Markov chains on $\Z^d\cup \{\partial\}$ described as follows.
\begin {itemize}
\item For $X$: at each time step, if $X$ is at $x$, it will jump to 
$\partial$ with probability $k (x) / (1+ k(x))$ (and then stay there forever). Otherwise it will choose one of the $2d$ neighbours of $x$ with equal probability and proceed from there. 
\item
For $Y$: on each edge $e$ of the graph, 
bells ring at a rate $1/(2d)$ (i.e., the gaps between each ring are independent exponential random variables with mean $2d$) 
and at each site $x$, a special bell rings at rate $k(x)$. Then, when $Y$ is at a site $x$,
it stays there until the first time at which either the bell of an adjacent edge rings (and then
$Y$ jumps along that edge and proceeds from there), or the special bell at $x$ rings (and then $Y$ jumps to the cemetery state $\partial$ and stays there forever). 
So, the time spent by $Y$ before jumping away from $x$ is an exponential variable with mean $1/ (1+ k(x))$. Moreover, if $\tau_n$ denotes the $n$-th jumping time of $Y$, then
the discrete chain $(X_n := Y_{\tau_n})$ is distributed as the walk $X$ described above, when both are stopped at their respective hitting times of $\partial$.
\end {itemize} 

We define $\tau$ and $\sigma$ to be the respective first 
times at which $X$ and $Y$ are not in $D$ (i.e. at which they either go to the cemetery state $\partial$ or to a point in $\partial D$). Then, we can define the massive Green's functions as follows. 
\begin {definition}[Massive Green's functions] 
\label {DefMassiveG}
The massive Green's function $G_{D,k,\text{discrete}}$ for the discrete-time random walk $X$, is the function on $D \times D$ defined by
$$
G_{D,k,\text {discrete}} (x,y) := E_x \Bigl[ \sum_{j=0}^{\tau-1} \mathbf{1}_{\{X_j = y\} } \Bigr]. 
$$ 
The massive Green's function $G_{D,k}$ for the continuous-time random walk $Y$ is the function on $D \times D$ defined by
$$
G_{D,k} (x,y) := 
E_x \Bigl[ \int_0^\sigma \mathbf{1}_{\{Y_s = y\} } ds \Bigr]. $$ 
\end {definition}
Of course, $$G_{D,k} (x, y) =\frac 1 {1+ k(y)} G_{D, k, \text {discrete}} (x,y),$$ which implies in particular that  
$$ (\prod_{j=1}^n (1+ k(x_j)) ) \det G_{D,k} = \det G_{D,k,\text{discrete}} .$$
Either directly (as in Remark \ref {Pn}) or using the Markov property (exactly as in the non-massive case), one sees that $(-U_{D,k}) \times G_{D,k}=\text{Id}$ so that $$\Sigma = G_{D, k}.$$

This description of the covariance function also allows us to generalise the definition of the massive GFF to infinite $D$. 
For instance, we can use it when $D$ is $\Z^d$ (even for $d=1, 2$), as long as $k$ 
is not identically $0$ (the case where $D$ is $\Z^d$ and $k\equiv m > 0$ is sometimes simply referred to the GFF with mass $m$ in the literature).

\subsection {GFF on electric networks} 

\label {elnet} 

For simplicity, we have focused so far on GFFs on subsets of $\Z^d$. However, the GFF can be naturally generalised to a broader class of weighted graphs that are 
often referred to as ``electric networks''. Let us now describe them.  

Let $V$ be a finite or countable set of vertices. We equip this vertex set with a function $c$ that assigns to each pair $\{ x , y \}$ of distinct vertices a conductance $c_{x,y} = c_{y,x}$ in $[0, \infty )$ 
(by convention $c_{x,x} =0$ for all $x$). 
We furthermore assume that for all $x \in V$, the quantity 
$
\lambda_x := \sum_{y: y\sim x} c_{y,x}  
$ 
is finite. This pair $(V,c)$ defines 
what is sometimes called an {\em electric network}. By convention, we say that in this electric network, there is an edge between $x$ and $y$ when $c_{x,y} > 0$.
This then defines an edge-set $E$. We 
will assume in the following that the graph $(V,E)$ is connected. 

On such electric networks, it is natural to define a discrete-time
random walk $(X_n)_{n \ge 0}$ in such a way that when it is at $x$, it chooses to visit a point $y$ at the next step with probability 
$c_{x,y} / \lambda_x$. It is also natural to consider the corresponding continuous-time Markov chain $Y$, that when at $x$, jumps along 
an edge $e$ at rate $c_e$. This means that for $Y$, the rate of jumping away from $x$ is $\lambda_x$ (i.e., the waiting time at $x$ before jumping is an exponential random variable with mean $1 / \lambda_x$). 
With this set-up, the measure that assigns the mass $\lambda_x$ to each site $x$ is then a 
reversible invariant measure for $X$ (although it is not necessarily finite if $V$ is infinite), and the  measure that assigns mass $1$ to each point 
of $V$ is a reversible invariant measure for $Y$.  

One example of such an electric network is $V=\Z^d$, with $c_{x,y}$ equal to $1/(2d)$ when $x$ and $y$ are neighbouring points, and $c_{x,y} = 0$ otherwise. In this case $\lambda_x =1$ for every $x\in \Z^d$.

All of the quantities that we will 
define in the coming paragraphs will implicitly depend on the function $c$, even if we omit this dependence in the notations. We suppose that $D$ is a finite subset of 
$V$. 
For a vector $(\gamma_x)_{x\in D}$, one can then define 
\[ \xi_D (\gamma):= \sum_{e \in E} ( c_{e} \times (\nabla \gamma (e))^2) \]
(with $c_e=c_{x,y}$ when $e$ joins $x$ and $y$), and with the convention that $\gamma_y=0$ for all $y \in V \setminus D$. 

\begin {definition}[GFF on electric networks]\label{def:gffelnet}
Assume that $V \setminus D$ is non-empty. We say that 
$(\Gamma (x))_{x\in D}$ is {\em a Gaussian Free Field} in $D$ (for the network $(V,c)$) with Dirichlet boundary conditions on $V \setminus D$ if its density 
is proportional to $e^{-\xi_D(\gamma)/2}$ at $(\gamma_x)_{x\in V}$ (again with the condition that $\gamma=0$ on $V\setminus D$).
\end {definition} 

\begin {remark} 
In order to define this GFF in $D$, it is actually sufficient that $\lambda_x < \infty$ for all $x \in D$ (it does not need to be finite for $x \in V \setminus D$); this comment is of course only relevant when 
$V$ is infinite.  
\end {remark}

Then the covariance function $\Sigma$ of the centred Gaussian process $(\Gamma (x))_{x \in D}$ 
can be described via another variant of the Green's function. To define this,
consider the discrete-time and continuous-time random walks $X$ and $Y$ described above,
write $\tau$ for the first time that $X$ reaches $V \setminus D$, and write $\sigma$ for the first time that $Y$ reaches $V \setminus D$. Then, defining the continuous-time 
Green's function by 
\begin {equation} 
\label {defGelnet}
G_{D,c} (x,y) 
:=\frac{1}{\lambda_y} {E}_x\Bigl[\sum_{n=0}^{\tau-1} \mathbf{1}_{\{X_n=y\}}\Bigr] = {E}_x \Bigl[ \int_0^\sigma \mathbf{1}_{\{Y_s =y\}} ds \Bigr]
\end {equation}
(in this notation, we omit the implicit dependence of $G_{D,c}$ on the larger graph $V$) we can prove -see Exercise \ref{ex:gff_en} below- 
that
$\Sigma =G_{D,c}$.

This last point again makes it possible to extend the definition of such a GFF to the case where $D$ is infinite, 
provided that the corresponding Green's function $G_{D,c}$ is finite. 

Note finally that 
the massive Green's function discussed in Section \ref {S51} is just a particular case of this more general set-up, where $V = \Z^d \cup \{ \partial \}$ (we add the cemetery point 
to our graph), $c_{x , \partial } = k (x)$ and $c_{x,y} = 1 / (2d)$ when $x$ and $y$ are in $\Z^d$.

\begin{exercise}\label{ex:gff_en}
Define the operator $ \Delta_{D,c}$ on functions $f:D \to \R$ by 
\[ \Delta_{D,c} f(x)=\sum_{y: y\sim x} c_{x,y} (f(y) - f(x)) \]
with the convention that $f =0$ outside of $D$. 
Let  $(\Gamma (x))_{x\in D}$ be a GFF in $D$ as in Definition \ref{def:gffelnet}.  
\begin{enumerate} \item Show the following resampling property:  for any $x\in D$, the conditional law of $ \Gamma (x)$ given $( \Gamma (y))_{ y \ne x}$ 
	is a Gaussian with mean $(\sum_{ y : y \sim x } c_{x,y}  \Gamma (y))/ \lambda_x$ and variance $1/\lambda_x$. 
	Deduce that \[ \sum_{y : y\sim x} {c_{x,y}} (\Gamma(x) - \Gamma(y))\] is a centred Gaussian with variance $ \lambda_x$, independent of $( \Gamma (y))_{y\ne x}$.
	\item Use this to prove that if 
	$\Sigma (x,y):=E[\Gamma (x) \Gamma (y)]$ is the covariance function of $ \Gamma$, then for $\Sigma_x(\cdot):= \Sigma (x,\cdot)$ we have $\Delta_{D,c}  \Sigma_x(y)=\1{y=x}$.
	\item Setting $g_x(y)=  G_{D,c}(x,y)$ for $x,y\in D$ show that $\Delta_{D,c}  g_x(y)=-\1{y=x}$. Deduce that $(\Gamma (x))_{x\in D}$ is the unique centred Gaussian process indexed by $D$, with covariance function $E[ \Gamma(x) \Gamma (y)]=G_{D,c}(x,y)$.	
\end{enumerate}
\end{exercise} 

\chapter {Loop-soups and the discrete GFF}
\label {Ch2}

\section {Uniform spanning trees and Wilson's algorithm}

We have already mentioned during our analysis of the determinant of the Laplacian that it was closely related to enumerations of spanning trees. The goal of this section is to describe this relation. 

Suppose that ${\mathcal D}$ is a finite connected graph, with  vertex set ${\mathcal V}$ and edge-set ${\mathcal E}$ (here we will allow the case of ``multiple edges'': where several edges of ${\mathcal E}$ join 
the same pair of points in ${\mathcal V}$). 

\begin {definition} [Spanning trees, uniform spanning trees]  
 A spanning tree in ${\mathcal D}$ is a subset $T$ of ${\mathcal E}$, such that the graph $({\mathcal V}, T)$ is a tree 
 (it does not contain a cycle that uses edges only once), and is connected (``spanning''). 
A uniform spanning tree (UST) in ${\mathcal D}$ is a random tree ${\mathcal T}$ that is chosen uniformly among all spanning trees of ${\mathcal D}$. 
\end {definition}

\subsection {In subsets of $\Z^d$}\label{sec::zd_ust}

In this section, we will study uniform spanning trees in particular graphs $\hat D$ that are defined as follows.
Let us work in the same setting as in the previous chapter ($D$ is a finite subset of $\Z^d$ with $n$ points, $\partial D$ is the set of points that are at distance $1$ from $D$, and $E=E_{\overline{D}}$ denotes the set 
of all edges of $\Z^d$ that have at least one extremity in $D$).
We define $x_0$ to be an abstract point obtained by the formal contraction of all the points in $\partial D$. We can then define $\hat D$ to be the graph with vertex set 
$D \cup \{ x_0 \}$ and with edge set $\hat E$ induced by $E$ (we keep the edges that join two points of $D$, and an edge from $x \in D$ to $y \in \partial D$ becomes an edge from $x$ to $x_0$). Note that it is possible for $x$ and $x_0$ to be joined by more than one edge in $\hat E$. We will denote by $\varphi$ the bijection taking edges in $E$ to their corresponding edges in $\hat E$.  

When $T$ is a spanning tree of $\hat D$, it can also be identified with the graph consisting of the vertices $D \cup \partial D$ and edges of $ \varphi^{-1} (T)$. 
This graph is now not necessarily connected any more (because there are no edges directly joining the various points of $\partial D$) but the set of edges
$\varphi^{-1}(T)$ is often referred to as a spanning tree of $D$ {\em with wired boundary conditions}. If ${\mathcal T}$ is a UST in $\hat D$, we therefore say that $\varphi^{-1} ( {\mathcal T})$
is a UST in $D$ with wired boundary conditions. 

\begin{figure}[h]
\centering
	\includegraphics[scale=0.2]{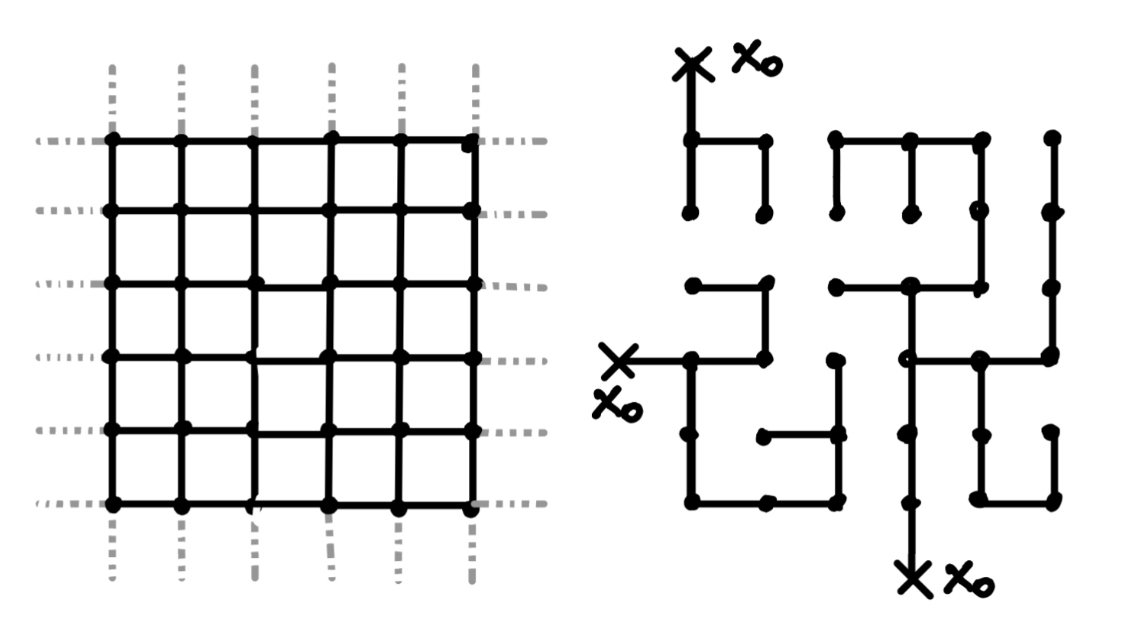}
	\caption{A subset $D\subset \Z^2$ (left) and a spanning tree in $\hat{D}$ (right)}
\end{figure}

Understanding uniform spanning trees in $\hat D$ is of course related to counting the number of spanning trees of $\hat D$. We are now going to describe an explicit procedure 
(known as Wilson's algorithm)
that constructs a random spanning tree ${\mathcal T}$ of $\hat D$, and we will show (even though this is far from obvious at first) that  the law of this tree  is actually uniform among all 
spanning trees. A by-product of the proof will be the following fact (recall that $n$ is the number of points in $D$):    

 \begin {proposition}
 \label {counting}
 The number of spanning trees in $\hat D$ is equal to $(2d)^n \det (-\Delta_D) $. 
 \end {proposition}
 
 Note that $2d \Delta_D$ is an integer-valued matrix, so it is not surprising that $(2d)^n \det (-\Delta_D) $ is an integer. 
\medbreak

Before describing Wilson's algorithm, 
we need to explain the notion of loop-erasure of a path, and the definition 
of loop-erased random walk.
\begin {itemize} 
\item 
Let us first clarify a little terminology issue that will be relevant throughout this chapter. We will often consider nearest-neighbour paths (or loops) in a graph. By this we will, loosely speaking, refer to  
a finite collection of points $Z=(Z_0, \ldots, Z_m)$ in the graph such that for each $j \in \{ 1, \ldots, m\}$, $Z_j$ and $Z_{j-1}$ are neighbours in the graph (and for loops, we will 
also require that $Z_0 = Z_m$). The quantity $m$ will denote the length of the path. However, we will always implicitly assume that the knowledge of such a path also includes the information about which edges $e_1, \ldots, e_m$ were used in the $m$ steps, so that in reality, a path $Z$ should be viewed as a collection 
$(Z_0, e_1, Z_1, \ldots, e_m , Z_m)$ where for each $j \in \{ 1, \ldots , m\}$, $e_j$ is an edge joining $Z_{j-1}$ and $Z_j$. This can be important when one enumerates 
paths because it could happen, for instance, that there are several edges joining $Z_{j-1}$ and $Z_j$, and this would mean that
$(Z_0, \ldots, Z_m)$ corresponds to several different possible paths. In the concrete setting of spanning trees in $\hat{D}$ as above, recall that interior points may be joined to $\{x_0\}$ via multiple edges. This is a situation where we should keep the preceding comment in mind. 
 \item 
We will now introduce the notion of loop-erasure of a path.   
For any path ${Z} = (z_0, \ldots, z_m) \in (\Z^d)^{m+1}$, we define the 
loop-erasure $$L({Z})=(L_0,\cdots, L_\sigma)(Z)$$ of ${Z}$ iteratively as follows:
we let $L_0 = z_0$, and then for each $j \ge 0$, we define  $r_j = \max \{ r \le m  \ : \ z_r = L_j \}$ and 
$
L_{j+1} = z_{1 + r_j} 
$ inductively,
until reaching $\sigma:=\min\{j:L_j=z_m\}$.
In other words, we have
 erased the loops of $Z$ in chronological order. The number of 
steps $\sigma$ of $L$ depends on $Z$; for instance, when $z_m = z_0$, then $\sigma = 0$. Again, the loop-erased path ``keeps track'' of the edges used by $Z$ that have not been erased. 
The edge from $L_j$ to $L_{j+1}$ is the edge from $Z_{r_j}$ to $Z_{r_j +1 }$ in $Z$. 
\item
Suppose that $Z= (Z_r, r \le \tau)$
 is a simple random walk started from $x_1 \in D$ and stopped at its first exit time $\tau$ of $D$. Let $L(Z)$ be its loop-erasure. This is now a nearest-neighbour path joining $x_1$ to a boundary point of $D$.    
Observe that for each such simple nearest neighbour path 
$y=\{y_0,\ldots, y_s\}$  from $x_1$ to $x \in \partial D$, when we decompose the probability that 
$L=y$ according to all possible $Z$'s with 
$L=y$, we have
\begin{eqnarray*}
\lefteqn{P (L = y)} \\ & &= (2d)^{-s}\prod_{i=1}^{s} \sum_{k_i\ge 0}  (2d)^{-k_i} \# \{\text{paths } y_{i-1} \to y_{i-1} \text{ with } k_i \text{ steps in } D\setminus \{y_0,\cdots, y_{i-2}\}  \}\\
& &=  (2d)^{-s} G_D (y_0, y_0) G_{D \setminus \{y_0\} } (y_1, y_1) \ldots G_{D \setminus \{y_0, \ldots 
, y_{s-2}\}} (y_{s-1}, y_{s-1}).
\end{eqnarray*}
In words: the term $(2d)^{-s}$ term corresponds to the jumps of the walk that are still present on the loop-erasure, and the other terms in the product correspond (for each given $j$) to the contributions of all possible paths that the random walk may perform between $r_{j}+1$ and $r_{j+1}$.  
\end {itemize} 

Now that we know how to define the loop-erased random walk from a point to the boundary of a domain, we are ready to describe Wilson's algorithm to construct a random ``tree'' in $D$ with wired boundary conditions: 

\begin{enumerate}
\item We order the $n$ points of $D$ as $x_1, \ldots ,x_n$.
\item We take a simple random walk $(Z_r)_{r\in \N}$ started from $x_1$ and stopped at its first exit of $D$. We consider its loop-erasure $X^{(1)}:=L(Z)$; here $X^{(1)}$ consists of the sites visited by this loop-erasure together with the edges along which the loop-erasure jumps. 
\item Iteratively, for each $ k \in \{ 2, \ldots ,  n \}$, we construct $X^{(k)}$ as follows. If $x_k \in X^{(k-1)}$, then we set $X^{(k)} = 
X^{(k-1)}$. On the other hand, if $x_k \notin X^{(k-1)}$, then we take a simple random walk started from 
$x_k$ and stopped at its first exit of  $D\backslash X^{(k-1)}$ and set $X^{(k)}$ to be the union/concatenation of its loop-erasure with $X^{(k-1)}$.  
\end{enumerate}
In this way, $\varphi (X^{(n)})$ is a tree in $\hat D$, and it contains all points of $D$. We have therefore defined a random 
spanning tree ${\mathcal T}$ of $\hat D$. 

\begin {proposition}[Wilson]
\label {Wilson}
The law of this tree ${\mathcal T}$ is that of a uniform spanning tree of $\hat D$. 
\end {proposition}

\begin {proof} 
If we are given a possible outcome $T$ for the tree $X^{(n)}$, then we re-label the points of $D$ as follows.
We denote by $y_1, \ldots, y_{s-1}$ the simple path (``branch'') in the tree going from $x_1$ to $\partial D$ (where $y_{s-1}$ is the last point in $D$ in this path). Then, we define $y_s$ to be the next $x_j$ in $D$ that is not in this already labelled set, and define $y_s, \ldots, y_{s'-1}$ to be the branch in $T$ that joins $x_j$ to $\partial D \cup \{ y_1 ,\ldots, y_{s-1} \}$. We proceed iteratively. This provides us with an ordering of the vertices of $D$ that is determined by the tree $T$.

Inductively, using the previously calculated probability for a \emph{single} branch, we see that the probability of this given tree $T$ being exactly the one constructed by Wilson's algorithm is
$$ P [ {\mathcal T} = T ] =  (2d)^{-n} \prod_{j=1}^n  G_{D \setminus \{ y_1, \ldots , y_{j-1} \} } (y_j, y_j).$$
However, we have seen in the previous chapter that this quantity is equal to $(2d)^{-n} \det G_D=((2d)^n \det(-\Delta_D))^{-1}$ and does not depend on the order of the points $y_1, \ldots, y_n$. 
It follows readily that the probability above does not depend on $T$ (hence, the algorithm samples a uniformly chosen spanning tree) and that this probability is
 $(2d)^{-n} \det G_D$. This implies both Proposition \ref {Wilson} and Proposition \ref {counting}.  
\end {proof}

As a warm-up to the considerations of Section \ref {SecLS} note that  this proof shows, in particular, that 
the probability of erasing no loop at all while performing Wilson's algorithm
is equal to $1/ \det G_D$, independently of the tree that one constructs. Indeed for any tree $T$,  given that $\mathcal{T}=T$, the probability that it was constructed without erasing any loops is equal to $(2d)^{-n}/((2d)^{-n}\det G_D)=1/\det G_D$.

The following remark will also be very useful: 
\begin {remark} 
\label {resamplingwarmup}
When one performs Wilson's algorithm as above, let us denote by $\lambda$ the (long, concatenated) loop from $x_1$ to $x_1$ that one erases when performing the algorithm. It is distributed 
like $(Z_0, \ldots, Z_\rho)$, where $Z$ is a simple random walk started from $x_1$ and $\rho$ denotes the last time at which it visits $x_1$ before hitting $x_0$ for the first time. The previous considerations 
then show that 
$$ P ( \lambda = (z_0, \ldots, z_r ) ) = (2d)^{-r} / G_D (x_1, x_1) $$
for each possible loop $(z_0, \ldots, z_r)$. By using the Markov property at the $j(\lambda) \ge 0$ successive return times to $x_1$ by $Z$, it is easy to see that the total number $j$ of returns to $x_1$ by $\lambda$  is geometric. The expectation of $j+1$, which is the mean number of visits of $x_1$ by the walk is $G_D( x_1, x_1)$. 
We see that conditionally on $j$, these $j$ excursions are independent and identically distributed (they each follow the law of a random walk started from $x_1$ up to its first return to $x_1$, conditioned 
to return to $x_1$ before hitting $x_0$). 
In particular, we see that reshuffling uniformly at random the order of these $j$ excursions, or reshuffling uniformly at random the order of the last $j-1$ excursions (keeping the first one fixed), 
or resampling the excursions themselves, will not change the law of $\lambda$. 
\end {remark}

\subsection {Some generalisations} 
\label {massivecase}

Let us now briefly mention uniform and weighted spanning trees in general finite graphs. The following remarks generalise our previous statements.

\begin {enumerate}

\item {\emph{The massive case.}}
Before turning to the general cases, let us first explain in the same set-up as Section \ref{sec::zd_ust} (with $D \subset \Z^d$) 
what sort of trees Wilson's algorithm constructs when we replace the simple random walks with ``massive'' ones. In this case, one adds the cemetery point $\partial$ to $\Z^d$, and joins each $x_j$ in $D$ to $\partial$ via an edge with non-negative conductance $k(x_j)$.
When we define $\hat D$, the vertex $x_0$ then corresponds to all sites outside $D$ (including $\partial$). We 
call ${\hat E}_\partial$ the set of edges of $\hat D$ that correspond to an edge from $D$ to $\partial$. 

We can then perform Wilson's algorithm ``rooted at $x_0$'' just as before, except that we now use the massive random walk on $\hat D$: when the walk is at $x_i \in D$, then it jumps to 
the cemetery point $\partial$ with probability $k(x_i) / (1 + k(x_i))$ and otherwise, it chooses uniformly one of the $2d$ neighbours of $x_i$. Then, for any given spanning tree $T$ of $\hat D$, 
we readily obtain that
$$ P [ {\mathcal T} = T ] =  ( \prod_{e \in T} c_e ) \times  \prod_{j=1}^n  G_{D \setminus \{ y_1, \ldots , y_{j-1} \}, k } (y_j, y_j) 
=  ({\prod_{e \in T} c_e }) \times { \det G_{D,k}},  $$ 
where $c_e = k(x_i) $ if the edge $e \in \hat E_\partial$ corresponds to the edge from $x_i$ to $\partial$, and $c_e = 1 / (2d)$ otherwise (here $G_{D,k}$ corresponds to the massive Green's function from 
Definition \ref {DefMassiveG}). 
The law of this random spanning tree ${\mathcal T}$ of $\hat D$ constructed by Wilson's algorithm is often referred to as a weighted spanning tree. 

 \item 
{\emph{In general graphs.}} 
Consider ${\mathcal D}=(\mathcal{V}, \mathcal{E})$ a finite connected graph. In order to be consistent with our previous study, we will assume that it has $n+1$ vertices that are labelled as 
$x_0, x_1, \ldots, x_n$. We denote by $d_i$ the number of neighbours of $x_i$ in ${\mathcal D}$ and we remark that (as opposed to the previous case), $d_i$ can vary from one point to another. 

We can then consider simple random walks on ${\mathcal D}$ and use them in order to construct a spanning tree ${\mathcal T}$ of ${\mathcal D}$ via Wilson's algorithm ``rooted at $\{x_0\}$''. We have:

\begin {proposition}[Wilson's algorithm, general case] 
\label{prop:wilson_gen}
The law of the random tree ${\mathcal T}$ constructed by Wilson's algorithm is that of a UST in ${\mathcal D}$. 
\end {proposition} 

The proof is essentially identical to that of Proposition \ref {Wilson} and left to the reader.  

\item {\emph{Weighted spanning trees.}}
Suppose we are also given a conductance function on ${\mathcal E}$; that is, for each edge $e \in {\mathcal E}$, we associate a positive conductance $c_e$. Then: 
\begin {definition}[Weighted spanning trees] 
A $c$-weighted spanning tree in ${\mathcal D}$ is a random spanning tree ${\mathcal T}$, chosen in such a way that for any spanning tree $T$,  
$$ P [ {\mathcal T} = T] = \frac {w_c(T)}{Z_c} $$
where 
$$ w_c (T) := \prod_{e \in T} c_e \hbox { and } Z_c := \sum_T w_c (T).$$ 
\end {definition}
So, a $c$-weighted spanning tree when $c$ is constant on ${\mathcal E}$ is just a uniform spanning tree. 

In order to be consistent with our previous study, we will again assume that ${\mathcal D}$ has $n+1$ vertices that are labelled as 
$x_0, x_1, \ldots, x_n$, and we denote the conductance of an edge between $x_i$ and $x_j$ by $c_{i,j}$. We let $D = \{ x_1, \ldots, x_n \}$. 

We can now define a random walk on ${\mathcal D}$ using these conductances: when the walk is at $x_i$ it jumps to $x_j$ with probability $c_{i,j} / \lambda_i$, where
$\lambda_i := \sum_{k \not= i} c_{i,k}$ (we assume that $c_{i,i}=0$). Then 
we can use this new random walk in order to construct a spanning tree ${\mathcal T}$ of ${\mathcal D}$ via Wilson's algorithm. 
Using almost exactly the same ideas, one can prove that:

\begin {proposition}[Wilson's algorithm, electric networks] \label{prop:wilson_en}
The law of the random tree ${\mathcal T}$ constructed by Wilson's algorithm is that of a $c$-weighted spanning tree, and 
$$ P [ {\mathcal T} = T ] = {w_c(T)} \times {\det G_{D,c}},$$  
where $G_{D,c}$ is the Green's function defined in Chapter \ref {Ch1}, Section \ref {elnet}
\end {proposition}

\begin{exercise}
	Prove Propositions \ref{prop:wilson_gen} and \ref{prop:wilson_en} using the same strategy used to prove Proposition \ref{Wilson}. 
\end{exercise}

\begin {remark}
Actually, the most general natural framework in which Wilson's algorithm can be made to work is that of Markov chains. The random trees that one constructs are then {\em oriented towards the chosen root}.
However, in the present notes, we will not treat this case (even if the generalisation is actually fairly immediate).  
\end {remark}
\end {enumerate} 

\section {The occupation time fields in Wilson's algorithm} 

We now come back to the setting of the UST in $D=\{x_1,\ldots, x_n\} \subset \Z^d$ with wired boundary conditions.  
Let us start this section with the following list of observations: 
\begin {itemize}
 \item When one performs a simple random walk starting from $x_1$ and stopped upon hitting $\partial D$, then by the strong Markov property, the number $N$ of returns that it makes to $x_1$ is clearly a geometric random variable. We also know that 
 $E[ N+ 1 ] = G_D (x_1, x_1)$. So, if we denote by $u=u(x_1, D)$ the probability
that the walk does not return to $x_1$ at all, we have that $P(N =  j)= (1-u)^j u$ and that $G_D(x_1, x_1) = E [ N+1 ] = 1 / u$. 
 \item If instead of the discrete-time simple random walk, we consider the continuous-time random walk that jumps with rate $1/(2d)$ along each edge, then we see that the total time $W(x_1)$ spent by this continuous-time 
 random walk at $x_1$ before exiting $D$ will be the sum of $N+1$ independent identically distributed exponential random variables with mean $1$, where $N$ is defined as before.
 But the sum of such a geometric number (plus one) of exponential random variables is also exponentially distributed, and we can conclude that $W(x_1)$ is distributed as an exponential 
 random variable with mean $G_D (x_1, x_1)$. 
\end {itemize} 
 
Motivated by this previous comment, we now also introduce a \emph{continuous-time} analogue of Wilson's algorithm. This algorithm is constructed in the same way as its discrete-time counterpart, except that one replaces the discrete-time simple random walks by continuous-time simple random walks with exponential waiting times of mean $1/ (2d)$. In this version of the algorithm, when a random walk hits the set of already discovered vertices, we instantaneously start the next random walk branch.

\begin {definition}[Occupation time fields]
Let us consider Wilson's algorithm constructing a UST of $\hat D$, in discrete or continuous time.
For all $x \in D$, we then define $V(x)$ (respectively, $W(x)$) to be the cumulative time spent at $x$ by all the discrete-time (resp. continuous-time)
random walks during the algorithm [by convention, for the discrete-time version, there is no time spent at the final vertex in each random walk ``branch'' (this vertex could be in $D$ if it is part of some previously discovered branch)]. The fields $V:= ( V(x))_{x \in D}$ and $W:= (W (x))_{x \in D}$ are called the (discrete-time and continuous-time) occupation time fields in Wilson's algorithm. 
\end {definition}

For a given $x$, by definition, the random variable $V(x)$ is a positive integer, while $W(x)$ is a positive real number. 
We have already seen that $N= V(x_1)-1$ is distributed like a geometric random variable with mean $G_D(x_1, x_1)-1$, 
and that $W(x_1)$ is distributed like an exponential random variable with mean $G_D (x_1, x_1)$. 
We will provide here a much more detailed description of the law of the fields $W$ and $V$.    

Note first that, by definition, there is an immediate relation between the law of $V$ and the law of $W$. If $(\xi_{i,j})_{j \le n, i \ge 1}$ are i.i.d. exponential random variables with 
mean $1$ that are independent of $V$, then the process $\tilde W$ defined by 
$$ \tilde W (x_j) = \sum_{i =1}^{V(x_j)} \xi_{i, j} $$ 
for $j \le n$ is distributed like $W$. In particular, this implies that the Laplace transform of $W$ can be determined easily from the Laplace transform of $V$ and vice versa: for all non-negative 
functions $k$ on $D$,  
 $$  E [ \exp ( - \sum_{j=1}^n k(x_j) W (x_j) ) ]   =   E [ \prod_{j=1}^n (1 + k (x_j))^{-V(x_j)} ].$$ 
 Indeed, 
\begin {eqnarray*}
\lefteqn {E [ \exp ( - \sum_{j=1}^n k (x_j) W (x_j) ) ] } \\
&& 
=  E [ \exp ( - \sum_{j=1}^n  \sum_{i=1}^{V(x_j) } k (x_j) \xi_{i,j}   ) ] 
=  E \Bigl[  E \bigl[ \prod_{j=1}^n \exp ( - \sum_{i=1}^{  V(x_j) } k (x_j)  \xi_{i,j} )    \big|  V \bigr] \Bigr]  \\
&& =  E \Bigl[   \prod_{j=1}^n   E [ \exp (- k (x_j) \xi_{1,j})]^{V(x_j)}   \Bigr] =  E \Bigl[ \prod_{j=1}^n (1 + k(x_j))^{-V(x_j)} \Bigr] 
.\end {eqnarray*} 
Now, these Laplace transforms turn out to have a nice compact expression in terms of determinants: 
\begin {proposition}[Laplace transforms of $V$ and $W$]\label{prop::laplace_waof}  For all non-negative functions $k$ on $D$, 
$$ E \Bigl[  \prod_{j=1}^n (1+k (x_j))^{-V (x_j) } \Bigr]  = E \Bigl[ \exp ( - \sum_{j=1}^n k (x_j) W (x_j) ) \Bigr] = \frac { \det ( - \Delta_D  )}{\det (- \Delta_D + I_k) }, $$
where as before, $I$ denotes the diagonal matrix $I (x_i, x_j) = k (x_i) 1_{i=j}$. 
\end {proposition}
\begin {remark} 
This shows that the laws of the fields $W$ and $V$ do not depend on the ordering of the $n$ points $x_1, \ldots, x_n$ that one uses when performing Wilson's algorithm. It does 
however not show yet (this will be derived in the next section) 
that $W$ and $V$ are actually independent of the constructed spanning tree ${\mathcal T}$.

A further observation is that in the special case where $k (x_j)= 0$ for all $j \ge 2$, when one develops $\det ( - \Delta_D + I_k )$ with respect to the first line of the matrix, one obtains that
$$ \det ( - \Delta_D + I_k ) = \det (- \Delta_D) + k (x_1)\det ( - \Delta_{D \setminus \{x_1 \}} ).$$ 
From here it follows that
$$
\frac { \det ( - \Delta_D  )}{\det (- \Delta_D + I_k) } = \frac {1}{ 1 + k (x_1)  (\det G_D / \det G_{D \setminus \{x_1 \}} )} = \frac {1}{1 +  k(x_1) G_D (x_1, x_1)} .
$$ 
This is as expected: it is consistent with the fact that $W(x_1)$ is distributed like an exponential variable with mean $G_D (x_1, x_1)$. However (and this is much less obvious directly) the same argument shows that $W(x_j)$ is an exponential 
random variable with mean $G_D (x_j, x_j)$ for $j \not= 1$. 
\end {remark}

\begin {proof}[Proof of Proposition \ref{prop::laplace_waof}]
The idea will be to couple two versions of Wilson's algorithm, constructing two different spanning trees: the UST ${\mathcal T}$ in the  graph $\hat D$ as described above, and the corresponding weighted massive spanning tree ${\mathcal T}'$ 
as described in Section \ref {massivecase} for mass function $k$.
 We choose to  work here with the continuous-time random walk, but the proof would also work for the discrete time version. {The continuous-time Wilson's algorithm in the massive case is defined just as in the discrete, but replacing the discrete-time random walk with killing rate $k$ by the continuous-time version: see Section \ref{S51}.}
 
We suppose that we construct the continuous-time walk using Poissonian bells on each edge: when the walk is at $x$ at a given time $t$, it waits there until 
the bell of an adjacent edge rings (and it then jumps along that edge). In our realisation of Wilson's algorithm, we construct the branches of the spanning tree ``one after the other'' (and we use the bells
 associated to the correct time intervals to define the branches). 
 
In this way, we can naturally couple a realisation of the algorithms for both the non-massive and massive cases:  the latter can hear the bells 
that ring on the additional edges in $\hat E_\partial$, while the former cannot. Note that by definition, in the latter case, when the walk is at $x$
and hears the bell on the additional edge $e' \in \hat E_{\partial}$ from $x$ to $x_0$, then it jumps along that edge. Since $x_0$ is the root point, 
this means that this edge $e'$ ends up being in the spanning tree ${\mathcal T}'$. In other words, if during the entire Wilson's algorithm used to construct ${\mathcal T}$, one did not mishear any bells on the additional edges 
in $\hat E_{\partial}$, then ${\mathcal T} = {\mathcal T'}$, and this happens if and only if ${\mathcal T}' \cap \hat E_{\partial} = \emptyset$.  
So, we can conclude that 
$$ E [ \exp ( - \sum_{j=1}^n k(x_j) W (x_j) ) ]  = P [ {\mathcal T}' \cap \hat E_{\partial} = \emptyset ];$$
the left-hand side being the probability that $n$ independent exponential random variables with means $(1/k(x_1),\ldots, 1/k(x_n))$ are greater than $(W(x_1),\ldots, W(x_n))$ respectively.
The term on the right-hand side can be directly calculated using the law of ${\mathcal T}'$: 
the weight of each spanning tree $T'$ containing no edge from the set {$\hat{E}_\partial$} is just $(2d)^{-n}$. Since the number of such spanning trees is $(2d)^n  / \det (G_D) $, we get that the term on the right-hand side is equal to
$$(\det G_D)^{-1} \times  \det G_{D,k}    = \frac {\det ( - \Delta_D )}{\det ( -\Delta_D + I_k )}.
$$
\end {proof}

Observe that this proposition fully describes the law of the processes $V$ and $W$ via their Laplace transforms.
Also recall Proposition \ref {LaplaceGFF}, which calculated the Laplace transform  of half the square of a GFF in $D$ 
to be equal to the square root of the same quantity. 
By comparison we immediately obtain the following relationship between these occupation time fields and the GFF:  
\begin {corollary} 
If $\Gamma_1$ and $\Gamma_2$ are two independent GFFs in $D$ with Dirichlet boundary conditions, then the field $(\Gamma_1)^2  / 2 + (\Gamma_2)^2 /2$ is distributed like $W$. 
\end {corollary}

In the next two sections we will discuss the occupation time fields of some different (but closely related) objects, known as loop-soups. In this case we will be able to obtain a connection similar to the above, but concerning the square of a \emph{single} GFF.

 \begin {remark} 
\label {remV-1}
In view of the loop-soup story that we will describe in the next sections, it is more natural to replace the discrete-time occupation field $V$ with the field  
$\tilde V := V -1$. We can for instance note that $V \ge 1$ anyway, because each site will be visited at least once during Wilson's algorithm. 
Then, instead of looking at $\Delta_D - I_k$, we should in fact consider the operator $- \Delta_{D, k}$, which is the Laplace operator associated to the random walk with killing (recall Chapter \ref {Ch1}, Section \ref{GFFonothergraphs}). One can obtain $\Delta_{D,k}$ from $\Delta_D-I_k$
by dividing the terms on each line $i$ by $1+ k(x_i)$. This means that $-\Delta_{D,k}$ has $1$'s on the diagonal and some terms of the form $1/(2d (1+ k(x_i)))$ on the 
off-diagonal (when $k$ is not constant, it is not a symmetric matrix any more). The formula in Proposition \ref{prop::laplace_waof} can then be rewritten as 
$$ 
 E \Bigl[  \prod_{j=1}^n (1+k (x_j))^{-\tilde V (x_j) } \Bigr]  = (\prod_{j=1}^n (1+ k(x_j))) \times \frac { \det ( - \Delta_D  )}{\det (- \Delta_{D}+I_k) } =  \frac { \det ( - \Delta_D  )}{\det (- \Delta_{D,k} ) }. 
$$
\end {remark}
\section {Discrete-time loop-soups and their occupation times} 
\label {SecLS}

\subsection {Some basic classical definitions and facts} 

It is worth first recalling here some basic properties of geometric random variables. Related ideas will be useful when trying to decompose the set of loops that 
have been erased in Wilson's algorithm into ``independent and identically distributed pieces''. For instance, one can keep in mind that in Wilson's algorithm started at $x_1$, the number of 
returns to $x_1$ before reaching $x_0$ is distributed like a geometric random variable. So, intuitively speaking, the long erased loop from $x_1$ to $x_1$ will consist of the concatenation 
of a geometric number of ``independent excursions'' away from $x_1$. 
\begin {reminder}[Infinite divisibility of geometric random variables]\label{rem:divis_geo}
Suppose that $K$ is a geometric random variable with distribution $P [ K=k] = (1-q) q^k$ for $k \in \N$ (for some given $q \in (0,1)$ that will be fixed in this reminder).
Then it is a classical fact that there exists a 
probability distribution on $\N$, such that the sum of two independent random variables with this law has the same law as $K$. 

Let us explain this in a 
more general setting.
First, when $a$ is a fixed positive 
integer, the sum $K_a$ of $a$ independent copies of $K$  satisfies 
\begin {equation}
 \label {lawKa}
P[ K_a = k ] =\frac { a (a+1) \ldots (a+k-1)  }{ k!}  q^k (1-q)^a\;\; (k\ge 1);\;\;\;\;\; P[K_a=0]=(1-q)^a
\end {equation}
(which can be seen by simply enumerating the number of possible choices for non-negative  $j_1, \ldots, j_a$ such that 
$j_1 + \ldots + j_a= k$). 
It also follows, since
$$ \frac {1}{(1-q)^a} =  \sum_{k \ge 0} \Bigl[  \frac {a (a +1) \ldots (a+k -1) }{k!} \times  q^k \Bigr],$$
that (\ref {lawKa}) defines the law of the random variable $K_a$ when $a$ takes \emph{non-integer} positive values. This distribution is known as the {\em negative binomial distribution},
and its Laplace transform is given by 
$$ E [ \exp (- \lambda K_a ) ] =  \Bigl[ \frac {1-q }{1 - q e^{- \lambda} } \Bigr]^a = E [ \exp (- \lambda K_1 ) ]^a$$
for all positive real $\lambda$. 
In particular, for any positive $a'$ and $a''$, if $K'$ and $K''$ are two independent random variables with the same laws as $K_{a'}$ and $K_{a''}$ respectively, then $K' + K''$ 
has the same law as $K_{a' + a''}$. 
As a consequence, one can decompose a geometric random variable $K_1$ into the sum of $m$ independent identically distributed random variables (each with law $K_{1/m}$) for any integer $m\ge 1$.
The example discussed at the beginning of this remark, and to be kept in mind, is simply the case $m=2$:
\[ K_1 \overset{(d)}{=} K'_{1/2} + K''_{1/2}, \] 
where $K'_{1/2}$ and $K''_{1/2}$ are independent copies of $K_{1/2}$.
\end {reminder}

The following observation may enlighten some of the combinatorics that will pop up in the next sections.

\begin {remark}
\label {lawofKm}
First notice (by considering its Laplace transform) that as $a \to 0+$, the law of $K_a$ converges to the Dirac mass at $0$. On the other hand, if 
one considers the conditional law of $K_a$ given $\{ K_a > 0 \}$, then it is straightforward to check that this converges to the probability distribution $\tilde \pi$ on $\N$ such that 
$ \tilde \pi (k) :=     {q^{k}} / ({uk}) $ for all integers $k\ge 1$ and $\tilde \pi(0)=0$; $u := \log (1 / (1-q))$. 

Writing $K_1$ as a sum of $m$ independent copies of $K_{1/m}$, we therefore have that as $m\to \infty$ each of these copies will be non-zero with probability roughly $(u/m)$, and the ones that are non-zero will have distribution given roughly by $\tilde{\pi}$. Keeping in mind the Poisson approximation of the binomial distribution and letting $m\to \infty$, we can deduce that 
$K_1$ may be decomposed as the sum of $R\sim \text{Poi}(u)$ independent random variables, each with distribution $\tilde \pi$. 
For instance, we obtain that 
$$ P [ K_1 = k ] = \sum_{r \ge 0}  \Bigl[ \frac {u^r e^{-u}}{r!} \sum_{(j_1, \ldots, j_r ) \in S_{k,r}} \tilde\pi(j_1) \ldots \tilde\pi (j_r) \Bigr], $$
 where $S_{k,r}$ denotes the set of positive $(j_1,\ldots, j_r)$ such that $j_1 +\cdots + j_r = k$.
Now, the probability on the left-hand side above is equal to $q^k (1-q)$, so by expanding the right-hand side, we obtain the following combinatorial identity: for each integer $k\ge 1$,
\begin {equation}
 \label {combinatorial}
 \sum_{r \le k} \sum_{(j_1, \ldots, j_r) \in S_{k,r} } \frac {1}{r! j_1 \ldots j_r } = 1.
 \end {equation}
Note that there are also much more direct ways to check (\ref {combinatorial}): one can for example look at the $z^k$ term in the power series expansion of 
$$\frac {1}{1-z} = \exp  \bigl[ - \log (1-z) \bigr]  = \exp \Bigl[ \sum_{j \ge 1} (z^j / j ) \Bigr]  = \sum_{r \ge 0}  \frac { \bigl[\sum_{j \ge 1} (z^j / j )   \bigr]^r}{r!} .$$
This classical identity is also related to the decomposition into cycles of uniformly chosen random permutations of a set with $k$ elements (this type of interpretation will show up when we consider
``resampling properties'' of loop-soups, see Section \ref{sec:sampling_ls}). 
\end {remark}

\begin{exercise}[Addition and splitting properties of Poisson random variables.] In this exercise we write $\mathrm{Poi}(\lambda)$  for the Poisson distribution with parameter $\lambda$. That is, if $X\sim \mathrm{Poi}(\lambda)$, we have $\mathbb{P}(X=n)=e^{-\lambda}\lambda^n/n!$ for all integers $n\ge 0$.

(i) Suppose that $(P_k; k\ge 1)$ is a sequence of independent random variables with $P_k\sim \mathrm{Poi}(\lambda_k)$ for each $k$, and $\lambda:=\sum_{k\ge 1} \lambda_k < \infty$. Show that $\sum_k P_k\sim \mathrm{Poi}(\lambda)$.

(ii) Let $M\sim \mathrm{Poi}(\lambda)$ and $(Y_j; j\ge 1)$ be a sequence of i.i.d. integer-valued random variables with $\mathbb{P}(Y_1=k)=p_k$ for $k\ge 1$.  Show that 
	$P_k :=\sum_{j=1}^M \1{Y_j=k} $
	defines a sequence of independent random variables, with $P_k\sim \mathrm{Poi}(\lambda p_k)$ for each $k$.
	\end{exercise}
	
\begin {reminder}[Basics on Poisson point processes]
\label {reminderPPP}
Suppose that $\mu$ is a $\sigma$-finite measure on some measurable space ${\mathcal M}$.  A Poisson point process with intensity $\mu$ can be loosely speaking thought of as  
a random cloud of points in ${\mathcal M}$ that somehow appear independently with ``intensity'' provided by $\mu$. Formally, it can be viewed as a random measure $N$ that assigns to each 
measurable $A$ in ${\mathcal M}$ the (integer) number of points $N(A) := N(\mathbf{1}_A)$ of the point process in $A$. The law of the process $(N(A))_A$ is characterised by the fact that 
for each $A$ such that $\mu (A) < \infty$, the variable $N(A)$ is a Poisson random variable with mean $\mu (A)$, together with the fact that for any disjoint measurable sets $A_1, \ldots, A_n$, 
the variables $(N(A_1), \ldots, N(A_n))$ are independent. 

In the special case where the space ${\mathcal M}$ is finite or countable (which we will be mostly dealing with), then for each $l \in {\mathcal M}$ 
we can define the number $N (l) := N(\mathbf{1}_{l})$ of points in the point process that are equal to $l$. The random variables $(N(l))_{l \in {\mathcal M}}$ are then 
simply independent Poisson random variable with respective means $\mu (\{l \})$.
\end{reminder}

\begin{exercise}
Suppose that $\mu$ is a finite measure on $\mathcal{M}$ with total mass $\lambda$. Show that for a suitable choice of distribution for $M$ and $Y_1$, letting $(Y_j;\, j\ge 1)$ be an i.i.d sequence, the process $N$ defined by 
\[N(A):= \sum_{j=1}^M \1{Y_j\in A} \]
for all measurable $A$ is a Poisson point process with intensity $\mu$.  How can you extend this construction to general $\sigma$-finite $\mu$?
\end{exercise}

Let us collect some useful features of Poisson point processes in the following exercises:
\begin{exercise}
	Show that if $N$ and $N'$ are independent Poisson point processes with intensity $\mu$ and $\mu'$ on ${\mathcal M}$, then $N+N'$ 
is a Poisson point process with intensity $\mu+\mu'$. 
\end{exercise}
This implies that Poisson point processes are infinitely divisible.  
For instance, if $N$ is a Poisson point process with intensity $\mu$, then it can be realised as the sum of two independent 
Poisson point processes with intensity $\mu /2$. 

\begin{exercise} \label{ex:lap_pp}
	Suppose that $F$ is a  non-negative measurable function on ${\mathcal M}$ such that $\int F d\mu < \infty$, and consider the random variable $N(F)$ that 
corresponds to the sum of the values of $F$ at all the points in a point process of intensity $\mu$. Show that the Laplace transform of $N(F)$ is given by  
$$ E [ \exp ( - \lambda  N(F) )] = \exp ( - \int_{\mathcal M} \mu (dx) (1- e^{- \lambda F(x)} ))
 . $$
 \end{exercise}

\emph{\textbf{Hint}: first use the expression for the Laplace transform of a Poisson random variable to deal with the case when $F$ takes only finitely many values.} 
\medbreak 

Finally, we note that when ${\mathcal M}$ is finite or countable, $N$ is a Poisson point process of intensity $\mu$, and when $A$ is a measurable set such that  $\mu (A)$ is finite, then   
conditionally on $N$ we can uniformly choose an order of the $N(A)$ elements of the process in $A$  
among all possible $N(A)!$ choices. In this way we obtain a finite ordered random family $U_1, \ldots, U_{N(A)}$ and we see that 
$$ P [ U_1 = u_1, \dots, U_n = u_n  | N(A) = n ] =\frac { \mu (u_1) \cdots \mu (u_n)} { \mu (A)^n} $$
and 

\begin {equation}
 \label {2.3.1}
P [ U_1 = u_1, \dots, U_n = u_n , N(A)= n ] = \frac { \mu (u_1) \cdots \mu (u_n)} { n!} e^ {- \mu (A)}.
\end {equation}

\subsection {Discrete-time loop-soups} 
Our goal is now to show that the family of loops erased during Wilson's algorithm is closely related to a Poisson point process of loops in $D$. We will 
call this process a loop-soup. Again in this section we work with $D\subset \Z^d$ and the standard (discrete-time) version of Wilson's algorithm.

Let us first provide some new basic definitions.
\begin {definition}[Rooted and unrooted loops]
We say that $l= (l_0, \ldots, l_m)$ {with $m>1$} is a rooted loop in $D$ if it is a nearest-neighbour sequence (in $D$) such that $l_0 = l_m$. 
An unrooted loop $L$ is an equivalence class of rooted loops under circular relabelling, i.e. $(l_1, \ldots, l_m , l_0)$ and $(l_0,\cdots, l_m)$ are equal as unrooted loops. 
\end {definition}

The length $m$ of a rooted  loop $l$  will be denoted by $|l|$ (this is the number of ``steps'' in the loop). When $l_0 = x$, we say that the loop $l$ is rooted at $x$. 
For all $y \in D$,  we denote by
$$j_l (y)  :=  \# \{ i \in \{1, \ldots, |l| \} , \ l_i = y \} $$
the number of visits of $y$ by the rooted loop $l$. Note that if the loop is rooted at $x$ and only returns to $x$ at the very end, then $j_l(x)$ is equal to $1$ and not $2$. 
We also denote by $|L|$ and $j_L (y)$ the corresponding quantities for unrooted loops. 
Again, keep in mind that the knowledge of $l$ (or $L$) also contains the information about the edges used in the loop (even though in the present case where $D\subset \Z^d$, there is always only one possible edge joining any two points). 

When $l$ and $l'$ are two rooted loops that are rooted at the same point, then we can define the rooted loop $l \odot l'$ of length $|l| + |l'|$ to be the concatenation of $l$ with $l'$ (the first $|l|$ steps are 
those of $l$ and the final $|l'|$ steps are those of $l'$).  
When $l$ is a rooted loop, we also define its multiplicity $J(l)$ to be the maximal integer $J$ such that $l$ can be written as the concatenation of $J$ identical rooted loops. 
It is easy to check that if $l$ and $l'$ are in the same equivalence class of unrooted loops, then $J(l) = J(l')$ -- we call this value the multiplicity $J(L(l))$ of the unrooted loop $L(l)$.

\begin {definition}[Rooted loop measure] 
For each $x \in D$, the rooted loop measure $\mu_D^x$ in $D$ rooted at $x$, is the measure on rooted loops from $x$ to $x$ in $D$ that assigns a mass $(2d)^{- |l|} / j_l(x) $ to each rooted loop $l$.
\end {definition} 

\begin {definition}[Unrooted loop measure] 
The unrooted loop measure $\mu_D$ in $D$ is the measure that assign a mass $(2d)^{-| L|} / J(L)$ to each unrooted loop $L$ in $D$. 
\end {definition}
For instance, an unrooted loop of the type $x,y,x,y,x,y,x,y,x$ where $x$ and $y$ are neighbours will have mass $(2d)^{-8} / 4$. It is a simple exercise to see that the total mass of $\mu_D$ on the set 
${\mathcal M}_D$ of all finite loops in $D$ is finite. If the reader would already like some concrete motivation for the precise nature of these definitions, see Lemma \ref{lem::relation_gf_ls}.
\begin{figure}[h]
\centering
\includegraphics[width=.6\textwidth]{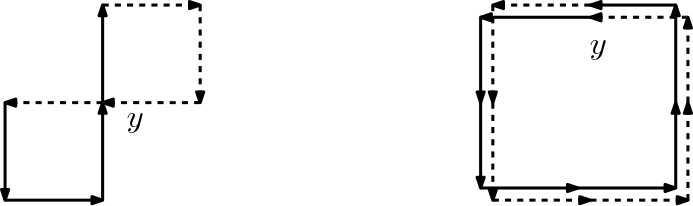}
\caption{Two loops rooted at $y$ (the right hand loop follows the same path twice). \textbf{Left}: $j_l(y)=2, J(l)=1$; \textbf{Right}: $j_l(y)=2, J(l)=2$.}
\end{figure}

We can also define the infinite measures $\mu$ and $\mu^x$ on unrooted and rooted loops in $\Z^d$ correspondingly (without the constraint that the loops remain in the finite set $D$). 

Let us now consider the set ${\mathcal M}_{D,x}$ of unrooted loops in $D$ that visit the point $x$. There is a close relationship between the restriction of $\mu_D$ to ${\mathcal M}_{D,x}$ and $\mu_D^x$. Indeed: 
\begin {itemize} 
 \item The measure $(\mu_D)|_{{\mathcal M}_{D,x}}$ is the image measure of $\mu^x_D$ under the map $l \mapsto L (l)$.
 \item In the other direction, for each $L \in {\mathcal M}_{D,x}$, one can choose a root for $L$ by picking uniformly one of the times that it visits $x$ (among the $j_{L} (x)$ possible choices). Then the image of $(\mu_D)|_{{\mathcal M}_{D,x}}$ after performing this operation is exactly $\mu^x_D$. 
\end {itemize}
When $j_l (x) = 1$, then it is clear that $\mu_D^x ( l) = (2d)^{-|l|} = \mu_D (L(l))$, so that in this case, the relations above are clear. 
When $j_l (x) = J(L(l)) \times k$ (i.e., $l$ is the concatenation of $J$ independent copies of the same rooted loop that visits $x$ exactly $k$ times), then there will be $k$
rooted loops (rooted at $x$) in the same equivalence class $L(l)$ of $l$, and indeed $\mu_D (L(l)) = (2d)^{-|l|} / J(L(l)) = k (2d)^{-|l|}/ j_l(x) = k \mu_D^x (l)$. On the other hand, when one chooses to randomly root the unrooted loop 
$L$ in the manner described above, the probability of ending up with the rooted loop $l$ is $1/k$. This shows the second statement.   

\begin {definition}[Discrete loop-soups] 
When $\alpha > 0$, a discrete loop-soup with intensity $\alpha$ in $D$ is a Poisson point process of unrooted loops in $D$ with intensity $ \alpha \mu_D$. 
\end {definition}

A loop-soup with intensity $\alpha$ can be thought of as  a random finite collection ${\mathcal L}$ of unrooted loops $(\Lambda_i)_{i \in I}$ in $D$, with cardinality given by a Poisson random variable of mean
$\alpha \mu ( {\mathcal M}_D)$.  
Note that if $D' \subset D$, then $\mu_D$ restricted to the set of loops that stay within $D'$ is exactly $\mu_{D'}$. It follows that if one considers the subset of all loops of ${\mathcal L}$ that stay within $D'$,
one gets exactly a sample of a loop-soup in $D'$. This is often referred to as the {\em restriction property} of loop-soups. 

By definition, we note that when ${\mathcal L}$ and ${\mathcal L}'$ are two independent loop-soups in $D$ with respective intensities $\alpha$ and $\alpha'$, then the union of these two loop-soups is 
a loop-soup with intensity $\alpha + \alpha'$. This will be important later on, but in the remainder of the present section we will mostly focus on the case where $\alpha =1 $ (this is the one 
that is directly related to Wilson's algorithm). So, for the rest of this section, ${\mathcal L}$ will denote a loop-soup with intensity $1$ in $D$. 

To each unrooted loop $L$ in the loop-soup ${\mathcal L}$, and each $x \in D$ that is visited by this loop, we can associate a rooted loop $l^x$ that starts and ends at $x$ as before, by choosing the starting point uniformly (and independently) at
random from one of the  $j_{L} (x)$ visits of the loop to $x$.

Let us now fix a point $x \in D$, and focus only on the collection ${\mathcal L}_x$ of loops in ${\mathcal L}$ that do go through the point $x$. This is a Poisson point process with intensity given by 
the measure $\mu_D$ restricted to the set of loops ${\mathcal M}_{D,x}$. So, the cardinality of ${\mathcal L}_x$ is a Poisson random variable with parameter $\mu_D ({\mathcal M}_{D,x})$.
\begin {lemma}\label{lem::relation_gf_ls}
One has $\exp (- \mu_D ( {\mathcal M}_{D,x}) ) = 1 / G_D (x,x)$. 
\end {lemma}
Note that this quantity is also equal to the probability that ${\mathcal L}_x$ is empty. 
\begin {proof}
Let us define $U$ to be the sum of $(2d)^{-|l|}$, over all rooted loops from $x$ to $x$ in $D$ that visit $x$ only once.
The quantity $U^2$ is therefore the sum of $(2d)^{-|l|}$ over all possible (rooted) loops that visit $x$ exactly twice, and similarly for higher powers of $U$. 
Hence the definition of $\mu_D$  shows that  
$$ \mu_D ({\mathcal M}_{D,x}) =  \sum_{j \ge 1} (U^j/ j) = - \ln (1-U).$$
On the other hand, we know from the definition of $G_D (x,x)$ (recall that we can express it as the sum over $k\ge 0$ of $(2d)^{-k}$ times the number of paths of length $k$ from $x$ to $x$ in $D$) that  
$$G_D (x, x)= 1 + \sum_{j \ge 1} U^j = 1/ (1-U),$$ which proves the lemma. 
\end {proof}
We can therefore deduce the following fact: 
\begin {corollary} 
The probability that the loop-soup ${\mathcal L}$ is empty is equal to $1 / \det G_D$.  
\end {corollary}
\begin {proof} 
The loop-soup is empty if and only if all the of the following occur: 
\begin {itemize} 
 \item There is no loop in ${\mathcal L}$ that goes through $x_1$. 
 \item There is no loop in ${\mathcal L}$ that stays in $D \setminus \{ x_1 \}$ and goes through $x_2$. 
 \item \ldots 
 \item There is no loop in ${\mathcal L}$ that stays in $D \setminus \{x_1, \ldots, x_{n-1} \}$ and goes through $x_n$. 
\end {itemize}
These $n$ events are independent (because the corresponding sets of loops are disjoint) and their probabilities are given by the previous lemma. We can conclude using the product formula for $\det G_D$. 
\end {proof}
\begin {remark} \label{loopsoupemptyalpha}
Note that the same proof shows that the probability that a loop-soup of intensity $\alpha$ is empty is equal to $1/ (\det G_D)^\alpha$.  
\end {remark}

Let us now turn back to the collection ${\mathcal L}_x$. Again for each unrooted loop $\Lambda_i$ in ${\mathcal L}_x$, we choose uniformly and independently
at random one of the $j_{\Lambda_i} (x)$ times that it visited $x$, and set $\lambda^x_i$ to be the corresponding rooted loop. In this way 
we obtain a Poisson point process of {\em rooted} loops (all rooted at $x$), with intensity measure $\mu_D^x$.  

Denote by $N$ the number of unrooted loops in ${\mathcal L}_x$ that visit $x$. They correspond (via the previous construction) to $N$ loops that are rooted at $x$, 
and  we can also choose uniformly at random (among all $N!$ choices) an order for these rooted loops. This defines an ordered collection of $N$ rooted loops which we call $\lambda_1$, \ldots, $\lambda_N$.
Finally, we can concatenate all these loops in the order we have chosen to form one single long rooted 
loop $\lambda$ from $x$ to $x$ in $D$ (when $N=0$, we just say that $\lambda$ is the empty loop).
\begin {lemma} 
For each given rooted loop $l$ from $x$ to $x$ in $D$ with $|l| \ge 1$ steps, $P [ \lambda = l ] =   {(2d)^{-|l|}}/ {G_D (x, x)}$. 
\end {lemma}
\begin {proof} 
Recall that $N$ is a Poisson random variable with mean $\mu_D ({\mathcal M}_{D,x})$.  
Also, by (\ref {2.3.1}), we know that for all $l^1, \ldots, l^r$ with concatenation $l$, it holds that
$$ P [ \lambda_1  = l^1 , \ldots, \lambda_N = l^N, N = r ] = (2d)^{- |l|} \times \frac {1}{j_1 \ldots j_r} \times \frac {1}{r! G_D (x, x)},$$ 
where the sequence $j_1$, \ldots,$j_r$ denotes the number of visits of $x$ by $l^1$, \ldots, $l^r$ respectively. 
Finally, for each given $l$ that visits $x$ exactly $k= j_l (x)$ times, each choice of $j^1, \ldots, j^r \ge 1$ with $j^1 + \ldots + j^r = k$ corresponds to exactly one 
decomposition of $l$ into a concatenation of $r$ rooted loops. Combining these observations, we obtain that 
$$ P [ \lambda = l ] = \frac {(2d)^{- |l|}} { G_D (x, x)}\sum_{r \ge 1} \sum_{(j_1, \ldots, j_r ) \in S_{k,r}  }  
 \frac {1}{ r! j_1 \ldots j_r} ,$$ 
and we can conclude by noting that (\ref {combinatorial}) holds for each fixed $k \ge 1$.
\end {proof}

Hence, if we now return to our description of Wilson's algorithm (rooted at $x_0$, in the case where $D=\{x_1,\cdots, x_n\}$, $x_1 = x$) we see that the loop $\lambda$ constructed above (allowing $\lambda$ to have zero length with probability $1/ G_D (x, x)$) 
is distributed in exactly the same way, by Remark \ref{resamplingwarmup}, as the \emph{long loop} from $x$ to $x$ that we erase before the last visit of the random walk to $x$.

Similarly, applying exactly the same reasoning to all steps of Wilson's algorithm (noting that the loop-soup restricted to those loops that do not go through $x_1$ is exactly a loop-soup in $D \setminus \{x_1\}$ and so on), 
we get the following result:
\begin {proposition}
\label {p1}
Sample an unrooted loop-soup (with intensity $1$) and an independent UST ${\mathcal T}$. 
Then, one can reconstruct a whole ``movie'' of Wilson's algorithm in $D$, with root $x_0$ and ordering $D=\{x_1,\cdots, x_n\}$ as follows.

Recall the notation $(y_1,\cdots, y_{s-1}) $ for the (loop-erased) path $X^{(1)}$ in $\mathcal{T}$ from $x_1$ to $x_0$. We construct the random walk \emph{with loops} along $X^{(1)}$ in the following manner:
\begin {itemize} 
 \item Consider the $N_1$ loops in the loop-soup that go through $y_1$. 
 Choose a root for each of them. independently and uniformly at random among the times they spend at $y_1$, and also choose an ordering of these $N_1$ rooted loops uniformly at random. Then, concatenate (i.e. trace one after the other) these loops into a single loop 
 from $y_1$ to $y_1$. 
 \item Jump from $y_1$ to $y_2$. 
 \item Consider the $N_2$ loops in the loop-soup that go through $y_2$ but not through $y_1$ and repeat the previous operation. 
 \item Proceed until reaching the root $x_0$. 
\end {itemize}
Then, we trace the other branches in an iterative fashion. 
\end {proposition}

Another way to describe this result is to start from Wilson's algorithm, and to read off the rooted loops $\lambda_i$ from $y_i$ to $y_i$ in $D \setminus \{ y_1, \ldots, y_{i-1}\}$ that are traced.
Then, for each $i$ independently, if $\lambda_i$ returns $k$ times to $y_i$, we choose to split this rooted loop into $r$ smaller loops with $j_1, \ldots, j_r$ returns to $y_i$ respectively, with a probability 
equal to $1 / ({r! j_1 \ldots j_r })$.
In this way we obtain a point process of rooted loops in $D$, which clearly also induces a point process of unrooted loops in $D$. Moreover, due to Proposition \ref{p1} we have the following key proposition.  
\begin {proposition} 
The obtained point process of unrooted loops is independent of the UST ${\mathcal T}$ that is constructed by Wilson's algorithm, and its law is that of a loop-soup in $D$. 
\end {proposition}

\subsection {Occupation times of these loop-soups}

Now, let us define ${\mathcal V} = ({\mathcal V}(x))_{x \in D}$ to be the occupation time field of a loop-soup ${\mathcal L}$ with intensity $1$ in $D$. That is, for any $x\in D$, we set $\mathcal{V}(x)$ to be the total number of visits to $x$ by all loops in $\mathcal{L}$. This is an integer-valued field, and as opposed to the 
field $V$ in Wilson's algorithm, it can take the value $0$ (it is possible that the loop-soup is empty -- while for the field $V$ we had to visit each point at least once). 
Actually, Proposition \ref{p1} implies the following:  
\begin {corollary} 
The law of $({\mathcal V}(x))_{x \in D}$ is identical to that of $(\tilde V (x))_{x \in D} := (V(x) -1 )_{x \in D}$. 
\end {corollary}

This already uncovers the following  feature: the occupation time fields $V$ and $\tilde V = V - 1$ in Wilson's algorithm are  infinitely divisible. Indeed, for each integer $k \ge 1$, we can consider 
$k$ independent realisations of the discrete loop-soup with intensity $1/k$ and their respective occupation time fields 
${\mathcal V}_{1, 1/k}, \ldots , {\mathcal V}_{k, 1/k}$. 
Then, since ${\mathcal V}_{1, 1/k} + \cdots + {\mathcal V}_{k, 1/k}$ is the occupation time field of a discrete loop-soup with intensity $1$, we get that
it is distributed like $(\tilde V(x))_{x \in D}$.

Furthermore, since we have determined the Laplace transform of the field $\tilde V$, we can deduce that of the occupation-time fields ${\mathcal V}_\alpha$: 
\begin {corollary} 
When $\alpha > 0$, and when ${\mathcal V}_{\alpha}$ is the occupation-time field of a loop-soup with intensity $\alpha \mu_D$, then 
$$ 
 E \Bigl[  \prod_{j=1}^n (1+k (x_j))^{-{\mathcal V}_\alpha (x_j) } \Bigr]  =\Bigl[  \frac { \det ( - \Delta_D  )}{\det (- \Delta_{D, k}) } \Bigr]^{\alpha}
 = \Bigl[ (\prod_{j=1}^n (1 + k(x_j)) ) \frac { \det ( - \Delta_D  )}{\det (- \Delta_{D} + I_k)  } \Bigr]^{\alpha}. $$
\end {corollary}
\begin {proof}
For each given non-negative $k$, the quantity on the left-hand side is a non-increasing function $\varphi (\alpha)$ of $\alpha$ on $\R_+$, and the infinite divisibility of loop-soup shows 
that $\varphi (\alpha) \varphi (\alpha') = \varphi (\alpha + \alpha')$. This implies that $\varphi (\alpha) = \varphi (1)^\alpha$, and we can conclude.
\end {proof}

\section {Continuous-time loop-soups and their occupation times} \label{ctdsloopsoup} 

We begin this section with two classical results, that can be viewed as the continuous limit of the infinite divisibility of geometric distributions, Reminder \ref{rem:divis_geo}
(recall that an exponential random variable can be viewed as an appropriately defined limit of geometric random 
variables). 

\begin {reminder}[The sum of two squared Gaussians is an exponential random variable] 
It is straightforward to check that when $Z$ and $Z'$ are two independent centred Gaussian random variables with variance $1$, then 
the law of $((Z)^2 + (Z')^2) /2 $ is exponential with mean $1$ (for example, by using the polar coordinate change of 
variables formula). Equivalently, this tells us that an exponential random variable of mean $1$ can be decomposed as a sum of two independent identically distributed 
random variables $Y_{1/2}$ and $Y_{1/2}'$, each having the law of $Z^2 /2 $. 
\end {reminder}

Actually, the previous fact can be considerably extended:

\begin{exercise}[Infinite divisibility of the exponential] 
\label {idexp}
	Consider a Poisson point process $N_1$ with intensity $ \pi (dt) := e^{-t} dt / t$ in $\R_+$. Note that there will be an infinite number of points in this point process, but they will 
	accumulate near $0$ (the total mass $\pi ( \R_+)$ is infinite, but $\pi ([\eps, \infty)) < \infty$ for all $\eps > 0$).
\begin{enumerate}
	\item Using Exercise \ref{ex:lap_pp}, show that the sum of the points $Y_1$ in $N_1$ is an exponential random variable with parameter $1$.
	\item For $\alpha>0$ let $N_\alpha$ be a Poisson point process with intensity $\alpha\pi$. Compute the Laplace transform of the sum of the points $Y_\alpha$ in this process.
	\item Deduce that the exponential distribution is infinitely divisible, and determine the law of each component when it is written as a sum of $n$ i.i.d. random variables.
\end{enumerate}
\end{exercise}

\begin {remark} 
In a way, the relation between (continuous-time) loop soups, both in the discrete and in 
the cable-graph setting that we will discuss later in these lectures, can be viewed as natural generalisations of this infinite divisibility property. Above, we represented both the square of a Gaussian random variable, and an exponential random variable (which can be viewed as the sum of squares of two independent Gaussian random variables), as sums of a Poissonian collection of smaller quantities. We will similarly decompose the square of the GFF (or alternatively the sum of two independent squares of the GFF) into a sum of a Poissonian collection of fields (which will be the occupation times of continuous loops). 
\end {remark}

We already defined the discrete loop-soup that appears in Wilson's algorithm, and we have noted the relation between the occupation time fields $V$ and ${\mathcal V}$. 
Our goal is now to look at the continuous-time counterpart of the loop-soup. This will lead us to discover an important relationship with the GFF: see Proposition \ref{W=Gamma2}.

Clearly, when one has an unrooted discrete loop $L$ in $D$, one can associate an unrooted continuous-time loop by sampling independent waiting times $\xi_1, \ldots, \xi_{|L|}$ for each 
of the $|L|$ steps of the loop. In this way, one obtains some continuous-time unrooted loop, with total time-length $\xi_1 + \ldots + \xi_{|L|}$.

Extending this, if we sample a discrete loop-soup with intensity $\alpha$, we can (independently for each loop in the soup) sample independent exponential waiting times with mean $1$ for each 
of the steps of the discrete loops. This yields a continuous-time loop-soup. Note that this loop-soup has the property that each loop visits at least two sites of $D$ (we need at least two sites to have a step).

Given such a continuous-time loop-soup, we can go on to define its cumulative occupation time field, and we denote this by $\hat {\mathcal W}_\alpha$ (we reserve the notation ${\mathcal W}_\alpha$ for another related object that we will define in a few lines). 
We note that the law of $\hat {\mathcal W}_\alpha$ is related to the law of ${\mathcal V}_\alpha$ just as $W$ was related to $V$: 
for all non-negative functions $k$ on $D$, 
$$ E \Bigl[  \prod_{j=1}^n (1+k (x_j))^{-{\mathcal V}_\alpha (x_j) } \Bigr]  = E \Bigl[ \exp ( - \sum_{j=1}^n k (x_j) \hat {\mathcal W}_\alpha (x_j) ) \Bigr].$$

In view of the various items that we presented so far, it is now very natural, for each $\alpha > 0$, to define a 
random field ${\mathcal Y}_\alpha = ( {\mathcal Y}_\alpha (x))_{x \in D}$ consisting of independent identically distributed random variables with the law of $Y_\alpha$
(as described in Exercise \ref {idexp}) for each $x$, 
and then to define the random field ${\mathcal W}_\alpha := \hat {\mathcal W}_\alpha + {\mathcal Y}_\alpha$, where 
${\mathcal Y}_\alpha$ and $\hat {\mathcal W}_\alpha$ are independent. 
Indeed, our description of the erased loops in Wilson's algorithm shows that the occupation time $W$ (in the continuous-time algorithm) is 
distributed exactly as ${\mathcal W}_1$. For general $\alpha$ the motivation will also soon become clear.

When $\alpha=1$,  ${\mathcal Y}_\alpha$ is defined by a collection of independent and identically distributed exponential random variables with mean 1: one for each $x$ in $D$. By the above discussion, there is also a concrete way to define $\mathcal{Y}_\alpha$ when  $\alpha = 1/2$: the independent identically distributed random variables 
$ {\mathcal Y}_{1/2} (x)$ have the law of $Z^2 /2$ where $Z$ is a centred standard Gaussian. 

This is actually very natural in the loop-soup framework as well, as $\mathcal{W}_\alpha$ corresponds to the occupation time field of a slightly different continuous-time loop-soup, including loops that visit only one point (they have a finite real life-time, but stay put at that point). This is defined by adding for each $x\in D$, to the continuous-time loop-soup with intensity $\alpha$ that we described above, an independent Poisson point process 
 of loops that visit only $x$. Such loops are described by their positive time-length, and the intensity of this time-length (in the Poisson point process that we add at each $x$) should be given by $\alpha \pi$. 

\begin {definition}[Continuous-time loop-measures and loop-soups]
In the sequel, we will refer to this Poisson point process of continuous-time loops with intensity $\alpha$ (including the Poisson point process of loops visiting only one point) as the continuous-time loop-soup (in $D$) 
with intensity $\alpha$. Its intensity measure (on the space on continuous-time unrooted loops) will be denoted by $\alpha \nu_D$.
\end {definition}

The definition of these loop-soups and their occupation-time fields are motivated by the following considerations: 
\begin {itemize} 
\item The occupation-time fields ${\mathcal W}_\alpha$ are infinitely divisible. For instance, ${\mathcal W}_{\alpha}$ has the same distribution as the sum of two independent copies of ${\mathcal W}_{\alpha /2 }$. 
Just as in the discrete case, we can note that this implies, for each given non-negative function $k$, that the function 
$$ \varphi: \alpha \mapsto E \Bigl[ \exp ( - \sum_{j=1}^n k(x_j) {\mathcal W}_\alpha (x_j) ) \Bigr] $$ 
is decreasing (if $k\ne 0$) and satisfies $\varphi (\alpha + \alpha' ) = \varphi (\alpha ) \varphi (\alpha') $. As before, this implies that $\varphi (\alpha) = \varphi(1)^\alpha$.
\item The occupation-time field ${\mathcal W}_1$ is distributed exactly like the continuous occupation-time field $W$ in Wilson's algorithm, so that for all non-negative $k_1, \ldots, k_n$, 
$$ \varphi (1)=E \Bigl[ \exp ( - \sum_{j=1}^n k_j {\mathcal W}_1 (x_j) ) \Bigr]
=E \Bigl[ \exp ( - \sum_{j=1}^n k_j W (x_j) ) \Bigr] 
= \frac { \det ( - \Delta_D  )}{\det (- \Delta_D + I_k) },$$
by Proposition \ref{prop::laplace_waof}.

\end {itemize}  
The above two points imply: 
\begin {proposition}
 For all $\alpha \ge 0$ and any non-negative function $k$ on $D$: 
$$E \Bigl[ \exp ( - \sum_{j=1}^n k (x_j)  {\mathcal W}_\alpha (x_j) ) \Bigr]
= \Bigl(  \frac { \det ( - \Delta_D  )}{\det (- \Delta_D + I_k) } \Bigr)^{\alpha}.$$
\end {proposition} 
Comparing this proposition for $\alpha = 1/2$ with the formula (Proposition \ref{LaplaceGFF}) for the Laplace transform of the GFF in $D$, we finally obtain: 
\begin {proposition}[How to construct the square of a discrete GFF via a Poisson cloud of continuous-time loops in $D$]
\label {W=Gamma2}
The field $({\mathcal W}_{1/2} (x))_{x \in D}$ is distributed like the process $(\Gamma^2(x) / 2)_{x \in D}$, where $\Gamma$ is a GFF in $D$. 
\end {proposition}

\vspace{0.1cm} 
Of course, this naturally raises the question of 
how one can construct the GFF itself out of a loop-soup, or how to describe the conditional law of the GFF given its square. The first question will 
be the motivation for the next section and the study of the GFF on cable systems, but we can already make a few comments related to the second question in 
order to illustrate the type of issue that arises. 

\begin {remark}[The sign of the GFF as an Ising model]
The only extra information one would need in order to recover $(\Gamma (x))_{x \in D}$ from $(\Gamma (x)^2)_{x \in D}$, is of course the sign $\sigma (x)$ of $\Gamma (x)$ at every point $x\in D$. The process 
$(\sigma (x))_{x \in D}$ takes its values in $\{-1, 1 \}^D$. In order to illustrate the nature of this question, it is worthwhile to make the following comment.

Given that all relevant joint distributions have smooth densities, it is easy to make sense of the conditional distribution of $\Gamma$, given 
that $(\Gamma(x)^2)_{x \in D} = (\gamma (x)^2)_{x \in D}$ for some non-negative function $(\gamma (x))_{x \in D}$. Recall that the density of $\Gamma$ at $(s (x) \gamma (x))_{x \in D}$ can be written as a multiple of  
$$ 
\exp \Bigl[  - \frac 1 {2 (2d)} \sum_{e=(x_e,y_e)} ( \1{s (x_e) = s(y_e)} (\gamma(x_e)- \gamma(y_e))^2 + \1{s(x_e) \not= s (y_e)} (\gamma(x_e) + \gamma (y_e))^2) \Bigr] . 
$$
Therefore, it follows that the conditional law of $\sigma$ given $(\Gamma(x)^2)_{x \in D} = (\gamma (x)^2)_{x \in D}$ has density of the type 
$$ P [ \sigma = s | (\Gamma(x)^2)_{x \in D} = (\gamma (x)^2)_{x \in D}] = \frac {1}{Z(\gamma)} \exp \Bigl[  - \sum_e  J(\gamma, e) \1{s (x_e) \not= s(y_e)}  \Bigr],$$ 
where  $ J(\gamma, e) = {\gamma(x_e) \gamma(y_e) }  / {d}$. 
This is nothing else than an Ising model on $D$ with  possibly inhomogeneous weights $J_e$  on different edges $e$. 
In the special case where $\gamma$ is constant, this is just the usual Ising model with inverse temperature $\gamma^2 / d$ in $D$. 
Hence, determining the conditional distribution 
of $\Gamma$ given its square is in fact an (inhomogeneous) Ising model question. 
\end {remark} 

\begin{exercise}[The square of the GFF does not satisfy the Markov property]  Consider a square $D=\{a,b,c,d\}=a+\{0,e_1,e_1+e_2,e_2\}\subset \mathbb{Z}^2 $, and let $\Gamma$ be a discrete Gaussian free field on $D$ 
with zero boundary conditions. We are going to compare the conditional distribution of $\Gamma^2 (a)$ given $(\Gamma^2(b), \Gamma^2 (c), \Gamma^2 (d)) = (1, 0, 1)$ and given 
$(\Gamma^2(b), \Gamma^2 (c), \Gamma^2 (d)) = (1, x, 1)$ for some large $x$. 

\begin{enumerate}
\item  What is the conditional law of $\Gamma (a)$ given $(\Gamma (b), \Gamma (c), \Gamma (d))$?

\item Determine the conditional law of $(\Gamma (b), \Gamma (d))$ given $(\Gamma^2 (b), \Gamma^2 (c), \Gamma^2 (d)) = (1, x, 1)$ for $x\ge 0$ (equivalently, the conditional probability that $\Gamma(b)$ and $\Gamma(d)$ have the same sign). What can you say when $x \to \infty$?  

\item Using (i) and (ii), determine the conditional law of $\Gamma^2 (a)$ given $$(\Gamma^2(b), \Gamma^2 (c), \Gamma^2 (d)) = (1, x, 1).$$ Does this depend on $x$?

\end{enumerate}
\end{exercise}

\begin{exercise}[Wilson's algorithm in one dimension and decomposition of (reflected) Brownian motion]
For $n\in \N$ let $D_n$ be the graph defined by taking the subgraph $\{0,1, \ldots ,n\}$ of $\Z$ together with an edge from the site $n$ to itself (so that the simple random walk on this graph has probability $1/2$ 
to stay at $n$ when it is at $n$, and probability $1/2$ to jump to $n-1$).
\begin{enumerate}
	\item Describe Wilson's algorithm (in continuous time) rooted at $0$ and starting at the vertex $n$? What is the law of the associated path (with and without the erased loops)?
	\item Describe how to form a collection of unrooted loops from the above erased loops, that has the law of a Poisson point process
	(\textbf{Hint}: although this setting is slightly different to that considered already in these lecture notes, you can check that the same arguments hold without modification, because the jump probabilities are $1/2$ 
	everywhere).
	\item Relate the law of a certain functional of this Poisson point process to the law of the time taken for a (continuous-time) random walk started at $n$ and reflected at $n$ to reach $0$.
	\item Prove that the time to hit $1$, for a reflected Brownian motion started from $0$ on $[0,1]$, has infinitely divisible law.
\end{enumerate}\end{exercise}

\textbf{Remark:} It is possible to relate the previous considerations to the Poisson point process of excursions above its future minimum of reflected random walk (or reflected Brownian motion) 
stopped at its first hitting time of some positive level. 

\section {Resampling and Markovian properties of unoriented loop-soups} \label{sec:sampling_ls}

In the previous section we saw that the loop-soups ${\mathcal L}_1$ and ${\mathcal L}_{1/2}$ were quite special, as they could be related 
to two nice probabilistic objects. Namely, ${\mathcal L}_1$ is related to uniform spanning trees (and Wilson's algorithm), and ${\mathcal{L}}_{1/2}$ to the Gaussian Free Field (via its square).

In the present section, we will discuss some properties of the loop-soups themselves, that are reminiscent of the resampling and Markovian properties of the Gaussian free field. They also
highlight that the relationship between ${\mathcal L}_{1/2}$ and the GFF arises from simple yet deep properties of this particular loop-soup. 
As we shall see, the features we are interested in are best expressed and understood when one considers the loops in ${\mathcal L}_{1/2}$ to be {\em unoriented}. In this setting, the 
loop-soup ${\mathcal L}_{1/2}$ has a unique special feature among all $({\mathcal L}_\alpha; \alpha>0)$.

It is worth stressing that
if we were to stick to oriented loops as in the previous sections, then we could derive very similar features for the loop-soup ${\mathcal L}_1$ (and these special properties give rise to the relationship with Wilson's algorithm) but we will not discuss this in detail here.  

\subsection {Unoriented loops and loop-soups}

In the previous sections, all our loops (rooted and unrooted) were oriented: the loops $xyztx$ and $xtzyx$ winding in different orientations around the same square 
were different (although both had  $\mu_D$ mass $1/(2d)^4$). Up to this point, it was important for us to consider oriented loops, in order to make 
sense of the concatenation of loops, and since loops in Wilson's algorithm are naturally oriented with respect to the time they appear. 
However, as we will explain in the present section, it is somehow more natural to consider unoriented loops when one is studying the relation with the discrete GFF. 

Suppose that $l= ( l_0, e_1, l_1 \ldots e_m, l_m=l_0)$ is a rooted loop as defined before. We can then define its  time-reversal
$r( l) := (l_m, e_m, l_{m-1}, \ldots, l_1, e_1, l_0)$. When the rooted loops $l$ and $l'$ are in the same equivalence class $L(l)=L(l')$ of unrooted loops, then $r(l)$ and $r(l')$ are clearly also in the 
same equivalence class of unrooted loops, and we call this class $r(L)$. We can then define a further equivalence relation, now on unrooted loops, that identifies $L_1$ and $L_2$ whenever $L_1=r(L_2)$. 

The \emph{unrooted unoriented} loop $U(L) = U(L(l))$ is the equivalence class of $L$ under this relation. Let us denote $\delta (L) = \delta (U)$ the cardinality of $\{ L, r(L) \}$, so this quantity is $1$ if $r(L)= L$, and is $2$ if $r(L) \not= L$.

\begin {definition}[Unoriented loop measure] 
The measure $\kappa_D$ on unoriented unrooted loops in $D$ is the measure assigning mass $(\delta (U) / 2 ) \times ((2d)^{-|U|} /  J(U))$ to each unrooted loop $U$ of length $m$ in $D$. 
\end {definition}

So, when $J(U) = 1$ and $\delta (U) = 2$, we see that $\kappa_D (U) = (2d)^{-|U|}$ (this will be typically the case for very long loops in a very large domain $D$). 

We can note that $\kappa_D$ is by definition the image measure of $\mu_D / 2$ under the map $L \mapsto U(L)$. Conversely, if one starts from an unrooted unoriented loop $U$, one can choose an orientation at random, 
in order to define an oriented unrooted loop $L$. The image measure of $\kappa_D$ under this operation will give rise to the measure $\mu_D / 2$, because 
$$ \mu_D (L) / 2  = \kappa_D (U) / \delta (U). $$ 

\begin {definition}[Unoriented loop-soup] 
An unoriented loop-soup in $D$ with intensity $c > 0$ is a Poisson point process of unoriented unrooted loops with intensity $c \kappa_D$. 
\end {definition}

By the above comments, when one samples a soup of (unrooted) oriented loops according to the 
loop measure $\alpha \mu_D$, and one forgets about the orientation of the loops, then one gets a soup of unrooted unoriented loops with intensity $2\alpha \kappa_D$.
Conversely, starting from a soup of unrooted unoriented loops with intensity $c\kappa_D$,  one can define a soup of oriented loops of intensity $c \mu_D / 2$
by choosing the orientation of each loop at random. 
In order to avoid confusions, we use the letters $\alpha$ to denote the intensity of soups of oriented loops (i.e. with intensity measure $\alpha \mu_D$) and
$c$ to denote the intensity of soups of unoriented loops (i.e. with intensity measure $c \kappa_D$). The natural relation between $c$ and $\alpha$ is therefore $c=2 \alpha $ and $\alpha = c/2$.  

In view of the previous sections, we can say that the  soup of oriented loops with intensity $\alpha=1$ is very closely related to Wilson's algorithm and that the 
soup of oriented loops with intensity $\alpha = 1/ 2$ is very closely related to the GFF. However, this relation to the GFF goes via the occupation time field, which 
can be just as well defined using the unoriented version of the loop-soup (the occupation time of an oriented loop does not depend on its orientation). So, we can say that 
 the soup of unoriented loops with intensity $c=1$ is very closely related to the square of the Gaussian Free Field. 

 We are now going to highlight some properties of this soup of unoriented loops with intensity $c=1$. 
Most of these properties have counterparts for the oriented loop-soup with intensity $\mu_D$ (i.e., for $\alpha = 1$).
However, since our prime motivation here is to discuss the relation to the GFF, we choose to focus solely on the case of unoriented loops with intensity $\kappa_D$.

To set up notation, suppose that we are given such a loop-soup. Then for each unoriented loop $U$, we can consider the number $N(U)$ of occurrences of $U$ in the loop-soup. By definition, these numbers $N(U)$ will be independent Poisson 
 random variables with respective means $\kappa_D (U)$.

\subsection {The law of the loop-soup given its occupation time measure on edges}

 \begin{exercise}[Warm-up to the resampling property of the loop soups.]
	Suppose that $D\subset \Z^2$ consists of the site $x_0$ and of the union $S_1\cup S_2$, where  $S_1$ is the unit square with bottom right hand corner $x_0$ and $S_2$ is the unit square with top left hand corner $x_0$.
	\begin{enumerate}
		\item
		Consider a (discrete-time) oriented loop-soup of intensity $\alpha$ in $D$, and the possible loop-soup configurations that give rise to exactly one jump along each of the eight edges of 
		$D$, in such a way that these jumps go clockwise around $S_1$ and around $S_2$. How many such loop-soup configurations are there?
		For each of them, determine the ratio of their probability with the probability that the loop-soup is empty. 
		What is special when $\alpha =1$? Can one interpret this in terms of ``resampling the connections at $x_0$''? 
		
		\item 
		Now consider (oriented) loop-soup configurations that give rise to exactly one jump along each of the eight edges of $D$, but with no constraint on their orientation. 
		Determine the number $K$ of such loop-soup configurations? For each of them, determine the ratio of their probability with the probability that the loop-soup is empty.
		If for each loop in the loop-soup, one forgets about its orientation, one gets a configuration of unoriented loops. Determine the number $K'$ of unoriented loops 
		there are that correspond to the previous $K$ configurations of oriented loops. Compare the probabilities of these $K'$ configurations. 
		What is special when $\alpha =1/2$? Can one interpret this in terms of ``resampling the connections of the unoriented loop-soup at $x_0$'' in that case? 
		\end {enumerate}\end{exercise}
	
In the previous sections, we studied the occupation time measure $V$ of a discrete-time oriented loop-soup with intensity $\alpha =1$, 
 which was defined on the sites of $D$. For each point in $D$, we counted how often this point has been visited by the loops in the loop-soup. However, the loop configuration also includes (extra) information about the number of jumps \emph{along each edge} between two points in $D$. 
For instance, if $D$ consists of a little square with four vertices abcd and if we know that the occupation time measure at each of the four points is $1$, then we do not know 
whether the loop-soup consists of two loops (aba and cdc, or ada and bcb) or of one loop (abcda). Conversely, if for each edge $e$ we know the total number of jumps along that edge, then 
we can deduce the occupation time at each site by summing the number of jumps along each edge adjacent to this site (and for all this, one does not need to know the orientation of the loop). 
Hence, as soon as $D$ contains a ``cycle'', we see that the occupation-time measure on edges contains \emph{more} information about the loop-soup than the occupation-time measure on sites. 

It turns out to be more convenient and natural to consider the trace of loops on edges rather than on sites, when one discusses resampling and Markovian properties of the loop-soups, or 
Markovian properties of their occupation time measures. 

In the remainder of this section, we will focus on a soup of unoriented unrooted loops in $D$ with intensity $\kappa_D$, and we will use in an essential way here that $c=1$. 
Such a soup defines an integer-valued occupation time field $(T(e))_{e \in E}$ on the set of edges $E = E_D$ (the set of edges with both end-points in $D$) -- mind that from now on in this chapter, the notation  
$T$ will be used solely for this field (we will
not discuss spanning trees in the remainder of this chapter). 
When $x \in D$, we denote by $E(x)$ the set of edges in $E_D$ that are adjacent to $x$ (i.e., one of their endpoints is $x$). Note that by definition, the occupation time field $(S(x))_{x \in D}$ on sites 
is related to the field 
$(T(e))_{e \in E}$ 
by 
$$ S(x) = \frac {1}{2} \sum_{e \in E(x)} T(e)$$ 
for all $x \in D$, so that the field $T$ contains at least as much information as $S$. (Typically, $T$ actually contains {\em strictly} more information than $S$: one can think for instance of a configuration in a square $D$ where each site 
has been visited once; then the occupation times of the four edges of the square could be $1, 1, 1, 1$ or $2, 0, 2, 0$ or $0, 2, 0 , 2$). 
Recall that the occupation time measure $S$ corresponds to the occupation time measure ${\mathcal V}_{1/2}$ of a loop-soup with intensity $\alpha = 1/2$, and that we have described its law 
via its Laplace transform: 
$$ 
E \Bigl[  \prod_{j=1}^n (1+k (x_j))^{-S (x_j) } \Bigr]  = \Bigl[  \frac { \det ( - \Delta_D  )}{\det (- \Delta_{D, k}) } \Bigr]^{1/2}. $$
for all non-negative functions $k$. 

We are going to address the following two questions: 
\begin {itemize} 
 \item Can we describe the law of the field $(T(e))_{e \in E}$? 
 \item What is the conditional law of the loop-soup, when one conditions it on its occupation time measure $T$ on the edges? That is, what can we say about the actual collection of loops?
\end {itemize}

Let us begin with the first question.
When $t=(t_e)_{e \in E}$ is a collection of integers defined on $E$, we define $s = (s(x))_{x \in D}$ 
to be the corresponding quantity on sites: $s(x)$ is half of the sum of $t(e)$ over the edges $e\in E(x)$.
A first remark is that the law of $T$ will be supported on the set ${\mathcal A}$ of (``admissible'') occupation measures $t$, such that $s(x)$ is an integer for each $x$. In other words, the sum of all $t(e)$ over $E(x)$ has to be even for each $x$. On the other hand, it is easy to see that for each $t \in {\mathcal A}$, the probability that $T=t$ will be positive. 
For such $t$ we define the total occupation time 
$$ | t | := \sum_e t(e).$$

For every even integer $2u$, we also set ${\mathcal P}(2u) = (2u)! / (2^u u!) = (2u-1)\times  (2u-3) \ldots 3 \times 1$ to be the number of possible ways of decomposing $\{1, \ldots, 2u \}$ into pairs (with the convention ${\mathcal P}(0) = 1$). 

\begin {proposition}[Law of the occupation time field] 
\label {proplawofoc}
For all $t \in {\mathcal A}$, 
$$ P [ T = t ] = \frac {1}{\sqrt { \det G_D}} \times (2d)^{-|t|} \times \Bigl[ \prod_{x} {\mathcal P} ( 2s(x))  \Bigr] \times \prod_e \frac 1 {t(e)!}. $$ 
\end {proposition}

We will prove this result together with the answer to the second question, that we now turn to. The following exercise serves as a warm up to these considerations.

\begin{exercise}[A renewal-type property (on edges) of the loop soups] 
	Suppose that for some subset $E'\subset E_D$ we condition on the event $\{T(e)=0 \;\forall e\in E'\}$. Using the definition of the loop soup as a Poisson point process of loops, describe the conditional law of the soup given this event.
	
	For $x$ a given site,  let $C_x$ denote the connected component containing $x$ of the graph defined by $D$ and by the set of edges $e$ such that $T(e) \not= 0$. 
	Show that conditionally on $C_x=C$, the law of $T^D$ restricted to the set of edges $E_{D \setminus C}$ is the law of $T^{D \setminus C}$. 
\end{exercise}

We note that the missing information (when one knows $T$) in order to determine the loop-soup,
is how to connect all these jumps along edges to each other in order to create loops. 
The intuitive way to think about it is as follows: 
\begin {itemize} 
 \item the occupation time measure on edges provides for each edge $e$, an integer number of straight pipes;
 \item in order to recover a loop-soup, one has at each site $x$, to decide which adjacent pipes will be paired and connected to each other. So, one has to choose a pairing of the (even) number $\sum_{e \in E(x)} T(e)$ of pipes.  
\end {itemize}
The answer to the second question is then the following.
\begin {proposition}[Resampling the connections at sites] 
\label{law_pipes}
The conditional law of the loop-soup given $T$ can be described as follows: 
at each site $x$, choose (independently for each $x$) a pairing of the $2S(x)$ adjacent pipes uniformly among all choices.
\end {proposition}
One striking feature of this proposition is that the choice of connections at different sites is made independently. This indicates that in some sense, the 
interesting ``long-range'' interaction properties of a loop-soup are already encapsulated by $T$.

Let us prove these two propositions together.
\begin {proof}
In order to prove these two facts, the following trick can be useful in order to avoid getting sidetracked into unnecessary combinatorial considerations. Instead of considering subgraphs of $\Z^d$, we will 
consider the graph that is obtained when each edge of $\Z^d$ appears $K$ times (i.e. there are $K$ identical copies of each edge). We can note that this change will not affect the definition of the random walk, of the Green's function, of the Laplacian or 
of the GFF. It also changes almost nothing about the loop-soup in $D$, except that if we want to keep track of the occupation times on edges, things will be a little different. Typically, when $K$ is very large 
for a given $D$, the loop-soup will tend to avoid using edges twice. More precisely, as $K \to \infty$, the probability that a loop-soup in $D$ (with these $K$-multiple edges) uses 
any edge twice goes to $0$. 

Let $U$ denote a configuration of the loop-soup that gives rise to the occupation-time measure $t$ on edges, and that uses no edge twice. By (3) of Reminder \ref {reminderPPP}, the probability 
of this configuration $U$ is simply proportional to 
$(2dK)^{-|t|}$. In order to determine the probability that we obtain $T = t$ (for the original loop-soup on $\Z^d$) it therefore suffices to enumerate the number of loop-soup configurations that give rise to $t$, multiply this by $(2dK)^{-|t|}$, renormalise to obtain a probability measure, and to finally
let $K \to \infty$. 

In order to perform this enumeration, we can first deal with the number of possibilities for each edge. One needs to choose $t(e)$ edges among the $K$ available ones, which gives rise to $K! / ( t(e)! (K-t(e))!)$ possibilities, and  this quantity behaves like $K^{t(e)}/t(e)!$ as $K \to \infty$. 
Then, once we know which edges are used, we need to pair them at each site in order to create the loop-soup configuration. Clearly, there are ${\mathcal P} (2s(x))$ choices at each site. 
Proceeding as described above, we therefore see that $P[T=t]$ must be proportional to $(2d)^{-|t|} \times \prod_{x} {\mathcal P} ( 2s(x))   \times \prod_e (  1 / {t(e)!} ) $. 
But we know by Remark \ref{loopsoupemptyalpha} that the probability that the loop-soup is empty (i.e., that $t$ is identically $0$) is equal to $1 / \sqrt{\text{det} G_D}$, which 
determines the renormalisation factor and concludes the proof of Proposition \ref {proplawofoc}.  

This argument also immediately  provides the resampling property at sites (Proposition \ref{law_pipes}).
\end {proof}

\subsection {Bridges and Markov property} 

In order to describe the Markov property of loop-soups and their occupation times, we need to make the following slight extensions and modifications to 
our set-up. Consider a finite subset $D$ of $\Z^d$ as before, and let $E=E_D$ again denote the set of edges of $\Z^d$ that have two endpoints in $D$.
We will also now consider a subset $E_1$ of $E$,
and consider the discrete-time random walk in $D$ that is killed at the first time $\sigma$ that it ``attempts'' to jump along an edge that is not in $E_1$. 
We can then define the corresponding Green's function $G_{E_1} (x,y)$, defined for $x,y$ with $E(x)\cap E_1\ne \emptyset, E(y)\cap E_1\ne \emptyset$, by 
$$ G_{E_1} (x, y) := E_x \Bigl[ \sum_{j=0}^{\sigma-1} \1{X_j = y} \Bigr], $$ 
and the corresponding Laplacian on $(D,E_1)$ (which is just obtained from $\Delta_D$ by replacing instances of $1/(2d)$ that correspond to edges in $E \setminus E_1$ by zeroes).
We say that a path from $x$ to $y$ stays in $(D, E_1)$ (or simply, in $E_1$) if it only uses edges in $E_1$. 

Consider two points $x$ and $y$ in $D$.
We say that a bridge $b$ from $x$ to $y$ in $E_1$ is a finite nearest-neighbour path in $(D,E_1)$ (as always, keeping track of the edges used) that starts at $x$ and finishes at $y$. 
We write $|b|$ for the length (number of jumps) of $b$, and a bridge from $x$ to $x$ is allowed to have zero length. 
By definition, $G_{E_1} (x,y)$ is then the sum over all bridges from $x$ to $y$ in $E_1$ of $(2d)^{-|b|}$. 

We can therefore define a probability measure  on bridges from $x$ to $y$ in $E_1$, that assigns a probability  $ (2d)^{-|b|}  / G_{E_1} (x,y)$ to each bridge $b$. 

\begin {definition}[Unordered unoriented bridges]\label{uub}
Suppose that $Z=(z_1, \ldots, z_{2N})$ are $2N$ points in $D$. An unordered unoriented $Z$-bridge in $E_1$ is a pairing $t$ of $\{1, \ldots , 2N \}$    
  (this is a 
  permutation  
  $(t^1_1, t^2_1) \ldots (t_N^1, t_N^2)$ of $\{1,\cdots, 2N\}$, where the transpositions are ordered according to some lexicographic rule), together with a collection of $N$ unoriented bridges joining the $N$ pairs $(z_{t_k^1}, z_{t_k^2})_{k\le N}$ in $E_1$.
\end {definition} 

\begin {definition}[Bridge measure]
   Suppose that $Z$ admits a pairing $t$ of $\{1,\ldots, 2N\}$ such that $G_{E_1} ( z_{t_k^1}, z_{t_k^2}) \not= 0$ for all $k \le N$. 
   
   Then   we define the measure $B_Z^{E_1}$ on unoriented unordered $Z$-bridges as follows: 
 \begin {enumerate}
  \item 
 we first sample a pairing $\tau$ in such a way that the probability of a given pairing $t$ is proportional to  $\prod_{k=1}^N G_{E_1} ( z_{t_k^1}, z_{t_k^2})$;
  \item
 given that $\tau=t$, we then sample $N$ independent (unoriented) bridges in $D$, from the measure described just before Definition \ref{uub}, joining the two points of each of the $N$ pairs $(z_{t_k^1}, z_{t_k^2})$.
 \end {enumerate}
 \end {definition}
 The definition basically means that we sample a $Z$-bridge in such a way that the probability of observing a given $Z$-bridge, whose $N$ sub-bridges have lengths summing to $K$, is just proportional to $(2d)^{-K}$.
 
We are now ready to describe the Markov property of loop-soups. Suppose that we consider an unoriented unrooted loop-soup in $D$, and that $E_1$ and $E_2$ form a partition of the set of edges that join two points of $D$.
Our goal is now to study the conditional law of the loop-soup in $E_1$ given ``its trace'' on $E_2$. For this, we write $\eta$ for the collection of all portions of loops or of entire loops in the loop-soup 
that use  edges in $E_2$. This consists of (counting all of these with their multiplicity): 
\begin {itemize} 
 \item The unoriented loops that use only edges of $E_2$. 
 \item Unoriented bridges in $E_2$ that are subsets of loops in the loop-soup, and that are maximal in the sense that they are contained in no longer bridge in $E_2$ in that loop. In other words, these are the collection of all excursions in $E_2$ 
 of the loops in the loop-soup. 
\end {itemize}
Each of these ${\mathcal N}$ unoriented bridges has two endpoints. We denote by ${\mathcal X}$ these $2{\mathcal N}$ endpoints of $\eta$. 
We finally define $\beta$ to be the ${\mathcal N}$ ``missing pieces'' that are needed to complete these ${\mathcal N}$ unoriented bridges in order to form the loops of the loop-soup. When one conditions on ${\mathcal X}$, 
then these missing pieces do form an unordered unoriented ${\mathcal X}$-bridge in $E_1$.

By adapting the previous ``multiplicity of edges $K$ tending to infinity'' idea, it is a simple exercise (that we leave to the reader) to prove the following fact: 

\begin {proposition}
\label {p111}
The conditional distribution of $\beta$ given $\eta$ is exactly the unordered unoriented bridge measure $B_{{\mathcal X}}^{E_1}$. 
\end {proposition}

It is worthwhile to highlight that this conditional law depends on $\eta$ only via the knowledge of $\mathcal{X}$. 
So, conditionally on these $2 {\mathcal N}$ endpoints, $\eta$ and $\beta$ are 
independent. Note however that the role of $E_1$ and of $E_2$ are not totally symmetric in this set-up For our definition, $\eta$ cannot have bridges of length $0$, while $\beta$ is allowed to have bridges of zero length.

\section {A quick survey of the GFF and loop-soups on cable graphs} 

In this section we survey, in a narrative and heuristic way without full proofs, aspects of the conditional law of the GFF given its square. 
This will enable us to not only define the square of the GFF from a loop-soup (as in the previous sections), but also how to define the GFF itself out of a (somewhat more complete) loop-soup. 

\subsection {The cable-graph GFF}

The  important new object is the so-called \emph{cable graph} associated to a subset $D\subset\Z^d$. This graph $C$ is defined to be the union of $D$ with all open edges (that we view as open intervals in $\R^d$) that are adjacent to the sites of $D$. The points of $\partial D$ are called the boundary points 
of the cable graph and will also be denoted by $\partial C$.
In view of our picture of the discrete GFF as being built of a collection of Gaussian springs, it is very natural to define the \emph{cable graph GFF} as a random continuous function on $C$ (or actually on $\overline C := C \cup \partial C$) using the following procedure. 

In the following construction, we consider all edges to have length $1$ -- i.e., we measure length using the parametrization of the edge between $x$ and $y$ by $t \mapsto t x + (1-t) y$ (so that all the Brownian bridges that we will refer to correspond to Brownian bridges on a time-interval of length $1$).

\begin {definition}[Cable graph GFF]
First sample a discrete GFF $\Gamma$ in $D$. Then, conditionally on $(\Gamma (x))_{x \in D}$, sample an independent Brownian bridge $\Gamma_e$ on each edge $e=(xy)$ of $C$, conditioned to satisfy $\Gamma_e (x) = \Gamma (x)$ and $\Gamma_e (y) = \Gamma (y)$, {with the convention that $\Gamma(y)=0$ for $y\in \partial D$.}
The obtained process $\Gamma$ defined on $C$ (equal to the discrete GFF $\Gamma$ on sites and to $\Gamma_e$ on edges) defines the cable-graph GFF on $C$. 
\end {definition} 

In other words, instead of interpolating the discrete GFF $\Gamma$ linearly on the edges (which also seems like a reasonable way to extend $\Gamma$ to the cable graph) we interpolate $\Gamma$ \emph{via Brownian bridges}. 
This is fairly natural, as one can then view the obtained function as the generalisation of one-dimensional Brownian motion, but when time is replaced by the whole cable graph $C$.

In the above construction, it is easy to see that the obtained cable-graph GFF $\Gamma$ is a continuous centred Gaussian process on $C$. 
Its law can therefore be described using only its covariance function $K_C (x, y)$ defined on $C \times C$.  By definition, 
this covariance function will coincide with $G_D$ when $x$ and $y$ are both sites of $D$. 
Furthermore, our definition also shows that for each $x$, the function $y \mapsto K_C (x, y)$ is actually harmonic (and therefore linear!) on each edge-portion that does not contain $x$. 

The first observation which  makes the cable-graph GFF so useful is the following. 
If one samples a cable-graph GFF, one can define its zero-set $Z$ to be the (closed) set of points $x$ in $C$ such that $\Gamma (x) = 0$. 
Then, one can define the ``excursions of $\Gamma$" to be the connected 
components of $C \setminus Z$. Since $\Gamma$ is a continuous function on the cable graph, its sign must be constant on each excursion.
Suppose now that one observes only the function $| \Gamma |$ (or equivalently the square of the function $\Gamma$). Then, one knows the set $Z$ (which is also the 
zero-set of $| \Gamma |$). The knowledge of the sign of the GFF on each excursion is the only missing information that would allow $\Gamma$ to be recovered. 

Now, the following fact is intuitively clear: 
\begin {proposition}[Law of $\Gamma$ given $\Gamma^2$ on the cable-graph] \label{P:g|g^2}
Conditionally on the zero-set $Z$, the law of the signs of $\Gamma$ on each of the excursions is given by a collection of independent fair coin tosses. 
\end {proposition} 

 Let us outline an elementary way to derive this fact. The idea is to verify that if we switch the sign of a single excursion of $\Gamma$, say containing a given point $x$, then the law of $\Gamma$ remains unchanged. 

\begin {itemize} 
 \item Donsker's invariance principle (based on the central limit theorem and a simple tightness argument) shows that suitably rescaled random walks on $\Z$ converge to one-dimensional Brownian motion. Similarly (using a ``local central limit theorem'' on each of the edges of the cable-graph) it is possible to show that the following processes will converge to the cable-graph GFF. For $n\ge 0$, divide each edge of the cable graph into $2n$ intervals of length $1/2n$. Consider the set of continuous functions on the cable graph that are equal to $0$ on $\partial D$, and that are differentiable with derivative equal to $+1$ or $-1$ on each of these 
length $(1/2n)$-intervals. For each $n$, let $f_n$ be a random function sampled from the \emph{uniform} measure on this set of functions: then, in the space of continuous functions on $C \cup \partial D$ endowed with the sup-norm, these $f_n$ converge weakly as $n\to \infty$ to the cable-graph GFF. For convenience, by Skorokhod's representation theorem, we can couple all the $f_n$ together with a realisation of the cable-graph GFF $\Gamma$, so that  $\sup | f_n - \Gamma | \to 0$ almost surely.
 
 \item Furthermore, if we define the mapping $f \mapsto \partial D \cup {\mathcal Z} (f)$ that associates to a continuous function on the cable-graph its zero-set (viewed as an element of the set of compact subsets of {$\overline{C}$} endowed with the Hausdorff distance) then with probability one, the cable-graph GFF $\Gamma$ is a continuity point of this map. This can be seen by considering sample path properties of one-dimensional Brownian bridges - for instance, using that any zero of such a bridge is almost surely an accumulation point of intervals where it is positive and an accumulation point of intervals where it is negative -  and applying this remark to the restriction of the GFF on each edge of the cable-graph.
 \item Finally, for each given $x$ and each continuous function $f$ on the cable-graph, we define $S_x(f)$ to be the function obtained from $f$ by just swapping the sign of the excursion of $f$ that contains $x$. Clearly, for all given $n$ and all given $x$, the law of $S_x (f_n)$ is the same as the law of $f_n$. 
 Furthermore, the fact that $\Gamma$ is a continuity point for the zero-set mapping readily implies that $ S_x (f_n) \to S_x (\Gamma) $ almost surely as $n \to \infty$. Hence, we conclude that $S_x (\Gamma)$ and $\Gamma$ do have the same distribution. 
\end {itemize}

\begin{exercise}
Turn the above outline into an actual proof.
\end{exercise}

\begin {remark} 
Note that we have now seen how to construct a cable graph GFF out of a discrete GFF (by adding Brownian bridges), and how to construct a cable graph GFF (and therefore a discrete GFF) out of the square of a cable graph GFF. However, this does not quite provide a recipe for how to construct a discrete GFF out of the square of a discrete GFF. In the next section, we will show how to  directly construct the square of a cable-graph GFF via other means (namely, using a Brownian loop-soup), from which one can then construct a discrete GFF. 
\end {remark} 

\begin {exercise}
1) Extend the definition of the cable graph and of the cable graph GFF to the case where the boundary points of $C$ are not necessarily points of $\Z^d$ but can lie anywhere on the edges of $\Z^d$. 

2) Extend the definition of local sets to the case of the GFF on cable graphs.

3) Show that for a given point $x$, the excursion of the cable-graph GFF that contains $x$ is a local set. 
\end {exercise}

\subsection {Brownian loop-soup on the cable graph} 

Just as in the case of the discrete graphs, it turns out that there is a direct way to construct the square of a GFF on the cable-graphs using loop-soups.  

Indeed, the covariance function $K_C$ of the GFF on the cable graph does correspond to an actual Green's function 
that can be interpreted in terms of Brownian motion on the cable-graph $C$.  It is 
easy to define such a Brownian motion, which is heuristically the trajectory of a particle moving at random on $C$ (mind that this Brownian motion is different from the GFF: it 
is parametrised by time and takes values in $C$). When this cable graph Brownian motion moves on an edge, 
it locally behaves like Brownian motion on that edge (i.e. like a one-dimensional Brownian motion), and each 
time it reaches a site $x \in D$ it chooses uniformly at random along which of the $2d$ adjacent edges to move along next. 
Once this Brownian motion is defined, we can choose to stop/kill it as soon as it reaches 
$\partial D$. Then, if it is started at $x\in C$ (at time $0$) the law of its position at time $t$ has a density $p^{\text{cable}}_t (x,y)$ with 
respect to one-dimensional Lebesgue measure on the cable graph. One can then define its Green's function 
$G_D^{\text{cable}} (x,y) = \int_0^\infty p_t^{\text{cable}} (x, y) dt$ and it turns out that $K_C$ is a constant multiple of $G_D$. 

If we divide each edges of the cable-graph into $2n$ pieces as before, and consider the random walk that at each step jumps to one of the neighbours at distance $1//2n$ of the actual position, then it is easy to check that (when suitably rescaled), this random walk converges to this cable-graph Brownian motion as $n \to \infty$. Actually, if we use the continuous-time random walks with suitably chosen exponential waiting times, then one can derive a simple Skorokhod-embedding type coupling between this random walk and the Brownian motion (which in turn, can actually be used to derive all the following results on cable-graphs from results on discrete graphs -- the motivated reader can do this as an instructive exercise!). 

Furthermore: 
\begin{itemize}
\item [-] One can define the analogue of the random walk loop-soup, using the above described Brownian motion on $C$ rather than random walks on $D$. This cable graph loop-soup will be a Poissonian cloud of (unrooted) Brownian loops on $C$.

\item[-] In a cable-graph loop-soup,  there will almost surely be only finitely many loops of diameter greater than $\eps$ for all $\eps$, but 
infinitely many loops of diameter smaller than $\eps$. An important observation is that if one is given a point $x \in C \setminus \partial D$, then the mass of the set 
of loops that pass through $x$ is infinite, so that $x$ will almost surely belong to infinitely many loops in the cable-graph loop-soup. Nevertheless, there will still exist 
{\em exceptional random} points in $C$ that belong to no loop of the cable-graph loop-soup: we denote this set of points by ${\mathcal Z}$. 

\item[-] Each loop $l_j$ in the loop-soup defines an associated occupation time measure density, which is a continuous function on the cable graph that measures ``how much time'' (or rather ``time-density'') it spent at 
each point. More precisely, this is a random continuous function $\omega_j (\cdot)$ on $C$, such that for all connected sets $I\subset C$ the time spent by $l_j$ in $I$ is given by
$\int_I \omega_j (x) dx$, where $dx$ is Lebesgue measure on the cable graph. This is the cable graph version of the local time process for a one-dimensional Brownian motion.

\item[-] One can then define the cumulative occupation time of the loop-soup to be the function 
$$F(x)  := \sum_j \omega_j (x); \;\;\; x\in C.$$ 

\end{itemize}

Not (so) surprisingly, the relation between the discrete loop-soup and 
the square of the GFF generalises nicely: 

\begin {theorem}[From the cable graph loop-soup to square of the GFF] \label{thm:cable_ls_gff_coupling}
The law of the process $F$ is exactly that of a constant times the square of a GFF on the cable graph $C$ (with boundary set $\partial D$). 
\end {theorem}

It is worth stressing that there is no longer any contribution of ``stationary'' loops here, as opposed to the discrete GFF case. In particular, one obtains a construction of the square of the GFF in $D$ directly from the 
cable-graph loop-soup. 
More precisely, the relation with discrete loop-soups in continuous time (as in Section \ref{ctdsloopsoup}) can be described as follows: 

\begin{itemize}
\item[-] The loops in the discrete (continuous-time) loop-soup in $D$ that visit more than one site of $D$ correspond exactly to the loops of the cable-graph loop-soup in $C$ that visit more than one site of $D$. 
\item[-] The Poisson point process of ``stationary'' discrete loops visiting only one point in $D$, that was used to construct $\mathcal{Y}_{1/2}$ in Section \ref{ctdsloopsoup}, 
correspond exactly to the 
loops in the cable-graph loop soup that visit only one point of $D$. 
\end{itemize}

Mind that in the cable-graph GFF, there are also loops (whose corresponding trajectories are contained strictly inside one edge) that visit no point of $D$. These clearly do not contribute to the occupation time density $F$ at sites $x$ of $D$.
\medbreak 

In the construction of the square of the GFF via the loop-soup  (Theorem \ref{thm:cable_ls_gff_coupling}), one can actually interpret the zero-set $Z$ of the GFF (or of its square) in terms of 
the loop-soup. Indeed: 
\begin {proposition} 
Almost surely, the set $Z$ of points at which $\Gamma^2 = 0$ 
is exactly the set ${\mathcal Z}$ of exceptional points on the cable-graph that belong to no cable-graph loop in the loop-soup. In other words, the excursions of $\Gamma$ are exactly the cable-graph loop-soup clusters. 
\end {proposition}
To prove this, one can use the following two features: 

- If $B$ is a one-dimensional Brownian motion, then for all time $t$, the local time 
of $B$ at time $t$ is strictly positive in the interior of its range $B[0,t]$. In other words (this follows for instance from the standard Ray-Knight theorems for one-dimensional Brownian motion, \cite {RY}). 
 It follows that if one considers one Brownian loop on a cable-graph, the local time will be strictly positive for all points that lie at the ``interior'' of its range. In particular, only finitely many points in the range of the Brownian loop will have a local time equal to $0$ (the ``boundary points'' of the range). 

- If one is given any (fixed) point on the cable graph, then almost surely, there will exist infinitely many Brownian loops in the loop-soup that will ``cover'' this point and have a positive local time at it. In particular, if we apply this iteratively (when discovering for instance the loops in the loop-soup one after the other in decreasing size), we see that almost surely, for all the loops in the loop-soup, the boundary points of that loop are ``covered'' by infinitely many smaller loops in the loop-soup. 

Combining these two facts, we see that if the cumulated local time at a point of the cable-graph is equal to $0$, it cannot be in the range of a Brownian loop in the loop-soup (if it was in the range of one of the loops, it would have to be at one of its finitely many ``boundary points'', but then the other loops in the loop-soup would actually cover it), so that $Z \subset {\mathcal Z}$.

Wrapping up, we conclude that we can construct the cable-graph GFF (and therefore the GFF on $D$) starting from a cable-graph loop-soup as follows: 
\begin {corollary}[Lupu's coupling: From the cable-graph loop-soup to the  GFF]
Consider a cable-graph loop-soup, and define its occupation-time density $F$ as before. Then, for each excursion $E_j$ of $F$, toss an independent fair coin to determine $\eps_j\in \{\pm 1\}$. When $x \in E_j$, define $\eps(x):= \eps_j$. 
Then the process $(\eps(x) \sqrt {F(x)} )_{x \in C}$ is a constant multiple of a cable-graph GFF. 
\end {corollary} 

This type of result turns out to be very useful when one wants to understand features of the continuum GFF. In some sense, the cable-graph GFF provides an interpolation between the discrete GFF and the continuum GFF, that
has the advantage (compared to both the discrete GFF and the continuum GFF) of being a continuous function. This means that it possesses additional properties, related to the reflection principle of one-dimensional Brownian motion, that make it very nice to work with.  

In order to illustrate the relation between loop-soups and the GFF on cable-graphs, it can be useful to revisit the definition of a square Bessel processes. This is directly related to the case where the underlying cable-graph is the positive half-line. 

\begin{exercise}For $d$ a positive integer, we define the law of the square Bessel process $X^{(d)}$ of dimension $d$, $\text{SQB}_d$, to be that of $|W_t|^2$, where $W$ is a $d$-dimensional Brownian motion started from the origin. What is the law of $|W_t|^2$ when $d=2$ for a given positive $t$?  
Can you see how to decompose $\text{SQB}_d$ as a sum of $d$ independent processes? 

In fact, although this is not totally obvious (one can use the so-called Yamada-Watanabe theorem, see \cite {RY}), the definition of $\text{SQB}_d$ can be extended to any positive real $d$ 
as a solution to the stochastic differential equation (and this coincides with the definition for $d\in \mathbb{Z}^+$): 
$$ dX_t^{(d)} = 2 \sqrt {X^{(d)}_t} dB_t + d \times dt, \;\;\;\; X^{(d)}_0=0.$$
Show that this process is infinitely divisible -- for instance, the sum of two independent copies of $X^{(d)}$ has the same law as $X^{(2d)}$. 
\end{exercise}

\begin {remark}
Pushing the previous exercise further, one can show $\text{SQB}_1$ can be defined via a Poisson point process of $\text{SQB}_{0+}$ excursions, that correspond exactly to the occupation times of a Poisson point process of Brownian loops in $\R_+$. Interested readers can then revisit the literature on Ray-Knight theorems (see for instance \cite {RY} and the references therein) with this perspective.  
\end {remark}

\section*{Bibliographical comments}

The relation between the number of spanning trees of a graph and the determinant of the Laplacian can of course be traced back all the way to Kirchoff \cite {Kirchoff} (or Tutte's PhD thesis one century later - 1948 - for the oriented 
version, see \cite {Tutte}). In a way, since electrons in a network will perform simple random walks and their loops can be traced in either direction, the basic relations between loop-erased walks and electric currents do also have a long history (see \cite {LyonsPeres} as a reference on electric networks).   
However, even if discrete loop-soups have been around implicitly for a long time, they have (to our knowledge) not been studied as such
until fairly recently. This study has been motivated by features of their continuous counterparts (the Brownian loop-soups introduced in \cite {LawlerWerner}) and the relation of these to SLE processes.

It is interesting and certainly not a mere coincidence that the loop-erased random walk was first introduced by Greg Lawler \cite {LawlerLERW} during his PhD under the supervision of Ed Nelson: one of the 
pioneers of the GFF's Markov property. The direct relation  between LERW and UST was unravelled in successive work by Lawler, Pemantle and then Wilson \cite {Wilson}. Somewhat surprisingly, the simple statement ``the loops erased during Wilson's algorithm can be viewed/rearranged as a sample of a random walk loop-soup'' had not been formalised until fairly recently, even though 
many earlier results were clearly quite closely related to it (see for instance, Wilson's proof with stacks \cite {Wilson}).

Some of the results presented in this chapter can be found in \cite {LJ,LJ2,LawlerLimic,CL,WernerMarkov,LupuWerner} and the references therein. See also the lecture notes by Lawler \cite {LawlerLN}.

The role of cable system loop-soups in order to construct a GFF (and not just its square)  was pointed out by Titus Lupu  \cite {Lupu1} (see also \cite{LupuWerner} for some of the results mentioned here).

\chapter {The continuum  GFF}
\label {Ch3}

We now turn to the continuum world, and we start our study of the {\em continuum GFF}. We will not rely on any of the results that we have derived in the discrete setting, but we will be guided by some 
of the features and intuitions that we have gathered so far.  

In this chapter, we will define and study the continuum GFF in open subsets of $\R^d$, $d \ge 2$. Recall from the warm-up chapter that for open subsets $I$ of $\R$ (i.e., when $d=1$), the 
GFF is nothing else than a collection of independent Brownian bridges in each of the connected components of $I$, which are well-known objects. 
Our goal here is to describe the analogue of the Brownian bridge when the parameter-space is higher dimensional.  
As we progress with this and obtain various results, the reader may find it interesting to draw analogies with (a) the corresponding properties of Brownian bridges and (b)
the corresponding properties of the discrete GFF.

Throughout of this chapter, $D$ will denote an open subset of $\R^d$, $d \ge 2$, satisfying the following properties: 
\begin {itemize}
\item if $d=2$, then $D \not= \R^2$ (together with the condition below, this will ensure that the Green's function in $D$ is finite);
\item if $\partial D \not= \emptyset$, then all boundary points $z\in \partial D$ are regular, meaning that 
for $B$ a $d$-dimensional Brownian motion started from $z$, we have $ \inf \{ t > 0 , \ B_t \notin D \} = 0 $ almost surely. This is a classical condition for existence of solutions to the Dirichlet problem in $D$ and is 
not very restrictive; for instance, it will be satisfied by any domain $D$ with a smooth boundary (it does however rule out domains $D$ such as $\R^d\setminus \{0\}$ for $d\ge 2$).
\end {itemize}
Sometimes, we will add further conditions on $D$, such as requiring it to be bounded or connected.  

\section {Definition of the continuum GFF} 

\subsection {Warm-up and heuristics} 

The object that we would like to define should be some sort of random function, or process, $(\Gamma (x))_{x \in D}$. The process $\Gamma$ should be a centred Gaussian process, and should correspond to the (appropriately normalised) 
limit of the discrete GFF on a lattice approximation $D_{\delta}\subset \delta \Z^d$ to $D$. Recall that the covariance function of the discrete GFF on $D_\delta$ 
(as in Section \ref{sec::informal_sl} of Chapter \ref {Ch1}, the discrete GFF on $D_\delta$ is defined by rescaling the discrete GFF on $\delta^{-1} D_\delta\subset \Z^d$) is the discrete Green's function on $D_\delta$. The only 
way to take $\delta \to 0$ in order to get a limiting process with some non-trivial correlation structure appears to be (as indicated in our warm-up chapter) to first normalise the discrete GFF in such a way that 
the discrete Green's functions converge to a non-trivial function in $D$. This limiting covariance function should still be harmonic away from the diagonal and positive, which essentially characterises it as the {\em continuum Green's function} $G_D$ in $D$. 

So, given that the weak limit of Gaussian processes is a Gaussian process, it looks like we are trying to define a centred Gaussian process $(\Gamma (x))_{x \in D}$ with covariance function 
$E[ \Gamma (x) \Gamma (y) ] = G_D(x, y)$.  As we have already pointed out in the warm-up chapter, this does not appear to be possible, due to the fact that $G_D (x, x) = \infty$.
Formally, this would mean that $\Gamma (x)$ is a Gaussian with infinite variance for every $x$. Note however that there are ways to heuristically interpret Gaussian random variables with infinite variance; for instance, as 
 formal sums of infinitely many independent Gaussian variables with variance $1$. 

In a different direction, if we suppose that $D$ is bounded, and $\Gamma$ has covariance structure as described in the previous paragraph, we could formally consider the ``integral'' $I_\Gamma(1)$ of $\Gamma (x)$ over $D$. Then, by Fubini, we would have 
$$ E [ I_\Gamma(1)^2 ] = \int_{D \times D} dx dy \, E [ \Gamma (x) \Gamma (y) ] = \int_{D \times D} dx dy \, G_D (x,y).$$ 
Now, as we will see in a moment, even if $G_D (x, y)$ explodes as $y \to x$, it is easy to see that for each given $x$, 
$\int_D G_D (x,y) dy$ is finite (one can write this in terms of the expected exit time of $D$ by a Brownian motion started from $x$). Thus the formal variance of $I_\Gamma(1)$ is actually 
finite. So even if for each given $x$, $\Gamma (x)$ does not make sense as a Gaussian random variable, it seems that $I_\Gamma(1)$ should be a Gaussian random variable with finite variance. 

More generally, for any given continuous test function $f$ with compact support in $D$, it turns out that the integral 
$$  G_D (f, f) := \int_{D \times D} dx dy  f(x) f(y) G_D (x, y)$$
is absolutely convergent. 
This in turn indicates that one should be able to define a quantity $I_\Gamma (f)$, that is a centred Gaussian random variable with variance given by $G_D(f,f)$, and can be formally interpreted as $\int_D f(x) \Gamma (x) dx$. 

Finally, if $f_1$ and $f_2$ are two continuous functions with compact support in $D$, then the same argument indicates (formally) that 
$$ E [ I_\Gamma (f_1) I_\Gamma (f_2) ] = \int_{D \times D} dx dy  f_1 (x) f_2 (y) G_D (x, y)=:G_D(f_1,f_2)<\infty.$$ 

In summary, it seems that it should be possible to define a family of random variables  $I_\Gamma (f)$ (indexed by the family of continuous functions $f$ with compact support in $D$) as a centred Gaussian process with covariance function $E[I_\Gamma(f_1)I_\Gamma(f_2)]=G_D (f_1, f_2)$.
\medskip

This formal heuristic conclusion will be the starting point of our definition of the continuum GFF. We will essentially {\em define} the GFF to be this Gaussian process $I_\Gamma$ (with the specified covariance structure). In fact, we will just use the notation $\Gamma (f)$ instead of $I_\Gamma (f)$. In other words, while the value of the continuum GFF at given points will not make sense, quantities that one can interpret as ``mean'' values of the GFF on bounded open domains $U$ (i.e., $\Gamma(\I_U)$) will be well-defined
Gaussian random variables.

\subsection {Basics on stochastic processes} \label{sec:sp}

We now quickly survey some basic results on stochastic processes and measure theory.

\begin{itemize}
\item A random real-valued process indexed by some set ${\mathcal A}$ is 
just a collection of random variables $(X_a)_{ a \in {\mathcal A}}$ defined on the same probability space. The {\em law of the process} is a measure on $\R^{{\mathcal A}}$ (endowed with the product $\sigma$-field)
and is characterised by its finite-dimensional distributions (i.e., 
the law of the finite-dimensional vector $(X(a_1), \ldots, X(a_n))$ for each $a_1, \ldots , a_n\in \mathcal{A}$). 

\item Conversely, if one is given a family of finite-dimensional distributions that is compatible (taking the marginal distribution of one of these distributions gives the correct corresponding finite dimensional distribution), then it is possible (this is Kolmogorov's extension theorem) to construct a probability space and process $(X_a, a \in {\mathcal A})$ on this probability space which has the given finite-dimensional distributions. 

\item When all finite-dimensional distributions are those of (centred) Gaussian vectors, we say that the process is a centred Gaussian process. In other words, a stochastic process $(X_a)_{ a \in {\mathcal A}}$ is a 
centred Gaussian process if and only if for any $n$, for any $a_1, \ldots, a_n$ in ${\mathcal A}$ and any real constants $\lambda_1, \ldots, \lambda_n$, the random variable 
$\lambda_1 X_{a_1} + \cdots + \lambda_n X_{a_n}$ is a centred Gaussian random variable. The law of a centred Gaussian process $(X_a)_{a \in {\mathcal A}}$ is fully described by 
its covariance function $\Sigma (a, a'):= E[ X_{a} X_{a'}]$ defined on $A \times A$. 

\item Combining the previous items shows that when ${\mathcal A}$ is a given set and $\Sigma$ is a real-valued symmetric function defined on $A \times A$ such that for all $n$, for all $a_1, \ldots, a_n$ in ${\mathcal A}$
and all $\lambda_1, \ldots, \lambda_n$ in $\R$, 
$$ \sum_{i, j \le n} \lambda_i \lambda_j \Sigma (a_i, a_j) \ge 0 ,$$
then it is possible to construct a probability space and a process $(X_a)_{a \in {\mathcal A}}$ on this probability space, such that $X$ is a centred Gaussian process with covariance function $\Sigma$. 
\end{itemize}

Note that in the above setting, it is generally not possible a priori to ``simultaneously observe'' more than a countable collection of the variables $X_a$. This is by definition of the product $\sigma$-field. So, a stochastic process does not define a measurable random {\em function} 
from ${\mathcal A}$ into $\R$. 

In order to make sense of certain concepts, for example, appropriate analogues of stopping times, it is useful to work with well-chosen realisations of a given stochastic process. That is, to work with processes having the prescribed 
finite-dimensional marginals, plus some additional features, such as being regular on a set of full probability. Brownian motion is of course a prominent example: one usually defines it as a centred Gaussian process indexed by some interval, and then works with a realisation that is actually a continuous function with probability one.   
More generally (for example, when working with the continuum GFF), it may be possible to construct a version of the process such that $a \mapsto X_a$ 
 {\em on some given subset of ${\mathcal A}$} is continuous on a set of full probability. 

We will see examples of this, in the context of the GFF, in Section \ref{sec:kol}.

\subsection {Basics on the continuum Green's function}

In order not to disrupt the flow of the presentation here, we quickly state without proofs some properties of the continuum Green's function. We will then provide a somewhat self-contained presentation, including proofs of these facts, in the next section. 

\emph{Suppose that $D \subset \R^d$ satisfies the conditions that we stated at the beginning of this chapter, and that it is connected.}

Let $y\in D$.
Then it is easy to see that, up to a multiplicative constant, there exists only one positive 
harmonic function $\tilde H_y$ in $D \setminus \{ y \}$ such that $\tilde H_y (x)$ tends to $0$ as $x \to \partial D$ or $x \to \infty$. For this function not to be identically zero, it has to tend to infinity when $x \to y$, and 
the way it does so depends on the dimension: when $d = 2$ it will explode like a constant times $\log (1 / |x-y|)$; and when $d >2$ it will explode like a constant times $|x-y|^{2-d}$. 
So for example, when $d \ge 3$ and $D = \R^d$, the function $\tilde H_y$ is in fact equal to a constant times $|x-y|^{2-d}$.  
The Green's function $G_D (x, y)$ is then defined to be this function $\tilde H_y (x)$, where the multiplicative constant is chosen in an appropriate way (in fact it is chosen to be equal to $a_d^{-1}$, where $a_d$ is the $(d-1)$-dimensional Lebesgue measure 
of the unit $(d-1)$-dimensional sphere in $\R^d$).  

This function then turns out to have all the properties that one would expect from the continuum analogue of the discrete Green's function. In particular (we will provide more details in the next section):  
\begin {enumerate} 
 \item it is a symmetric function, i.e., $G_D (x, y)= G_D(y,x)$;
 \item it can be interpreted as an integral operator that turns out to be the inverse of $- \Delta$ (or rather of a multiple of $- \Delta$), for $\Delta$ the continuum Laplacian;
 \item it can be interpreted in terms of expected occupation times by Brownian motion stopped when exiting $D$. 
\end {enumerate}

We define ${\mathcal M}_D^+$ to be the set of finite measures that are supported in $D$, and such that 
$$ \int_{D\times D} G_D (x,y) d\mu (x) d \mu (y) < \infty.$$
We also define ${\mathcal M}_D$ to be the vector space of signed measures $\mu^+ - \mu^-$, where $\mu^+$ and $\mu^-$ are in ${\mathcal M}^+_D$. We will omit the subscript $D$ and simply 
write $\{\mathcal M,\mathcal M^+\}$, unless we are discussing various domains simultaneously. 

One class ${\mathcal M}_c$ of measures $\mu\in \mathcal M$ is given by the set of measures of the form $f(x) dx$, when $f$ is continuous with compact support in $D$ (and $dx$ denotes the Lebesgue measure). 
It is easy to check that such measures do lie in $\mathcal{M}$, simply by considering the rate at which $G_D(x,y)$ explodes as $|x-y|\to 0$.

\subsection {Definition} 

We are now ready to define the continuum GFF.
 As in the discrete case, there are several possible ways to do this. We choose here to first define it as a random process. This will be quite useful, for instance, when we consider local sets in the next chapter.

\begin {definition}[Continuum GFF]
\label {GFFdef}
We say that the process $(\Gamma (\mu))_{ \mu \in {\mathcal M}}$ is a Gaussian Free Field in $D$ if it is a centred Gaussian process with covariance function 
$$ \Sigma (\mu, \nu)  := \int_{D \times D} G_D (x,y) d \mu (x) d \nu (y).
$$
\end {definition} 

In order to check that this definition makes sense, it suffices to check that 
this  function $\Sigma (\mu, \nu)$ is indeed a well-defined covariance function, i.e., that  for any $\mu_1, \ldots, \mu_n$ in ${\mathcal M}$ and any real $\lambda_1, \ldots,\lambda_n$, 
$$ \sum_{i, j \le n} \lambda_i \lambda_j \Sigma (\mu_i, \mu_j )  \ge 0.$$ 
Note that the left-hand side is equal to 
$ \Sigma (\mu, \mu) $ 
for $\mu = \lambda_1 \mu_1 + \ldots + \lambda_n \mu_n$, and the fact that this quantity is non-negative will 
 follow from elementary features of the Green's function  (see (\ref {appro}) in the next section). So, the GFF in $D$ does indeed exist. 

When $f$ is a measurable function in $D$ such that 
$$ \int_{D \times D} | f(x) f(y) | G_D (x, y) dx dy < \infty,$$ 
then $\mu_f = f(x) dx$ lies in $\mathcal{M}$, and we will often write $\Gamma (f)$ as a shorthand notation for $\Gamma(\mu_f)$.

We can immediately note that for all $\lambda \in \R$ and all $\mu$ and $\nu$ in ${\mathcal M}$, one has 
$$ \Gamma (\lambda \mu) = \lambda \Gamma (\mu) \hbox { and }
\Gamma (\mu + \nu) = \Gamma (\mu) + \Gamma (\nu) \text{ almost surely }$$
(by simply noticing that in both cases the second moment of the differences between left and right-hand sides vanish).  
It readily follows that the law of the GFF is characterised by this linearity relation 
and the fact that for each given $\mu$ in ${\mathcal M}$, $\Gamma (\mu)$ is a centred Gaussian random variable with variance  $\Sigma (\mu, \mu)$.
 
Let us stress once again that when dealing with processes with uncountable index sets (such as $\mathcal{M}$) one has to pay close attention to the positioning of ``for all $\mu\in \mathcal{M}$'' or ``for each given $\mu\in \mathcal{M}$'' in statements such as the above. As an illustration, we would like to 
mention straight away that for general {\em random} measures $\mu \in {\mathcal M}$ (coupled to the GFF), 
$\Gamma (\mu)$ is not necessarily a well-defined random variable (i.e. it is not necessarily measurable), and may not make any sense at all. This contrasts with Brownian motion, where 
$B_t$ can be defined for all $t$ simultaneously, because one can choose it to be a continuous function). So, the GFF 
cannot be viewed as a random {\em function} from ${\mathcal M}$ into $\R$: this would mean that one is able to define  $\Gamma (\mu)$ for all $\mu \in {\mathcal M}$ ``simultaneously'', which turns out 
not to be possible. On the other hand, we will be able to define it simultaneously for all $\mu$ in certain nice subsets of $\mathcal{M}$.

\subsection {Other boundary conditions} 

Suppose now that $D$ is bounded and that 
$H$ is a given harmonic function in $D$. As opposed to the discrete case where harmonic functions were necessarily bounded in the neighbourhood of the boundary points, and (in the case where $D$ is bounded) are the harmonic extensions of their finite values on $\partial D$, in this continuum case, the function $H$ may be unbounded near $\partial D$ (this will typically happen when we will discuss the Markov property of the continuum GFF). However, the harmonic function is fully determined by its value on any neighbourhood of $\partial D$, so that one can (at least informally) think of it as the harmonic extension of ``its trace on $\partial D$''. In other words, we can view the information of a ``boundary conditions'' as the same information as the knowledge of the entire harmonic function. 
This (and of course the corresponding features of the discrete GFF) does lead to the following definition: 

\begin {definition}[GFF with non-constant boundary conditions] 
\label {GFFnonconstant}
We say that $\hat \Gamma$ is a GFF in $D$ with boundary conditions given by $H$ if $\hat \Gamma = H + \Gamma$, where $\Gamma$ is a Dirichlet GFF in $D$.
\end {definition}

The equation $\hat \Gamma = H + \Gamma$ should be understood in the sense that 
$$ \hat \Gamma (\mu) = \int H(x) \mu (dx) + \Gamma (\mu).$$
If $H$ is unbounded in the neighbourhood of $\partial D$, one can restrict the definition to the set of measures $\mu$ in ${\mathcal M}$ with compact support in $D$ to be on the safe side (in order to be sure that $\int H(x) \mu (dx)$ is well-defined).

\begin {remark}
We will briefly discuss other type of boundary conditions (Neumann, or periodic) in Chapter~\ref{Ch6}. 
\end {remark}

\section {A closer look at the continuum Green's function}     \label{sec::greensfunction}

In order to help readers who are not so familiar with potential theory, 
let us now quickly review a few properties of the Green's function in a $d$-dimensional domain $D$ (satisfying the conditions described at the beginning of this chapter). Of course, the Green's function can be defined and studied in many different 
ways, and the following approach is just one of many possibilities. We will mostly try to highlight ideas, and therefore leave out some of the classical details as exercises. 

\emph{For presentation purposes only, we will assume in addition to the usual hypotheses, that $D$ is a bounded and connected domain of $\R^d$.}

 However, it is 
easy to adapt all the results in this section to unbounded and/or non-connected domains (except those related to the 
eigenfunction decomposition of the Laplacian, see Remark \ref{EV}). We will briefly comment on this later (see point (vii) at the end of this section). Of course $D=\R^d$ for $d \ge 3$ (where $G_D$ is defined to be $a_d^{-1} |x-y|^{2-d}$) is an important example of a domain that satisfies the usual hypotheses but is not bounded. 

Some basic terminology first. We will say that a function $f$ is harmonic on an open subset $O$ of $\R^d$ if it is continuous and if for any closed ball
$\overline {B}(z, \epsilon):=\{x\in \R^d: |x-z|<\eps\}$ contained in $O$, the mean value $\overline f^\epsilon (z)$ of $f$ on the boundary $\partial \overline{B}(z,\epsilon)$ is equal to $f(z)$.

\begin{exercise}
	Show that the above definition of harmonicity is equivalent to the fact that $f$ is smooth with $\Delta f=0$ in $O$.
\end{exercise} 

We will denote by $P^x$ (and $E^x$) the probability measure (and the corresponding expectation) under which $B$ is a $d$-dimensional Brownian motion started from 
$x$. For any open set $O$, we write $\tau_O$ for the exit time of $O$ by $B$, and simply write $\tau$ for $\tau_D$. 

For each fixed $y \in D$, the Green's function $x \mapsto G_D (x, y)$ is going to be a positive harmonic 
function in $D \setminus \{ y \}$ that vanishes on $\partial D$. As we shall see, this in fact already characterises the function up to a multiplicative constant (when $D$ is unbounded, one would have to add the condition that 
the function tends to $0$ as $x \to \infty$ in $D$). 

Let us now go about constructing such a function. It is easy to check that for any $y\in \R^2$, the function $x \mapsto \log (1/|x-y|)$ is harmonic in $\R^2 \setminus \{ y \}$, and that for any $y\in \R^d$ with $d>2$, the function $x \mapsto |x-y|^{2-d}$ is harmonic in $\R^d \setminus \{ y \}$. We use the general notation $H_y$ for this function \emph{divided by} the $(d-1)$-dimensional measure $a_d$ of the boundary of the $d$-dimensional unit ball.  For instance, when $d=2$, we have
$$ H_y (x) := \frac 1 {2\pi} \times \log \frac 1 {|x-y|}, $$ 
and when $d=3$, 
$$ H_y (x) := \frac 1 {4 \pi} \times \frac 1 {|x-y|}.$$  
The choice of normalisation (dividing by $a_d$) will ensure that the Laplacian of $-H_y$ can be viewed as the Dirac mass at $y$ (and consequently that $G_D$ will be the inverse of $- \Delta$: see Lemma \ref{lap_green_delta}).

However, the function $x\mapsto H_y(x)$ is not yet the correct choice for the Green's function, because it does not vanish on $\partial D$. Nonetheless, we can simply subtract from $H_y$ the harmonic function in $D$ with the same boundary values as $H_y$ on $\partial D$.  More precisely, we can define 
for each $x, y \in D$, 
$$ h_{y,D} (x) := E^x [ H_y (B_\tau) ],$$ 
which is the unique solution to the Dirichlet problem in $D$ with boundary conditions $H_y$ on $\partial D$. This function 
is continuous on $\overline D$ and equal to $H_y$ on $\partial D$, by the assumption that $\partial D$ is regular. 
\begin{exercise}\label{ex:harm}
	Suppose that $f$ is a continuous function with compact support in $D$. Show that the function $x \mapsto \int_D f(y) h_{y,D} (x) dy$ is bounded and continuous on $\overline D$, and harmonic in $D$.
\end{exercise}

Then we define the Green's function as follows.

\begin {definition}[Green's function] 
For $x \not= y$ in $D$, we set
$$
 G_D (x,y) :=  H_y (x) - h_{y,D} (x).
 $$
\end {definition} 

Alternatively, the Green's function can be characterised in the following manner. 

\begin {lemma}
\label {charactlemma}
For each given $y$ in $D$, the function $x \mapsto G_D (x, y)$ is the unique continuous function defined on $\overline D \setminus \{ y \}$ such that: 
(i) it is equal to $0$ on $\partial D$, (ii) it is harmonic in $D \setminus \{ y \}$, and (iii) the function $ x \mapsto G_D (x,y) - H_y (x)$ remains bounded in a neighbourhood of $y$.   
\end {lemma}
\begin {proof} 
The function $x \mapsto G_D (x, y)$ clearly satisfies these three conditions, and if $g(x)$ is another function satisfying these conditions, then $F(x):=g(x)-G_D (x,y)$ is 
harmonic on $D \setminus \{ y \}$, bounded on $\overline{D}$ and vanishes on $\partial D$. It is easy to see that such a function must be identically $0$ (for example, one can show that $F(B_t)$ is a martingale when $B$ is a Brownian motion started from $x\in D$, and then apply the optional stopping theorem at time $\tau_D$).
\end {proof}

Just as in the discrete case, one can interpret
$G_D$ as the inverse of $-\Delta $, when we define the continuum Laplacian $\Delta$ of a twice differentiable function $f$ by
$$ \Delta f (x) = \sum_{j=1}^d \frac {\partial^2 f}{\partial x_j^2}(x)$$
in the standard way. 
 Before explaining this, let us make a very brief side-remark.
On discrete regular lattices such as $\Z^d$, 
we chose to define the discrete Laplacian of a function $f$  to be $\overline f - f$, where $\overline f(x)$ was the mean-value of $f$ on the neighbours of $x$. However, we cautioned
that this was not really the standard definition of the discrete Laplacian; in the literature it is often instead defined by the operator $f\mapsto 2d(\overline f - f)$.
 The connection between the discrete Laplacian and the continuum Laplacian $\Delta$ goes as follows: consider the difference $\overline f^\eps- f$ where $\overline f^\eps(x)$ is the average of $f$ on $\partial \overline B(x,\eps)$ as before. 
 Then it follows from Taylor's expansion of $f$ that 
\begin{equation*} \overline f^\eps (x) - f(x) =\frac { \eps^2}{2d} \Delta f(x) + o (\eps^2). \end{equation*} as $\eps\to 0$.
That is, it is actually $\Delta f(x) / (2d) $ that describes the behaviour of $(\overline f^\eps (x)-f(x))$ as $\eps \to 0$.

\vspace{0.1cm} 
Let us now move on to the connection between $G_D$ and $-\Delta$.

\begin {lemma}[The Green's function is the inverse of $-\Delta$] \label{lap_green_delta}
When $f$ is a continuous function with compact support in $D$, we define 
$$ F(x) = \int_D f(y) G_D (x, y) dy=: G_D(f)(x).$$
Then $F$ is continuous on $\overline{D}$, smooth in $D$, vanishes on $\partial D$ and satisfies $ - \Delta F = f$.
\end {lemma} 

\begin {proof} 
Given Exercise \ref{ex:harm}, it remains only to check 
that $x \mapsto \int_D f(y) H_y(x) dy$ is continuous on $\overline D$ and that its Laplacian is equal to $-f$. We leave this as a further exercise for the reader.  
\end {proof}

Note that this also shows the inequality
\begin {equation}
 \label {posit}
\int_{D \times D}  f(x) G_D (x, y) f(y) dx dy
= - \int_D F(x) \Delta F (x) dx   = 
 \int_D | \nabla F (x) |^2 dx \ge 0
\end {equation}
for any such function $f$. 
Also, for any two such functions $f_1$ and $ f_2$ 
$$ 
\int_{D \times D}  f_1(x) G_D (x, y) f_2(y) dx dy
=   \int_D (\nabla F_1 (x)  \cdot \nabla  F_2 (x))  dx ,$$ 
which provides one way, among many, to see that $G_D (x,y) = G_D (y,x)$. 
\medbreak

\begin {remark}[Dirichlet Laplacian eigenfunction decomposition]
\label {EV}
Recall that when $D$ is bounded, it is possible to find an orthonormal basis $(\varphi_j)_{j \ge 1}$ of $L^2(D)$ consisting of eigenfunctions of $-\Delta$
that vanish on $\partial D$ (these are sometimes referred to as the eigenfunctions of the Dirichlet Laplacian). 
We denote by $(\lambda_j)_{j \ge 1}$ 
the corresponding eigenvalues (that are all positive), so $- \Delta \varphi_j  = \lambda_j \varphi_j$. The existence of such a basis actually follows immediately from the fact that when $\varphi$ is an eigenfunction of the Laplacian, the space of functions that 
are orthogonal to $\varphi$ (with respect to the $L^2$ inner product) is stable under $-\Delta$. 

Decomposing according to this orthonormal basis, one can define $f_j = \int_D f(x) \varphi_j (x) dx$ for any $f\in L^2$, and then one has $f(\cdot) = \sum_{j \ge 1} f_j \varphi_j (\cdot)$, where the infinite sum is viewed as a limit in $L^2$.  One also has, by Weyl's law, control on the asymptotic behaviour as $\lambda\to \infty$ of the number of eigenvalues $\lambda_j$ that are smaller than $\lambda$ (we will come back to this later). 

Given that $G_D$ can be viewed as the inverse of $- \Delta$, we can therefore also use this basis of eigenfunctions to describe $G_D$. Indeed, for all $x$ and $y$ in $D$, one has
\begin {equation} \label{Geig}
G_D (x, y) = \sum_{j \ge 1} \frac {1}{\lambda_j} \varphi_j (x) \varphi_j (y), 
\end {equation}
where the sum on the right-hand side is a convergent series in $L^2 (D \times D)$. Then, one has 
$$ \int_D \varphi_i (x) G_D (x, y) dx = \sum_{j \ge 1} \lambda_j^{-1} \varphi_j (y) \int_{D} \varphi_i (x) \varphi_j (x) dx = \lambda_j^{-1} \varphi_j (y), $$ 
so that $\varphi_j$ can be viewed as an eigenfunction of the operator $f\mapsto G_D(f)$, with associated eigenvalue $1 / \lambda_j$. 

We see that the two natural Hilbert spaces to consider here are on the one hand the Sobolev space of functions $f (\cdot) = \sum \lambda_j^{-1} f_j \varphi_j(\cdot)$ with $\sum f_j^2 < \infty$ (which is also the closure of the set of smooth functions that vanish on the boundary, with respect to the norm 
$\int | \nabla f (x) |^2 dx$) and on the other hand the space $L^2 (D)$. 
Then, $-\Delta$ is a bijection from the former onto the latter, and its inverse $G_D$ is a bijection from the latter onto the former. 
\end {remark}

Recall that the value of the discrete Green's function at $(x,y)$ had a natural interpretation in terms of the expected number of visits to $y$ by simple random walk started at $x$. The continuum counterpart 
of this interpretation will relate the continuum Green's function to expected occupation times of Brownian motion. Before stating this, let us again take a few lines to clarify how normalisations differ between the discrete and continuum cases (i.e., in the end,  an extra factor of $2$ will appear in the continuum setting): 

Since the co-ordinates of
a $d$-dimensional Brownian motion $(B_t)_{t \ge 0}$ started from the origin are $d$ independent standard one-dimensional Brownian motions, we have 
$E [ \| B_t \|^2 ] = d \times t$. 
Moreover, using the classical fact that $\| B_t \|^2 - d \times t$ is actually a martingale and that the exit time $\tau_\eps$ of $B(0,\eps)$ by $B$ is almost surely finite, one can deduce via the optional stopping theorem that $E [ \tau_\eps] = { \eps^2 }  / d$. Note that this $1/d$ term does not appear for a simple random walk $X_n$ on a lattice that spends time $1$ between each jump: this is because only one co-ordinate moves at each time for $X_n$, so we simply have $E(\| X_n \|^2)=n$. 
On the other hand, we noticed just a few paragraphs above that the continuum Laplacian $\Delta$ really comes with an additional factor of $2d$ compared to our definition of the discrete Laplacian. The outcome of these two differences will that  (compared to the discrete case) there will be a factor of $2$ appearing in the relationship between the continuum Green's function and occupation times of Brownian motion.

Let us now actually derive this relationship. When $f$ is smooth with compact support in $D$, we define the function
$$ J (x) = J_f (x) :=   \frac{1}{2} E^x \Bigl[ \int_0^\tau f(B_t) dt \Bigr].$$ 
It is easy to see that this function is continuous and that it vanishes on $\partial D$. Furthermore, by the strong Markov property we have 
$$ \overline {J}^\eps (x) - J (x) = - \frac {1}{2} E^x \Bigl[ \int_0^{\tau_\eps} f(B_t) dt \Bigr], $$
and by the previous estimates this behaves like $-f(x) \times \eps^2 / (2d)$ as $\eps \to 0$. 
On the other hand, it is not difficult to check that $J$ is smooth, so that $\overline {J}^\eps (z) - (J)(z) \sim \eps^2 \Delta (J) (z)/ (2d)$ as $\eps \to 0$.
We can therefore conclude that 
$$ \Delta J = -f,$$ 
and by applying the maximum principle to $J - G_D(f)$, that we must have $J= G_D (f)$.

By dominated convergence, approximating a bounded, measurable function by smooth functions), the same result will hold true for the following class of functions $f$: 
\begin {lemma}[Green's function as an expected occupation time of Brownian motion]
For all bounded, measurable functions $f$ with compact support in $D$, 
$$\int_D f(y) G_D (x, y) dy = \frac 1 2  E^x \Bigl[ \int_0^\tau f(B_t) dt    \Bigr].$$ 
\end {lemma} 
In particular, when $f$ is the indicator function of an open subset $A$ of $D$, we obtain that 
$$ \frac 1 2 E^x \Bigl[  \int_0^\tau \1{B_t \in A} dt  \Bigr]  = \int_A G_D (x, y) dy .$$
Let us stress the normalising factor $2$ appearing here (that does not appear in the discrete setting). 
In some sense, the more natural operator associated to Brownian motion in the continuum is $\Delta/2$ (rather than $\Delta$, which is the inverse of $G_D$).

\begin {remark}\label{G_eig}
It is also possible and not very difficult to see that for each $t$, 
the law of $B_{\min (t, \tau)}$ has a density $p_{D,t} (x,y)$ in $D$ that is jointly continuous for $(t,x,y)\in(0,\infty) \times D \times D$. 
By Fubini, we then see that for any bounded measurable $f$: 
\begin {eqnarray*}
 \lefteqn{ 
E^x \Bigl[ \int_0^\tau f(B_t) dt    \Bigr] = \int_0^\infty E^x [f(B_t) \1{t < \tau} ] dt}
\\
&& = 
\int_0^\infty \int_D p_{D,t} (x, y) f(y) dy dt = \int_D (\int_0^\infty p_{D,t}(x,y) dt)  f(y) dy.
\end {eqnarray*}
This gives us another useful expression  
$$ G_D (x,y) = \frac 1 2  \int_0^\infty p_{D,t} (x,y) dt $$
for the Green's function.

We can then relate this to the discussion about the eigenfunction decomposition of $- \Delta$ when $D$ is bounded. 
Indeed, for each $y$, $(x,t) \mapsto p_t (x, y) := P_t(x)$ solves the heat equation 
$\partial_t P = \Delta P / 2$,  and it can be written (for $t >0$) as 
 $$ p_t (x, y) =  \sum_j e^{- \lambda_j t /2} \varphi_j (x) \varphi_j (y).$$ 
 Integrating this with respect to $t$ gives 
 $$ \frac 1 2  \int_0^\infty p_{D,t} (x,y) dt = 
 \sum_{j} \lambda_j^{-1} \varphi_j (x) \varphi_j (y), $$ 
 which is consistent with the expression \eqref{Geig} that we saw previously.
 \end {remark}

\noindent Let us conclude this subsection with some final remarks on the Green's function. 

\medbreak (i)
Suppose that $d=2$ and $D\subset \R^2$ is a domain satisfying our usual conditions. Then it is a rather simple exercise to check  that discrete Green's functions in suitable fine-mesh approximations to $D$ converge to the continuous Green's function $G_D$ as the mesh-size goes to $0$. One approach to this is to note that 
the appropriate discrete random walks converge to Brownian motion in a strong sense, and to then use (minding the two factors of $2$ that appear) the description of the Green's functions as occupation time densities. 

\medbreak (ii) Let $d\ge 2$ again be general, and suppose that $O$ is an open subset of $D$ such that all boundary points of $O$ are regular. Define for any $x \not= y$ in $D$, 
$$ H_{D,O} (x,y) = G_D (x,y) - G_{O} (x,y)$$
Note that this function is equal to $G_D (x,y)$ as soon as either $x$ or $y$ are not in $O$, and that it can be extended continuously to the points $\{(x,x): x\in O \}$ (by part (iii) of Lemma \ref{charactlemma}). Moreover, when $x \in O$, the function 
$ y \mapsto H_{D,O} (x,y)$ is harmonic in $O\setminus \{ x \}$ and continuous at $x$, so therefore also harmonic at $x$. Hence $ y \mapsto H_{D,O} (x,y)$ is the (unique) harmonic extension in $O$ of the function that is equal to $G_D (x,y)$ for all $y \in \partial O$. 

We further remark that for $A \subset O$, one has 
$$
\int_A H_{D, O} (x, y)  dy = E^x \Bigl[ \int_{\tau'}^\tau \1{B_t \in A} dt  \Bigr]
$$ where $\tau'$ is the exit time of $O$ by $B$,
so that one can interpret the quantity $H_{D, O} (x,y)$ as the density of the cumulative occupation time in a neighbourhood of $y$, 
for a Brownian motion started from $x$ and restricted to those times between $\tau'$ and $\tau$. 

\medbreak (iii) With the same notation as in (ii), when $x\in O$ we define  $\nu_{x, \partial O}$ to be  the law of $B_{\tau'}$. This measure is often called the harmonic measure on $\partial O$ (seen from $x$) and is an element of ${\mathcal M}_D$ for any $x$. To see this, observe that for all {$z \in \partial O$} the function {$z' \mapsto G_D( z, z')$} is harmonic in $O$, from which it follows that 
$$ \int_{\partial O} d\nu_{x, \partial O} (z') G_D (z, z') = G_{D} (x, z),$$
and so $\nu_{x, \partial O} \in \mathcal{M}_D$ (because $\partial O$ is at positive distance from $x$). 

A similar argument shows that 
for all $x \not= y \in O$,
\begin {eqnarray*}
H_{D, O} (x,y) &=& G_D (x,y) - G_{O} (x,y) = \int d\nu_{x, \partial O} (dz) d\nu_{y, \partial O} (dz') (G_D (z,z')- G_{O}(z,z')) 
\\
&=& \int d\nu_{x, \partial O} (dz) d\nu_{y, \partial O} (dz') G_D (z,z') .
\end {eqnarray*}

\medbreak (iv) Now let us comment on the positive definiteness of $G_D$. We note (and leave as a simple exercise to check) that if $\mu \in {\mathcal M}$,
and if for all $\epsilon >0$, one defines $\mu^\epsilon$ to be the convolution of $\mu$ with a well-chosen smooth, positive test function supported on the ball of radius $\epsilon$ around the origin (for example, if $\varphi$ is a smooth positive function supported in the unit ball of $\R^d$ and with $\int_{\R^d} \varphi(x) dx =1$, then we could consider $\mu^\eps = \mu * \varphi^\eps$ where $\varphi^\eps(x)=\frac{1}{\eps^d}\varphi(\frac{x}{\eps^d})$),
then by dominated convergence one can conclude that 
$$
\int_{D \times D} G_D (x,y) d\mu^\epsilon (x) d \mu^\epsilon (y) 
\to 
\int_{D \times D} G_D (x,y) d\mu (x) d \mu(y).
$$
In particular, in view of (\ref {posit}) - recall that this held for continuous functions, we see that
\begin {equation} 
 \label {appro} 
\int_{D \times D} G_D (x,y) d\mu (x) d \mu(y) \ge 0
\end {equation}
and similarly
$$
 \int_D G_D (x,y) (d\mu^\epsilon (x)- d \mu (x)) (d \mu^\epsilon (y) - d \mu (y)) \to 0$$
 as $\eps\to 0$. 
Recall that (\ref {appro}) is what ensures that the covariance structure of the GFF is licit. 

\medbreak (v)
We note that if $D_n$ is an increasing sequence of open sets such that $\cup D_n = D$, then for all $x \not= y$ in $D$ (by applying dominated convergence and using that $\tau_{D_n}\to \tau$ almost surely for a Brownian motion started at $x$) we have 
$$ G_{D_n} (x,y) \to G_D (x,y) $$ as $n \to \infty$. 
For instance, and we shall be using the following notation throughout this section, this will hold if we define $A_n$ to be the union of all closed $2^{-n}$-dyadic hypercubes $$S_{j_1,\cdots, j_d}^n:= [ j_12^{-n}, (j_1+1) 2^{-n}] \times ... \times  [j_d 2^{-n}, (j_d+1) 2^{-n} ]\; , \;\;\;\;( j_1,\cdots, j_d)\in \mathbb{Z}^d$$ that intersect the complement of $D$, and then set $D_n := D \setminus A_n$.

\medbreak (vi)
Moreover, the same result as (iv)  holds if for every $n$ we subdivide $D$ into dyadic hypercubes of the 
type $S^n_{j_1,\cdots, j_d}$, 
and define the measure $\mu_n= 
 \sum_{j_1,\cdots, j_d} \mu (S^n_{j_1,\dots, j_d}) \lambda_{j_1,\cdots, j_d}$, for $\lambda_{j_1,\cdots, j_d}$ the uniform (i.e. multiple of Lebesgue) probability measure on $S^n_{j_1,\cdots, j_d}$. 
Then, as before, it is easy to check that 
$$ \int_D G_D (x,y) (d\mu_n (x)- d \mu (x)) (d \mu_n (y) - d \mu (y))
\to 0  $$
as $n\to \infty$.

\medbreak (vii) As mentioned at the beginning of this section, the definition of Green's functions and the derivation of most their properties do not actually require the additional assumptions of connectedness or boundedness on $D$. 
For non-connected domains $D$, the definition is trivial: one decomposes $D$ into its connected components $(D_i)_{i\in I}$, and sets $G_D (x,y)= G_{D_i}(x,y)$ when $x$ and $y$ are in the same 
connected component $D_i$; $G_D(x,y)=0$ otherwise.

For unbounded connected domains $D$ (with the condition that $D \not= \R^2$ if $d = 2$, and that the boundary of $D$ is regular), all results can be easily extended, except the ones 
that involve the spectrum of the Laplacian. One just has to be a little careful, and slightly modify some of the statements. For instance, one needs to add an extra condition on the behaviour at infinity when characterising the function $x \mapsto G_D (x,y)$ as the only harmonic function in $D$ with certain properties (as in Lemma \ref {charactlemma}).

Recall that if $D = \R^d$ for $d \ge 3$, then the Green's function is given by
$$ G_{\R^d} (x,y) = \frac 1 {a_d |x-y|^{d-2}}.$$ 

\medbreak (viii) 
The Green's functions have a simple scaling property: when $\lambda > 0$ we can relate the Green's function in $D$ with the Green's function in $\lambda D$ for any $D$, by 
$$ G_{\lambda D} (\lambda x, \lambda y) = \lambda^{2-d} G_D (x,y ).$$ 
Note that $d=2$ plays a special role here.

\medbreak (ix) 
Given the above, let us now restrict to the case $d=2$. Take $D\subsetneq \R^2$ and consider a conformal transformation $\Phi$ from $D$ to $\Phi(D)$. Then the characterisation of $G_D ( x ,\cdot)$ 
as the unique harmonic function in $D \setminus \{ x \}$ with prescribed boundary conditions immediately implies that 
$$ G_{\Phi (D)}  ( \Phi (x), \Phi (y)) =  G_D (x,y)$$ for all $x\ne y$ in $D$.
Thus, one can obtain expressions for $G_D$ with $D$ arbitrary if we know the Green's function in some reference domain. 

One can for instance recall or check that 
$$G_\U(x,y)= \frac 1 {2 \pi} \log \frac {|1 - x   \overline {y}|} {|y-x|}$$ 
for the unit disc $\U\subset \C=\R^2$, 
and that  the Green's function in the upper half-plane $\HH= \{ z \ : \ \Im (z) > 0  \}\subset \C$ is given by 
$$ G_{\HH} (x, y) = \frac 1 {2 \pi} \Bigl( \log \frac 1 {|x-y|}  - \log \frac 1 {|x- \overline y|}  \Bigr). $$
 
Suppose finally that $D\subset \R^2$ is simply connected, and that $\Phi$ is the unique conformal transformation of the unit disc onto $D$ with $\Phi ( 0) = x$ and $\Phi' (0) \in \R_+$. This derivative $\Phi'(0)$ is sometimes called the conformal radius of $D$ at $x$. Then, using the explicit expression for $G_\U$, we get that as $x \to y$,  
$$G_D (x,y) = \frac 1 {2 \pi} \log \frac 1  {|x-y|}  + \frac 1 {2 \pi}  \log \Phi' (0) + o (1).$$
More generally, when $d > 2$, the coefficient of the constant term in the Laurent expansion of $G_D( x,y)$ as $x \to y$ provides some information about how ``close'' $y$ is to the boundary of $D$.

\section {First comments on the regularity of the GFF}

After having collected all these basic facts about the Green's function, we can proceed to study the GFF as introduced in Definition \ref {GFFdef}. 
\emph{In this section, and the rest of the chapter, we continue to assume that $D$ is a bounded and connected domain of $\R^d$ with regular boundary (and $D\ne \R^2$ when $d=2$).}
\subsection {Approximation via mean values on dyadic hypercubes}

 One standard way to construct a stochastic process $(X_a)_{a \in {\mathcal A}}$ indexed by a large set ${\mathcal A}$ is to first define the random variables $X_{a'}$ for all $a'$ in a countable 
subset of ${\mathcal A}$, and then for each individual $a \in {\mathcal A}$, to define $X_a$ as a limit of some sequence $X_{a_n}$, where $a_n$ is a well-chosen sequence in the countable set. 

In the present case, one natural countable ``dense'' subset of ${\mathcal M}$ to consider is the set of measures 
$\mu_{j_1,\cdots, j_d}^n$ that are the uniform distributions on the 
dyadic hypercubes $S^n_{j_1,\cdots, j_d}$ contained in $D$. 

Suppose that one knows the countable collection of 
random variables $(\Gamma (\mu_{j_1,\cdots, j_d}^n))$. 
Then using Remark (vi) of the previous section, 
we know that given any $\mu \in \mathcal{M}$ it is possible to find a sequence $\mu_n$ of linear combinations of the $\mu_{j_1,\cdots, j_d}^n$'s 
so that $ \Gamma (\mu_n)$ converges in $L^2$ to $\Gamma ( \mu)$
(and therefore choosing some appropriate deterministic subsequence, converges almost surely). 
Hence, the knowledge of all these $(\Gamma (\mu_{j_1,\cdots, j_d}^n))$'s enables one to recover each $\Gamma (\mu)$ individually 
(i.e., for each $\mu\in \mathcal{M}$, one can almost surely recover $\Gamma (\mu)$).

In the sequel, an important role will be played by the $\sigma$-field ${\mathcal F}_A$ generated by all the random variables $\Gamma (\mu)$, where $\mu$ ranges over elements of $ {\mathcal M}$ that are supported in some compact set $A$.

\subsection {The GFF as a random Fourier series, the GFF as a random generalised function} 
\label {S.Fourier}
Recall in the case where $D$ is bounded (see Remark \ref {EV}) the existence of an orthonormal basis $(\varphi_j)_{j \ge 1}$ of $L^2(D)$ that 
consists of the eigenfunctions of $-\Delta $ that vanish on the boundary of $D$. We denote the associated eigenvalues by $(\lambda_j)_{j \ge 1}$.

Weyl's law tells us that the number $N(\lambda)$ of eigenvalues smaller than $\lambda$ satisfies 
$$\lim_{\lambda\to \infty} \frac{N(\lambda)}{\lambda^{d/2}}=c_d\text{vol}(D)$$
for some finite dimension-dependent constant $c_d$. 

Also recall that (if $\Sigma$ denotes the covariance of the GFF in $D$, as in Definition \ref{GFFdef}) one has $\Sigma(\varphi_i,\varphi_j)=\lambda_i^{-1}\mathbf{1}_{i=j}$. This means that  
$( {\mathcal N}_j :=  \sqrt {\lambda_j} \, \Gamma (\varphi_j) )_{j \ge 1}$ is a sequence of independent standard Gaussian random variables.

Conversely, one can actually start from such a family of i.i.d. centred normal variables $({\mathcal N}_j)_{j \ge 1}$ and (re)construct the GFF. 
 For instance, for any given 
$L^2$ function $f$ with compact support in $D$, we can 
decompose $f$ using the orthonormal basis $(\varphi_j)_j$ as $f(\cdot)= \sum_{j \ge 1} f_j \varphi_j ( \cdot)$, where
$ f_j := \int_D f(x)\varphi_j(x) \, dx$ 
and the sum is converging in $L^2$. 
Then, we can simply define 
\begin {equation} 
\label {RFS}
  \Gamma (f) := \ \sum_{j \ge 1}    \frac{{\mathcal N}_j}{\sqrt{\lambda_j}} f_j
  \end {equation} 
(this sum converges in $L^2$ as $\lambda_j \to \infty$ and $\sum_{j} f_j^2 < \infty$). 

In fact, for any fixed $\mu\in \mathcal{M}$, if we set $\mu_j:=\int_D \varphi_j(x) d\mu(x)$ then the defining property of $\mathcal{M}$ implies that $\sum_{j\ge 1}  \lambda_j^{-1} \mu_j^2 <\infty$. Thus we can set
$\Gamma(\mu):=\sum_{j \ge 1} \lambda_j^{-1/2}{\mu_j} \mathcal{N}_j   $, where the sum also converges in $L^2$. The obtained process $(\Gamma (\mu))_{\mu \in {\mathcal M}}$ is easily seen to be a GFF.

Conversely, we can note that if we start with a GFF $\Gamma$, we can also recover the variables ${\mathcal N}_j = \Gamma (\varphi_j)/ \sqrt {\lambda_j}$.  

We may wonder whether it is actually possible to use this description of the GFF to define $(\Gamma(f))_{f \in {\mathcal S}}$ for all $f$ in some class 
${\mathcal S}$ of smooth functions {\em simultaneously}. Recall that we formally interpret $\Gamma(f)$ as ``$\int f(x) \Gamma(x) \, dx$''.  Then the above expressions suggest that (formally),
$$ \Gamma (f) = \sum_{j \ge 1}  \Bigl[   \frac{{\mathcal N}_j}{\sqrt{\lambda_j}} \int_D f(x) \varphi_j (x) dx  \Bigr] 
= \int_D f(x) \Bigl[\sum_{j \ge 1} \frac{{\mathcal N}_j} {\sqrt{\lambda_j}}  \varphi_j (x) \Bigr] dx $$ 
so that one could try to say, in some appropriate space of generalised functions, that
\begin {equation}
 \label {randomFourierSeries}
 \Gamma (\cdot) = \sum_{j \ge 1}   \frac{{\mathcal N}_j} {\sqrt{\lambda_j}}  \varphi_j (\cdot).
\end {equation}  
To make sense of this, let us take some $s>d/2-1$ and consider the set of functions $f\in L^2(D)$ that satisfy 
$$ \| f \|_{{\mathcal H}^s}^2 := \sum_{j \ge 1} \lambda_j^s f_j^2 <\infty$$ 
(note that the set ${\mathcal H}^s$ of such functions equipped with the corresponding inner product is a Hilbert space). 
Defining $\Gamma(f)$ by (\ref {RFS}) as above, we see, using Cauchy--Schwarz that 
\[  \sum_{j \ge 1}  | \frac{f_j}{\sqrt{\lambda_j}} {\mathcal N}_j  |  
\le \Bigl[ \sum_{j\ge 1} \lambda_j^s f_j^2\Bigr]^{1/2} \times  \Bigl[ \sum_{j \ge 1}\frac{{\mathcal N}_j^2}{\lambda_j^{1+s}}\Bigr]^{1/2}. \] 
Moreover, by Weyl's law we know that 
$ \sum_{j \ge 1} \lambda_j^{-\beta}$ is finite as soon as $\beta > d/2$. 
We can therefore deduce that almost surely
$$  C(s) := \sum_{j \ge 1}  \frac{{\mathcal N}_j^2}{\lambda_j^{1+s}} < \infty,$$
because 
$$  \sum_{j \ge 1}  \frac{E [ {\mathcal N}_j^2] }{|\lambda_j^{1+s}|} = \sum_{ j \ge 1} \frac 1 {\lambda_j^{1+s}} < \infty.$$    
Hence, we can control the absolute convergence 
of the sum in (\ref {RFS}) for all $f \in {\mathcal H}^s$ simultaneously. In other words, we can almost surely define $\Gamma (f)$ for all $f \in {\mathcal H}^s$ at once. 
Furthermore, we see that for all $f, g $ in ${\mathcal H}^s$, 
$$
\bigl| \Gamma (f) - \Gamma (g) \bigr| = \bigl| \Gamma (f - g ) \bigr| 
  \le  \sum_{j \ge 1}     
 \Bigl| \frac{f_j - g_j}{\sqrt{\lambda_j}} {\mathcal N}_j  \Bigr| 
\le C(s)^{1/2} \times \| f -g \|_{{\mathcal H}^s}
$$
(with the obvious definition for $g_j$). This shows that $\Gamma$ can be viewed as a random generalised function, when acting on the space $\mathcal{H}^s$ of test functions, 
and the map $f \mapsto \Gamma (f)$ is then continuous on ${\mathcal H}^s$. In fact, this exactly says that $\Gamma$ can be viewed as a random element of a Sobolev space of negative exponent. 
 
 \subsection{The Cameron-Martin space of the GFF} \label{acbasics}
 The above description of the GFF as a random Fourier series provides a particularly nice framework for discussing absolute continuity relations, and describing the Cameron-Martin space (which is a general concept for Gaussian processes) in the particular case of the GFF. 
 More specifically, we will describe here the class of deterministic functions that one can add to the GFF, such that the sum 
 of the GFF with that function remains absolutely continuous with respect to the GFF itself.
 
 Let us recall the following elementary fact: when $X$ is a standard Gaussian and $a$ is some positive constant, then the Radon-Nikodym derivative of the 
 law of $X-a$ with respect to the law of $X$, at the point $x$, is just $\exp (ax - a^2/2)$. This can be seen by simply writing down the ratio of densities of the two laws. 
 Also recall that if $(a_i)_{i \ge 1}$ satisfies $\sum_{i \ge 1} a_i^2 < \infty$ and $({\mathcal N}_i)_{i \ge 1}$ is a 
 sequence of independent standard Gaussians, then by standard results on series of independent random variables, the series $\sum_{i=1}^n a_i {\mathcal N}_i$ converges
 almost surely, and the limit is a Gaussian random variable with variance $\sum_i a_i^2$.
 
 From these two ingredients, one can easily deduce the following classical fact:
 \begin{lemma}[Cameron-Martin space]\label{acgaussians}
 	Suppose that $(\mathcal{N}_i)_{i \ge 1}$ is a sequence of independent standard Gaussian variables and that $(a_i)_{i\ge 1}$ 
 	is a deterministic sequence of real numbers with $\sum_i a_i^2 <\infty$. Then the law of the process $({\mathcal N}_i + a_i )_{i \ge 1}$ is absolutely continuous with 
 	the law of $({\mathcal N}_i)_{ i \ge 1}$ if and only if $\sum_i a_i^2 < \infty$. Furthermore, the Radon-Nikodym derivative of the former law with respect to the latter, at the point $(x_i)_{i \ge 1}$ 
 	in the support of the law of $({\mathcal N}_i)_{i \ge 1}$, is given by 
 	$ \exp ( \sum_i ( a_ix_i)  -  \sum_i (a_i^2/2))$ 
 	(where the sum $\sum_i (a_i x_i)$ is defined to be the limit as $n \to \infty$ of $\sum_{i=1}^n a_i x_i$).
 \end{lemma}

\begin{remark}
When one decomposes Brownian motion $(B_t, t \in [0, 1])$ (or rather its generalised derivative) using an orthonormal basis of $L^2([0,1])$, then this lemma allows one to describe the space of continuous functions $f$ for which the law of $(B_t + f(t), t \in [0,1])$ is absolutely continuous with respect to that $(B_t, t \in [0, 1])$. This is the Cameron-Martin space of Brownian motion, and consists of the set of functions that can be written as integrals of $L^2$ functions. 
\end{remark}

Let us now turn to the particular case of the GFF. Consider a bounded domain $D$, an orthonormal basis $(\varphi_i)_{i \ge 1}$ of $L^2(D)$, and  $\Gamma$ a GFF in $D$. We denote by $\mathcal{N}_1,\mathcal{N}_2,\cdots $ the sequence of  i.i.d.\ standard Gaussians $(\mathcal{N}_j=\Gamma(\varphi_j)/\sqrt{\lambda_j})$ that corresponds to the orthonormal decomposition of $\Gamma$ described in the previous subsection. 

Let us now define the space ${\mathcal F}$ of functions $f$ in $D$
that can be written as $f (\cdot) = \sum_i f_i \varphi_i (\cdot)$ for some $(f_i)_{i \ge 1}$ with 
$\sum_i \lambda_i f_i^2 < \infty$. 
	This space of functions is a certain Sobolev space, often denoted by $\mathcal{H}_0^1(D)$, which is
	the set of functions in $D$ with finite Dirichlet energy and zero boundary values (more precisely, zero ``trace'') on $\partial D$. Indeed, 
	the quantity $\sum_i \lambda_i f_i^2$ is the $L^2$ scalar product of $f$ with $-\Delta f$, which (because $f$ 
	has zero boundary conditions) is the same as the integral over $D$ of $| \nabla f |^2$.

Translating Lemma \ref{acgaussians} into the language of the GFF (simply writing $a_i = \sqrt {\lambda_i} f_i$ and $x_i = \sqrt {\lambda_i} g_i$) , we get that:  
\begin{lemma}[Cameron-Martin space of the GFF]
\label{cmspace}
	The law of $\Gamma+f$ is absolutely continuous with respect to the law of $\Gamma$ (as stochastic processes indexed by $\mathcal{M}$) if and only if $f \in {\mathcal F}$. 
	Furthermore, the Radon-Nikodym derivative between these two laws at the generalised function
	$g = \sum_i g_i \varphi_i (\cdot)$ (in the support of the law of $\Gamma$) is equal to 
	$\exp ( \sum_i  (\lambda_i g_i f_i) -  \sum_i \lambda_i f_i^2 / 2 ))$ 
	(where $\sum_i  (\lambda_i g_i f_i) := \lim_{n \to \infty} \sum_{i=1}^n  (\lambda_i g_i f_i) $). 
\end{lemma} 

Let us now briefly describe some consequences of this result, that will turn out to be very useful later on. 
Suppose that $D$ is a bounded two-dimensional domain with a smooth boundary, and that $L$ is a finite union of smooth loops within $D$ and smooth curves that remain within $D$ except at their endpoints. 
Suppose that $h_0$ is some Lipschitz function on $L$, that tends to 0 at the finitely many intersection points of $L$ and $\partial \U$. Then one can define the harmonic extension $h$ of $h_0$ to $D$; the value of $h$ at $x$ is simply the expected value of $h_0 (B_T)1_{\{B_T\in L\}}$ when $B$ is a Brownian motion started from $x$ and $T$ its exit time from $D \setminus L$.
\begin {corollary} 
\label {CMcorollary}
The function $h$ is in the Cameron-Martin space of the GFF 
(with Dirichlet boundary conditions in $D$). 
\end {corollary}
\begin {remark} 
It is possible to relax the conditions on $h_0$ quite a bit, but this result will be sufficient for our later purposes.
\end {remark}
\begin {proof}[Sketch]
Suppose that $x \in D \setminus L$, and let $r := d ( x, \partial (D \setminus L))$. We want to evaluate $|h(y)- h(x)|$ when $y \to x$ (in order to bound $| \nabla h |$). For this, we consider $y$ to be very close to $x$ and we can use the mirror coupling between 
two Brownian motions $\{B,B'\}$ started from $x$ and $y$. Write $S$ for the time at which the first of them reaches the circle of radius $r/2$ around the midpoint between $x$ and $y$. The probability that these two Brownian motions do not couple before time $S$  is bounded 
by a constant times $|y-x|/r$. On the event $E$ where they do not couple, we let them run in parallel (instead of being mirror-coupled) after time $S$. It is then a simple exercise (using the smoothness of 
the boundary) to see that conditionally on $E$, $|h(B_S) - h(B_S')|$ is bounded by a constant times $r \log (1/r)$. We can therefore conclude that 
$| \nabla h (x) |$ is bounded by a constant times $\log (1/r)$, which in turn implies that the integral of $|\nabla h|^2$ is finite and that $h$ is indeed in the Cameron-Martin space of the GFF.  
\end {proof}

\subsection {Circular/Spherical averages} \label{circav}
Suppose that $z_0$ is fixed and that $r_0$ is smaller than  $d(z_0, \partial D)$. {We use the notation $\lambda_{z_0,r}$ for the uniform, i.e., multiple of Lebesgue, probability measure on the boundary of $B(z_0,r)$.} Then we define, for all $r\le r_0$, the  average (often referred to as the \emph{circle average} when $d=2$, for obvious reasons)
$$ \gamma (z_0, r) = \Gamma (\lambda_{z_0, r}),$$
which makes sense since $\lambda_{z_0,r}\in \mathcal{M}$ (see the discussion in Section \ref{sec::greensfunction}). 
Suppose that $\mu\in \mathcal{M}$ is another probability measure in $D$ that is supported in $D \setminus \overline{B}(z_0, r_0)$.
Then, we observe that for $r \le r_0$, by the harmonicity properties of the Green's function,
$$ E [ \gamma (z_0, r)  \Gamma (\mu)] 
= \int  d\mu (x) G_D (x,y) d\lambda_{z_0, r}  (dy)
= \int d\mu (x) G_D (x,z_0),$$
and consequently 
$$ E [ ( \gamma (z_0, r) - \gamma (z_0, r_0))  \Gamma (\mu) ] =0.
$$
Hence, the process $(r \mapsto \gamma (z_0, r) - \gamma (z_0, r_0))_{r\in (0,r_0]}$ is independent of any $\Gamma (\mu)$ with $\mu\in \mathcal{M}$ supported outside of $\overline{B} ( z_0, r_0)$.

Similarly, we obtain that for all $r < r' \le r_0$,
$$
E [ (\gamma (z_0, r) - \gamma (z_0, r') )^2 ] = 
\int  d\lambda_{z_0, r} (x) G_D (x,z_0)
-   \int d\lambda_{z_0, r'} (x) G_D (x,z_0) $$ 
which can be shown (using the Markov property) to be equal to
$$ \int  d\lambda_{z_0, r} (x) G_{B(z_0, r')} (x,z_0).$$
This then shows that  
\begin {equation}
 \label {bmcov}
 E [ (\gamma (z_0, r) - \gamma (z_0, r') )^2 ] = \begin{cases} \log(r'/r)   & d=2 \\ (r)^{2-d}-(r')^{2-d} \; & d>2.\end{cases} 
 \end {equation}
It follows (recall that two Gaussian random variables are independent if and only if they have covariance zero) that the process 
$$(b_{z_0, r_0}(u) , u\ge 0):= \begin{cases} (\gamma (z_0, r_0 e^{-u}) - \gamma(z_0, r_0) , u \ge 0)  & d=2 \\ (\gamma(z_0, (u+r_0^{2-d})^{1/(2-d)})-\gamma(z_0,r_0) , u\ge 0 ) \; & d>2 \end{cases} $$ has the same finite dimensional distributions as a one-dimensional Brownian motion, and is independent of the $\sigma$-field generated by $\{\Gamma (\mu): \text{supp}(\mu) \cap \overline{B} ( z_0, r_0)=\emptyset\}$. In particular: 

\begin {proposition}[Spherical averages as independent Brownian motions] 
If we are given a  countable collection of disjoint open balls
$B(z_j, r_j)_{j \ge 1}$ in $D$, then the  
processes $(b_{z_j, r_j})_{j \ge 1 }$  are independent Brownian motions.
\end {proposition} 
{\begin{remark}
	\label{thickpoints}
	Since for any given $z_0$, $b_{z_0,r_0}(u)$ has the same finite dimensional distributions as a one-dimensional Brownian motion, one easily deduces that $b_{z_0,r_0}(u)/u\to 0$ almost surely as $u\to \infty$. However, this does not rule out the existence of \emph{exceptional} points $z_0$ such that this limit is actually something non-zero. For instance, one can ask whether for given $\alpha>0$, there exist exceptional points where this limit is equal to $\alpha$. In fact, it follows from a simple first moment argument that such points do not exist when $\alpha$ is greater than $\sqrt{2d}$. A more refined analysis shows that they do exist when $\alpha\le \sqrt{2d}$, and furthermore, the Hausdorff dimension of the set of exceptional points is strictly positive if and only if $\alpha<\sqrt{2d}$. These points are often referred to as $\alpha$-thick points of the field.
\end{remark}}

\subsection {Kolmogorov's criterion and a first application} \label{sec:kol}
One classical tool to construct ``continuous modifications'' of a stochastic process $ (X_a)_{a \in {\mathcal A}}$ is Kolmogorov's criterion. 
Suppose that ${\mathcal A}$ is a subset of $\R^d$ and that the law of the process $X=(X_a, a \in {\mathcal A})$ is such that there exists a even integer $2N$ and positive constants $\delta$ and $C$ so that for all $a, a' \in {\mathcal A}$, 
\begin{equation}
\label{eqn:kolgen}
 E  \Bigl[ ( X(a) - X(a'))^{2N} \Bigr]   \le C | a - a'|^{d+\delta}.\end{equation}
Then Kolmogorov's criterion guarantees the existence of a modification $X'$ of $X$ such that $a \mapsto X'(a)$ is continuous on ${\mathcal A}$ on some set of probability one (there is a little subtlety here: the event that 
$\{a \mapsto X(a)\}$ is continuous is not measurable, but one can construct a measurable set with probability one that is contained in it). By a modification $X'$, we mean another process defined 
on the same probability space, such that for any given $a \in {\mathcal A}$, one has $X_a = X_a'$ almost surely (note the order of ``for any given $a$'' and ``almost surely'' here). 

This idea of the proof of this is as follows. One first defines $X'$ to be equal almost surely to $X$ on a countable dense subset of ${\mathcal A}$, and then shows using 
the Borel-Cantelli lemma that this process can be extended into a continuous function $a \mapsto X(a)$ on ${\mathcal A}$. Finally, one checks that for any given $a \in {\mathcal A}$, 
$X_a = X_a'$ almost surely. 

In the special setting of Gaussian process, we can make use of the following trivial fact. When $X$ is a centred Gaussian random variable with variance $\sigma^2$ then (for some universal constants $c_N$), one has $E ( X^{2N} ) = c_N \sigma^{2N}$.
From this it follows immediately (applying Kolmogorov's usual criterion with $N$ chosen so that $N \times  \epsilon > d$) that:

\begin {lemma}[Kolmogorov's criterion, Gaussian case] 
\label {kolmo}
If $(X_a, a \in {\mathcal A})$ is a Gaussian process indexed by ${\mathcal A} \subset \R^d$, and there exists positive $\epsilon$ and $C$ such that for all $a, a'\in 
{\mathcal A}$, 
$$ E \Bigl[ ( X(a) - X(a'))^2 \Bigr]   \le C | a - a'|^\epsilon ,$$
then there exists a modification of $X$ such that on a set of probability one, the map $a \mapsto X(a)$ is 
continuous on ${\mathcal A}$.
\end {lemma}

One can of course apply this criterion to show existence of a modification of Brownian motion that is almost surely continuous since Brownian motion satisfies the 
inequality with $C=1$ and $\epsilon =1$.

One example of how to apply this criterion for the GFF goes as follows: 
we have already seen that setting $\gamma(z,r)=\Gamma(\lambda_{z,r})$, when $ r < r' < d(z_0, \partial D)$, the quantity 
$E [ (\gamma (z_0, r) - \gamma (z_0, r') )^2 ] $ is given by (\ref {bmcov}).  
 Moreover, 
 similar considerations show that when $d( z , z') < r_0  $ and $d( z , \partial D) > r_0, d(z',\partial D)>r_0$,  
$E[ (\gamma (z, r) - \gamma (z', r))^2 ]$ is bounded by a constant (depending on $r_0$) times $| z -z'|$. 

Hence, we get that for any given $r_0$, there exists a constant $C(r_0, D)$ such that  
$$E [ (\gamma (z, r) - \gamma (z', r') )^2 ]  \le  C(r_0, D) \times ( | z-z'| + |r-r'| ) $$ 
for all $r, r' > r_0$ and $z, z'$ in $D$ with $d(z, z') \le r_0/2$, $d(z, \partial D) > r_0$. 
Applying Lemma \ref {kolmo} in this setting, one can therefore deduce that there exists a version of $\gamma (z, r)$ such that $(z, r) \mapsto \gamma (z,r)$ is 
continuous on the set of $(z,r)$ such that $d(z, \partial D) > r_0 > r$. Since this is true for all rational $r_0$, one can readily conclude the following. 
\begin {proposition}[Spherical averages as a bi-continuous function]
 There exists a version of the process $(z,r) \mapsto \gamma (z,r)$ that is continuous on 
$\{ (z,r) \in D \times (0,\infty), \ r < d (z, \partial D) \}$.
\end {proposition}
This shows, for instance, that it is possible to construct a modification of the process $\Gamma$ such that, almost surely, {\em all} the $b_{z,r}$ are (simultaneously) continuous Brownian motions.

Finally, if $\overline \gamma (z,r)$ denotes the value of $\Gamma(\mu)$ with $\mu$ equal to Lebesgue measure on the ball $B(z,r)$, then it is possible to define a version of this process that is jointly continuous in $z$ and $r$. This time $r=0$ is allowed, because the Lebesgue measure has vanishing total mass as $r \to 0$, and this readily implies that if we set $\overline \gamma (z, 0) = 0$ then the process will still be continuous. Note  that $r \mapsto \overline \gamma (z,r)$ is in fact differentiable and that its derivative is related to the spherical averages.

\subsection {Translation/scale/conformal invariance of the  GFF} 
\label {CIGFF}
When $d=2$, the GFF inherits a conformal invariance property from the conformal invariance of the Green's function. 
More precisely, suppose that $D$ and $\tilde D$ are two conformally equivalent domains in the plane (i.e. there exists an angle-preserving bijection $\Phi$ from $D$ onto $\tilde D$). Then, we have seen that $G_D (x,y) = G_{\tilde D} ( \Phi (x) , \Phi (y))$. Hence, if the GFF were an actual function, then the law of this function would be conformally invariant. In reality, it is conformally invariant ``as a generalised function'', which means that for any $\mu \in \mathcal{M}$, $\Gamma_D (\mu)$ is distributed like $\tilde \Gamma_{\tilde D} ( \tilde \mu)$ (where this $\tilde \Gamma_{\tilde D}$ is a GFF in $\tilde D$), for $\tilde \mu$ the push-forward measure defined by
$$ \tilde \mu (\Phi (A)) := \int_{A} \mu (dx) |\Phi' (x)|^2.$$
In other words, if $\Gamma$ is a GFF in $D$ and if we {\em define} for each $\tilde \mu$ in ${\mathcal M}_{\tilde D}$ 
the random variable $$ \tilde \Gamma (\tilde \mu)  = \Gamma_D (\mu), \;\; \mu (A) := \int_{\Phi (A)} \tilde \mu (dy)  |(\Phi^{-1})' (y)|^2, $$
then $\tilde \Gamma$ is a GFF in $\tilde D$. In the sequel, we will simply refer to this GFF $\tilde \Gamma$ as {\em the image of $\Gamma$ under the conformal map $\Phi$}, and denote it by $\Gamma \circ \Phi^{-1}$. 

 In a similar fashion, when $d> 2$, since for any connected $D$, $r>0$ and $x,y$ in $D$ we have $G_{rD}(rx,ry)=r^{2-d}G_D(x,y)$, it follows that $r^{d/2-1}\Gamma_{rD}$ is equal in law to $\Gamma_{D}$. Moreover, for any $a\in \R^d$, it is clear that $G_{D+a}(x+a,y+a)=G_D(x,y)$, and so $\Gamma_{D+a}$ is equal in distribution to $\Gamma_D$. Observe here that $d=2$ plays a special role: it is the only dimension in which the Gaussian free field is scale invariant.

\section {Relation with Brownian loop-soups (a non-rigorous warm-up)} 

In view of the relation between the discrete GFF and loop-soups, it is natural to wonder if some analogous results might hold in the continuum setting. The answer is that (at least in dimension $2$) most of these results do indeed hold, and there are actually some rather nice additional features (for instance, due to the role played by conformal invariance). In the present section, we will survey (without proper proofs) some results in this direction, in order to provide 
some motivation and guiding principles for the next chapters.

\subsection {The Brownian loop-soup} 

It is rather easy to guess how the definitions of discrete loop-soups and cable-graph loop-soups should be extended to the continuum. In this case, the loops will be described by closed trajectories of $d$-dimensional Brownian motion. 
Let us define it step by step. 
\begin {itemize}
\item Recall that one can define a one-dimensional Brownian bridge of time length $T$ $$\beta=(\beta_t, t \in [0, T])$$ from $0$ to $0$, to be the process $\beta_t = B_t - (tB_T/T)$, where $B$ is a one-dimensional Brownian motion (there are actually several equivalent definitions, for example using Fourier decomposition, or ``conditioning'' Brownian motion to be at $0$ at time $T$). 
\item Similarly, if one uses the same definition but replaces $B$ by a $d$-dimensional Brownian motion, 
then $\beta$ is a $d$-dimensional Brownian loop from $0$ to $0$ of time-length $T$. Its $d$ coordinates are then $d$ independent one-dimensional Brownian bridges. Let us denote the law of this $d$-dimensional Brownian loop by $P_{0 \to 0, T}$.
\item The density at the origin (with respect to the Lebesgue measure) for $B_T$ is given by $(2 \pi T)^{-d/2}$. This makes it natural, in view of the definition of the discrete loop-soup and of the Green's function, to define a Brownian loop-measure rooted at the origin by
$$ \mu_{0 \to 0} := \int_0^\infty \frac {dt}{(2 \pi t)^{d/2}} P_{0 \to 0, t}.$$
\item We would now like to define an \emph{unrooted} Brownian loop-measure, which should be invariant under translations of $\R^d$. 
To do so, we first note that an unrooted Brownian loop of time-length $T$ will (heuristically) have $T$ times 
more possible starting points than an unrooted Brownian loop of time-length $1$. This suggests that it is better to start with the modified rooted measure on Brownian loops defined by 
$$ \tilde \mu_{0 \to 0} := \int_0^\infty \frac {dt}{t \times (2 \pi t)^{d/2}} P_{0 \to 0, t} . $$
We can then define $\tilde \mu_{x \to x}$ to be the image of $\tilde \mu_{0 \to 0}$ by $z \mapsto z + x$ (so that $\tilde \mu_{x \to x}$ is a measure on Brownian loops from $x$ to $x$), and set
$$ \tilde \mu := \int_{\R^d} dx \tilde \mu_{x \to x}.$$ 
Finally, we define the unrooted Brownian loop-measure $\mu$ to be measure induced by $\tilde \mu$ on the set of equivalence classes of unrooted loops (i.e., when one erases the information about where the root $x$ was on the loop). 
\end {itemize}
\begin {definition}[Brownian loop-measure] 
This measure $\mu$ defines the Brownian loop-measure in $\R^d$. For any open subset $D$  of  $\R^d$, we define $\mu_D$ to be the restriction of $\mu$ to the set of loops that remain entirely in $D$. 
\end {definition}
Note that these are measures on oriented (unrooted) loops.  It is then easy to check the following properties: 
\begin {enumerate}
 \item if one starts with the unrooted loop measure $\mu$, and for each loop, chooses its root uniformly at random on the loop, then one obtains the measure $\tilde \mu$;
 \item the measure $\mu$ is invariant under translations $\gamma \mapsto \gamma + x$ and under multiplications $\gamma ( \cdot) \mapsto \lambda \gamma (\cdot / \lambda^2)$;
 \item the measure $\mu$ is ``locally finite.'' For instance, the $\mu$-mass of the set of loops intersecting the cube $[0,1]^d$, with time-length between $1$ and $4$ (say), is finite.
\end {enumerate} 
If we now consider the measure $\mu_D$, for $D$ a bounded domain, we see that for all $\eps>0$ the $\mu_D$-mass of the set of loops of diameter greater than $\eps$ is finite. 
This mass goes to $\infty$, asymptotically as $\eps \to 0$, like a constant times $\eps^{-d}$. Equivalently, the $\mu_D$-mass of the set of loops of time-length greater than $u$ behaves asymptotically 
like a (different) constant times $u^{-d/2}$. 

\begin {definition}[Brownian loop-soup]
 For $\alpha>0$, we define the oriented Brownian loop-soup with intensity $\alpha$ in $D$ to be a Poisson point process of unrooted oriented Brownian loops, with intensity $\alpha \mu_D$. 
\end {definition}

Given that random walks converge to Brownian motion, it should not be surprising that in some appropriate sense, the Brownian loop-measure and the Brownian loop-soup can be viewed as 
limits of random walk loop-measures and random walk loop-soups on fine-mesh approximations of $D$ (or on its cable-graph).

\begin {remark} 
It is possible to show that the resampling properties and the Markov-type properties of the loop-soups all have analogues for these Brownian loop-soups. In particular, this will happen for the soup of oriented Brownian loops with intensity $\alpha = 1$. It is also possible to use the soup of oriented loops corresponding to $\alpha=1$ to understand the ``scaling limit'' of Wilson's algorithm (we will briefly discuss this in Chapter \ref{Ch6}, Section \ref {SLE2}).

One can also define the measure $\kappa_D$ to be the image on the set of {\em unoriented} Brownian loops of $\mu_D$, and construct for each $c > 0$, the {\em unoriented} loop-soup with intensity $c$ to be a Poisson point process with intensity $c \kappa_D$. Then, the unoriented loop-soup with intensity $c=1$ will be the one with the special resampling property. This will also be the one that is going to be discussed in the next section, as it can be related to the GFF via occupation times. 
\end {remark}

\subsection {Loop-soup occupation time and square of the GFF} 
\label {S.BLSGFF}

Suppose that one samples a Brownian loop-soup of positive intensity $\alpha$ in a bounded domain $D$. Let us first make some back-of-the-envelope a priori estimates. 

- The number of loops in the loop-soup with time-length between $2^{-n}$ and $2^{-n +1}$ will be a Poisson random variable with mean of order $2^{nd/2}$ as $n \to \infty$. 

- Hence the cumulated 
time-length $T_n$ of all such loops (i.e., the sum of their individual time-lengths)  will have expectation of order $2^{n(d/ 2 -1)}$ and variance of order $2^{n(d/2 -2)}$. 

- The sum over $n \ge  0$ of these variances therefore converges if and only if $d<4$. On the other hand, the sum of the expectations diverges as soon as $d \ge 2$.
 \medbreak 
 
 This leads to the following feature of the loop-soups when $d=2$ or $d=3$:
Consider an oriented Brownian loop-soup ${\mathcal L}_\alpha$ with intensity $\alpha$ in $D$ (or equivalently, an unoriented Brownian loop-soup with intensity $c= 2 \alpha$ as they will define the same occupation times). Then: 
\begin {proposition}[Renormalised occupation time measure] 
For each open set $O \subset D$, let us define ${\mathcal T}_n (O)$ to be the total time spent in $O$ by all loops in the loop-soup with time-length at least $2^{-n}$. Then the sequence 
${\mathcal T}_n(O) - E [ {\mathcal T}_n (O) ]$ converges in $L^2$, as $n\to \infty$, to a finite random variable ${\mathcal T}_\infty (O)$. 
\end {proposition}
It should be stressed that ${\mathcal T}_\infty (O)$ has zero expectation and can therefore take negative values - so it is not really an occupation time!
Intuitively, ${\mathcal T}_n (O)$ blows up as $n \to \infty$ because of the very large number of very small loops. On the other hand, there is some law of large numbers behaviour occurring, which means that ${\mathcal T}_n (O)$ 
actually stays quite close to its expectation at first order. This is why the variance of ${\mathcal T}_n (O)$ remains finite.

Motivated by the discrete result relating the square of the discrete GFF to the random walk loop-soup with $\alpha=1/2$, we now discuss how one can try to define the square of the continuum GFF $\Gamma$. This is of course non-trivial, because 
$\Gamma$ is not a proper function, so one cannot a priori take its square. Let us briefly and heuristically explain one way to proceed. 
One natural option  is to use the decomposition of $\Gamma$ on an $L^2$-basis $(\varphi_n)_{n \ge 1}$ of eigenfunctions of the Laplacian, and to view $\Gamma$ (as in Section \ref {S.Fourier}) as the limit 
when $N \to \infty$ of 
$$ \sum_{1\le n \le N} \frac {{\mathcal N}_n}{\sqrt {\lambda_n}} \varphi_n (\cdot),$$
where $({\mathcal N}_n)_{n \ge 1}$ is a sequence of independent identically distributed standard Gaussian random variables.
With this in mind it is tempting to investigate the behaviour as $N \to \infty$ of the function
$$ \Lambda_N (x) := \Bigl( \sum_{1 \le n \le N} \frac {{\mathcal N}_n}{\sqrt {\lambda_n}} \varphi_n (x) \Bigr)^2. $$
For instance, one can look at the integral $ (\Lambda_N ,1) := \int_D \Lambda_N(x) dx$ of $\Lambda_N$ with respect to the Lebesgue measure on $D$. 
Expanding the square, and using the fact that $(\varphi_n)_{n \ge 1}$ is an orthonormal basis of $L^2 (D)$, we get
$$  (\Lambda_N ,1) = 
 \Bigl( 2 \sum_{1 \le n<m\le N} \frac { {\mathcal N}_n {\mathcal N}_m}{\sqrt{\lambda_n \lambda_m}} \int_D \varphi_n(x) \varphi_m (x) dx \Bigr) + 
\sum_{1 \le n \le N} \frac { {\mathcal N}_n^2}{\lambda_n} \int_D \varphi_n^2 (x) dx
= \sum_{1 \le n \le N} \frac { {\mathcal N}_n^2}{\lambda_n}.$$
Hence, we see using Weyl's law (this implies that $\lambda_n$ behaves like a constant times $n^{2/d}$) that:
\begin{itemize}
	\item on the one hand 
	$ E[ (\Lambda_N ,1) ] = \sum_{1 \le n \le N} ( 1  / {\lambda_n}) $
	goes to infinity as $N \to \infty$ when $d \ge 2$;
\item  on the other hand, the variance
$$ E [( (\Lambda_N ,1) - E [ (\Lambda_N ,1) ] )^2 ] = \sum_{1 \le n \le N} \frac {1}{\lambda_n^2} E [ ( {\mathcal N}_n^2 - 1 )^2 ]=  \sum_{1 \le n \le N} \frac {2}{\lambda_n^2}$$
actually converges as $N \to \infty$, when $d=2$ or $d=3$. 
\end{itemize} We can therefore interpret the limit in $L^2$ of 
$ \sum_{1 \le n \le N} ({ {\mathcal N}_n^2 - 1}) / {\lambda_n}$ as the 
integral over $D$ of the ``renormalised square'' of the GFF $\Gamma$. We will denote this limit as $(\Lambda, 1)$ (even though at this point, we have not shown that it corresponds to the integral of a generalised function $\Lambda$. 

A first statement that relates this squared GFF to loop-soup occupation times goes as follows: 
\begin {proposition} 
The law of $(\Lambda , 1)$ is identical to that of a constant multiple of the renormalised occupation time ${\mathcal T}_{\infty}(D)$ of a Brownian loop-soup of intensity $\alpha=1/2$. 
\end {proposition}

This result can be upgraded into a stronger statement relating the two processes ${\mathcal T}_{\infty}(O)$ and $(\Lambda, 1_O)$ indexed by the collection of open subsets $O$ of $D$. 
But in order to state this properly, one first needs to make sense of the latter random variables. For that, one can heuristically use the decomposition $\Gamma = \Gamma_A + \Gamma^A$ that will be described in the next chapter.  for $A= D \setminus O$, and the fact that
$\Gamma^A$ and $\Gamma_A$ have zero mean and are independent, so that can expand the sum of $h_A$ and $\Gamma^A$. In other words, the ``square of $\Gamma$'' integrated on $O$ (that we denote by $(\Lambda, 1_O)$ would be the square of $\Gamma^A$ (as defined above) plus the integral of the square of 
$h_A(x)^2 - E [ h_A(x)^2] $ over $O$, plus twice $\Gamma^A (h_A)$. 
Again, this process has zero expectation, so it is 
not necessarily positive -- it is not really a square! 

It turns out that the relationship between occupation times of discrete loop-soups and the square of the discrete GFF has a natural continuum counterpart:   
\begin {theorem}[Renormalised loop soup occupation time and square of the GFF] 
\label {BLSGFF}For $d=2$ and $d=3$, 
the two processes $(\Lambda (\mathbf{1}_O))_{O \in {\mathcal O}}$ and  $({\mathcal T}_\infty (O))_{O \in {\mathcal O}}$ have the same law. 
\end {theorem}

\subsection {The excursion decomposition of the continuum 2D GFF}

We have just seen that when $d=2$ and $d=3$, there is still a natural coupling between a Brownian loop soup and the
(renormalised) square of the GFF.
One may next wonder if there is a relation between the GFF itself and the Brownian loop-soup.
In view of the discrete GFF results (and the fact that the signs of the GFF on the cable graphs can be chosen to be independently for each cluster), it seem natural to guess that the sign of the GFF will be chosen independently for each cluster of Brownian loops.
The goal of this section is to briefly survey without proof some recent results in this direction, restricted to the two-dimensional case. Even though the statements will probably be surprising (and perhaps confusing) on a first reading, we hope that having them  in mind already will help to guide the intuition in the next chapters. 

Recall the construction given earlier, of a GFF starting from a 
cable graph loop-soup in the discrete setting. One alternative way to describe this can be summarised as follows. 
When one samples a cable-graph loop-soup for $c=1$, one obtains a partition of the set of vertices of the graph into clusters (we can call them $C_i$). One can also associate to each cluster $C_i$ a non-negative function $\alpha_i$ with support exactly $C_i$, and a sign $\eps_i \in \{-1, 1\}$, in such a way that: 

(a) The GFF is equal to $\sum_i \eps_i \alpha_i (\cdot)$. 

(b) Conditionally on the collection $(C_i, \alpha_i)$, the signs 
$\eps_i$ are chosen to be independent with $P [\eps_i = 1 ] = 
P [ \eps_i = -1 ] = 1/2 $. 

It turns out that when $d=2$, a similar result can be shown to hold for the continuum GFF. Let us describe this without proof (some of the results that we are going to state in the coming paragraphs would require very long proofs, that build among other things on ideas we will present in the next chapters!).

Let us take $D$ to be the unit disk in the plane, and sample an intensity $c=1$ Brownian loop-soup in $D$. It is useful to keep in mind that for every fixed $x \in D$, one can almost surely find infinitely many (small) Brownian loops in the loop-soup that do surround $x$ (and disconnect it from $\partial D$), but it turns out that there are also some random exceptional points that are surrounded by no Brownian loop in the loop-soup. 
One can partition the set of loops in this loop-soup into
clusters: any two Brownian loops $l$ and $l'$ in the loop-soup will be in the same cluster as soon as one can find a finite chain of Brownian loops $l_0=l, \ldots, l_m=l'$ in the loop-soup such that $l_j \cap l_{j-1} \not= \emptyset$ for all $j < m$. It turns out that (again, we warn the reader that 
none of these statements is easy to prove...): 

\begin {enumerate}
 \item 
 Almost surely, there exist infinitely many clusters of loops. This is actually not obvious -- indeed one can show that as soon as $c$ is greater than $1$, all loops are in the same cluster! Let us call the family of clusters $(C_i)_{i \in I}$. It will actually be convenient to define the $C_i$ to be closures of the unions of all Brownian loops in the same cluster, so that each $C_i$ is a compact set. 

\item The sets $C_i$ are almost surely all disjoint (this is not obvious because one defined each $C_i$ as the closure of the union of Brownian loops) and they are all almost surely at positive distance from $\partial D$. 

\item Actually, when one looks at the geometry of each of the $C_i$, one can observe that the boundary consists of the 
union of countably many disjoint simple loops. In particular (as this will be relevant later), the outer boundary of each $C_i$ is a simple loop. 

\item Note that each $C_i$ contains the union of countably many Brownian loops, so that it will be a rather ``fat'' random object (in particular, its Hausdorff dimension is going to be equal to $2$). 

\item To each $C_i$, it is possible to associate in a measurable (deterministic) way a measure $\alpha_i$ supported on $C_i$, so that if one introduces a family of random variables $(\eps_i)_{i\in I}$ with 
$P [\eps_i = 1 ] = 
P [ \eps_i = -1 ] = 1/2 $ that are conditionally independent given the $(C_i)$, then the field 
$ \sum_i \eps_i \alpha_i (\cdot) $
is a GFF. Here, the sum over $i$ should be viewed in $L^2$. One way to phrase this statement in more detail, is that for any given order $i_1, i_2, \ldots$ of the clusters, then for each finite family of smooth functions $(f_1, \ldots, f_k)$, the random variable
$$ (\sum_{j=1}^n \eps_{i_j} \alpha_{i_j} (f_1) , 
\ldots , \sum_{j=1}^n \eps_{i_j} \alpha_{i_j} (f_k) )$$ 
converges in $L^2$ to a limit $(\Gamma (f_1), \ldots, \Gamma (f_k))$, 
where $\Gamma$ is a GFF. 

 In this decomposition, the sets $C_i$ and the signs $\eps_i$ are in fact deterministic functions of the GFF $\Gamma$. So, one can interpret the formal identity 
$$ \Gamma = \sum_i \eps_i \alpha_i $$ 
as an ``excursion decomposition'' of the continuum GFF. 

\end {enumerate}

In this list of results, we can already observe one feature that appears a little surprising, because it indicates (in a way) that something stronger happens in the continuum setting than in the discrete. Here, the measures $\alpha_i$ are deterministic functions of the clusters $C_i$ -- while in the discrete setting, the knowledge of the clusters was clearly not sufficient to recover the actual random functions $\alpha_i$. 

We are now going to describe another striking feature of this 
relation between the two-dimensional GFF and the Brownian loop-soup. Let us do this in very loose terms. Suppose that one discovers a loop-soup cluster from the outside. For instance, one can define the 
outermost cluster $C_{i_0}$ that disconnects the origin from the unit circle. This outermost cluster has an outer boundary which is a simple loop $\eta_i$ around the origin. The fact that we chose the ``outermost cluster'' indicates that it is possible to somehow algorithmically discover the loop $\eta_{i_0}$ ``from the outside''.  
Now, here is the question which will turn out to have a surprising answer: 
{\em Conditionally on $\eta_{i_0}$, what is the conditional law of the restriction of $\Gamma$ to $O_{i_0}$?} 
The answer in the analogous question in the case of the cable system GFF would be something like a GFF in $O_{i_0}$ conditioned by the event that all the points of $\eta_{i_0}$ belong to the same cluster. 
But in this two-dimensional continuum case, it turns out that this  complicated conditioned GFF is very easy to describe. The answer to the previous question is the following:  {\em the conditional law is that of a GFF in $O_{i_0}$ with boundary conditions $2 \eps_{i_0} \lambda$ for some universal constant $\lambda$.}
Intuitively, the $\eps_{i_0}$ would correspond to the coin tossing that decides about the sign of the GFF on the cluster $C_{i_0}$.

This suggests that the complement of $O_{i_0}$ would be some sort of 
local set for the continuum GFF, where the harmonic function associated to it would be the constant function $2 \eps_{i_0} \lambda$. To understand this, it is therefore useful to understand how to make sense of the Markov property and of local sets for the continuum GFF, which will be the topic of the next chapter. 

Furthermore, while the GFF would have boundary condition $2 \eps_{i_0} \lambda$ on the inside of the loop $\eta_{i_0}$, it is not difficult to work out that, conditionally on $\eta_{i_0}$, the law of the loop-soup in the complement of $O_{i_0}$ is just a loop-soup in the domain $D \setminus \overline {O_{i_0}}$ conditioned to have no cluster that surrounds $O_{i_0}$, which is an event of positive probability. Hence, this suggests that one the ``outside'' of the loop, the GFF $\Gamma$ looks like a GFF with zero boundary conditions. In other words, $\eta_{i_0}$ is a little bit like a cliff -- with $0$ boundary conditions on one side, and $2 \eps_{i_0} \lambda$ on the other side. 
The chapter on the Schramm-Loewner Evolution SLE$_4$ will make sense of these cliff-lines.

\section*{Bibliographical comments} 
The content of the first three sections of this chapter is rather classical. For thick points of the GFF, we refer to \cite {HMP}. The results presented (mostly without proof) in the final section are more recent. The Brownian loop-soup was introduced in \cite {LawlerWerner}, the relation between the square of the GFF and the loop-soup occupation time is due to Le Jan \cite {LJ2}, and the excursion decomposition of the continuum GFF is due to Aru, Lupu and Sep\'ulveda (see \cite {ALS1,ALS2} and the references therein -- it also builds on the corresponding results 
on cable systems and on SLE-type considerations). The cliff-lines of the GFF will be discussed and commented on in Chapter \ref {Ch5} (see the bibliographical comments there).

\chapter {The Markov property and local sets in the continuum}
\label {Ch4}

In this chapter, we will discuss in rather abstract terms what the analogue of local sets are for the continuum GFF.

\section {The Markov property} 

The goal of this section is to describe the continuum analogue of the Markov property for the discrete GFF. There are as usual several ways 
to tackle this, and we will present one route, which is possibly not the most elegant one! 

Let us fix some compact subset $A$ of $\overline D$, such that 
the boundary of $O := D \setminus A$ is regular as well. 
Let $\Gamma$ denote a continuum GFF in $D$. Our goal, inspired by the corresponding results in the discrete case, is to decompose $\Gamma$ into the sum 
of two independent processes $\Gamma_A$ and $\Gamma^A$, i.e., $\Gamma (\mu) = \Gamma_A (\mu) + \Gamma^A (\mu)$ for all $\mu \in {\mathcal M}_D$, where:
\begin {itemize} 
 \item the process $\Gamma^A$ is a continuum GFF in $O= D \setminus A$ (with zero boundary conditions);
 \item the field $\Gamma_A$ should be thought of as ``equal to $\Gamma$ in $A$''  and to be defined in $O$ as the harmonic extension $h_A$ of the ``values of $\Gamma$ on $\partial O$''.
\end {itemize}
We note that, just as in the discrete case, the decomposition $\Gamma=\Gamma_A+\Gamma^A$ together with the first bullet point implies that $\Gamma_A$ must be a centred Gaussian process indexed by $\mathcal{M}_D$, with covariance kernel given by $H_{D, O}= G_D - G_O$. Explicitly,  for all $\mu$ and $\mu'$ in ${\mathcal M}_D$, it must be that
$$ E [ \Gamma_A (\mu) \Gamma_A (\mu') ] = \int_{D \times D} d\mu (x) d\mu'(x') (G_D (x, x') - G_O (x, x') ) .$$ 

Let us first comment on why it is possible  to make sense of the harmonic function $h_A$ described above, even though $\Gamma$ is not defined pointwise (so it is a priori not so 
clear what this harmonic extension should mean). Consider the example where $O$ is a ball in $\R^d$, centred at $z_0$ and of radius $r_0$.
In this case, when one is given a continuous function $g$ on $\partial O$, the unique harmonic function in $O$ with boundary value $g$ is equal at $z_0$ to the average value of $g$ on the sphere $\partial O$.
This suggests that it is natural 
to set $h_A (z_0) := \Gamma (\lambda_{z_0, r_0})$ (which is indeed a well-defined random variable: see Section \ref{circav}). 

More generally, for each $z \in O$, recall that $\nu(z, \partial O)$ denotes the law of the first point on $\partial O$ that a Brownian motion started from $z$ hits. We can then set  

\begin{equation}
\label{hA}
 h_A (z):= \Gamma (\nu_{z, \partial O})
 \end{equation}
(recall that we know by (iii) of Section \ref {sec::greensfunction} that $\nu_{z, \partial O} \in {\mathcal M}_D$). 
By (iii) of Section \ref {sec::greensfunction} again, we then know that 
$$E [ h_A (z) h_A (z') ] = H_{D, O} (z, z')$$ 
as expected. Note that this definition immediately implies that $(h_A(z))_{z \in O}$ is a centred Gaussian process. 

One would now like to say that $h_A$ can actually be realised as a harmonic function. Naturally, the first step is to show that $\nu_{z, \partial O}$ and $\nu_{z', \partial O}$ are close when $z$ and $z'$ are close. For this, one uses the usual ``mirror coupling'' between two Brownian motions starting from $z$ and $z'$: until the paths meet, the increments of one are the reflection of the increments of the other in the hyperplane bisecting $[z,z']$, and after this time they coincide. It can be shown that for given $\eps>0$ there exists a constant $C(\eps)<\infty$ such that under this coupling, the probability that the Brownian motions do not coincide before hitting the boundary is less than $C(\eps)|z-z'|$, uniformly in $z,z'$ with $d(z,\partial O),d(z',\partial O)\ge \eps$. Hence the total variation of the measure $\nu_{z, \partial O} - \nu_{z', \partial O}$ is bounded by $C( \eps) | z-z'|$ for all such $z,z'$.
From this, it follows that if the process $h_A$ is defined as in \eqref{hA},
then 
\begin{equation}\label{eqn:varhA}E [ (h_A (z) - h_{A} (z'))^2]  = \int_{\partial O}(G_D(z,y) - G_D(z',y)) (\nu_{z, \partial O}-\nu_{z', \partial O}) (dy) \end{equation} is less than or equal to $C'(\eps) | z -z'|$, where $C'(\eps)=2C(\eps)\times G_{\R^d}(0,\eps)$. Applying Kolmogorov's criterion then readily shows that: 
\begin {lemma}[Defining the harmonic extension]  
There exists a continuous version of the process $(h_A(z))_{z \in O}$, and this continuous version is a harmonic function in $O$. 
\end {lemma}
For the last part of the statement, one only needs to verify that this continuous version satisfies the mean-value property in $O$, i.e. that for all $z \in O$ and all $r < d(z, \partial O)$, 
the mean-value of $h_A$ on the sphere $S(z,r)$ of radius $r$ around $z$ is equal to 
$h_A(z)$. Indeed, we have already mentioned that a continuous function satisfying the mean value property is actually smooth with vanishing Laplacian. However, the mean-value property of $h_A$ follows directly from the fact that $\nu_{z, \partial O}$ is the mean-value of $\nu_{z', \partial O}$ with $z'$ ranging over $S(z,r)$: a consequence of the strong Markov property of Brownian motion. So, $h_A$ is (on a set of full probability) indeed a harmonic function. 

On the other hand, it is important to note that when $\partial O \subset D$ is a deterministic set, then the harmonic function $h_A$ will not be bounded in any neighbourhood of $\partial O$. It will typically
start oscillating pretty wildly: this corresponds 
to the fact that $\Gamma$ is not defined pointwise on $\partial O$.

The next step in our quest for the Markov property
is to define the random variable $\Gamma_A (\mu)$ for $\mu \in {\mathcal M}$. The first idea (keeping in mind the Markov decomposition of the discrete GFF) would be to define it as 
$$  \Gamma ( \mu \mathbf{1}_A ) + \int_{O} h_A (x) \mu (dx).$$
However, care is required, because it is not clear whether the integral of $h_A$ is well-defined in the usual sense. As we have already mentioned, the function $h_A$ will 
not be bounded near $\partial A$ and indeed, in general, it might happen 
that $\int |h_A|  \mu (dx) = \infty$ for some measure $\mu$. 
One way around this is to instead define another measure 
$\nu_{\mu, \partial O}$ which is the integral with respect to $\mu (dx) 1_O$ of $\nu_{x, \partial O}$, and to then  define  
$$ \Gamma_A ( \mu)  := \Gamma ( \mu \mathbf{1}_A ) + \Gamma (\nu_{\mu, \partial O}),$$
(which seems a good alternative to $\Gamma ( \mu \mathbf{1}_A ) + \int_{O} h_A (x) \mu (dx)$ given that $h_A(z)=\Gamma(\nu_{z, \partial O})$).

To justify this definition, we need to explain why $\nu_{\mu, \partial O} \in {\mathcal M}_D$. For this, first assume that $\mu$ is non-negative, and note (using (iii) of Section \ref {sec::greensfunction} again) that
\begin {eqnarray*}
 \lefteqn {\int_{\partial O \times \partial O} \nu_{\mu, \partial O} (z) \nu_{\mu, \partial O} (z') G_D (z, z')} \\
 &=& \int_{O \times O} \mu (dx) \mu( dx') \, [ 
\int_{\partial O \times \partial O} \nu_{x, \partial O} (z) \nu_{x', \partial O} (z')  G_D (z, z')  ]  \\
& 
=&  \int_{O \times O} \mu (dx) \mu( dx') H_{D, O} (x, x') \le \int_{D \times D } \mu (dx) \mu (dx') G_D (x, x') < \infty.
\end {eqnarray*} 
The justification for general $\mu \in \mathcal{M}$ follows by splitting $\mu$ into positive and negative parts.

Since $\mu \mathbf{1}_A$ and $\nu_{\mu, \partial O}$ are deterministic measures (i.e., deterministic functions of $\mu$), 
it follows that the process $\Gamma_A$ (equal to $\Gamma(\mu \mathbf{1}_A)+\Gamma(\nu_{\mu, \partial O})$ for each $A$) is a centred Gaussian process. By (iii) of Section \ref {sec::greensfunction}, one easily checks that  
$$ E [ \Gamma_A ( \mu) \Gamma_A (\mu') ] 
= \int_{D \times D} (G_D (x,y) - G_{O} (x,y)) \mu (dx) \mu' (dy)$$ 
as expected. Hence,  $\Gamma_A$ fulfils the properties that we are looking for: it 
is nothing else than $\Gamma$ when restricted to measures supported on $A$, and  when restricted to measures $\mu$ with {\em compact support} in $O$, it is exactly the integral $\int h_A(z) d\mu (z)$. 

Note also that by definition, $\Gamma_A (\mu)$ is equal to $\Gamma ( \overline \mu)$ for  $\overline \mu = \mu\mathbf{1}_A + \nu_{\mu, \partial O}$, which is a measure supported in $A$. Therefore, the process 
$\Gamma_A$ is ${\mathcal F}_A$ measurable (recall that $\mathcal{F}_A$ is the $\sigma$-field generated by all the $\Gamma(\mu)$ for $\mu\in \mathcal{M}$ supported in $A$).
Conversely, it is clear that ${\mathcal F}_A \subset \sigma (\Gamma_A)$ (because $\Gamma (\mu) = \Gamma_A (\mu)$ for any measure $\mu$ supported in $A$). 
Hence,  $\mathcal{F}_A$ is exactly the $\sigma$-field generated by the process $\Gamma_A$.
 
Finally, we define, as in the discrete case, 
$$ \Gamma^A := \Gamma - \Gamma_A .$$ 
Again using the harmonicity properties of the Green's function, we see that the processes $\Gamma^A$ and $\Gamma_A$ are independent.
Indeed, first note that $\Gamma^A$ vanishes on all measures supported in $A$. Then, observe that for all $\nu$ supported in $A$ and all $\mu$ supported in $O$, 
we have
\begin {eqnarray*}
 \lefteqn { 
E [ \Gamma (\nu) \Gamma_A (\mu)] = \int_{A \times O} d\nu (x) d\nu_{\mu, \partial O} (y) G_D (x, y) }
\\
&&
= \int_{A \times O} d\nu (x) d\mu (y) G_D (x, y) 
= E[ \Gamma (\nu) \Gamma (\mu) ]. 
\end {eqnarray*} 
This implies the independence between $\Gamma (\nu)$ and $\Gamma^A (\mu)$, and therefore between $\Gamma^A (\mu)$ and ${\mathcal F}_A=\sigma(\Gamma_A)$. 

The covariance function of $\Gamma^A$ is thus given by the difference between that of $\Gamma$ and that of $\Gamma_A$, so that for any $\mu \in \mathcal{M}$,
$$ E [ \Gamma^A (\mu)^2 ]  = \iint d\mu (x) d\mu(y) G_{O} (x,y).$$
In other words, the process $\Gamma^A$ is a GFF in $O$. 

\begin {remark} One can reformulate the decomposition of $\Gamma_A + \Gamma^A$ in terms of conditional expectations. We have just proved that for any $\mu \in {\mathcal M}$, 
$$\Gamma_A ( \mu) = E [ \Gamma (\mu) | {\mathcal F}_A ]$$
almost surely.
\end {remark}

Let us summarise the above discussion with the following proposition.
\begin{proposition}
	[Markov property of the continuum GFF] 
	Let $D$ and $A$ satisfy the assumptions stated at the beginning of this section, and $\Gamma$ be a GFF in $D$. Then defining 
	\[\Gamma_A(\mu)= \Gamma ( \mu \mathbf{1}_A ) + \Gamma (\nu_{\mu, \partial O}); \;\;\; \Gamma^A(\mu)=\Gamma(\mu)-\Gamma_A(\mu)\] for all $\mu\in \mathcal{M}$ as above, one has that: 
	\begin{itemize}
		\item $\Gamma^A$ and $\Gamma_A$ are independent Gaussian processes;
		\item $\Gamma^A$ has the law of a Gaussian free field in $O$;
		\item there exists a version of $\Gamma_A$ such that $\Gamma_A|_{O}$ is almost surely equal to a harmonic function $h_A$ in $O$.	
	\end{itemize}	
\end{proposition}

\begin{figure}[h]
\centering
	\includegraphics[width=0.7\textwidth]{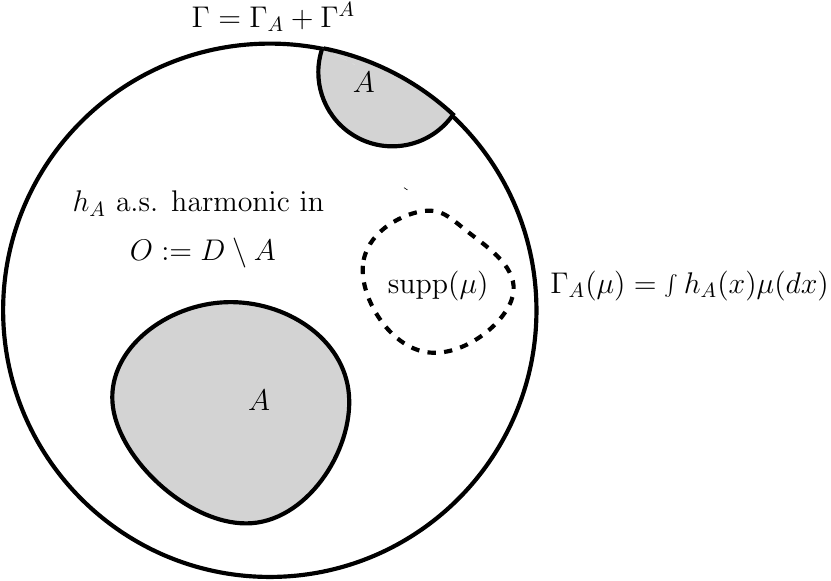}
	\caption {An illustration of the Markov property for a set $A\subset \overline D\subset \R^2$. The process $\Gamma^A$ is a GFF in $O=D\setminus A$, independent of $\Gamma_A$. }
\end{figure}

\begin {remark}
 Suppose now that $A$ and $A'$ are two compact subsets of $\overline D$ such that  $A \subset A'$, and such that both $D \setminus A$ and $D \setminus A'$ have a regular boundary. 
 Then it is a simple exercise, that we safely leave to the reader, to show that (almost surely)
$$ \Gamma^{A'} = (\Gamma^A)^{A'} \hbox { and } \Gamma_{A'} = \Gamma_A + (\Gamma^A)_{A'}.$$
In other words, $\Gamma$ can be decomposed into the sum of the three independent Gaussian processes:
$$\Gamma = \Gamma_A + (\Gamma^A)_{A'} + \Gamma^{A'}.$$
\end {remark}

\section {Local sets of the continuum GFF}

\subsection {Warm-up}
We are now ready to define {\em local couplings of random sets to the continuum GFF} (or in short, local sets of the GFF). 
Our random sets will be random compact subsets of $\overline D$, and we endow this space of compact sets with the usual Hausdorff metric. 
As in the discrete case, these will be random sets for which a strong version of the Markov property of the GFF can be made sense of.

\begin {definition}[Strong Markov property]\label{smp} 
Suppose that $A \subset \overline D$ is a random compact set such that $D \setminus A$ has a regular boundary. $A$ is said to 
satisfy the strong Markov property for the GFF $\Gamma$ 
if there exists a pair $(\Gamma_A, \Gamma^A)$ of processes such that: 
\begin {itemize}
\item the GFF $\Gamma$ is the sum of $\Gamma_A$ and $\Gamma^A$;
 \item the process $(\Gamma_A (\mu))_{\mu \in {\mathcal M}}$ is linear in with respect to $\mu$, and  
 {there exists a function $h_A$ in the complement of $A$ that is harmonic with probability one, and such that for any $\nu\in \mathcal{M}$ the equality $\Gamma_A(\nu)=\nu(h_A)$ holds on the event that the support of $\nu$ is contained in $D\setminus A$};
 \item conditionally on $(A,\Gamma_A)$, the law of the process $\Gamma^A$ is a GFF in $D \setminus A$. 
\end {itemize} 
\end {definition} 

\begin {remark}
 There is no difficulty in making sense of the last statement in this definition: the conditional law of $\Gamma^A$ given $(A,\Gamma_A)$ here is described explicitly in
terms of the random compact set $A$ (in a way that depends in a measurable way on $A$ via the Green's function in $D \setminus A$).
Also note that the conditional law of $\Gamma^A$ given $(A , \Gamma_A)$ is a function of $A$ only, so that 
conditionally on $A$, the fields $\Gamma_A$ and $\Gamma^A$ are independent. 
\end {remark}

\begin {remark}
Note that if $A$ satisfies the strong Markov property for $\Gamma$, then all the information about the joint distribution of $(A, \Gamma_A, \Gamma^A, \Gamma)$ is encapsulated in the joint distribution of $(A, \Gamma_A)$ 
(because we know the conditional distribution of $\Gamma^A$ given $(A, \Gamma_A)$).
\end {remark}

In the next paragraph, we will define (via approximation) what we will call \emph{local sets} of the GFF $\Gamma$, and then show that a set is local if and only if it satisfies the strong Markov property as in Definition \ref{smp}.
Having the equivalence between these two notions is useful. 
For instance, proving that the union of two conditionally independent local sets is a local set (recall Lemma \ref{disc_unionls}) is easier with the approximation approach. On the other hand, it is often easier to use the strong Markov property in order to actually show that a set is local.  

As a motivation for this section on local sets, and also as a warm-up for the coming chapters, we would like to stress that it is actually possible to define local sets (at least when $d=2$) such 
that the harmonic function $h_A$ is identically zero. The boundaries of such local sets are somewhat special, because seen from $O=D \setminus A$, $\Gamma$ is somehow 
``equal to $0$'' on $\partial O$ (if one keeps in mind the heuristic interpretation of $h_A$ as the harmonic extension of the values of $\Gamma$ on $\partial O$).

\subsection {Definition via dyadic approximations} 

We are going to define local sets in two stages: first we will define a notion of local sets that concerns random unions of $2^{-n}$-dyadic cubes, and then we will define (by approximation) the general notion of local sets. 

We will use the following terminology in the present section:  a set of the type $[i_1 2^{-n}, (i_1 +1) 2^{-n} ]  \times \cdots \times [i_d 2^{-n}, (i_d +1) 2^{-n} ]$ in $\R^d$ for $i_1, \ldots , i_d \in \Z$
will be called a $2^{-n}$ closed dyadic cube. 
Deterministic finite unions of closed dyadic cubes will be denoted by small letters $a$, $b$ etc., and we will be able to apply the results of the previous section (construction of $\Gamma_a$ etc.) to those. 
Capital letters, e.g. $A$, will be used to denote {\em random} closed subsets of the unit cube.

\begin {definition}[Dyadic local sets]
Suppose that $D$ is as in the previous section, so that one can define a GFF $\Gamma$ in $D$, and that $n \ge 1$ is fixed. 
We say that the random compact set $A \subset \R^d$ (defined on the same probability space as the GFF $\Gamma$) is a $2^{-n}$-dyadic local set if:
\begin {itemize}
\item it is the intersection of $D$ with a random finite union of closed $2^{-n}$-dyadic cubes (it can therefore take only countably many values);
\item for any deterministic finite union $a$ of closed $2^{-n}$-dyadic cubes, the GFF $\Gamma^a$ (in $D \setminus a$) is independent of 
$\sigma ( {\mathcal F}_a, \{ A=a \} ) = \sigma ( \Gamma_a, \{ A =a \})$.
\end {itemize}   
\end {definition}

This definition is of course reminiscent of the definition of local sets for the discrete GFF. Mind that (as in the discrete setting), this is a property of the joint distribution of $(A, \Gamma)$. 

Given that {any random set $A$ that is a finite union of $2^{-n}$-dyadic cubes} can take only countably many values, we can define without any problem 
\begin{equation} \label{eqn::2ncube_decomp} \Gamma_A := \sum_a \1{A=a} \Gamma_a, \ 
\Gamma^A := \sum_a \1{A=a} \Gamma^a, \ h_A := \sum_a \1{A=a} h_a.
\end{equation}
and note that with probability one, $h_A$ is in fact 
a random harmonic function in the random set $D \setminus A$. 
Furthermore, {if $A$ is a $2^{-n}$-dyadic local set}, then conditionally on the random set $A$ (and on $\Gamma_A$), $\Gamma^A$ is a GFF in $D \setminus A$. 

\begin{remark}\label{rmk:equivlscubes}
Conversely, if $A$ is a random union of $2^{-n}$-dyadic cubes that satisfies the strong Markov property (with associated processes $\tilde \Gamma_A, \tilde \Gamma^A$) 
then it follows from the definition that $A$ is a $2^{-n}$-dyadic local set of $\Gamma$. Moreover, $\tilde \Gamma_A, \tilde \Gamma^A$ must be equal to $\Gamma_A, \Gamma^A$ defined by \eqref{eqn::2ncube_decomp} (indeed, on the event $\{A=a\}$ one must have 
	$\tilde \Gamma_A(\mu)=\tilde \Gamma_A(\mu\mathbf{1}_{a})+\tilde  \Gamma_A(\mu\mathbf{1}_{D\setminus a})=\tilde \Gamma_A(\mu\mathbf{1}_{a})+\tilde \Gamma_A(\nu_{\mu, \partial (D\setminus a)})=\Gamma(\mu\mathbf{1}_{a})+\Gamma(\nu_{\mu, \partial (D\setminus a)})=\Gamma_a(\mu)$). 
	
In other words, a random union of $2^{-n}$-dyadic cubes is a local set if and only it satisfies the strong Markov property, and the decomposition $(\Gamma_A,\Gamma^A)$ in the strong Markov property is uniquely defined by \eqref{eqn::2ncube_decomp}.
\end{remark}

One can also define the $\sigma$-field ${\mathcal F}_A$ to be the $\sigma$-field generated by all events $U$ such that  
$ U \cap \{ A =a \} \in {\mathcal F}_a$ for all $a$ (where $\mathcal{F}_a$ is as defined as before, since $a$ is deterministic) -- this is very natural, given how one defines the stopped $\sigma$-algebra for Brownian motion. In other words, this is the set of events that can be decomposed as $\cup_a ( \{ A = a \} \cap U_a )$ with 
$U_a \in {\mathcal F}_a$ for each $a$. It is then immediate that $\Gamma_A(\mu)= \sum_a \1{A=a} \Gamma_a(\mu)$ is ${\mathcal F}_A$
measurable, and that $\Gamma_A  (\mu) = E ( \Gamma (\mu) | {\mathcal F}_A )$ almost surely.

Exactly as in the discrete setting, one can easily prove that:
\begin {enumerate}
\item 
if $A$ is a random union of $2^{-n}$ dyadic sets that is independent of $\Gamma$, then it is a $2^{-n}$-dyadic local set (of course these are not particularly interesting examples);
\item 
the definition of dyadic local sets is unchanged if we replace
$\sigma ( {\mathcal F}_a, \{ A=a \} )$ by $\sigma ( {\mathcal F}_a, \{ A \subset a \} )$ in the last line of the definition. To see this, simply use the fact that $\{ A \subset a \} =\cup_{a'\subset a} \{A = a'\}$ and recall that 
for $a' \subset a$, one can decompose $\Gamma^{a'}$ into the sum of the two independent processes $\Gamma^{a}$ and $(\Gamma^{a'})_a$.
\end {enumerate}

We now  define general local couplings. For any compact subset $A \subset \overline D$, we define its $2^{-n}$-dyadic approximation $A_n$ 
to be the intersection of $\overline D$ with the union of all closed $2^{-n}$ dyadic cubes that intersect $A$.  We then  let $O_n = D \setminus A_n$. Note that the sets $A_n$ are decreasing with $\cap A_n=A$, and that the sets $O_n$ are therefore increasing with $\cup O_n=O=D\setminus A$.

\begin{figure}[h]
\centering
	\includegraphics[width=0.35\textwidth]{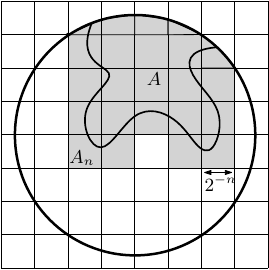}
	\caption {The dyadic approximation $A_n$ (shaded) to a compact set $A$.}
\end{figure} 

\begin {definition}[Local sets]
Let $A$ be a random compact set defined on the same probability space as a GFF $\Gamma$, and such that $D \setminus A$ has a regular boundary. 
We say that $A$ is locally coupled to $\Gamma$ (or equivalently that it is a local set of the GFF $\Gamma$) if for every $n \ge 1$, the set $A_n$ is a $2^{-n}$-dyadic local set of $\Gamma$.
\end {definition}
 
Of course, a rather unexciting class of of local sets are given by the random compact sets with regular boundary that are independent of $\Gamma$.

\subsection {Equivalence between the two notions} \label{subsec:equiv_markov_nots}

The goal of this section is to derive the following fact: 
\begin {proposition} 
\label {propequi}
$A$ is a local set for $\Gamma$ if and only if it satisfies the strong Markov property for $\Gamma$. 
\end {proposition}

\begin {proof}
Let us first assume that $A$ is a local set for $\Gamma$, and go about constructing $\Gamma_A$ (hence also $\Gamma^A:=\Gamma - \Gamma_A$). 
Notice that if one knows $A$, then one knows all the $A_m$ for $m\ge 1$, and conversely, if one knows all the $A_m$ for all  $m$ larger than some given $n$, 
then one also knows $A = \cap_{ m \ge n} A_m$. Hence, the $\sigma$-field generated by $A$ and the $\sigma$-field generated by $(A_m , m \ge n)$ coincide for all $n$. 

We now define for each $n$, the $\sigma$-field
$$ {\mathcal G}_n = \sigma ( A, \Gamma_{A_n}) = \sigma (A_n, A_{n+1}, A_{n+2}, \ldots, \Gamma_{A_n}).$$ 
Note that if one knows that $A_n=a$, $A_{n+1}=a'$ and $\Gamma_{a}$, then one also knows $\Gamma_{a'}$. Therefore we can equivalently write
$$ {\mathcal G}_n = \sigma ( A_n, A_{n+1}, A_{m+2}, \ldots, \Gamma_{A_n}, \Gamma_{A_{n+1}}, \ldots),$$
and see that ${\mathcal G}_{n+1} \subset {\mathcal G}_n$ is decreasing.
In fact, the decomposition   $\Gamma_{a} = \Gamma_{a'} + (\Gamma^{a'})_a$ shows that 
$$ \Gamma_{A_{n+1}} ( \mu) =  E ( \Gamma_{A_n} ( \mu ) | {\mathcal G}_{n+1} )$$
for any $\mu\in \mathcal{M}$, 
or in other words, that for any $n \ge 1$, 
$$ \Gamma_{A_n} (\mu) = E ( \Gamma (\mu ) | {\mathcal G}_n ). $$
This implies that 
 $(\Gamma_{A_n} (\mu), n \ge 0)$ is an inverse martingale for the inverse filtration $({\mathcal G}_n)_{n \ge 0}$.
 Note that the random variable $\Gamma (\mu)$ is Gaussian, and consequently lies in $L^p$ for any $1 \le p < \infty$. The inverse martingale therefore converges almost surely and in  $L^p$ (for any such $p$) to 
 $ E [ \Gamma_{A_1} ( \mu ) | {\mathcal G}_{\infty} ] $, where 
${\mathcal G}_\infty := \cap_n {\mathcal G}_n$.
We now {\em define} $\Gamma_A$ to be this limit: 
$$\Gamma_A (\mu) :=  E [ \Gamma_{A_1} ( \mu ) | {\mathcal G}_{\infty} ] \text{ for all } \mu \in \mathcal{M},$$
and also denote $\mathcal{F}_A:= \mathcal{G}_\infty$.
Finally, we define $\Gamma^A (\mu) := \Gamma (\mu) - \Gamma_A (\mu)$, which is therefore also the limit as $n \to \infty$ (almost surely and in any $L^p$) of $\Gamma^{A_n} (\mu)= \Gamma (\mu) - \Gamma_{A_n} (\mu)$. 

There is now a bit of work to be done to argue that $\Gamma_A$ restricted to the complement of $A$ corresponds to a harmonic function $h_A$.  Here is one way to proceed.
For any $z \in D$ and  $r< d (z, \partial D)$, we first define for each $n$ the average $\gamma_n (z,r)$ of $\Gamma_{A_n}$ on the sphere of radius $r$ around $z$. 
Note that when the the closure of the ball $B (z,r)$ is a subset of $O_n$, this spherical average is equal to $h_{A_n}(z)$ almost surely (because we know that $h_{A_n}$ exists and is almost surely harmonic). 

We have just seen (taking $\mu$ to be uniform measure on $\partial B(z,r)$) that for fixed $(z,r)$, $\gamma_n (z,r)$ converges almost surely and in any $L^p$ 
to the spherical average $\gamma_\infty (z,r)$ of $\Gamma_{A}$.
Hence, for any fixed $(z, r)$, 
if the define the event $E_r(z)$ that the closure of the ball $B(z,r)$ is in $D \setminus A$,  the sequence $\mathbf{1}_{E_r(z)} h_{A_n} (z)$ converges almost surely, and  in any $L^p$.
We then define $h_A (z) := \1{z \notin A}  \lim_{n \to \infty} h_{A_n} (z)$ (where the limit is therefore both an almost sure limit and a limit in any $L^p$). 
Observe that $\Gamma_A(\nu)=\nu(h_A)$ almost surely on the event that the support of $\nu$ is contained in $D\setminus A$, 
since $\Gamma_A(\nu)=\lim_{n \to \infty}\Gamma_{A_n}(\nu)=\lim_{n \to \infty} \nu(h_{A_n})$ on this event.

Let us now argue that this process $z \mapsto h_A (z)$ can be modified into a continuous and harmonic function on $D \setminus A$.    
We know that for all $z$, $z'$ that are at distance at least, say, $2r$ from the boundary of $D$ (using the $L^6$ convergence
of $\gamma_n (z,r)$ to $\gamma_\infty (z,r)$ and then the conditional Jensen inequality) that
\begin {eqnarray*}
 \lefteqn { 
 E \Bigl[ (\gamma_\infty (z,r) - \gamma_\infty (z',r))^6 \Bigr] 
 = \lim_{n \to \infty}   E \Bigl[ (\gamma_n (z,r) - \gamma_n (z',r))^6 \Bigr]  }\\
&& = \lim_{n \to \infty} E \Bigl[  E [\gamma (z,r) - \gamma (z', r)    | {\mathcal G}_n ]^6 \Bigr]   
 \le E \Bigl[ ( \gamma (z,r) - \gamma (z', r))^6  \Bigr]   \le C(r) |z-z'|^3.
\end {eqnarray*}
From this we deduce using Kolmogorov's criterion (the general case \eqref{eqn:kolgen}; here we do not have that $\gamma_{\infty}$ is Gaussian, which is why we use 6th rather than 2nd moments), that there exists (for each given $r$) a continuous modification $\tilde \gamma_\infty$ of $z \mapsto \gamma_\infty (z,r)$.
For all given $z \in D$ and $r > 0$,  one has $h_A (z) = \tilde \gamma_\infty (z,r)$ almost surely on the event where $d(z, A) > 2r$. 
Since this is true for all arbitrarily small rational $r$,  we deduce that there exists a version of $z \mapsto h_A (z)$ that is continuous on $D \setminus A$.

We now need to show that this function $h_A$ is necessarily harmonic on $D \setminus A$. One way to see this is as follows. The goal is to check that for a given $z \in D$ and $r >0$, almost surely on the event $E(z,r)$ where 
$z \in O$ and $3r < d(z, \partial O)$ (so that $d(z',\partial O)>2r$ for all $z'\in \overline{B}(z,r)$), $h_A(z)$ is equal to the mean value of $h_A$ on the sphere $\partial B (z,r)$. Since we know from the definition that 
$$h_A(z)= \gamma_\infty (z,r) = \lim_{n \to \infty} \gamma_n (z,r) = \lim_{n \to \infty} \int d\lambda_{z,r} (y) h_{A_n} (y) $$ 
almost surely on this event, it is sufficient to have that 
\begin {equation}
\label {CVp}
  \mathbf{1}_{E(z,r)} \int d\lambda_{z,r} (y) h_{A_n} (y)  \to  \mathbf{1}_{E(z,r)} \int d\lambda_{z,r} (y) h_{A} (y) 
  \end {equation}
in probability. 
On the one hand, we have seen that almost surely with respect to the product probability measure $P \otimes \lambda_{z,r}$, 
$\mathbf{1}_{E(z,r)} h_{A_n} (y) \to \mathbf{1}_{E(z,r)} h_A (y)$. 
But we also have the a priori bound, 
$$ E [ \mathbf{1}_{E(z,r)} \int d\lambda_{z,r} (y) (h_{A_n} (y))^2 ] 
\le E [   \int d\lambda_{z,r} (y) (\gamma (y,r)^2) ] \le  \int d\lambda_{z,r}(y) E [ \Gamma (\lambda_{y,r})^2 ]$$ 
which is finite. Hence, since a sequence of random variables that is bounded in $L^2$ and that converges almost surely also converges in $L^1$, 
we can conclude that the convergence (\ref {CVp}) holds in $L^1$, and therefore also in probability.

 Let us now  argue that the conditional law of $\Gamma^A$ given $(A, \Gamma_A)$ is that of a GFF in $O$. 
 This follows from the fact that process $\Gamma^A$ is the almost sure limit of $\Gamma^{A_n}$ as $n \to \infty$, where conditionally on the event $\{A_n = a\}$, 
 $\Gamma^{A_n}$ is a GFF in $D\setminus A_n$ that is independent of {$\mathcal{G}_n\supseteq \mathcal{G}_\infty \supseteq \sigma(A, \Gamma_A)$}. To see the independence between $\Gamma^{A_n}$ and $A_{m}$ for $m\ge n$, note that conditionally on $\{A_{m}=\tilde{a}\}$, $\Gamma^{A_m}$ is independent of $A_m$ and $\Gamma^{A_n}=\Gamma^a=(\Gamma^{\tilde{a}})^a$ is measurable with respect to $\Gamma^{A_m}$.
Hence, the local set $A$ satisfies the strong Markov property.

Finally, it remains to show the converse statement, namely that if $A$ satisfies the strong Markov property,
 then $A_n$ is a $2^{-n}$-dyadic local set for each $n\ge 1$. For this, we use that conditionally on $(\Gamma_A, A)$, $\Gamma^A$ has the law of a GFF in $O=D \setminus A$, and that for each $n$, the 
set $A_n$ is a deterministic function of $A$. This means that conditionally on $(\Gamma_A,A)$, $\Gamma^A$ can be further decomposed as
$$ (\Gamma^A)^{A_n} + (\Gamma^A)_{A_n},$$
using the usual Markov property (and the corresponding notation) for the GFF $\Gamma^A$. Then, we see that conditionally on $(\Gamma_A, A)$ and $(\Gamma^A)_{A_n}$, the process
$(\Gamma^A)^{A_n}$ is a GFF in $O_n$, and that $\Gamma_A + (\Gamma^A)_{A_n}$ when restricted to $O_n$ is a harmonic function (as the sum of two harmonic functions). {Remark \ref{rmk:equivlscubes} then implies that $A_n$ is a $2^{-n}$ dyadic local set.}
\end {proof}

{\begin{remark}\label{dmpunique}The above proof along with Remark \ref{rmk:equivlscubes} implies that the strong Markov decomposition in Definition \ref{smp} is actually unique.
	\end{remark}}

Let us note the following immediate facts.
\begin {itemize}
 \item Suppose that $A$ is a local set that is coupled with a GFF $\Gamma$. If we observe $A$ and $\Gamma$, then (since $\Gamma_A$ is a limit of the $\Gamma_{A_n}$) we know everything about the process $\Gamma_A$, and therefore also about the harmonic function $h_A$. 
 \item If $A$ and $B$ are two local sets, and $z$ is some fixed point in $D$, let $O_A (z)$ and $O_B (z)$ denote
 the connected components of $D \setminus A$ and $D \setminus B$ (if they exist) containing $z$.
 Then, almost surely on the event that $\{ O_A (z) = O_B (z) \}$, one has $h_A (z) = h_B (z)$ (again just using the definition of $h_A(z)$ as a limit of $h_{A_n}(z)$).
\end {itemize}

\subsection {Unions of conditionally independent local sets} 

Let us now explain how to derive the following continuum counterpart of Lemma \ref{disc_unionls} for the discrete GFF.

\begin {proposition}[Unions of local sets] \label{prop:union}
If $A$ and $A'$ are both locally coupled to the GFF $\Gamma$, and if they are conditionally independent given $\Gamma$, then $A \cup A'$ is also locally coupled to $\Gamma$. 
\end {proposition}

\begin {remark} 
The proof of this result turns out to be rather easy using the definition of local sets via dyadic approximation. It would have been more of a challenge to derive this directly using the ``strong Markov property'' characterisation of local sets. 
\end {remark}

\begin {proof} 
Let us first show the result when $A$ and $A'$ are two dyadic local sets; this part of the  proof is almost a copy-and-paste of the proof in the discrete case. 

Let $b$ be a finite, deterministic union of $2^{-n}$-dyadic cubes. For any given measures $\mu_1, \ldots, \mu_m$, $ \nu_1, \ldots, \nu_{m'}$ and open sets $U_1, \ldots, U_m$ and $ V_1, \ldots, V_{m'}$ in $\R$, we define the $\sigma (\Gamma^b)$ and $\sigma (\Gamma_b)$  measurable events 
$$ 
U^b = \{ \forall j \le m, \ \Gamma^b (\mu_j) \in U_j \} \hbox { and } 
V_b = \{ \forall j \le m', \ \Gamma_b (\nu_j) \in V_j \}.
$$
Note that the set of events $U^b$ (with $\mu,\nu,U,V$ varying) is stable under finite intersections and generates $\sigma (\Gamma^b)$, and that the family of events $V_b$ 
is stable under finite intersections and generates $\sigma ( \Gamma_b)$.
Since it is clear that $A\cup A'$ is a finite union of $2^{-n}$ dyadic cubes, to prove the lemma it is sufficient to show that
\begin{equation}\label{unionls}
{ \P \Bigl[ {U^b}, \ {V_b }, \ A \cup A' = b \Bigr] } 
= \P ( U^b ) \times  \P \Bigl[ V_b , \ A  \cup A' = b \Bigr] .
\end{equation}
Note that for any $a\subset b$ the family of events $V_b$, $V_b \cap \{ A =a \}$ 
is stable under finite intersections and generates $\sigma ( \Gamma_b,\{ A =a \} )$
and a similar statement holds replacing $A$ by $A'$.
 Then, for all $a$ and $a'$, with $a \cup a'=b$, we have that
\begin{eqnarray*}
\lefteqn { \P \Bigl[ U^b , \ V_b , \ A = a, \ A' = a' \Bigr] 
 =  \E \Bigl[ \P ( U^b , \ V_b , \ A = a, \ A' = a' \mid \Gamma )  \Bigr]  } \\
&&=   \E \Bigl[  \mathbf{1}_{U^b, V_b } \P(  A  = a , \ A' = a' \mid \Gamma )  \Bigr] =   \E \Bigl[  \mathbf{1}_{U^b} \mathbf{1}_{V_b} \P(  A  = a \mid \Gamma) P( A' = a' \mid \Gamma )  \Bigr].
\end{eqnarray*}
However we know that $\Gamma^b$ is independent of $\sigma ( \Gamma_b, \1{ A=a })$ (since $a \subset b$), from which it follows that 
$$ \P(  A  = a \mid \Gamma) = \P ( A  = a \mid \Gamma_b )$$
is measurable with respect to  $\sigma (\Gamma_b$), and that the same is true for $\P ( A' = a'  \mid \Gamma)$. 
Hence, since $\Gamma_b$ and $\Gamma^b$ are independent, it follows that 
\begin {eqnarray*}
 { \P \Bigl[ U^b, \ V_b , \ A  = a, \ A'  = a' \Bigr] } 
&=& \P [ U^b ] \times  \E \Bigl[ \mathbf{1}_{V_b } \P(  A  = a \mid \Gamma) P( A' = a' \mid \Gamma )  \Bigr]\\
&=& \P [ U^b ] \times  \P \Bigl[ {V_b }, \ A = a , \ A' = a' \Bigr]  
\end {eqnarray*}
Summing over all $a$ and $a'$ such that $a \cup a' =b$ we can deduce \eqref{unionls}, and therefore that $A \cup A'$ is a $2^{-n}$-dyadic local set. 

Suppose now that $A$ and $A'$ are general local sets and that $n \ge 1$. We have just proved (because $A_n$ and $A_n'$ are $2^{-n}$-dyadic local sets that are conditionally independent given $\Gamma$) that $A_n \cup A_n'$ is also a $2^{-n}$-dyadic local set. Since $(A \cup A')_n= A_n \cup A_n'$ by definition, we therefore have that $(A\cup A')_n$ is a $2^{-n}$-dyadic local set, and this concludes the proof. 
\end {proof} 
\begin{remark}\label{rmkforbcs}
		Let $A,A'$ be as in Proposition \ref{prop:union}. The above proof, together with the same reasoning as in the penultimate paragraph of the proof of Proposition \ref{propequi}, shows that even conditionally on $(A,A',\Gamma_{A\cup A'})$ - rather than just conditionally on $(A\cup A', \Gamma_{A\cup A'})$ - $\Gamma^{A\cup A'}$ is a GFF in $D\setminus \{A\cup A'\}$.
\end{remark}
Note that in the previous proposition,  we did not describe the harmonic function $h_{A \cup A'}$ in terms of $h_A$ and $h_{A'}$. This is a trickier issue than it appears at first glance. 
Intuitively, one would like to say that the ``boundary conditions'' of $h_{A \cup A'}$ are just given by those of $h_A$ on $\partial A$, and by those of $h_{A'}$ on $\partial A'$, but putting this 
on a rigorous footing is delicate.
We will do this in a special case, when $d=2$, in Section \ref {furthersection}.

\subsection {Thin local sets} 

In this section we describe a particularly useful class of local sets.
\begin {definition}[Thin local sets]
The local set $A$ is said to be thin, if for all bounded functions $f$ with compact support in $D$, 
$\Gamma_{A_n} (f \mathbf{1}_{A_n})$ converges to $0$  in probability as $n \to \infty$. 
\end {definition}
This means in particular that for all such $f$, 
$$\Gamma_A (f) = \lim_{n \to \infty}  \Gamma_{A_n} (f\mathbf{1}_{O_n})$$ 
in probability. 
But it is easy to see that conditionally on $A$, 
$$ \Gamma_A (f \mathbf{1}_{O_n}) - \Gamma_{A_n} (f\mathbf{1}_{O_n}) 
=  -\Gamma^A ( f \mathbf{1}_{O_n}) + \Gamma^{A_n} (f\mathbf{1}_{O_n}) 
= (\Gamma^A)_{A_n} (f\mathbf{1}_{O_n}) \to 0 $$ 
in probability as $n \to \infty$. 
In particular, this shows that 
$$\Gamma_A (f)=\lim_{n \to \infty}\Gamma_A (f \mathbf{1}_{O_n})  = \lim_{n \to \infty}\int_{O_n} h_A (x) f(x) dx, $$ which 
is a function  of $A$ and $h_A$. Hence, all $\Gamma_A (f)$ (and therefore the whole process $\Gamma_A$) can be recovered from the knowledge of $(A, h_A)$.  
So, for thin local sets, $\Gamma_A$ carries no more information than $h_A$.

Note that a deterministic compact set of zero Lebesgue measure is always a thin local set. 
One may wonder whether a local set is thin as soon as it has zero Lebesgue measure.
As we shall
mention later in these lecture notes, this turns out not to be the case: 
there exists local sets $A$ that are not thin, but have Lebesgue measure zero almost surely. 
The goal of the next few paragraphs is to describe  simple criteria ensuring that a local set is thin.  
\medskip

We start with the following criterion when $d=2$: 

\begin {proposition}[Small local sets are thin ($d=2$ case)]
\label{prop:thincrit_d2}
Suppose that $D\subset \R^2$ is bounded, and that $A$ is a local set of a GFF $\Gamma$ in $D$.
Define $|A_n|$ to be the Euclidean area of $A_n$ (which is $4^{-n}$ times the number of closed
$2^{-n}$-dyadic squares that $A$ intersects). If there exists a sequence $n_k \to \infty$ such that 
almost surely, $|A_{n_k}|= o(1/n_k)$, then the local set $A$ in thin.
\end {proposition}

\begin {proof} 
Let $S$ denote any closed $2^{-n}$-dyadic square contained in $D$. 
We can bound $E [  \Gamma (f \mathbf{1}_S)^2 ]$ (via the double integral of $G_D (x,y)$) and see that there exists $C$ such that for any $n$, and for any $S$, 
$$ E [  \Gamma_{A_n} (f\mathbf{1}_S)^2 ]  \le E [  \Gamma (f \mathbf{1}_S)^2 ] \le C^2   \| f \|_\infty^2 n 4^{-2n}.$$
But, using the fact that $\Gamma_{A_n} (f\mathbf{1}_S)$ is a Gaussian random variable, we get the tail estimate  
$$ 
P \Bigl[ |\Gamma_{A_n} ( f\mathbf{1}_S)^2 | > C \| f \|_\infty  M \sqrt {n} 4^{-n} \Bigr] \le \exp ( - M^2 / 2)
$$
for all large enough $M$. 
Summing this over all $O(4^n)$ of the  $2^{-n}$-dyadic squares in $D$, we see that for each $n$, 
the probability that there exists one or more such squares $S$  for which $|\Gamma (f\mathbf{1}_S)| > CM \sqrt {n} 4^{-n}$ is bounded by $4^n \exp (-M^2 / 2)$. 
If we choose $M= M(n)= x\sqrt {n}$, then  for some fixed large enough $x$, this bound decays exponentially in $n$. Hence by Borel--Cantelli, we know that almost surely, for all large enough $n$, 
for all the $2^{-n} $ dyadic squares $S$, 
$$| \Gamma_{A_n} (f\mathbf{1}_S) | \le Cx n 4^{-n}.$$  
We then conclude using the fact that  $A_{n_k}$ is almost surely the union of $o( 4^{n_k} / n_k)$ such squares. 
\end {proof} 

The following corollary provides a condition that is easier to check in practice.
\begin {corollary} 
If $A$ is a local set such that $E[ | A_n | ] = o( 1/n)$, then it is a thin local set.
\end {corollary}

\begin {proof}
We can find $n_k \to \infty$ such that $E [ |A_{n_k} | ] \le 1/(k^2 n_k)$, and by the Borel-Cantelli lemma, we see that the previous criterion is satisfied for this choice of the sequence $n_k$.  
\end {proof}

It can also be useful to define thin local sets when $D \subset \R^2$ is not bounded. In this case, we can just use conformal invariance. For instance, we will say that 
a local set $A$ in the upper half-plane is a thin local set if its image under the map $z \mapsto (z-i)/(z+i)$ is a {thin} local set in the unit disk. 

Later on we shall be interested in some very particular thin local sets of the GFF, for which  $h_A$ is actually a function of $A$. This property is satisfied, for instance, by some local sets defined by SLE-type curves in two dimensions. 

\begin {remark}
\label {confradremark}
In Chapter \ref {Ch5} we will discuss one special thin local set of the GFF, for which the harmonic function $h_A$ can only take
the values $a$ and $-a$ (where $a$ is a fixed constant). This means in particular that $h_A$ is constant in each of the connected components of $D \setminus A$. Moreover, for each fixed $z$ with the property that $z \notin A$ almost surely, one has that $P [ h_A (z) = a ] = P [ h_A (z) = -a] = 1/2$. Indeed, this follows by considering the expectation of the spherical average $\Gamma(\lambda_{z,\eps})$ at radius $\eps$ around $z$, and letting $\eps \to 0$.

In fact, we can further note that when $\eps$ is very small, the difference between the variances of the two Gaussian random variables 
$\Gamma (\lambda_{z,\eps})$ and $\Gamma^A ( \lambda_{z, \eps})$ is equal to the limit as $y \to z$ of 
$G_D (z,y) - G_{D \setminus A} (z,y)$. If we denote this quantity by $C(z,A,D)$, and note that $\Gamma (\lambda_{z,\eps})-\Gamma^A ( \lambda_{z, \eps})$ is bounded and converges almost surely to $h_A(z)$ as $\eps\to 0$, it follows readily that $C(z,A,D)$ has the law of the exit time from $[-a,a]$ by a standard one-dimensional Brownian motion.
In the two-dimensional case this means that the expected value of the difference between the log-conformal radius of $D$ and of $D \setminus A$ at $z$ is equal to $2\pi a^2$.  (This difference is exactly equal to $2\pi C(z,A,D)$ by definition of $G_D$ and of the log conformal radius -- see Chapter 3.)
\end {remark}

\medbreak 
Finally, let us state the corresponding criterion for a local set to be thin, when one considers a GFF in dimension $d \ge 3$: 
\begin {proposition}[Small sets are thin ($d \ge 3$)]
Suppose that $D\subset \R^3$ is bounded and that $A$ is a local set of a GFF $\Gamma$ in $D$.
Define $|A_n|$ to be the $d$-dimensional Lebesgue measure of $A_n$ (which is $2^{-nd}$ times the number of closed $2^{-n}$-dyadic cubes that $A$ intersects). If almost surely, 
$|A_n|= o(4^{-n})$, then the local set $A$ is thin.
\end {proposition}

The proof is almost identical to the case $d=2$ and left to the reader. It says in particular that if the Minkovski dimension of a local set $A$ is smaller than $d-2$, then it is a thin local set.  

\subsection {Some further features of local sets in two dimensions} 
\label {furthersection}

We are now going to derive some further results for local sets and unions of (conditionally) independent local sets.
We choose to describe only the results that will be actually used later on in these notes, and do not strive for the most general statements.
So, even though parts of this section would also work in some $d$-dimensional domains, we will restrict ourselves 
to the two-dimensional setting. The main goal here is to derive Proposition~\ref {harmfunctionunion} below, that describes the harmonic function $h_{A \cup A'}$ of Proposition~\ref {prop:union} 
in certain special cases.

\medbreak

(1) Let us first consider the case when $D$ is a connected open subset of the unit disc $\U$, such that $D$ contains a neighbourhood (in $\U$) of some point on $\partial \U$ (in 
particular $\partial D$ contains some open arc of the unit circle).
Of course, one particular case is when $D= \U$. 
Suppose that $o$ is a deterministic open subset of $D$ and let $a$ be the closure of $D \setminus o$. 
For each $\eps >0$, we define the following sets: 
$\partial_\eps$ will denote the part of $\partial o \cap \partial \U$ that is at distance greater than $\eps$ from $a \cup (\U \setminus D)$;  $u_\eps$ will denote the subset of points in $o$ 
that are at distance smaller than $\eps /4$ from $\partial_\eps$; and finally, $v_\eps$ will denote the union $u_\eps \cup \partial_\eps$. 

Let $\Gamma$ be a zero boundary Gaussian free field in $D$.
Recall that $h_a$ is the harmonic function in $o$ defined by $h_a(z) = \Gamma ( \nu_{z, \partial o})$. Our first observation is the following:
\begin {lemma} \label{harmoniccontinuity_det}
Almost surely, the function $h_a$ can be extended by continuity to be equal to $0$ on $\cup_{\eps >0} \partial_\eps$.
\end {lemma}
In other words, for all $\eps > 0$, the function $h_a(z)$ almost surely tends to $0$, uniformly as the distance between $z$ and $\partial_\eps$ tends to $0$.
\begin {proof}
Let us fix an arbitrary positive $\eps$. It is sufficient to check that if we define $h_a$ to be equal to $0$ on $\partial_\eps$, then there exists a version of $h_a$ that is continuous on $v_\eps$
(because we already know continuity in $o$).  
We will show this using the version of Kolmogorov's criterion for Gaussian processes, Lemma \ref{kolmo}, simply using the explicit expression for $E [ (h_a(z)- h_a (z'))^2]$. 

Recall from \eqref{eqn:varhA}, by harmonicity of the Green's function, that
 $$ E[(h_a(z)-h_a(z'))^2]  =  \int_{\partial o}(\nu_{z,\partial o}(dy)-\nu_{z',\partial o}(dy))(G_D(z,y)-G_D(z',y))$$ 
 for all $z,z'$ in $o$ (recall that $G_D (z,y) = 0$ as soon as $z$ or $y$ are on $\partial \U$). Note in particular that this identity is still valid 
 when $z$ and/or $z'$ are on $\partial \U$. 
 We therefore have the bound 
 $$ E[(h_a(z)-h_a(z'))^2] 
 	\le  (p(z,o)+p(z', o)) \times \sup_{y\in a} |G_D(z,y)-G_D(z',y)| $$
 	where $p(z,o)$ denotes the probability that a Brownian motion starting at $z$ exits $o$ strictly before exiting $\U$.
 
 When $z \in u_\eps$,  this probability $p(z,o)$ can be bounded  by the probability that a Brownian motion starting from $z$ reaches the circle of radius $\eps / 2$ around $z$ before hitting some given line
 (tangent to $\partial \U$) that is at distance $d(z, \partial \U)$ from $z$ (recall that this distance is smaller than $\eps/4$ by definition of $u_\eps$). 
 This easily implies that $p(z,o)$ is  bounded above by a universal constant times $d(z, \partial \U)  / \eps$. 

 Let us now show that for some constant $C(\eps)$ depending only on $\eps$, 
 \begin {equation}
 \label {GDz}
 |G_D(z,y)-G_D(z',y)| \le C(\eps) |z-z'| / d(z, \partial \U)\end {equation}
  whenever $z,z' \in u_\eps$ and $y \in a$ (so that $y$ is at distance at least $\eps/2$ from $z$ and $z'$). 
 To see this, recall the mirror coupling of two Brownian motions started from $z$ and $z'$ (see the paragraph preceding \eqref{eqn:varhA}), and denote such a coupling by $(B^1,B^2)$. Then by harmonicity of the Green's function, one can rewrite the difference $G_D(z,y)-G_D(z',y)$ as the expectation $E[G_D(B_\tau^1,y)-G_D(B_\tau^2,y)]$, where $\tau$ is the first time that either $B^1$ or $B^2$ leaves the ball of radius $d(z,\partial \U)$ around its starting point. Observe that, on the event that the two Brownian motions couple before time $\tau$, the quantity in the expectation is 0. Moreover, one can show that the probability of this \emph{not} occurring is bounded above by  $|z-z'| / d(z, \partial \U)$. Finally, to deal with this complementary event, we can note that $G_D (z'',y)\le  G_\U (z'', y)$ and that $G_\U (z'' , y)$ is bounded by some absolute constant whenever $| y - z'' | \ge \eps/2$. This yields (\ref {GDz}).

Combining these estimates, we obtain that  
	\[E[(h_a(z)-h_a(z'))^2]\le C'(\epsilon)|z-z'|\]
	for all $z,z'$ in ${v_\eps}$ (the case where $z$ and/or $z'$ is on $\partial_\eps$ is also easily treated), for some constant $C'(\epsilon)$ depending only on $\epsilon$.
{Since $h_a$ is a Gaussian process, we can then conclude using the continuity criterion for Gaussian processes (Lemma \ref{kolmo}). }
\end {proof}

\medbreak

(2) We now consider a variant of the above where $a$ is replaced by a local set $A$ of the GFF $\Gamma$ in $D$ (with the same conditions on $D$ as before). 
We define  $O := D \setminus A$, and for each $\eps > 0$, we now define  
$\partial_\eps= \{ z \in \partial \U , \ d(z, A \cup (\U \setminus D)) > \eps \}$, we let $U_\eps$ be the set of points in $O$ that lie at distance 
smaller than $\eps /4$ of $\partial_\eps$, and set $V_\eps := U_\eps \cup \partial_\eps$. We also assume that $O$ 
is almost surely connected. 
Then, we have the following generalisation of the previous lemma for local sets:

\begin {lemma} 
\label {harmoniccontinuity}
Almost surely, the  harmonic function $h_A$ can be extended by continuity to be equal to $0$ on $\cup_{\eps >0} \partial_\eps$.
\end {lemma}
In other words, for all $\eps > 0$, the function $h_A(z)$ almost surely tends (uniformly) to $0$ as the distance between $z$ and $\partial_\eps$ tends to $0$. 
\begin {proof} 
We fix $\eps$ and choose $m$ with $2^{-m} < \eps /8$. 
 First, we note that Lemma \ref{harmoniccontinuity_det} applied to each of the finitely many possible options for $A_m$ 
 shows that for each $m$, the dyadic local set $A_m$ does satisfy the conclusions of Lemma \ref{harmoniccontinuity}. This means that for each $m$, we get the existence of an almost surely continuous extension of $h_{A_m}$ to $V_\eps$. 

 The idea is now to construct the Markovian decomposition $(\Gamma_{A_m}, \Gamma^{A_m})$ (and in particular the function $h_{A_m}$) in two steps. First we discover $A$ and $\Gamma_A$; 
  we know that conditionally on $(A,\Gamma_A)$, $\Gamma^A=\Gamma-\Gamma_A$ is a GFF in the complement $O$ of $A$, that is also independent of $A_m$ (since $A_m$ is a deterministic function of $A$). This means that we can decompose $\Gamma^A=(\Gamma^A)_{A_m}+(\Gamma^A)^{A_m}$ on the set $A_m$, and for this decomposition, the conclusion of Lemma~\ref{harmoniccontinuity_det} does still hold. In particular, writing $h^A_{A_m}$ for the restriction of $(\Gamma^A)_{A_m}$ to the complement of $A_m$, we have that $h_{A_m}^A$ extends continuously to $V_\eps$. Note that on the event that $O$ does not contain a neighbourhood of any point in $\partial \U$, the conclusion of this lemma is trivial.
  
Moreover, by uniqueness of the Markov decomposition, it must be that $\Gamma_{A_m}=\Gamma_A+(\Gamma^A)_{A_m}$ and so one can almost surely write $h_{A_m}=h_A+h_{A_m}^A$.
 In other words, when restricted to the complement of $A_m$ (which in particular contains $U_\eps$), we have $h_A = h_{A_m} - h_{A_m}^A$. Since $h_{A_m}$ and $h_{A_m}^A$ can almost surely be extended by continuity to 
 $V_\eps$, this shows that $h_A$ can almost surely be extended by continuity as well.
\end {proof}

\medbreak 

(3) We now finally turn our attention to the description of $h_{A \cup A'}$ when $A$ and $A'$ are two conditionally independent local sets. 
The previous results allow us to derive the following useful fact. 
\begin {proposition}
\label {harmfunctionunion}
Consider a GFF in the unit disk $\U$. Suppose that $A$ and $A'$ are two conditionally independent local sets of $\Gamma$, such that all connected components of 
$O=\U \setminus A$ and $O' = \U \setminus A'$ 
are simply connected. We know from Proposition \ref{prop:union} that $A'':= A \cup A'$ is a local set of $\Gamma$.  Then the harmonic function $h_{A \cup A'}=h_{A''}$ defined on
$O'':= {\U} \setminus  A''$  almost surely satisfies that for all $\eps > 0$, 
\begin {itemize}
\item $(h_{A''} - h_A)(z)$ goes uniformly to $0$ when $z \to \partial O''$ with $d(z,A') >  \eps $. 
\item $(h_{A''} - h_{A'})(z)$ goes uniformly to $0$ when $z \to \partial O''$ with $d(z,A) >  \eps $.
\end {itemize}
\end {proposition}
\begin {remark}
Note that on the event where the two local sets $A$ and $A'$ are disjoint, then these two conditions do fully characterise $h_{A''}$. Indeed, if $h$ was another harmonic function in $O''$ 
that satisfied these conditions,
then $h(z)- h_{A''} (z)$ would tend uniformly to $0$ as $z \to \partial O''$, which in turn would imply by the maximum principle that $h= h_{A''}$. This will be very useful later on.  
\end {remark}

\begin{figure}[h]
\centering
	\includegraphics[width=0.45\textwidth]{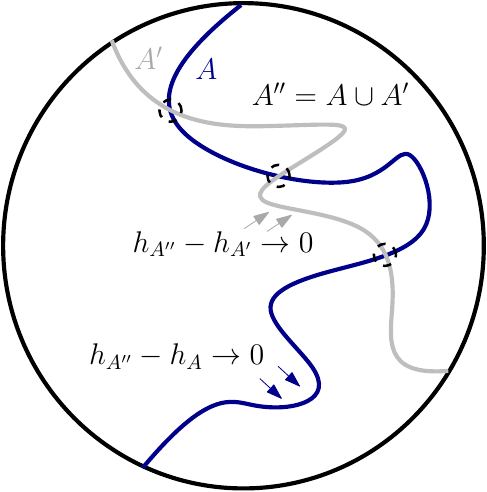}
	\caption {Two local sets $A$ and $A'$ in $\U$ given by two curves, with associated harmonic functions $h_A$ and $h_{A'}$. The harmonic function $h_{A''}$ is associated with the local set $A''=A\cup A'$. Proposition \ref{harmfunctionunion} provides information about the behaviour of $h_{A''}$ near the boundary of $O''=\U\setminus A''$ in terms of $h_A$ and $h_{A'}$.}
\end{figure} 

\begin {proof}
Since $A$ and $A'$ play symmetric roles, it suffices to prove the statement about $h_{A''}- h_{A}$.
The rough strategy will be to show that the closure of $A' \setminus A$ is in fact a local set for the GFF $\Gamma^A$, and to then apply Lemma~\ref {harmoniccontinuity} and conformal invariance to conclude. 
However, we will circumvent the question of making sense of local sets in random domains, by using conformal invariance.

 We first note that for any $\eps >0$, there almost surely exist only finitely many connected components of $O''$ that contain a ball of radius $\eps/2$.
Therefore, it suffices to show the uniform convergence for any single connected component of $O''$. For given $z_0\in \U$, write $O''(z_0)$ for the connected component of $O''$ (when it exists) that contains $z_0$.
Rephrasing the previous statement, it is sufficient to show that for any given $z_0 \in \U$, almost surely, the function 
$\1{z_0 \in O''} ( h_{A''} (z)- h_A (z))$ goes uniformly to $0$ when $z \to \partial O''$ with $d(z,A') >  \eps $ and $z\in O''(z_0)$.

Let us now fix $z_0 \in \U$. We denote by $O (z_0)$ the connected component of $O= \U \setminus A$ that contains 
$z_0$. We will assume that the probability of $O(z_0 )$ being non-empty is positive (otherwise, the statement is obvious). 
 When $O(z_0) \not= \emptyset$, let us define the conformal map $\Phi = \Phi_{A, z_0}$ from $O(z_0)$ onto the unit disk such that $\Phi (z_0) = 0$ and $\Phi' (z_0) \in \R_+$. The fact that $A$ is a local set shows that 
 conditionally on the event $O( z_0) \not= \emptyset$, the image under $\Phi$ of the field $\Gamma^A$ (as described in Section \ref {CIGFF}), restricted to $O(z_0)$, is a GFF in $\U$ that is (conditionally) independent of $(A, \Gamma_A)$. We will denote this GFF by $\tilde{\Gamma}$.
By possibly extending the probability space, we can also define $\tilde \Gamma$ to be a GFF that is conditionally independent of $A$ and $\Gamma$ on the event that $O(z_0) = \emptyset$. Finally, we define $\tilde{O}$ to be the set $\Phi (O' \cap O  (z_0) )$ on the event where $z_0 \in O$, and when $z_0 \notin O$, we just set $\tilde O := \emptyset$. We write
$\tilde{A}$ for closure of $\U \setminus \tilde{O}$.   

So, now we have a GFF $\tilde \Gamma$ defined on the unit disk, and a random subset $\tilde A$ of the unit disc. 
The next step, which is the key to the proof is to show that $\tilde A$ is actually a local set of $\tilde \Gamma$, and that on the event where $\tilde O$ is not empty,
 $h_{\tilde A} = (h_{A''} - h_A) \circ \Phi^{-1}$ in $\tilde O$. 

The proposition will then follow because:
	\begin{itemize}
		\item we can apply Lemma \ref{harmoniccontinuity} to the local set of $\tilde{\Gamma}$ given by the complement of the connected component of $\tilde{O}$ containing $0$;
		\item by Koebe's quarter theorem (for instance), we almost surely have that $d(z,\partial \U)\to 0$ with $z\in \tilde{O}$ if and only if $d( \Phi^{-1} (z), \partial O(z_0))$ tends to $0$ with $\Phi^{-1}(z)\in O'$. 
	\end{itemize}
	
Let us now finally show that $\tilde A$ is indeed a local set of $\tilde \Gamma$ by giving its Markovian decomposition. 
First note that when restricted to $O(z_0)$, the field $\Gamma - \Gamma^A$ is the harmonic function $h_A$, and that when restricted to $O''$, the field 
$\Gamma - \Gamma^{A''}$ is the harmonic function $h_{A''}$. Hence, when restricted to $O''$, the field $\Gamma^{A''}- \Gamma^A$ is the 
harmonic function $h_{A} - h_{A''}$. By taking the image under $\Phi$, we then get that when restricted to $\tilde O$, the 
field $\tilde \Gamma - ( \Gamma^{A''} \circ \Phi^{-1} ) $ is the harmonic function $H:= ( h_{A''} - h_A) \circ \Phi^{-1}$.

We can recall from Remark \ref {rmkforbcs} that conditionally on $(A, A', \Gamma_{A''})$, the GFF $\Gamma^{A''}$ is a GFF in the random set $O''$. Since $\Phi$ is measurable with respect to $A$, we get that 
on the event where $O''$ is not empty, the field $\Gamma^{A''} \circ \Phi^{-1} $ is a GFF in the random set $\tilde O$ that is (conditionally) independent of $(A, A')$ and $\Gamma_{A''}$.

We therefore get the Markovian decomposition of $\tilde \Gamma$ into the sum of a GFF in $\tilde O$ with a field that coincides with the harmonic function $H$ in $\tilde O$, so that $\tilde A$ is a 
local set of $\tilde \Gamma$ with $h_{\tilde A} = H$. 
\end {proof}

\section*{Bibliographical comments} 
The definition and main regularity properties of the continuum GFF are rather classical facts. The spatial Markov property of the GFF was of course also pointed out early on (see e.g., \cite {Nelson,Rozanov}). 
The very closely related notion of local sets of the GFF was coined by Schramm and Sheffield \cite {SchrammSheffield} for the two-dimensional GFF (in relation to SLE curves), and then extensively used in the work of Miller and Sheffield (see e.g. \cite {IMG1} and the references therein). See \cite {ASW} for some features of thin local sets.

\chapter {Topography of the continuum Gaussian Free Field}
\label{ch:sle4gff}
\label{Ch5}

\section {Warm-up and overview} 

We will now focus on the Gaussian Free Field in two dimensions: more specifically, we will work with the GFF in a simply connected domain $D \not= \R^2$. As we have already mentioned, conformal invariance shows 
that the particular choice of $D$ does not really matter; so, let us discuss the case where $D$ is the unit disc $\U$.

In the sequel, $\partial_+$ and $\partial_-$ will denote the top and bottom half-circles of the unit circle respectively, that join $-1$ to $1$. We denote by $\mathbf{h}^+_0$ the bounded harmonic function in $\U$ that extends continuously 
to $\partial_+$ and to $\partial_-$, and is equal to $1$ on $\partial_+$ and to $0$ on $\partial_-$. That is,  $\mathbf{h}_0^+ (z)$ is the probability that a Brownian motion started from $z$ exits $\U$ through $\partial_+$.

A special role will be played in this chapter by the GFF in $\U$ with boundary conditions given by the harmonic function $2 \lambda \mathbf{h}^+_0$, that we will often just describe as the GFF with boundary conditions $2 \lambda$ on $\partial_+$ and $0$ on $\partial_-$ (where $\lambda$ is some positive constant). Recall that this GFF is just the sum of $2 \lambda \mathbf {h}^+_0$ and a GFF with Dirichlet boundary conditions in $\U$.

We are going to describe a particular random continuous curve $\gamma$ from $-1$ to $1$ in $\overline \U$, that is simple (i.e., non self-intersecting) and does not intersect $\partial \U \setminus \{ -1, 1 \}$. Such a curve
divides $\U$ into the two connected components of $\U \setminus \gamma$ that we denote by $U^+$ and $U^-$ ($U^+$ being the one which has $i$ on its boundary). 
We denote by $\mathbf{h}^+_\infty (z)$ the function $\1{ z\in U^+}$.

\begin{figure}[h]
\centering
\includegraphics[scale=1]{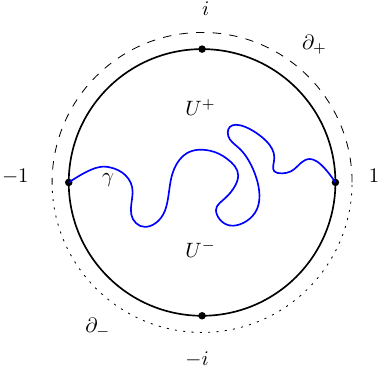}
\caption {The domains $U^+$ and $U^-$}
\end{figure}

 Conditionally on $\gamma$, we then consider two independent GFFs $\Gamma^+$ and $\Gamma^-$ (with zero boundary conditions) in $U^+$ and $U^-$ respectively. We also define $\Gamma^{++}$ to be the sum of $\Gamma^+$ with the constant function $2 \lambda$ in $U^+$ (so $\Gamma^{++}$ is a GFF in $U^+$ with constant boundary conditions $2 \lambda$). Finally, we define a
process $\tilde \Gamma$, indexed by the set ${\mathcal S}$ of smooth functions $f$ with compact support in $\U$, by: 
$$ \tilde \Gamma (f) := \Gamma^{++} (f \mathbf{1}_{U^+}) + \Gamma^- ( f \mathbf{1}_{U^-}),$$
where (here and throughout this chapter) $\lambda$ is some well-chosen (positive) constant. 
In other words, $\tilde \Gamma$ restricted to $U^+$ is a GFF in $U^+$ with $2 \lambda$-boundary conditions, and $\tilde{\Gamma}$ restricted to $U^-$ is a GFF in $U^-$ with zero boundary conditions. One main result of this chapter is the following: 

\begin {theorem}
\label {SLE4thm}
There exists a constant $\lambda>0$ and a simple random curve $\gamma$ from $-1$ to $1$ in $\U$, called the Schramm-Loewner Evolution SLE$_4$, such that the process $\tilde \Gamma$ described above has the law of a GFF in $\U$ with boundary conditions $2 \lambda$ on $\partial_+$ and $0$ on $\partial_-$. 
\end {theorem}

We can note that this theorem implies that the curve $\gamma$ is a local set of this GFF $\tilde \Gamma$ (mind that $\tilde \Gamma$ and the Dirichlet GFF $\tilde \Gamma - 2 \lambda {\mathbf h}_0^+$ 
generate the same $\sigma$-fields -- so we just define local sets for $\tilde \Gamma$ to be the local sets of $\tilde \Gamma - 2 \lambda {\mathbf h}_0^+$ ). 

There is only one value of $\lambda$ for which this will work. The quantity $2 \lambda$ is called the {\em natural height-gap of the two-dimensional GFF} (see discussion below).

Because the random curve $\gamma$ given by this theorem turns out to satisfy the criterion (Proposition \ref{prop:thincrit_d2}) about how the areas of its neighbourhoods decay, it will be a {\em thin local set} of the GFF as discussed in the previous chapter.

\begin {remark}
\label {R53}
In fact, we will see that the curve $\gamma$ in Theorem \ref {SLE4thm}  is a deterministic function of the GFF $\Gamma$. The proof of this fact  will require additional (non-trivial) considerations and will 
be dealt with in a subsequent section.
\end {remark}

One way to think of Theorem \ref{SLE4thm} (and the fact that $\gamma$ is a deterministic function of $\Gamma$) is to view $\gamma$ as a natural ``cliff-line'' of the field $\tilde \Gamma$, or equivalently of the field $\hat \Gamma := \tilde \Gamma - \lambda$.  
One  starts with $\hat \Gamma$, which is a GFF with boundary conditions $\lambda$ on $\partial_+$ and $- \lambda$ on $\partial_-$. Then, Theorem \ref{SLE4thm} says that $\hat \Gamma$ possesses a ``cliff-line'' $\gamma$ from $-1$ to $1$ in 
$\U$, such that the GFF has boundary conditions $+ \lambda$ on the top side of $\gamma$ and $-\lambda$ on the bottom side of $\gamma$ (and therefore the boundary condition of $\hat \Gamma$ when 
restricted to $U^+$ is $+\lambda $ on the whole of $\partial U^+$, and the boundary condition of $\hat \Gamma$ when restricted to $U^-$ is $-\lambda$ on the whole of $\partial U^-$). 

\begin{figure}[h]
	\centering
	\includegraphics[scale=1]{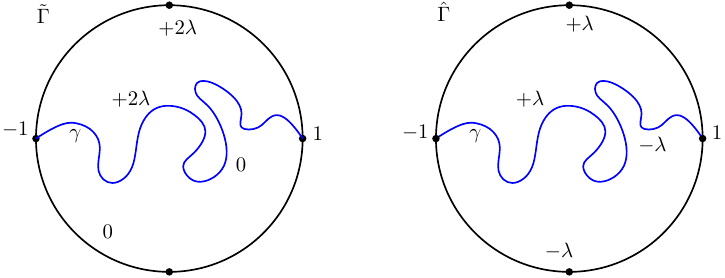}
	\caption {Constructing $\tilde \Gamma$ and $\hat \Gamma$}
\end{figure}

This suggests that it may be possible, given $\hat \Gamma$, to explore $\gamma$ progressively from $-1$ to $1$ (similarly to how, given a continuous function on $\U$ with boundary conditions $-\lambda$ on $\partial_-$ and $+\lambda$ on $\partial_+$, one could explore the line from $-1$ to $1$ on which it takes the value $0$). If one were to stop such an exploration at a time $t$, then one would expect that the boundary condition of 
$\hat \Gamma$ on the boundary of $U_t := \U \setminus \gamma [0,t]$ would be $+\lambda$ on the top side of $\gamma [0,t]$, and $- \lambda$ on the bottom side of $\gamma[0,t]$. So, if one applied the conformal transformation 
$\phi_t$ mapping $U_t$ onto $\U$ with $\phi_t (\gamma_t) = -1$, $\phi_t (1) = 1$ and $\phi_t' (1) = 1$, then the image of the restriction of $\hat \Gamma$ to $U_t$ under $\phi_t$ should be distributed like 
$\hat \Gamma$ itself (it should be a GFF with boundary conditions $+ \lambda$ on $\partial_+$ and $-\lambda$ on $\partial_-$). This remark will actually lie at the root of the definition of $\gamma$. Indeed, it indicates 
that $\gamma$ can be constructed as an iteration of random conformal maps (the law of $\phi_{2t}$ will be distributed like the composition of two independent copies of $\phi_t$). 

More precisely, the random curve $\gamma$ that we are going to construct will have the following property.
For each time $t$ and $z \in U_t$, if we define the harmonic measure $\mathbf{h}_t^+(z)$ of $\partial_+ \cup (\gamma[0,t])_+$ in $U_t$ at $z$, 
where $(\gamma [0,t])_+$ denotes the ``top side'' of $\gamma [0,t]$, then:
\begin{itemize}
	\item for each $z \in \U$, the process $(\mathbf{h}_t^+ (z), t \ge 0)$ will be a martingale with respect to the filtration $({\mathcal F}_t = \sigma (\gamma [0,t]))_{t \ge 0}$.
\end{itemize}
In fact, we can say exactly what the law of $\gamma$ has to be, if we want it to satisfy this property. Namely, it has to be a so-called SLE$_4$ curve from $-1$ to $1$ in $\U$:
\begin{itemize}
\item the law of an SLE$_4$ curve from $-1$ to $1$ in $\U$ is, up to time-change, the unique law on random curves with this martingale property.
\end{itemize}

When we set $\gamma$ to be an SLE$_4$, the fact that this process is a martingale will actually enable us to derive Theorem \ref{SLE4thm}.

\begin {remark}[Level-lines]
\label {LevelLines}
This cliff-line $\gamma$ is often also referred to as a {\em level-line} of the continuum GFF. So for instance, in the case that we have just discussed ($-\lambda$ on one side 
of the curve and $+\lambda$ on the other side), it would be $0$-level line of $\hat \Gamma$. The reason for this comes from the following interesting feature: Suppose that one considers an approximation of $\U$ by a triangular lattice with width $\delta$, and that one considers a discrete GFF $\hat \Gamma^\delta$ on this graph with boundary conditions $+ \lambda$ and $-\lambda$ on the two discrete approximations of $\partial_+$ and $\partial_-$. Then, the lower boundary of the cluster of sites where the GFF is positive and that contains $\partial_+$ coincides with the upper boundary of the cluster of sites where the GFF is negative, and that contains $\partial_-$ -- it is a simple curve $\gamma^\delta$ drawn on the hexagonal lattice (dual to the triangular lattice) such that the discrete GFF is positive on its neighbouring sites  one of its side and negative on its neighbours on the other side (so in this discrete case, it can really be viewed as a ``level line''). It turns out (this is a highly non-trivial result by Schramm and Sheffield that we will not discuss here) that as $\delta \to 0$, the joint law of $(\hat \Gamma^\delta, \gamma^\delta)$ converges to the joint law of $(\hat \Gamma, \gamma)$ that we have just described. So, in this sense, $\gamma$ can be viewed as a level-line itself (as scaling limit of level lines).
\end {remark}

\section {Deterministic Loewner chains background}

We quickly review without proofs some basic facts about deterministic Loewner chains and simple curves in the upper half-plane. 
Notation-wise, we will always use $\U$ and $\HH$ to denote the open unit disc and the open upper half-plane in the complex plane i.e. $\overline \HH = \{ x +iy \in \C \ , \ y \ge 0 \}$. $\Im (z)$ (resp. $\Re (z)$) will denote the imaginary (resp. real) part of a complex number $z$. The Loewner chain set-up is easier 
to first describe in the upper half-plane (although Theorem \ref{SLE4thm} involves an SLE$_4$ in $\U$, this can be obtained from an SLE$_4$ defined in $\HH$ by conformal mapping).

\begin {itemize}
 \item 
Suppose that $(\gamma(u), u \in [0, \tau))$  is a continuous (deterministic) simple curve in 
$\overline \HH $ such that 
$ \gamma (0) = 0 $ and $\gamma (0, \tau) \subset \HH$.
Then, for each $u<\tau$, by Riemann's mapping theorem, one can uniquely define the two conformal transformations $\tilde g_u$ and $\tilde f_u$ from $\HH \setminus \gamma (0,u]$ into $\HH$, that are chosen to satisfy 
$$\tilde f_u (\gamma(u)) = 0, \text{ and as } z\to \infty: \, \tilde f_u(z) \sim z; \  \tilde g_u (z) = z + o(1).$$
Both of these functions have a Laurent series expansion near $\infty$ with real-valued coefficients. That is,
$$ \tilde f_u (z) = z - \tilde W_u + \tilde a(u)z^{-1}+ o(|z|^{-1}) \text{ as } |z|\to \infty$$ 
for some $\tilde W_u\in \R$ 
(and one then has $\tilde g_u (z) = \tilde f_u (z) + \tilde W_u$ and $\tilde g_u (\gamma (u))= \tilde W_u$). Note that this defines real-valued functions $\tilde a$ and $\tilde W$ from $\gamma$.

It is easy to see that the mapping $u \mapsto \tilde a(u)/2$ is an increasing continuous function (that converges to some $\sigma \in (0, \infty]$ as $u \to \tau-$).
We can therefore define the reparametrised continuous curve $\eta: [0, \sigma) \to \overline \HH$ such that for all $u < \tau$, $\eta ( \tilde a(u)/2 ) = \gamma (u)$. 
From now on we work with $f,g$ and $W$ defined  by $f_{\tilde a(u)/2} := \tilde f_{u}$, $g_{\tilde a(u)/2} := \tilde g_{u}$ and $W_{\tilde a(u)/2} := \tilde W_u$, so that for all $u<\tau$ the maps $g_u$ and $f_u$ are conformal transformations from $\HH\setminus \eta(0,u]$ onto $\HH$, and 
$$f_u(\eta(u))=0\; \;\;\; g_u(\eta(u))=W_u$$ $$f_u(z)=z-W_u+2u z^{-1}+o(|z|^{-1}) \text{ as } z\to \infty; \;\;\; g_u=f_u+W_u.$$ 

In summary, after performing a deterministic simple time-change, we have obtained a path $(\eta(t), \,t\in [0,\sigma))$ such that 
for all $t$ in $[ 0, \sigma)$, one has $f_t(z) + W_t = g_t (z) = z + 2t/z + o(1/z)$ as $z\to \infty$. Note that if $|\Im (\eta (t))|$ is unbounded in $t$, then one necessarily has
 $\sigma = \infty$. 

\item
Loewner's equation provides a recipe to recover $\eta$ from the function $t \mapsto W_t = g_t (\eta(t))$ (which in particular shows that the curve $\eta$ is fully determined by $W$). Indeed, for all $t \ge 0$, when $z \in \HH \setminus \eta (0, t]$, it turns out that 
\begin {equation}
\label {Loewner}
  \partial_t g_t (z) = 2 / ( g_t (z) - W_t ).
  \end {equation}
In particular, this enables one (via the ``reverse flow''), for each $y \in \HH$ and each $T \ge 0$, to construct $g_T^{-1} (y)$ as the value at time $T$ of the function $y(\cdot)$ with $y(0)=y$ and $\partial_t y(t) = -2 / (y(t)- W(T-t))$. Then, one can recover $\eta (0,T]$ as $\HH \setminus g_T^{-1} (\HH)$. 

Let us emphasise that for each simple curve $\eta$, there exists a continuous function $W$ from which one can recover $\eta$ uniquely using this procedure, but that if we are given an arbitrary continuous $W$, it may happen that it does not correspond to a continuous curve $\eta$. 

\item
Let us summarise a few trivial properties of the Loewner flow. If we fix $z \in \HH$ and define $Z_t = X_t + i Y_t := f_t (z)$ and $\theta_t := \arg  (f_t (z)) \in (0, \pi )$, then (as long as $z \notin \eta [0,t]$),
$$ Y_t - Y_0 = \int_0^t  \Im ( 2/Z_s) ds, \quad X_t = Y_t / \tan (\theta_t), \quad W_t = -X_t + X_0 + \int_0^t \Re (2 /Z_s) ds .
$$
Hence, we see that $W$ can obtained from $X$ and $Y$ by a simple transformation starting from $\theta$ and involving only addition or compositions with smooth functions. In particular, if 
we happen to know that $(\theta_t)$ is a semi-martingale with respect to some filtration (and that $\eta$ is also adapted to this filtration, so $X$ and $Y$ are too), then it follows immediately that $(W_t)_t$ is also a 
semi-martingale with respect to the same filtration.

\item 
At each time $t$, the domain $H_t := \HH \setminus \eta (0, t]$ is simply connected. Clearly, $H_t$ is decreasing with $t$, so that the functions $t \mapsto G_{H_t} (x,y)$ are non-increasing. Recall that the Green's function in $\HH$ is given by 
$$ G_\HH (x, y)= \frac 1 {2 \pi}\log \frac { | x - \overline y | }{|x - y |}=\frac 1 {2 \pi} \Re ( \log (x - \overline y) - \log (x-y)), $$
and so by conformal invariance, we get 
$$ G_{\HH \setminus \eta (0, t]} (x,y) = G_\HH ( f_t (x), f_t (y)) = G_\HH ( g_t (x), g_t (y)). $$
For any $x\in \HH$, the function $t \mapsto g_t (x)$ is smooth up until the possibly finite time at which $\eta$ hits $x$. 
Differentiating the previous expression with respect to $t$ shows immediately that (for $x,y \in H_t$),
$$ \partial_t  G_{\HH \setminus \eta (0, t]} (x,y) = - \frac {1}{2 \pi} I_t (x) I_t (y),$$
where here and in the sequel, $I_t (x) = \Im (-2 / f_t (x))$. 

This shows in particular that for all smooth test functions $\varphi$ (continuous with compact support), 
\begin{equation}
\label{eqn:ghgh}
\iint 
\varphi (x) \varphi (y) ( G_\HH (x,y) - G_{H_\infty} (x,y) )  dx dy 
= \frac {1}{2 \pi} 
\int_0^\infty \iint I_t (x) I_t (y) \varphi(x) \varphi (y) dt dx  dy.
\end{equation} 
Note that the right-hand sided is therefore bounded by $\iint \varphi (x) \varphi (y) G_\HH (x,y) dx dy$, uniformly in the curve $\eta$. 

\end {itemize}
 
\section {SLE$_4$, harmonic measure martingales and coupling with the GFF}

Here we will not give a detailed construction of SLE$_4$, only a brief summary of some of its features. Let us state without proof the following result about SLE that is essentially due to Rohde and Schramm:
\begin {proposition} \label{prop:slecurve}
For all $\kappa \le 4$, there exists a random continuous simple curve $\eta$ such that its corresponding driving function $(W_t)_{t \ge 0}$ is a one-dimensional Brownian motion running at speed $\sqrt {\kappa}$ (so $\beta_t := W_t / \sqrt {\kappa}$ is a standard Brownian motion).
Furthermore, $ \Im (\eta(t)) $ is unbounded and $| \eta (t) | \to \infty$ as $t \to \infty$ almost surely.

Finally, if one considers the image $\tilde \eta$ of the $\eta$ via the map $z \mapsto (z-i)/(z+i)$ from $\HH$ onto $\U$, then the expected area of the $\eps$ neighbourhood of $\tilde \eta$ is bounded by a power of $\eps$ 
as $\eps \to 0$. 
\end {proposition}

\begin{remark}
	\label{rmk:sleDab}
	In fact, if $\phi$ is a conformal map from $\HH$ to itself that fixes $0$ and $\infty$ (i.e. a scaling map), then it follows from the scale invariance of Brownian motion that the image of an SLE$_\kappa$ under $\phi$ again has the law of an SLE$_\kappa$. This gives us a way to define SLE$_\kappa$ in any simply-connected domain $D$ between two marked boundary points unambiguously: we let it be the curve whose law is obtained by taking the image of an SLE$_\kappa$ in $\HH$, under a conformal transform that sends $\HH$ to $D$, and maps $0$ (resp. $\infty$) to the chosen starting (resp. ending) point on $\partial D$.
\end{remark}
Note that the statements of Proposition \ref{prop:slecurve} are non-trivial. The final point essentially says that the Hausdorff dimensions of $\tilde \eta$ and of $\eta$ are strictly smaller than $2$: this will ensure 
that when we take $\gamma$ to be an SLE$_4$ in Theorem \ref{SLE4thm}, it will indeed be a thin local set of the GFF. 
In fact, it can be proved that the Hausdorff dimension of an SLE$_\kappa$ curve with $\kappa \le 4$ is almost surely equal to $1 + ( \kappa /8)$. 

Another equivalent way to state this proposition is the following. Start with a one-dimensional Brownian motion $\beta$, and define $W_t = \sqrt {\kappa}  \times \beta_t$. For each 
$z \in \HH$, let $g_t(z)$ be the solution to Loewner's equation (\ref {Loewner}), started from $g_0 (z) = z$ and defined up to time 
$$ T (t) := \inf \{ t \ge 0 , \ \inf_{s \in [0,t)} | g_s (z) - W_s | = 0 \}. $$
Then, denoting $$K_t := \{ 0 \} \cup \{ z \in \HH , \ T(z) \le t \},$$ there almost surely exists a continuous simple curve $\eta$ such that for all $t$, 
$K_t = \eta(0, t]$, and $g_t$ is the unique conformal map from $\HH \setminus K_t $ to $\HH$ that is normalised to satisfy $g_t (z) = z + 2t/z + o(1)$ as $|z|\to \infty$.
In other words, starting from the Brownian motion $\beta$, one has a concrete recipe to construct the curve $\eta$, by solving Loewner's equation. The non-trivial part of the proposition is  to prove that the set $K_t$ defined above is described by a random simple curve when $\kappa \le 4$ (this actually fails to be true when $\kappa > 4$).   

Applying Proposition \ref{prop:slecurve}, we see that if $\eta$ is an SLE$_\kappa$ and $f_t=g_t-W_t$ is the associated Loewner flow, then for each $z \in \HH$ one has that $f_0 (z) = z$ and
$$ df_t (z) = - \sqrt {\kappa} d \beta_t + \frac 2 {f_t (z)} dt $$
as long as $f_t(z)$ does not hit $0$. The only point $z$ for which $f_t (z)$ hits $0$ at time $t$ is $z=\eta(t)$. In particular, for every point $z$ that is not on the curve $\eta$, $f_t(z)$ is defined for all time.

For each $z \in H_t$, we let $\theta_t (z) \in (0, \pi)$ be the argument of $f_t (z)$ as before. This is clearly a continuous function of $t$ (as long as $f_t (z)$ stays away from $0$) and applying It\^o's formula, we get that
\begin {equation}
 \label {} 
d\theta_t (z) = \Im ( d \log (f_t(z))) = \Im ( - \sqrt {\kappa} / f_t(z))   d\beta_t + (2 - \kappa /2) \Im (1/ f_t(z)^2) dt.
\end {equation}
The fact that the drift term disappears at $\kappa =4$ is one of the things that makes this value very special. Let us highlight this simple fact as a proposition:

\begin {proposition}[The nice SLE$_4$ martingales]\label{sle4mgale}
When $\eta$ is a SLE$_4$, then for all fixed $z \in \HH$, $(\theta_t (z), t \ge 0)$ is a continuous martingale with respect to the filtration $(\sigma ( \eta[0,t]))_{t\ge 0}$ (note that since $\eta$ and $W$ encapsulate ``the same'' information when run up to any given time ($\eta$ can be recovered from $W$ via Loewner's equation and $W$ from $\eta$ since  $W_t=g_t(\eta(t))$ for all $t$) the filtrations generated by $\eta$ and $W$ coincide).
\end {proposition}

Note that for every given (deterministic) $z\in \HH$, $\eta$ almost surely does not hit $z$. This means that the martingale $\theta_t (z)$ is almost surely defined for all time, and stays in $(0, \pi)$ (this explains why it is actually a martingale, and not only a local martingale). Of course, there do exist random exceptional points that end up being on the curve, and for which the corresponding $\theta_t$ is only defined up to a finite time.

From now on, we will assume that $\kappa =4$ and that $\eta$ is an SLE$_4$. Then, for each given deterministic $z$, $t \mapsto \theta_t (z)$ is a martingale, and as $t \to \infty$ it converges to $\theta_\infty(z)$. Note that $\theta_\infty(z)$ is either $0$ or $\pi$ depending 
on whether $\eta$ passes to the left of $z$ or to the right of $z$. We see for instance (viewing $\theta_t$ as a time-changed Brownian motion, see below) that the quantity 
$ \int_0^\infty I_t (z)^2 dt $ can be interpreted (recall that $I_t (z) := \Im (-2 / f_t (z))$) as the time at which a Brownian motion starting from $\theta_0(z)$ exits the interval $(0,\pi)$.

\begin{remark}\label{rmk:bmmgale}
	Here we record the following classical fact: if $\beta$ is a standard Brownian motion 
and $H=(H_t)$ is a continuous process adapted to the filtration of $\beta$, 
then the process $M_t = \int_0^t H_s d \beta_s$ is a local martingale that can be interpreted as a time-changed Brownian motion. Indeed, if one sets $U_t = \int_0^t H_s^2 ds$ for all $t$, defines $\tau$ to be the inverse of $U$ and takes ${\mathcal G}_u = {\mathcal F}_{\tau (u)}$, then 
$ B_u := M_{\tau (u)}$ is a Brownian motion with respect to the filtration ${\mathcal G}_u$ (possibly stopped at the stopping time $U_\infty$ if this quantity is finite).

We will use this result in the particular case where there exists a deterministic $U_0$ such that $U_\infty < U_0$ almost surely. 
In this case, if we condition on $(B_t; t<U_\infty)$ and then add to
$B_{U_\infty} = M _\infty$ a (conditionally) independent random Gaussian variable with mean $0$ and variance $V:=U_0 - U_\infty$, we obtain (via the strong Markov property) a random variable that is distributed like a Brownian motion at time $U_0$ i.e. a Gaussian random variable with mean $0$ and variance $U_0$. 
\end{remark}

We will use this line of reasoning very shortly in the proof of Proposition \ref{prop:wholecurvelocal}.

Our goal is to construct a coupling of $\eta$ with a GFF in $\HH$, so that $\eta$ will be a local set of the GFF (we will then apply a conformal map from $\HH$ to $\U$ to obtain Theorem \ref{SLE4thm}). As the Hausdorff dimension of $\eta$ is almost surely strictly smaller than 2, $\eta$ will necessarily be a thin local set, and so the law of the coupling will be totally described by the associated harmonic function defined in $\HH \setminus \eta$. In fact, the harmonic function will turn out to be given by 
\begin{equation}
\label{h}
h(z) = \frac {1}{\sqrt {2\pi}} (\theta_\infty  (z) - \theta_0 (z)),\end{equation}
This will be a very special example of a local set, because the harmonic function is actually deterministic given $\eta$ (in fact it is not at all clear a priori that such local sets should even exist).

Note that if $\eta$ is an SLE$_4$ curve, then $\HH \setminus \eta$ consists of two open simply connected domains $H_-$ and $H_+$ that lie respectively to the left and to the right of $\eta$. Therefore, one can first sample $\eta$ and second (conditionally on $\eta$) sample two independent (zero boundary condition) GFFs $\Gamma^-$ and $\Gamma^+$ in the domains $H_-$ and $H_+$. Finally, one can add the harmonic function $h$ described above, and ask if the resulting field $\Gamma$ has the law of a GFF. The answer is that it does, and we summarise this in the following proposition.

\begin {proposition}[The whole curve is a local set]\label{prop:wholecurvelocal} 
This construction provides a local coupling of an SLE$_4$ $\eta$ and a GFF $\Gamma$ with zero boundary conditions. In this coupling, $\eta$ is a thin local set of $\Gamma$ and the associated harmonic function is the function $h$ of \eqref{h}.
\end {proposition}

\begin {proof}
Let $\overline{\Gamma}=\Gamma^+ + \Gamma^-$ where $\Gamma^{\pm}$ are as described in the paragraph above.  We would like to show that $\Gamma:=\overline{\Gamma}+h$ has the law of a (zero boundary condition) GFF in $\HH$. In order to prove this, it suffices to show that for any fixed function $\varphi$ on $\HH$, that is continuous with compact support, the random variable $\overline\Gamma (\varphi ) + \int \varphi (z) h(z) d z $ is a Gaussian random variable with mean $0$ and variance 
$$ 
U_0 := \iint G_\HH (x,y) \varphi (x) \varphi (y) dx d y
$$
(indeed, this determines the characteristic function of the random vector $(\Gamma (\varphi_1), \ldots , \Gamma (\varphi_k))$ for any finite family $\varphi_1, \ldots , \varphi_k$, and therefore 
the finite-dimensional distributions of $\Gamma$). 
 
By definition, the conditional distribution of 
$\overline \Gamma (\varphi )$ given $\eta$ is a centred normal with variance 
$$ 
V := \iint G_{\HH \setminus \eta} (x,y) \varphi (x) \varphi (y) dx d y.
$$ 
Let us now define 
$$ M_t := \frac {1}{\sqrt {2\pi}}\int \varphi (z) \theta_t (z) d z .$$ 
By (stochastic) Fubini, we get that 
$$ 
M_t = \int_0^t H_s d \beta_s, 
\hbox { where } 
 H_s := \frac {1}{\sqrt {2\pi}}\int \varphi (z) I_s (z) d z,$$
and hence one can view $M_\infty$ as a Brownian motion $B$, stopped at the time 
$U_\infty := \int_0^\infty H_s^2 ds$.

Integrating this with respect to $\varphi(x) \varphi (y) d x d y d t$ over $\HH\times \HH\times [0,\infty)$, and then applying \eqref{eqn:ghgh}, we obtain that
$$ U_0 -  V 
= \iint (G_\HH (x,y) - G_{\HH \setminus \eta} (x,y))  \varphi (x) \varphi (y) dx d y 
= \int_0^\infty H_s^2 ds = U_\infty.$$ 
By Remark \ref{rmk:bmmgale}, this concludes the proof. 
\end {proof} 

\begin {proposition}[The SLE trace until a fixed finite time $t$ is also a local set]\label{prop:partcurvelocal}
Suppose that one couples an SLE$_4$ $\eta$ with a GFF $\Gamma$ as in Proposition \ref{prop:wholecurvelocal}. Then, for any fixed $t$, the curve $\eta [0,t]$ is also a thin local set with respect to the same GFF $\Gamma$. The harmonic function associated to $\eta [0, t]$ is 
$$ h_t (z) = \frac {1}{\sqrt {2\pi}} (\theta_t  (z) - \theta_0 (z)).$$
\end {proposition}

This result shows that $\eta$, viewed as a growing curve, is in fact a ``continuously increasing family of local sets''. 
Heuristically, it is a continuum counterpart of the iterative procedures  considered in Chapter \ref{Ch1}, for discovering local sets of the discrete GFF by ``uncovering" its values at sites one by one.

In order to prove this proposition, we need the following property, that follows immediately from the Markov property of Brownian motion. 
Suppose that $t >0$ is fixed, and let us sample $\eta [0,t]$. Then, the random path $\eta^t := (f_t ( \eta (t+s)), s \ge 0) $ is an SLE$_4$ that is independent of $\eta [0,t]$. 

\begin{figure}[h]
\centering
	\includegraphics[scale=1]{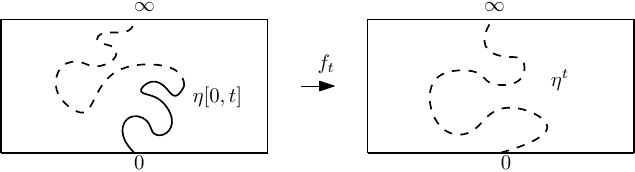}
	\caption{The function $\theta_0$ is the harmonic function in $\HH$ that is equal to $\pi$ on $\R_-$ and $0$ on $\R_-$. $\theta_t$ is the harmonic function in $\HH\setminus \eta[0,t]$ that is equal to $\pi$ on $\R_-$ and on the left-hand side of $\eta[0,t]$, and equal to $0$ on $\R_+$ and on the right-hand side of $\eta[0,t]$. A consequence of Proposition \ref{prop:partcurvelocal} is that if one samples $\eta[0,t]$, and then defines $\Gamma$ to be the sum of a GFF in $\HH\setminus \eta[0,t]$ plus the function  $h_t=(2\pi)^{-1/2}(\theta_t-\theta_0)$, then $\Gamma$ has the law of a GFF in $\HH$.} 
\end{figure}

In particular, we can apply Proposition \ref{prop:wholecurvelocal} to the random path $\eta^t$. This allows us to prove Proposition \ref{prop:partcurvelocal} as follows. 

\begin{proof}[Proof of Proposition \ref{prop:partcurvelocal}]
Observe that by conformal invariance of the GFF, we can construct the coupling $(\eta,\Gamma)$ of Proposition \ref{prop:wholecurvelocal} by:
\begin{itemize} 
	\item first sampling $\eta[0,t]$; 
	\item conditionally on $\eta[0,t]$, sampling $\eta^t$ in $\HH$; 
	\item defining a field in $\HH$ by sampling independent GFFs on either side of $\eta^t$ and adding to them the function $h$; 
	\item finally setting $\Gamma$ to be the image of this field in $\HH$ under $f_t^{-1}$ (and $\eta=\eta[0,t]\cup f_t^{-1}(\eta^t)$).
\end{itemize}
The observation that $\eta^t$ is an SLE$_4$ in $\HH$ that is independent of $\eta[0,t]$, together with Proposition \ref{prop:wholecurvelocal} applied to $\eta^t$, means that conditionally on $\eta[0,t]$:
\begin{itemize}
	\item $\Gamma$ is exactly a GFF + $h_t$ in the complement of $\eta[0,t]$, that is independent of $\eta[0,t]$.
\end{itemize} This proves Proposition \ref{prop:partcurvelocal}. 
\end{proof}

Theorem \ref{SLE4thm} also follows from this, with $\lambda=\sqrt{\pi/8}$, after mapping everything to the unit disc and using conformal invariance of both the GFF and SLE$_4$ (recall that we define SLE$_4$ from $-1$ to $1$ in $\U$ as the conformal image of SLE$_4$ in $\HH$, Remark \ref{rmk:sleDab}). Indeed, for this value of $\lambda$ the harmonic function $h$ of the proposition (see \eqref{h}) is the function that is equal to $2\lambda$ on the left of the curve and $0$ on the right of the curve, minus the harmonic function equal to $2\lambda$ on the negative real line and $0$ on the positive real line.

\medbreak

\begin{remark}
	\label{rmk:othercoupling}
Finally, let us notice that if we replace $h$ by $-h$ in Proposition \ref{prop:wholecurvelocal}, then by symmetry this defines another local coupling between SLE$_4$ and the GFF.  This works because the obtained field is then distributed like $-1$ times a GFF, which of course just has the law of a GFF. 
\end{remark}

\section {SLE$_4$ is a deterministic function of the GFF}

We now explain the following important property of the above-described coupling of SLE$_4$ and the GFF, as promised in Remark \ref {R53}.

\begin {proposition}[SLE$_4$ is determined by the GFF] 
\label {p11}
In the  coupling between SLE$_4$ and the GFF described by Theorem \ref{SLE4thm}, the SLE$_4$ is a deterministic function of the field.
\end {proposition}

As a by-product of the proof, one also obtains the following non-trivial property of SLE$_4$:

\begin {proposition}[SLE$_4$ reversibility] \label{prop:sle4rev}
The law of $\eta$ is reversible in the following sense: up to time-reparametrisation, an SLE$_4$ from $-1$ to $1$ in $\U$ is distributed like {the time reversal of} an SLE$_4$ from $1$ to $-1$ in $\U$.
\end {proposition}

\medbreak

Let us first summarise the outline of the proof (of Proposition \ref{p11}). We denote by $r$ a conformal map from $\U$ to itself, that maps $-1$ to $1$ and $1$ to $-1$. We define on the same probability space a triple $(\eta, \hat \eta, \Gamma)$, such that $\eta$ is an SLE$_4$ from $-1$ to $1$ in $\U$, $\hat \eta$ is an SLE$_4$ from $1$ to $-1$ in $\U$, $\Gamma$ is a GFF in $\U$, and:  
\begin {itemize}
\item the joint law of $(\eta, \Gamma)$ is that of the coupling between SLE$_4$ and the GFF in Theorem \ref{SLE4thm} (if $\tilde{\Gamma}$ is as in the theorem then we take $\Gamma=\tilde{\Gamma}-2\lambda \mathbf{h}_0^+$ the associated Dirichlet GFF);
\item the joint law of $(\hat \eta,  \Gamma)$ is such that $(r(\hat{\eta}), -r(\Gamma))$ is coupled as in Theorem \ref{SLE4thm} (this corresponds to the coupling described in Remark \ref{rmk:othercoupling});
\item conditionally on $\Gamma$, the two curves $\eta$ and $\hat \eta$ are independent.
\end {itemize} 
	Since in the coupling of Theorem \ref{SLE4thm}, $\eta$ and $\Gamma$ are both random variables taking values in Polish spaces, the regular conditional distribution $\eta$ given $\Gamma$ exists. In particular, the coupling described above is well-defined.
	
The following lemma will imply the proposition: 
\begin {lemma}
\label {l4}
For such a triple $(\eta,\hat{\eta},\Gamma)$, the trace of $\eta$ is 
almost surely equal to the trace of $\hat \eta$.
\end {lemma}
Indeed, if the lemma is true, then the trace of $\eta$ is conditionally independent of itself given $\Gamma$ (because by assumption it is conditionally independent of the trace of $\hat \eta$). It is therefore a deterministic measurable function of $\Gamma$. Furthermore, it is clear that $\eta$ is (up to time reparametrisation) the time-reversal of $\hat \eta$ (since they are two simple curves with the same trace, and one goes from $-1$ to $1$ while the other goes from $1$ to $-1$). This implies that the time-reversal of an SLE$_4$ is indeed an SLE$_4$ (Proposition \ref{prop:sle4rev}).

\medbreak
Let us now explain our previous remark that SLE$_4$ (in $\HH$) is actually characterised by the fact that the processes $( \theta_t(z)=\arg(f_t(z)), t \ge 0)$ for $z\in\HH$ are martingales. This will be our main tool for  the proof of Lemma \ref{l4}.

Suppose that $(\gamma_u, u < \tau)=(\gamma(u), u < \tau)$ is a random simple curve in $\overline \HH$ (with $\gamma (0)=0$ and $\gamma (0,\tau) \subset \HH)$)  that is adapted with respect 
to some filtration $({\mathcal F}_u)=({\mathcal F}_u)_{u\ge 0}$, and such that $\tau$ is a possibly infinite $({\mathcal F}_u)$-stopping time. Define for each $z \in \HH$, the conformal transformation $\tilde f_u$ from $\HH \setminus \gamma (0,u)$ onto $\HH$ with $\tilde f_u (\gamma_u) = 0$ and $\tilde f_u (z) \sim z$ as $z \to \infty$. Define $\varphi_u (z)$ to be the argument of $\tilde f_u(z)$. We will also assume that as $u \to \tau-$, either $d(\gamma_u, \R) \to 0$ or $\Im (\tilde f_u(z))$ is unbounded.  

\begin {lemma}
\label {l5}
Under these conditions, if  $\varphi_u (z)$ is a martingale in the filtration $({\mathcal F}_u)$ for each $z\in U$, then $\gamma$ is distributed like a time-changed SLE$_4$ from $0$ to infinity in $\HH$. It follows in particular that as $u \to \tau-$, $|\gamma(u) | \to \infty$.
\end {lemma}

\begin {proof}[Proof of Lemma \ref {l5}]The goal is to show that a certain Loewner chain is an SLE$_4$, i.e., that its driving function $W$ is distributed like a Brownian motion with speed $4$. The assumptions of the lemma will tell us that certain explicit functions are local martingales with respect to the filtration generated by the curve (or equivalently by $W$), and using It\^o calculus, one can then deduce that both $W_t$ and $(W_t^2 - 4t)$ are local martingales. This is enough to deduce that $W_{t/4}$ is a Brownian motion. We remark that this type of argument is now quite standard, and has been used in many other instances to identify the scaling limit of some lattice model interface in terms of an SLE curve. The argument goes as follows:
\begin {itemize}
 \item First, time-change $\gamma$ into $\eta$ as in Section $1$, and for all $t \ge 0$, set ${\mathcal G}_t = {\mathcal F}_{u(t)}$ where $u(t)$ is  the inverse of $v \mapsto \tilde a(v)/2$.  
 Since $u(t)$ is a stopping time for $({\mathcal F}_u)$, it follows that for each fixed $z$, $\theta_t := \varphi_{u(t)}$ is a martingale (stopped at the stopping time $\sigma (\tau)$) for the filtration $({\mathcal G}_t)$.
 \item Then, the previous considerations imply that $(W_t)$ (stopped at $\sigma$) is a semi-martingale with respect to the filtration $({\mathcal G}_t)$. This semi-martingale can be decomposed into its local martingale term $M_t$ and its finite variation term $V_t$. Using the fact that 
 $\theta_t = \Im ( \log f_t (z))$ (where $f_t=\tilde{f}_{u(t)}$), we can apply It\^o's formula for semi-martingales, from which it follows that 
 $$ d\theta_t =   \Im ( - dM_t / f_t(z)) + \Im ( - 1 / f_t(z)) dV_t + (2 dt  -  d \langle M \rangle_t /2) \Im (1/ f_t(z)^2) dt .$$ 
As the finite variation part of this semi-martingale is zero, it follows that for each fixed $z$, for all $t \ge 0$, 
 $$ \int_0^t \Im (1/f_s (z)) dV_s + \Im (1/f_s^2 (z)) (2 ds - d \langle M \rangle_s /2 ) = 0 $$
 almost surely. 
 By looking at $(1/n)$ times this quantity when $z = \mathrm{i} n$ and letting $n \to \infty$ (using the fact that almost surely, $f_s (\mathrm{i}n) / \mathrm{i}n$ converges uniformly to $1$ in the time-interval $[0,t]$ as $n \to \infty$), we see that $V_t=0$.
By continuity of $V$ with respect to $t$, we get that $V=0$ at all times. 
 Hence, 
we know that for any given $z$ and all $t \ge 0$,
$$\int_0^t  \Im (1/f_s^2 (z)) (2 ds - d \langle M \rangle_s /2 ) = 0. $$
Moreover, almost surely for each $t$, one can find an $n$ large enough such that 
$\Im (1/ f_s^2 (n+\mathrm{i}))$ remains positive on $ s \in [0,t]$, and therefore conclude that  $ \langle M \rangle_s - 4s =0 $ almost surely on $[0,t]$. By continuity with respect to $t$, it finally follows that $W_t / 2$ behaves like a Brownian motion up to time $\sigma$.
 \item Hence, we may conclude that the path $\eta$ is a SLE$_4$ up to the time $\sigma$. But our conditions on the behaviour of $\gamma$ near time $\tau$ then imply that $\sigma = \infty$.
\end {itemize}
\end {proof}

We can now finally prove Lemma \ref {l4}.
\begin {proof}[Proof of Lemma \ref {l4}]
We know that almost surely, both $\eta$ and $\hat \eta$ are simple curves. 
It will therefore be sufficient to prove that for any given $s \ge 0$, the point $\hat \eta (s) $ on $\hat{\eta}$ is also almost surely on $\eta$. Then it follows that on a set of full probability this holds simultaneously for a countable dense set of times $s$, and therefore for all times by continuity. 

We will prove the claim for fixed $s$ by showing that conditionally on $\hat \eta (0, s ]$, the path $\eta$ up to the time $T$ at which it hits $\hat \eta [0,s ]$, is distributed like an SLE$_4$ from $-1$ to $\hat \eta(s)$ in $\U \setminus \hat \eta(s)$. In particular, $\eta$ does indeed hit $\hat \eta (s)$. 

In order to prove this statement, it suffices the check that the characterisation of SLE$_4$ (in terms of martingales, Lemma \ref{l5}) holds.

So let us fix $s \ge 0$, define $T$ as above and let $\eta^T (t) := \eta (\min (t, T))$. Define for each $t > 0$, the $\sigma$-field ${\mathcal F}_t = \sigma ( \eta^T [0,t], \hat \eta[0,s])$. 
We know that $\eta [0,t] \cup \hat \eta [0,s]$ is a local set of $\Gamma$, as it is a union of conditionally independent local sets. Let us denote by $h_{t,s}$ the corresponding harmonic function (so that $h_t:=h_{t,0}$ is the harmonic function associated to $\eta[0,t]$ and $\hat h_s := h_{0,s}$ is the harmonic function associated to 
$\hat \eta [0,s]$). 
In fact, by Proposition \ref {harmfunctionunion}, we actually know what $h_{t,s}$ is as long as $t < T$: it is equal to the harmonic function with boundary values provided by those of $\hat h_s$ and $h_t$ respectively on 
the two sides of $\hat \eta$ and $\eta$ (and equal to $0$ on $\partial \U$).
	
\begin{figure}[h]
\centering 
	\includegraphics[scale=1]{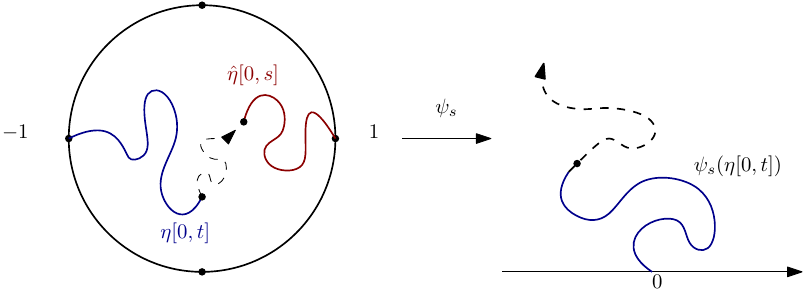}
	 \caption {The curve $\eta$ in the complement of $\hat \eta [0,s]$}
\end{figure} 
It follows from this that conditionally on $\hat \eta [0,s]$, if we apply a conformal map $\psi_s$ from $\U\setminus \hat\eta[0,s]$ to $\HH$ sending $-1$ to $0$ and $\hat{\eta}(s)$ to infinity, then $\psi_s(\eta(t))_{t\le T}$ satisfies the hypotheses of Lemma~\ref{l5}. 

Indeed, for $z\in \U\setminus \hat{\eta}[0,s]$, $t \mapsto h_{\min (t,T),s}(z)$ is a martingale with respect to the filtration $(\mathcal{F}_t)$. One can see this, for instance, using that if $z$ is at distance greater than $\eps$ from $\hat{\eta}[0,s]$ and $T_\eps$ is the first time that $d(z,\eta[0,t])\le \eps$, then $h_{t,s}(z)\mathbf{1}_{\{t<T_\eps\}}$ is just the conditional expectation $E[\Gamma(\lambda_{z,\eps})|\mathcal{F}_t]\mathbf{1}_{\{t<T_\eps\}}$. Moreover, if $f_t$ is the conformal map from $\HH\setminus \psi_s(\eta[0,t])$ to $\HH$ with $f_t(\psi_s(\eta(t)))=0$ and $f_t(z)\sim z$ as $|z|\to \infty$, then $\sqrt{2\pi} h_{t,s}(z)+\arg(z)=\arg(f_t(z))$.

 Hence, the curve $(\eta(t))_{t\in [0,T)}$ is an SLE$_4$ in $\U\setminus \hat{\eta}[0,s]$, which concludes the proof of Lemma~\ref {l4}.
\end {proof}

\section {Variants of this coupling result} 

\subsection {Other boundary conditions} 

Let us now explain how the fact that the SLE$_4$ is a deterministic function of the GFF (as discussed in the previous section), combined with absolute continuity considerations (as discussed in Chapter \ref {Ch3}), makes it 
possible to define level-lines for GFFs with some other boundary conditions and derive some of their properties.  
The goal of this section is to explain the specific statements that will be used in 
the next sections, and we will therefore not strive for generality.

We will first reformulate the results of the previous section in a very slightly different way. Suppose that $\hat \Gamma$ is a GFF in $\U$ with boundary conditions $\lambda$ and $-\lambda$ on  $\partial_+$ and $\partial_-\subset \partial \U$ respectively. 
Suppose that $O$ is some open subset of $\U$ that contains a neighbourhood of $-1$ in $\U$ (i.e., the distance between $\U \setminus O$ and $-1$ is positive). 
Then, our previous proofs imply  that 
there exists a unique random simple curve $(\eta_t)_{t \le \sigma}$ (up to time-reparametrisation)  that joins $-1$ to $ \eta_\sigma \in \partial O$, and is such that for all $t$, when $t < \sigma$, the conditional law of $\hat \Gamma$ in $\U \setminus \eta [0,t]$ given $\eta [0,t]$ is that of a GFF with boundary conditions $+ \lambda$ (resp. $-\lambda$) on the upper semi-circle and on the left-hand side of $\eta$ (resp. on the lower semi-circle and on the right-hand side of $\eta$).

{The existence statement is clear, and to see the uniqueness, one can first complete any such curve after time $\sigma$ by coupling an SLE$_4$ from $\eta_\sigma$ to $1$ with the conditional law of the field in $\U\setminus \eta[0,\sigma]$ as in Theorem \ref{SLE4thm}, and then applying our previous uniqueness result}. This shows in particular that 
the law of $(\eta_t)_{t \le \sigma}$ must be that of an SLE$_4$ up to its exit time of $O$, and that the curve $\eta$ is a deterministic function of $\hat \Gamma$. Recall that SLE$_4$  almost surely does not hit $\partial \U \setminus \{ -1, 1 \}$, so in particular, we have that $\eta (\sigma) \notin \partial \U \setminus \{1\}$. 

We next note that the Markov property of the GFF shows that $\hat \Gamma$ can be decomposed by first discovering $\hat \Gamma_{\overline O}$ and then $\hat \Gamma^{\overline O}$ (the latter process being defined to be a GFF with Dirichlet boundary conditions in $\U \setminus \overline O$). On the other hand, 
the Markov property and the fact that $\eta[0,\sigma]$ is a local set that is contained in $\overline O$ ensures that $\hat \Gamma^{\overline O}$ is a GFF in the complement of $O$, that is independent of $\eta [0,\sigma]$. This implies 
that $\eta [0,\sigma]$ is actually a deterministic function $F$ of $\hat \Gamma_{\overline O}$. 

Let us now consider a given harmonic function $h^*$ in $\U$. We now denote the GFF with boundary conditions given by $h^*$ by $\Gamma^*$. We make the assumption that $O$ is an open subset of $\U$ that contains a neighbourhood of $-1$ in $\U$ and such that the law of 
$(\Gamma^*)_{\overline O}$ is absolutely continuous with respect to the law of $\hat \Gamma_{\overline O}$. 
This makes it possible to define the random curve $\eta^* := F (( \Gamma^*)_{\overline O})$. 
Combining the previous statement with the absolute continuity result then readily shows that $ \eta^* = (\eta_t)_{t \le \sigma^*}$ is the unique random curve 
(up to time-reparametrisation) joining $-1$ to (some other point of) $\partial O$, such that for each $t  \ge 0$, the conditional law (given $\gamma^* [ 0, t]$ and the fact that $t < \sigma$) of $\Gamma^*$ in the complement of $\eta^* [0,t]$ is a GFF with boundary conditions given by those of $h^*$ on $\partial \U$ and by $\pm \lambda$ on the two sides of $\eta^*$. In other words, the boundary conditions are given by the harmonic function $h_t^*$, such that $h_t^* - h^*$ goes to $0$ on $\partial \U$, remains bounded in the neighbourhood of $\eta^* [0,t]$, and tends to $\pm \lambda$ on the two sides of the curve. 

By definition, the law of $\eta^*$ up to its first exit time of $O$ is then absolutely continuous with respect to that of $\eta$, from which we get that $ \eta^* (\sigma^*) \not\in \partial \U \setminus \{ 1 \}$ almost surely.

\begin{figure}[h]
\centering
	\includegraphics[scale=1]{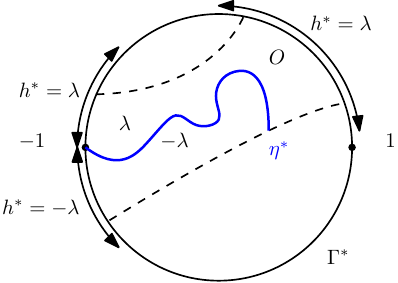}
	\caption{Suppose that $h^*$ is a harmonic function in $\U$, having boundary values equal to $\pm \lambda$ on the indicated regions of $\partial \U$. Denote the open set bounded by the two dashed lines in the figure and $\partial \U$ by $O$. By Corollary \ref{CMcorollary},  a GFF $\Gamma^*$ in $\U$ with boundary values given by $h^*$ is absolutely continuous, when restricted to $O$, with respect to $\hat{\Gamma}_{\overline{O}}$. The above discussion then gives the existence and uniqueness of the curve $\eta^*$, which run up to any time less than $\sigma^*$ is a local set of $\Gamma^*$ with boundary values as marked on the figure (equal to $h^*$ on $\partial \U$). In particular, the curve $\eta^*$ does not exit $O$ through $\partial \U \setminus \{ 1 \}$.}
\end{figure}

\subsection {A global characterisation} 

We now provide a slightly stronger statement than Proposition \ref {p11} about the uniqueness of the coupling of SLE$_4$ with the GFF. One motivation for explaining this here is that this proof strategy can be adapted to explain properties of an important coupling between the GFF and CLE$_4$, that we will mention in the next section. 

\begin {proposition}[Coupling characterisation via the whole curve]\label{prop:sle_determ} 
Consider a GFF $\hat \Gamma$ in $\U$ with boundary conditions $\lambda$ and $-\lambda$ on $\partial_+$ and $\partial_-$ respectively. Then there exists a unique random simple curve $\eta$ joining $-1$ to $1$ in $\U$ with the property that conditionally on the entire curve $\eta$, $\hat{\Gamma}$ restricted to the two connected components of $\U \setminus \eta$ are two independent GFFs with respective boundary conditions $\lambda$ and $-\lambda$.
\end {proposition}

In other words, it is sufficient to know that the entire curve is a local set with the prescribed harmonic function (i.e., one does not need to know that each $\eta [0,t]$ is a local set) to conclude that it is the SLE$_4$ ``level line''.

\medbreak

\begin {proof}[Outline of the proof]
We already know the existence of such a curve $\eta$ -- the SLE$_4$ that was described in the previous sections. It therefore remains to show the uniqueness statement. So, suppose that $\tilde \eta$ is another such curve.

We are going to look at the continuously growing family of local sets $K_{t}^\#:= \tilde \eta \cup \eta [0, t]$ (since $\eta$ is a deterministic function of $\hat{\Gamma}$, we  know that  
the two local sets $\tilde \eta$ and $\eta [0,t]$ are conditionally independent given $\hat{\Gamma}$ so that $K_t^\#$ is a local set of $\hat{\Gamma}$, but we could actually have otherwise chosen $\tilde \eta$ to be conditionally independent of $\eta$ given $\hat{\Gamma}$, and not relied on this fact). We let $\tilde h$, $h_t$ and $h_t^\#$ denote 
the corresponding harmonic functions (in the complements of $\tilde \eta$, in the complement of $\eta [0,t]$ and in the complement of $K_t^\#$). 

If the probability that $\eta \not\subset \tilde \eta$ is not zero, then for some deterministic rational time $t_0$,
the probability that $\eta (t_0) \notin \tilde \eta$ {and $\eta(t_0)$ lies on the boundary of a connected component of $\U\setminus K_{t_0}^\#$ that also has $1$ as a boundary point}, is positive. By symmetry, it will suffice to consider the case where $\eta (t_0)$ lies on the ``top side'' of $\tilde \eta$ (and to show that this leads to a contradiction). 

We will now work with a fixed $t_0$, condition on $\tilde \eta$ and $\eta[0,t_0]$, and define $U$ to be the connected component of $\U \setminus K_{t_0}^\#$ containing $1$. Since $\eta$ is a simple curve from $-1$ to $1$ in $\U$, there are only two possibilities: either it will exit $U$ at a point of $\tilde \eta \setminus \eta [0, t_0 ]$ after time $t_0$, or it will go all the way to $1$ without hitting 
$\tilde \eta$ again. 

\begin{figure}[h]
\centering
	\includegraphics[scale=1]{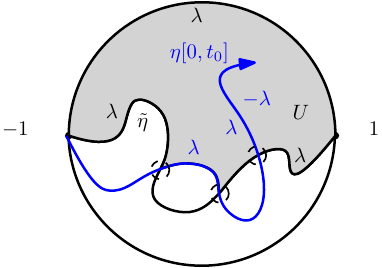}
	\caption{The component $U$ is shaded. The harmonic function $h_{t_0}^\#$ is equal to $\pm \lambda$ on $\partial U$ as indicated on the figure (mind however that we do not have 
	proved anything about its behaviour near the intersection points of $\eta$ and $\tilde{\eta}$).}
\end{figure}

Let us consider the harmonic function $h_{t_0}^\#$ restricted to $U$. 
We do not know exactly what the boundary values of $h_{t_0}^\#$ are near the points where $\eta$ and $\tilde \eta$ 
intersect, 
but as a consequence of Proposition \ref {harmfunctionunion}, we know that 
on the part of $\tilde \eta$ that is at positive distance from $\eta$, it behaves like $\tilde h$ i.e., the boundary condition is $+ \lambda$. Hence,  the boundary condition is in fact $+ \lambda$ in the neighbourhood of all points of $\partial U$ that lie at positive distance of $\eta [0, t_0]$, including the target point $1$.

Then, using the final statement in the previous section (mapping $U$ back onto $\U$), {together with Corollary \ref{CMcorollary}}, we conclude that $\eta$ can exit $U$ only through a point of $\eta [0,t_0]$. This contradicts the fact that $\eta$ is a simple curve. 
 \end {proof}

\section {CLE$_4$ and the GFF} \label{cle4gff}

In the previous sections, we have explained how it is possible to find ``level-lines'' of a GFF with well-chosen boundary conditions. 
These are random curves along which the GFF will be equal to $0$ on one side, and to $2 \lambda$ 
on the other side (or $\pm \lambda$ on the two sides if one subtracts $\lambda$ from the GFF). One may wonder if it is also possible to find such lines in the case of a GFF with zero boundary conditions. Due to the fact that the boundary conditions of the GFF are $0$, one can guess that all these curves should actually be closed loops.

In this case, for reasons that will become clearer later, it is actually a good idea to simultaneously  look for level-lines where the GFF is $0$ on one side and $- 2 \lambda$ on the other side $(\lambda=\sqrt{\pi/8})$.

A variation of the coupling construction of SLE$_4$ with the GFF that we described in the previous chapter actually makes it possible to derive the following result. 
\begin {theorem}[CLE$_4$]
\label {CLE4}
It is possible to construct a thin local set of a GFF $\Gamma$ with Dirichlet boundary conditions, such that $h_A$ takes its values in $\{ - 2\lambda, 2 \lambda \}$. 
\end {theorem}

The set $A$ (or the collection loops that form the outer boundaries of the $O_j$) is called a CLE$_4$ (conformal loop ensemble with parameter $4$). 
The way in which this random fractal set is constructed uses a variant of SLE$_4$. More precisely, one chooses a boundary point $x_0$ (say, $-1$ on $\partial \U$), and for each $x \in \U$, one defines the law of a 
random path $\eta_x$ from $x_0$ to $x$, that is described by some variant of SLE$_4$. For each $x$, one can define a particular stopping time $\tau_x$ on $\eta_x$, at which $\eta_x$ creates a particular loop $l_x$ around $x$. 
The law of $\eta_x$ has the following two important properties:

\begin{itemize}
	\item [-]it is conformally invariant, in the sense that if $\phi$ is a Moebius transformation of the unit disk with $\phi (x) = x'$ and $\phi (x_0)= x_0$, then the law of $\phi(\eta_x)$ is the law of $\eta_{x'}$, and 
	the image of the stopping time $\tau_x$ is $\tau_{x'}$;
	\item [-] It is also target-invariant, in the sense that it is possible for $x\ne x'$ to couple $\eta_x$ and $\eta_{x'}$ in such a way that $\eta_x$ and $\eta_{x'}$ coincide until the first time at which $\eta$
	disconnects $x$ from $x'$ in $\U$. Furthermore, in this coupling if $x'$ is surrounded by $l_x$, then $l_{x'} = l_x$.
\end{itemize}

This makes it possible to choose a countable dense collection of points $(x_n)_{n\ge 1}$ in $\U$, and to couple all the $\eta_{x_n}$, so that the union of these curves create a so-called CLE$_4$ exploration tree. The union of all 
the interiors of the $l_{x_n}$ is then the complement of $A$, and the domain encircled by a loop $l_{x_n}$ will correspond to a connected component of the complement of $A$.
The corresponding function $h_A$  will end up being equal to  $2 \lambda$ or $-2 \lambda$ in that domain, depending on whether
the loop $l_{x_n}$ has been traced clockwise or anti-clockwise by the exploration tree. So, we see that the outer boundaries of the connected components 
of $\U \setminus A$ consist of SLE$_4$-type loops (in particular, their Hausdorff dimension is equal to that of SLE$_4$, namely $3/2$). 

To check that this local set is thin, one can directly see from the construction of $\eta_x$ that the probability that it creates a loop $l_x$ before reaching distance $\eps$ from $x$ will decay at least like a power of $\eps$ as $\eps \to 0$. 
\begin {remark}
\label {confradremark2}
Remark \ref {confradremark} actually implies that for each given $z$, the decrease in log conformal radius from $\U$ to $\U\setminus A$ with respect to $z$ is distributed like the hitting time of $\sqrt{2\pi} \times 2\lambda=\pi$ by a one-dimensional Brownian motion (which is one way to justify the power-law decay).
\end {remark}

\begin {remark} 
It is possible to derive further properties of this random set $A$:

- The conformal invariance of $\eta_x$ can be used to see that the law of $A$ is invariant under any 
given Moebius transformation of the unit disc.

- The Hausdorff dimension of $A$ turns out to be almost surely equal to $15/8$. Note that this value is strictly larger than the Hausdorff dimension $3/2$ of the loops.
\end {remark}

\medbreak 
The next important feature of CLE$_4$ is that: 
\begin {proposition}[CLE$_4$ and its labels are determined by the GFF] 
This set $A$ is the unique thin local set of a GFF with Dirichlet boundary conditions in the unit disc such that $h_A \in \{-2 \lambda, 2 \lambda \}$. 
This implies in particular that $A$ and $h_A$ are deterministic functions of the GFF. 
\end {proposition}

To show this proposition, one uses a similar strategy as
the one that we outlined for Proposition \ref {l4}. 
Suppose that $\tilde A$ and $h_{\tilde A}$ do satisfy the conditions of the proposition for a GFF $\Lambda$. Consider (possibly extending the probability space) a CLE$_4$ that is coupled with $\Lambda$ as in Theorem \ref {CLE4}, such that $A$ and $\tilde A$ are conditionally independent given $\Lambda$. Recall that $A$ is constructed via a branching tree of SLE$_4$-type processes. 

For each given $z$, we define $O(z)$ and $\tilde O(z)$ the connected components of $\U \setminus A$ and $\U \setminus \tilde A$ that contain $z$. The key lemma is the following: 

\begin {lemma} 
For each $z$, with probability one, $O(z) \subset \tilde O(z)$ or $\tilde O(z) \subset O(z)$. 
\end {lemma}

Let us first explain how this lemma can be proven: If $O(z) \not\subset \tilde O(z)$, then the branch $\eta_z$ enters $\tilde O(z)$ before creating the CLE$_4$ loop around $z$. But then a similar argument to that used in the proof of Proposition \ref {l4}
shows that the path $\eta_z$ can not exit the domain $\tilde O(z)$ before completing the loop around $z$. This then implies that $O(z) \subset \tilde O(z)$. 

To conclude this section, let us explain how to deduce the proposition from the lemma. The lemma implies that $\hat A := A \cup \tilde A$ is a thin local set which also has the property that the harmonic function $h_{\hat A}$ takes its values in $\{-2 \lambda, 2 \lambda\}$. 
Remark \ref {confradremark} and Remark \ref {confradremark2} then in turn imply that (for a given dense set of points $z$ in $\U$), the law of the decrease of the log-conformal radius of $\hat O(z)$ at $z$ (compared to that of $D$) is described in terms of the hitting time of $\{-2 \pi, 2 \pi \}$ by a one-dimensional Brownian motion, in exactly the same way as that of $O(z)$. In other words, the log-conformal radius of $O(z)$ and of $\hat O(z)$ have the same law. But if $\hat O(z) \subset O(z)$ almost surely, it means that $O(z) = \hat O(z)$. The same argument (applied to the log-conformal radius of $\tilde O(z)$) shows that $\tilde O(z) = \hat O(z)$ almost surely. 
This in turn finally implies that $\hat A =\tilde A = A$ almost surely.     

\medbreak 

Finally, this local set has the following striking feature:
\begin {proposition}
Conditionally on $A$, the values $2 \eps_j \lambda$ of $h_A$ in the connected components of $\U \setminus A$ are obtained by independent identically distributed with 
$P [ \eps_j = 1 ] = P [ \eps_j = -1 ] = 1/2$.
\end {proposition}

One way to explain this fact will come from another construction and interpretation of CLE$_4$ via Brownian loop-soups (see Section \ref {S.last}).

\section {Constructing a GFF as a nested CLE$_4$}

Given that the two-dimensional GFF is conformally invariant, Theorem \ref{CLE4} has a counterpart in any bounded simply connected domain $U$. One chooses a conformal map $\psi$ from $\U$ onto $U$, and defines a CLE$_4$ in $U$ to be the 
conformal image of a CLE$_4$ in $\U$ under $\psi$. Given such a CLE$_4$ $A$ in $U$, one can define just as in $\U$, independent GFF's $\Gamma^j_U$ in each connected component of the complement of $A$. By putting these together with a collection 
of independent fair coin tosses in each component, one constructs a GFF in $U$. 

Suppose that one now starts with a CLE$_4$ in the unit disc $\U$, and has coin tosses $\eps_j$ and the GFFs $\Gamma^j$ as in Theorem \ref {CLE4}. In order to construct each of the $\Gamma^j$, one can actually use the same procedure: 
first sample an independent CLE$_4$ $A_j$ in each of the $O_j$, plus more coin tosses, and then sample independent GFFs in the connected components of each $O_j\setminus A_j$. 
In this way we get a decomposition of $\Gamma$ as $ \Gamma = h_2 + \Gamma_2$, where:
\begin{itemize} \item $h_2$ is a function that is constant in each of the connected components of $\U\setminus A^2$; $A^2 := \cup_j A_j$, and takes values in $\{ -4 \lambda, 0, 4 \lambda \}$; 
	\item $\Gamma_2$ consists of independent GFFs in each of the connected components of $\U \setminus A^2$. 
\end{itemize}

Iterating this procedure, we get that for each $k \ge 2$, one can define a local set $A^k$ of $\Gamma$ such that $A^{k-1} \subset A^k$, and such that
$\Gamma = h_k + \Gamma_k$, where $h_k$ is constant in each connected component of $\U \setminus A^k$ with values in $\{ -2k\lambda,  (-2 k +4 ) \lambda , \ldots, 2 k \lambda \}$, 
and where $\Gamma_k$ consists of independent GFFs in each of these connected components. Theorem \ref {CLE4} actually shows (by induction) that all the sets $A^k$ are in fact deterministic functions
of the GFF $\Gamma$ in this construction. 

For a given $x \in \U$, it is not difficult to see that the diameter of the connected component of $\U \setminus A^k$ that contains $x$ almost surely tends to $0$ as $k \to \infty$. By dominated convergence, one 
can then deduce that for any bounded measurable function $f$, 
$ E[ \Gamma_k (f)^2 ] \to 0 $ as $k \to \infty$. Since $\Gamma(f) = \int_{\U} f(x) h_k (x) dx + \Gamma_k (f)$, one concludes that: 
\begin {proposition} 
For each bounded measurable function $f$, $\Gamma (f)$ is the limit in probability of $\int_{\U \setminus A_k} h_k (x) f(x) dx$ as $k \to \infty$.
\end {proposition}
Hence, we can recover the GFF $\Gamma$ from the knowledge of its ``topographic'' map (provided by the nested CLE$_4$), together with the coin tosses associated to each of the loops.

\begin {remark}
In this set-up, the value of the GFF at a point $z$ heuristically appears as the limit of the 
simple random walk (with step-size $2 \lambda$) $h_k (z)$, which is of course not a well-defined random variable in the $k \to \infty$ limit (this is similar to the fact that $\Gamma(z)$ can be viewed as the limit of circle averages as their radii tend to $0$). However, we can see that for $z \not= z'$, the correlation between 
$h_k(z)$ and $h_k(z')$ will come from the loops that surround both $z$ in $z'$. In particular, one sees that that expected number of nested CLE$_4$ loops that surround both $z$ and $z'$ will be equal to $G_D (z,z') / (2 \lambda)^2$ -- which sheds some simple interpretation on the covariance structure of the Gaussian generalised function $\Gamma$.  
\end {remark}

\section {Brownian loop-soup and CLE$_4$ -- the three couplings are the same}
\label {S.last}

One may wonder if the previous CLE$_4$ is not related in some way to the Brownian loop-soup that can be coupled to the GFF as explained in Section \ref {S.BLSGFF}. This indeed turns out to be the case.
When one considers a Brownian loop-soup (with $c=1$) in the unit disk, one can define (as in the discrete GFF) the loop-soup clusters. More precisely, one says that two Brownian loops $\gamma$ and $\gamma'$ in the loop-soup are 
in the same cluster if there exists a finite chain of loops $\gamma_0, \ldots, \gamma_k$ in the loop-soup such that $\gamma_0 = \gamma$, $\gamma_k = \gamma'$, and $\gamma_i \cap \gamma_{i+1} \not= \emptyset$ for $i =0, \ldots, n-1$. 

One says that a cluster is an outermost cluster if it is surrounded (i.e. disconnected from $\partial U$) by no other cluster, and we denote by $(L_j)_{j \in J}$ the collection of all outer boundaries of the outermost clusters. In fact, these happen to themselves be simple disjoint loops.
We have the following: 
\begin {theorem}[CLE$_4$ via Brownian loop-soup clusters]  
\label {CLE4BLS}
The collection of loops $(L_j)_{j \in J}$ is distributed exactly like the loops of a CLE$_4$. 
\end {theorem}

Furthermore, the three couplings (GFF-CLE, CLE-loop soup, loop soup-GFF) can be made to coincide. Intuitively speaking, this corresponds to the fact that for every $j$, $\eps_j$ is the sign of the GFF on the loop-soup cluster that $\partial O_j$ is the 
outer boundary of. 
More precisely, we have the following.
\begin {theorem}[The three couplings commute]
\label {thmcommute}
It is possible to couple a CLE$_4$, a GFF $\Gamma$ and a Brownian loop-soup in such a way that: 

- The CLE$_4$ and the GFF are coupled as in Theorem \ref{CLE4}; 

- The CLE$_4$ and the Brownian loop-soup are coupled as in Theorem \ref {CLE4BLS};

- The square of the GFF $\Gamma$ and the Brownian loop-soup are coupled as in Theorem \ref {BLSGFF}. 

\end {theorem}

One way to heuristically understand all of this goes as follows. If one starts with the Brownian loop-soup, one can define its loop-soup clusters $(K_i)_{i \in I}$. The (renormalised) square of the Gaussian free field can then be constructed via the (renormalised) occupation time density of the loop-soup. The GFF itself can then still be obtained by tossing some independent fair coins $(\eps_i)_{i \in I}$ to 
decide the ``sign'' of the GFF on each cluster $K_i$ (of course, all this is just heuristic and the actual statements are more subtle). Then, one can look at the outermost clusters $(K_i)_{i \in I_o}$ 
and at their outer boundaries $(\gamma_i)_{i \in I_o}$. These boundaries form a CLE$_4$. 

Now, when one zooms in on a portion of some $\gamma_i$ for $i \in I_o$, because it is the outer boundary of some loop-soup cluster, we see an asymmetry: no Brownian loop ``touches'' $\gamma_i$ from the outside (which corresponds to the fact that the GFF has boundary values $0$ when seen from the outside), but a number of Brownian loops do touch $\gamma_j$ from the inside, and they therefore 
contribute to the values of the square of the GFF 
on the ``interior side'' of $\gamma_i$. What this theorem says, is that the effect of these ``inside-touching loops'' is 
to create a deterministic shift of exactly $\pm 2 \lambda$ for the boundary conditions of the GFF (where the sign is given by $\eps_i$).

\section {Some related couplings} \label{sec:relatedcouplings}

We have just described natural couplings between the GFF, CLE$_\kappa$ for $\kappa=4$ and Brownian loop-soups with intensity $c=1$. 
There are extensions of these couplings to other values of $\kappa$ and $c$. 

\begin {itemize} 
 \item {[Imaginary Geometry: Coupling the GFF with other SLE$_\kappa$]}
 The coupling of SLE$_4$ with the GFF, where the former is viewed as a ``level-line'' of the latter, can be extended into a coupling of 
 any SLE$_\kappa$ with the GFF. In these couplings, the SLE$_\kappa$ will also be a local set of the GFF.  However,  the boundary conditions of the GFF on the two sides of the SLE$_\kappa$ will not just be constant and equal to $\pm \lambda$, but will involve an additional (unbounded) ``twist'' term. 
 
\item {[CLE$_\kappa$ as clusters of Brownian loops]} 
The construction of CLE$_4$ as clusters of Brownian loops in a loop-soup with intensity $c=1$ can be generalised as follows. If one 
considers a Brownian loop-soup of intensity $c < 1$, then the collection of outermost cluster boundaries defined in the very same way
turn out to be a Conformal Loop Ensemble CLE$_\kappa$ where $c \in (0,1]$ and $\kappa \in (8/3, 4]$ are related by the formula 
\begin {equation}
 \label {ckappa} c = \frac {(6- \kappa) (3 \kappa - 8)}{2 \kappa}.
 \end {equation} 
The boundaries of the loop-soup clusters (in a loop-soup of intensity $c$) are therefore SLE$_\kappa$-type loops for this value $\kappa(c)$. Note that the derivation of this result in \cite {SheffieldWerner} is not simpler for $c=1$ than for the other values of $c$. 

\item {[SLE duality]} 
There is always a natural coupling between SLE$_\kappa$ type curves for $\kappa < 4$ (these are almost surely simple curves) and SLE$_{\kappa'}$ type curves for $\kappa' = 16 / \kappa > 4$ (these are non-simple curves) because the former can be viewed as ``outer boundaries'' of the latter. One instance of this will arise in our discussion of the UST scaling limit. 
\end {itemize}

\section*{Bibliographical comments}
SLE was introduced by Oded Schramm in his seminal paper \cite {Schramm}. The continuity of the trace of SLE (and in particular of SLE$_4$) is due to Rohde and Schramm \cite {RohdeSchramm} (see also Lawler's book \cite {Lawlerbook} for a general introduction and the main properties of SLE).  

The relation between SLE$_4$ and the GFF originated in a series of papers by Schramm and Sheffield (see \cite {SchrammSheffield} for the particular coupling presented in Theorem \ref {SLE4thm} and the references therein -- see also Dub\'edat \cite {Dubedat}), and then further developed (to couplings with other SLEs) by Miller and Sheffield in the ``imaginary geometry'' series of papers (see \cite {IMG1} and the references therein). These are the papers for which the notion of local sets of the continuum GFF is instrumental. For the result mentioned in Remark \ref {LevelLines}, see \cite {SchrammSheffield2}.

For the construction of conformal loop ensembles and their properties, see \cite {Sheffield, SheffieldWerner} -- see also \cite {ASW} for a brief survey of the branching tree construction of CLE$_4$. The coupling between CLE$_4$ and the GFF is due to Miller and Sheffield (unpublished manuscript) -- a derivation with the fact that the CLE$_4$ is a deterministic function of the GFF can be found in \cite {ASW}. The fact that GFF is then a deterministic function of the (labelled) nested 
CLE$_4$ follows rather directly from this coupling. 

The relation between Brownian loop-soup clusters and CLE$_4$
has been derived in \cite {SheffieldWerner} (interestingly, this paper pre-dated and motivated many of the developments of discrete and cable system loop soups presented in Chapter \ref {Ch2}). This result was then used by Lupu \cite {Lupu2} to show that discrete loop-soup clusters converge (in the scaling limit) to CLE$_4$, which in turn, was one instrumental input in the proof of Theorem \ref {thmcommute} in \cite {QianWerner}. 

The general SLE duality results are due to Zhan and Dub\'edat \cite {ZhanDuality, DubedatDuality}. General SLE reversibility results (for other SLE$_\kappa$ than SLE$_4$) are somewhat more difficult to derive -- see \cite {ZhanReversibility,IMG3}. For the Hausdorff dimension of SLE curves and the dimension $3/2$ of SLE$_4$, see \cite {RohdeSchramm,Beffara} and for that of CLE carpets and the dimension $15/8$ of CLE$_4$, see \cite {SSW,NW}.

\chapter {Quick review of further related results} 
\label {Ch6}
In this chapter, we review without proofs (but trying to provide some general ideas) some results that are 
related to items described so far in these notes. We will focus here mainly on the two-dimensional case.

\section {Liouville quantum gravity} 
\label{S.LQG}

We are now going to describe one important instance where the two-dimensional continuum GFF is instrumental. 
Suppose that one wants to find a natural (and hopefully also physically relevant) way
to define a random area measure in a domain $D$. For instance, something that could be interpreted as a canonical 
perturbation of the Euclidean area measure.

In earlier chapters, we explained how the continuum GFF  (loosely speaking) describes some kind of canonical random fluctuation away from a constant function on the disk.  It therefore makes sense to take the GFF as a basic building block for our measures. Moreover, 
in view of the Markov property, it is actually very natural to try and define random measures with densities given by constant multiples of  ``exponentials'' of the GFF. Of course, it is not clear a priori what this should mean, because the continuum GFF is not a proper function. So, one must proceed with caution. (As we shall see, due to the roughness of the GFF, these constant multiples will in some sense have to be $0$ -- it is the usual feature of randomness where infinities cancel out and their differences/ratios give rise to interesting random objects).

One way to go about this is the following. The aim is to give a rigorous meaning to measures having Radon-Nikodym derivatives, with respect to the Lebesgue measure, equal to constant multiples of $\exp ( \gamma  \Gamma)$, for $\Gamma$  a GFF and $\gamma>0$ a real parameter.  We can think of $\gamma$ as controlling how wild a fluctuation away from Lebesgue measure we would like to construct. Let us assume that $\Gamma$ is a GFF on the unit disk, i.e. $D=\U$. The strategy is then to:   
\begin{itemize}
    \item choose some natural approximation $\Gamma_n$ of $\Gamma$, where $\Gamma_n$ is a proper function for every $n$;
	\item define a measure $\mu_n$ for every $n$, whose density with respect to the Lebesgue measure is proportional to $\exp(\gamma \Gamma_n(x))$;
	\item show that as $n \to \infty$, the measures $\mu_n$ converge to some non-trivial random measure $\mu$ in the unit disk; 
	\item finally check that this limit $\mu$ is independent of the specific choice of approximations $\Gamma_n$ being used (among a wide class of natural possible choices). 
\end{itemize}

\begin {remark}[Area measures vs. distances]
When one has an area measure with continuous positive density $f$ with respect to the Lebesgue measure in a domain $D$, then it immediately defines
a metric $d$ in $D$. Namely, the distance $d(x,y)$ can be set equal to the minimum of the integral of $f$ over all smooth paths (parametrised by arc-length) from $x$ to $y$.

This raises the question of whether the measures we discuss here can actually be used, in some analogous way, to define a (physically relevant) random metric in $D$. At first glance this appears to be very difficult, because the measures in question -- like the GFF -- are extremely rough. However, it has been shown in a recent series of papers that such a definition is possible.
\end {remark} 

Let us describe heuristically one way in which this strategy can be implemented.
Recall first that when $(\beta_u)_{u \ge 0}$ is a one-dimensional Brownian motion, then (for any $\gamma \not= 0$), the process 
$\exp ( \gamma \beta_u - \gamma^2 u /2 )$ is a positive martingale started from $1$, that tends almost surely to $0$ as $u \to \infty$ (and is therefore not uniformly integrable).  

As we have seen in Chapter \ref{Ch3}, Section \ref {circav}, one natural way to approximate the two-dimensional GFF by proper functions is via its circle averages. Let us denote 
the average of $\Gamma$ on the circle of radius $r$ around $z$ by $A(z,r)$ (here we do not use the letter $\gamma$, as it is already used for other purposes in this section). 
We also
let $a(z,r) = E [ A(z, r)^2 ]$, so that $A(z, r)$ can be viewed as Brownian motion at time $a(z,r)$ (recall that $a(z,r_0 e^{-u})$ grows linearly in $u$ when $u \ge 0$ and $r_0 < d(z, \partial \U)$). 

This makes it natural to consider for each $r$, the measure $\mu_r$ in $\U$ with density
$$ M(z,r) :=  \exp ( \gamma A(z,r) - \frac {\gamma^2}2 a(z,r) ) $$ 
with respect to the uniform measure $dz / \pi$. We would like to see what happens to  $\mu_r$ as $r \to 0$. 

As we have just explained, for each fixed $z$, the density $M(z,r)$ tends almost surely to $0$ as $r \to 0$, while $E [ M(z,r)]=1$. We can now instead study the limiting behaviour of $\mu_r ( \U )$. By Fubini, 
$$ E [ \mu_r ( \U ) ] = \int_{\U}  E [ M(z,r)] \frac {dz}{\pi} = 1,$$
but it could still be that $\mu_r ( \U) \to 0$. 
Let us now explain, heuristically, what is going on here. 

- For a fixed $z$, the exponential martingale $M(z,e^{-u})$ is not uniformly integrable and the main contribution 
to $E[ M(z,e^{-u}) ]$ for large $u$ will come from an event $E_u$ that has  probability going to $0$ as $u \to \infty$. This means that it is not ``seen'' by any sample path when $u \to \infty$ (explaining, roughly, the fact that $M(z,r) \to 0$ almost surely while $E [ M(z,r)] = 1$). We can also note that the larger $\gamma$ is, the smaller the probability of this exceptional event $E_u$ becomes. 

- However, when we are looking at $\mu_r ( \U)$, we are looking at the mean value of $M(z,r)$ over 
$z$ in $\U$. In this set-up, one therefore has ``more chances'' to capture some exceptional large values for $M(z,r)$, that will (via a simple averaging out effect) turn $\mu_r ( \U)$ into a uniformly integrable family. So, $\mu_r (\U)$ could still converge to a 
random variable with expectation $1$. Whether this is the case or not will depend how large $\gamma$ is. 

\begin {remark}
One useful analogy to have in mind here is with branching Brownian motion, or the branching random walk. Let us briefly review aspects of this classical theory:  
Suppose that one takes a regular $k$-ary tree, and considers a set of i.i.d. standard Gaussian random variables $({\mathcal N}_u)$ indexed 
by the nodes $u=(u_1, \ldots, u_n) \in \{ 1, \ldots, k\}^n$ of the tree. 
Define the random walk indexed by the tree to be
$$S_u= \sum_{j=1}^n {\mathcal N}_{u_1, \ldots, u_j}.$$
If one follows $S_u$ along each given infinite branch of the tree, then $S$ is a simple random walk with Gaussian increments, and the corresponding exponential martingale goes to $0$ almost surely. 
However, if one defines 
$M_n$ to be $k^{-n}$ times the sum of $\exp ( \gamma S_u - \gamma^2 n / 2 )$ over all the $k^n$ $n$-th generation nodes of the tree, then depending on $k$ and $\gamma$, the situation will be different. 

Indeed, it is easy to find a constant $C=C(k)$ so that the maximum of $S_u$ over all $n$-th generation nodes of the tree will almost surely grow slower than $C\times n$ (for example, one can use Markov's inequality to see that the probability of this maximum being greater than $C\times n$ for fixed $n$ is less than $\exp(n(\log k + 1/2 -C))$).  This implies in particular that when $\gamma$ is very large, $M_n$ will still go to $0$ almost surely. On the other hand, when $\gamma$ is very small, an elementary calculation for jointly Gaussian random variables gives that $M_n$ is uniformly bounded in $L^2$. Therefore the limit - which exists because $M_n$ is a martingale - cannot be $0$ almost surely (in fact, a zero-one argument allows one to conclude that it is actually strictly positive almost surely). Slightly more delicate arguments provide an explicit critical value of $\gamma$, depending on $k$, that separates these two regimes.

The connection with GFF circle averages, and the corresponding approximate LQG measures $\mu_r$, comes through the following observation. For $z,w$ distinct, the circle average processes $A(z,r)$ and $A(w,r)$  will look quite similar up to $r$ of the order of $|z-w|$, and after that time they will evolve essentially independently. We also note that there is a nested CLE$_4$ approach to approximating the LQG measures - to be discussed in Remark \ref {CLE4approach} - which is related even more closely to the above discussion.
\end {remark} 

It should therefore be no surprise that the following result holds:  
\begin{proposition}\label{lqgmeas}
Let $\gamma>0$. As $n \to \infty$, almost surely,  the measure $\mu_{2^{-n}}$ converges weakly to a limiting measure $\mu$. This measure $\mu$ is almost surely equal to $0$ when $\gamma\ge 2$, and when $\gamma<2$, it is almost surely a finite measure in $\U$ with positive total mass, such that $E [ \mu ( \U)] = 1$.
\end{proposition}

This measure $\mu$ (when $\gamma < 2$) is often referred to in the literature as a Liouville quantum gravity area measure.
It turns out that the same conclusion as in Proposition \ref{lqgmeas} can be reached if one uses convolutions of the GFF with smooth approximations to the identity in place of circle averages. Furthermore, one has the desired property that the limiting measures obtained in this way agree, regardless of the specific convolution that one chose. They also agree with the measures obtained via other natural approximation schemes (for instance, taking partial sums in the Fourier decomposition of the GFF -- see equation \eqref{randomFourierSeries}). 

{\begin{remark} There is a very close relation between Proposition \ref {lqgmeas} with the notion of thick points of the GFF that was discussed in Remark \ref {thickpoints}. Indeed, it roughly speaking turns out (and this is not very difficult to prove using an argument based on Girsanov Theorem) that  the measure $\mu$ will be supported on the set of $\gamma$-thick points of the GFF $\Gamma$.  
 This explains, at least heuristically, why the measures $\mu$ are non-trivial only when $\gamma <2$, which is 
 the range for which the set of such thick points has strictly positive Hausdorff dimension. The value $\gamma=2$ is a borderline case, that needs to be treated a little differently (and we will not expand on this here) \end{remark} }

We can remark that the conformal invariance property of the GFF gives rise to a similar feature for these LQG measures. Suppose that $f:\U\to \U$ is a conformal transformation, $\Gamma$ is a GFF in $\U$ and denote $Q:=2/\gamma+\gamma/2$. It then follows readily from the construction using circle averages, and considering how circles transform under conformal maps, that the image under $f$ of the $\gamma$-measure associated to $\Gamma$ is the $\gamma$-measure for the field given by the image of $\Gamma$ under $f$ \emph{minus} the function $Q \log |f'(z)|$ on $\U$. Note that although we have so far only defined these LQG measures for Gaussian free fields, there is no problem extending the definition to a field given by a GFF plus a function that is continuous away from the boundary. The result is the following:

\begin{lemma}[Change of coordinates]
\label{lem:coc}
Suppose that $\Gamma$ is a GFF on a domain $D\subset \C$ and $f: \U\to D$ is a conformal map, so that $\Gamma \circ f$ (viewed as a generalised function) is a GFF in $\U$. Let $\mu$ be the LQG measure associated with $\Gamma\circ f+Q\log |f'|$. Then we can unambiguously define the LQG measure associated with $\Gamma$ to be the image of $\mu$ under $f$.
Note that this definition is independent of the choice of function $f$.
\end{lemma}

\begin {remark}[Construction via CLE$_4$] 
\label {CLE4approach} 
There exists another natural way to construct these same LQG area measures $\mu$, arising from the construction of the GFF using CLE$_4$ (see the previous chapter). Indeed this construction gives rise to some natural filtrations, that make it possible to view the approximating measures as martingales.  

Recall that if we consider a nested CLE$_4$ (together with i.i.d. fair coins tosses for each loop) in the unit disc, and define 
the function $h_n$ (with values in $\{ -2 n \lambda, \ldots, 2 n \lambda \}$) in the interior of the $n$-th level CLE loops as at the end of the previous chapter, then
if one adds to $h_n$ a field $\Gamma_n$ that corresponds to a collection of independent zero boundary condition GFFs in the interior of each $n$th level loop, then the obtained field is a zero-boundary GFF. 
This implies that as $n\to \infty$, the function $h_n$ converges (in probability) to a zero-boundary GFF.

It will be useful to define for each $n \ge 0$ and $z \in \U$, the quantity $C_n(z)$ such that that $\exp ( - C_n(z))$ is the conformal radius seen from $z$ 
of the interior of the $n$-th level loop that contains $z$. (The quantity $\exp (- C_0(z))$ is then the conformal radius of $\U$ seen from $z$.) 

The main underlying feature that allows one to construct the CLE$_4$/GFF coupling, is that when one considers the exploration tree 
(i.e., all the curves $\eta_x$) stopped when it constructs the CLE$_4$ loops, then one has a natural filtration for which the values of 
the conditional expectation of the field (i.e., of the harmonic functions) at all the points are continuous martingales (with respect to the same filtration). 
In this set-up, the value $h_{1} (z) - h_0 (z)$ (and iteratively each $h_{n+1} (z) - h_n (z)$) is interpreted as the value of a Brownian motion at the first time at which it exits $(-2 \lambda, 2 \lambda)$, and this time is exactly $C_1 (z) - C_0(z)$ (respectively $C_{n+1}(z) - C_n (z)$). 

This makes it natural to consider the measures $\mu_n$ with density  
$$  \exp(\gamma h_n(z)-\frac{\gamma^2}{2} (C_n(z) - C_0 (z)) $$ 
with respect to the Lebesgue measure on $\U$. 
Now, almost by construction, for every open set $O$, the sequence $\mu_n(O)$ will actually be a positive martingale with respect to the filtration $({\mathcal F}_n)_{n \ge 0}$, where ${\mathcal F}_n$ is the $\sigma$-field generated by all nested CLE$_4$ loops up to level $n$ and the functions $h_1, \ldots, h_n$. It therefore follows immediately that $\mu_n$ almost surely converges to a limiting measure $\mu$ on $\U$. The question is then to decide when $\mu$ is almost surely trivial. In this set-up, it turns out to be quite easy to see directly via Girsanov-type arguments that the measure $\mu$ is non-trivial if and only if $\gamma<2$ (it actually turns out that  when $\gamma <2$, $\mu_n (O)$ will be bounded in $L^p$ for some $p(\gamma) > 1$, but this value $p(\gamma)$ will tend to $1$ as $\gamma \to 2$). 
It is then possible to check that (for $\gamma<2$), the limiting measure $\mu$ does agree almost surely with the measure constructed using the circle averages.
\end {remark}

\section {GFF with Neumann boundary conditions and on compact surfaces} 

There are other important and natural versions of the GFF than the GFF ``with Dirichlet boundary conditions'' that we have focused on so far in these notes. 
The following two variants (in their continuum versions) actually appear to be 
particularly relevant in relation to LQG measures. 

\subsection {GFF on compact surfaces} 

For convenience of exposition, we will describe this in some detail in the case of discrete and continuum tori. Let us start with the discrete case. For some dimension $d \ge 2$ and some given positive integers $w_1, \ldots, w_d$, we consider the discrete torus $T:= \Z^d / ( w_1 \Z, \ldots, w_d \Z)$. In other words, we identify all points of the type $(n_1 + a_1 w_1, \ldots, n_d + a_d w_d)$ with 
$(n_1, \ldots, n_d)$ (when $(a_1, \ldots, a_d) \in \Z^d)$). 
This is now a graph with $\# T = w_1 \times \cdots \times w_d$ sites and $d\times \#T$ edges. Each of these sites has exactly $2d$ neighbours, but as opposed to the finite graphs that we considered in Chapter \ref {chapter:discrete GFF}, there are no boundary sites. 

For each function $\gamma$ defined on (the sites of) $T$, we can still define its Dirichlet energy ${\mathcal E}_T(\gamma)$, which we set to be the sum of $|\nabla \gamma (e) |^2$ over all edges $e$ of the graph. We note that adding the same constant everywhere to any function $\gamma$ does not change its Dirichlet energy. 

This makes it natural to define the GFF on $T$ as a random function defined {\em up to additive constants}. In other words, we consider the quotient space obtained from the set of functions $\R^T$ on $T$, when one identifies any $\gamma$ and $\gamma'$ such that $\gamma'-\gamma$ is a constant function. The space $\overline {\mathcal F}$ of such equivalence classes of functions is therefore a $(\# T - 1)$-dimensional space, and it is possible to unambiguously define the Dirichlet energy ${\mathcal E}_T (\overline \gamma)$ of any $\overline \gamma$ in ${\overline {\mathcal F}}$. There are several natural ways to choose one element in each equivalence class $\overline \gamma$ of functions. These include:
(a) picking the function $\gamma_0 := \gamma_0 ({\overline \gamma})$ that takes the value $0$ at some given point $x_0 \in T$, and (b)picking the function $\gamma_1$ such that $\sum_{x \in T} \gamma_1(T) = 0$. 

The GFF on $T$ is then a random function in $\overline {\mathcal F}$ with density given by a multiple of 
$\exp ( - {\mathcal E}_T ({\overline \gamma}) / (2 \times 2d))$ 
in this $(\# T-1)$-dimensional space. If we choose the representative of $\overline \gamma$ via (a), then we get a proper random function with density given by its Dirichlet energy, in the space of functions defined on $T$ which take the value $0$ at $x_0$. In other words, we are in almost the same set-up as for the definition of the Dirichlet GFF, but now $x_0$ plays the role of the boundary point. 

One can easily generalise to this compact setting most of the results that we described for the Dirichlet GFF. Let us mention a few of the little  tweaks that are needed to make things work smoothly. 

\begin{itemize}
	\setlength \itemsep{0em}
\item[-]  The main issue to be dealt with is that the Green's function on $T$ is infinite, because the random walk on $D$ is {recurrent}. However, if we focus on the representative of the GFF that takes the value $0$ at a given point $x_0$, then this proper function will then be a centred Gaussian vector on $T \setminus \{ x_0 \}$ with covariance given by the Green's function on $T \setminus \{ x_0 \}$ (corresponding to the random walk on $T$ killed when it reaches $x_0$). We note that the {determinant of the Laplacian corresponding to this Green's function} does actually not depend on the choice of $x_0$, by transitivity of the graph $T$. 

\item[-]  One can also  define the random walk loop-soup on $T$ in discrete and continuous times (and the corresponding cable-graph Brownian motion), but given that the random walk on $T$ is recurrent, the mass of the large loops is infinite, so that there will be infinitely many very long loops in this loop-soup (this mirrors exactly the fact that the Green's function on $T$ is infinite). On the other hand, if we kill all the loops that go through a given point $x_0$, then one has a ``usual'' loop-soup in $T \setminus \{ x_0 \}$ with only finitely many large loops. 

\item[-] Since the GFF on $T$ is defined ``up to constants'', its square is obviously not well-defined, but it is still possible to relate the square of the version of the GFF that is equal to $0$ at $x_0$, to the occupation time of a loop-soup in $T \setminus \{x_0\}$. 

\item[-] There is no real difficulty in extending the notion of local sets for this GFF. In the ``strong Markov property'' decomposition (see Definition \ref{smp} of Chapter \ref {Ch4}) one will instead have a zero boundary GFF in the complement of the local set, plus a function (the ``extension of the boundary values'') that is defined up to additive constants too. 
\end{itemize}

One useful and canonical way to think of the GFF on $T$ is to take the ``dual'' perspective. Define ${\mathcal F}_0$ to be the vector space of functions defined on $T$ such that $\sum_x f(x) = 0$. One can then view an element $\overline \gamma$ in ${\overline {\mathcal F}}$ as a linear function that associates to each function $f \in {\mathcal F}_0$ the quantity $\overline{\gamma}(f_0):=\sum_x  \gamma (x) f_0 (x)$ (the fact that $f_0 \in {\mathcal F}_0$ shows that the choice of the representative $\gamma$ of $\overline \gamma$ does not matter). Then, we see that this process $\overline \gamma (f_0)$ indexed by ${\mathcal F}_0$ is just a centred Gaussian process, with variance given by 
$$ E [ \overline \gamma (f_0)^2 ] = E [ {\sum_x (\gamma(x) f_0(x))^2} ] $$

\medbreak 
Let us now turn to the continuum case. We consider the torus $T=\R^d / ( w_1 \Z, \ldots, w_d \Z)$ where $w_1, \ldots, w_d$ are now any positive real numbers. 
Intuitively, the GFF in $T$ will be a random generalised function defined up to an additive constant. Given that we defined the continuum Dirichlet GFF as a process indexed by a set of measures or functions, one natural choice in the present setting is to define the continuum GFF on $T$ to be a centred Gaussian process $(\Gamma(f))_{f \in {\mathcal F}_T}$. Here ${\mathcal F}_T$ can be (for instance) the set of bounded measurable functions such that 
$\int_T f(x) dx = 0$, where $dx$ denotes  the Lebesgue measure on $T$. 

One option to define its covariance function is to consider an {orthonormal} basis of $L^2(T)$ made up of eigenfunctions of (minus) the Laplacian. Note that (as opposed to the case of the Laplacian on a domain with boundary, and its eigenvectors with Dirichlet boundary conditions), the constant function $\varphi_0$ will be an eigenfunction with eigenvalue $0$. In the particular case of the torus, the eigenfunctions and eigenvectors are explicit. One can for instance take eigenfunctions that are multiples  of $\cos(2\pi \sum_1^d (m_i x_i/w_i))$ and $\sin(2\pi \sum_1^d (n_i x_i/w_i))$ to form an orthonormal basis of $L^2(T)$ (the corresponding eigenvalues are then $4\pi (\sum (m_i/w_i)^2)$ and $4\pi (\sum (n_i/w_i)^2)$ respectively). We can then order the eigenfunctions $(\varphi_n)$ in some way (for instance according to increasing eigenvalue $\lambda_n$), and formally define 
$$ \Gamma = \sum_{n \ge 1 } \frac {{\mathcal N}_n}{\sqrt {\lambda_n}} \varphi_n ( \cdot)$$ 
(mind that we omit $n=0$ in the sum here). 
We can note that this corresponds to the choice (b) of representative of the GFF in the discrete case, since the integral of $\Gamma$ (as above) on the torus is equal to $0$.

\subsection {Neumann GFF} 

Suppose now that $D$ is some finite subset of $\Z^d$, and consider a set of edges joining neighbouring points of $D$, in such a way that the obtained graph is connected. For the other edges of $\Z^d$ that have at least one endpoint in $D$, we imagine that we cut it in the middle, creating two half-edges. 
In this way, one creates a connected cable-system, so that each point of $D$ has $2d$ neighbouring edges or half-edges. 
On $D$, one can then define a random walk, that at each step chooses one of the $2d$ possible directions. If the direction is that of a full edge, the walk jump along that edge and lands on its other end, which is a point of $D$. If this direction is a half-edge, and presents a dead-end, then the walk bounces back from that wall, and just decides to stay put. 
So, for instance, if the walk is at a site $x$ that has only $2$ full outgoing edges (and so $2d-2$ half edges), then it will stay put with probability $1 - (1/d)$.

One can now make sense of {\em harmonic functions} for this graph -- those are the functions $h$ on $D$ such that for each $x \in D$, if the walk $X$ is at $x$ at time $0$, then the expected value $E_x [h (X_1)]$ is equal to $h(x)$. The corresponding operator (that associates to a function $f$ the function 
$ \big( x \mapsto E_x [ f(X_1)] - f(x)\big) $ is a discrete analogue of the Laplacian with Neumann boundary conditions. 

The situation is now very similar to the case of random walks in tori; the random walk is recurrent and the Green's function is infinite. 
On the other hand, it is still possible to define the Dirichlet energy of a function $f$ on $D$ as the sum of 
$|\nabla f (e) |^2$ over all full edges $e$ in the graph. 

The GFF can then, just as in the discrete torus, be defined in either of the following ways: (a) as a random function on $D$ such that the sum of its values is $0$; (b) as a random function on $D$ that is equal to $0$ at some given point; (c) as a  linear form acting on the space of functions on $D$ with zero mean. 

In the continuum, let us first consider the case where $d \ge 3$. In that case, in order to make sense of the GFF with Neumann boundary conditions, one needs some regularity conditions on the boundary. One can for instance assume that 
$\partial D$ is a $C^2$ hyper surface, so that one can define the normal vector to $\partial D$ at each point of $\partial D$. 
In that case, it is well-known that there exists an orthonormal basis $(\psi_n)_{n \ge 0}$ of $L^2 (D)$ consisting of eigenfunctions of the Laplacian with the property that on $\partial D$, the normal derivative of $\psi_n$ vanishes (we choose the first eigenfunction $\psi_0$ to be a constant function). Let us denote the corresponding family of eigenvalues by $\mu_n$. The Laplacian is then a bijection from the completion of the space 
of $C^2$ functions on $\overline D$ with vanishing normal derivatives (with respect to the Dirichlet energy) into $L^2 (D)$. 
Then, we can again simply define the Neumann GFF
$$ \Gamma :=  \sum_{n \ge 1} \frac {{\mathcal N}_n}{\sqrt {\mu_n}} \psi_n(\cdot)$$
as a random generalised function with zero mean (or view it as one representative of a linear form acting on the space of functions with zero mean). 

In the two-dimensional case, things are simplified by conformal invariance. Indeed, the previous definition of the Neumann GFF will be conformally invariant (in the same sense as the Dirichlet GFF is conformally invariant), 
so when the domain $D$ is conformally equivalent to some domain $D'$ with a smooth boundary via a conformal map $\Phi$, then one can {\em define} the Neumann GFF in $D$ as the image under $\Phi^{-1}$ of the GFF in $D'$. Using Koebe's uniformisation theorem, this allows for instance to define the GFF in any finitely connected subset of the plane.

\begin {remark}
One can decompose a Neumann GFF as the sum of a Dirichlet GFF in $D$ with a random harmonic function (defined modulo additive constants) that intuitively corresponds to the harmonic extension in $D$ of the values of the Neumann GFF on $\partial D$. This shows in particular that (up to the additive constant issue) the Dirichlet GFF and the Neumann GFF are absolutely continuous with respect to each other when looking at a piece of $D$ that is at positive distance from $\partial D$. In particular, one will be able to also associate LQG-type area measures to Neumann GFFs.  
\end {remark}

\begin {remark} 
\label {NeumannExample}
One concrete example that is useful to have in mind is the case of the Neumann GFF defined in the upper half-plane (even if it is unbounded, we can use conformal invariance and first define it in $\U$). Recall that the Dirichlet Green's function in $\HH$ is 
$$ G (x, y) = \frac {1}{2 \pi} \left( \log \frac {1}{|y-x|} - \log \frac {1}{|y- {\overline x}|} \right) . $$ 
A representative of the Neumann GFF can be obtained as the centred Gaussian generalised function with covariance function obtained by simply changing the minus sign into a plus in the previous expression: 
$$ G^* (x, y) = \frac {1}{2 \pi} \left( \log \frac {1}{|x-y|} + \log \frac {1}{|x- {\overline y}|} \right).$$ 
In particular, we observe that in this case, when $x$ is on the real line, 
$$G^* ( x, y) = \frac 1{\pi} \log \frac {1}{|x-y|}$$
blows up twice faster as $y \to x$ (compared to when $\Im (x)>0$).  
\end {remark}

\section {Quantum zipper and LQG}
\label {S.zipper}
\subsection {Slicing open an LQG surface along an independent SLE}
In this section, we will very briefly describe a further coupling between simple SLE curves and the GFF, that is particularly relevant in the context of Liouville quantum gravity. Recall from Proposition \ref{sle4mgale}
of the previous chapter that if one draws an SLE$_4$ curve from $0$ to $\infty$ in the upper half plane $\HH$, and fixes some point $z\in \HH$, then the process $\arg(f_t(z))$ is a martingale for any $z$ (one can recall that $f_t$ is the conformal map from the slit domain $\HH$ minus the curve up to time $t$, that sends the tip of the curve to $0$ and behaves like the identity at infinity).
This observation leads to a coupling between SLE$_4$ and the GFF with Dirichlet boundary conditions, in which the curve drawn up to any time is a local set of the field. 
As we hinted at previously, if one looks at $\arg(f_t(z))$ plus a well-chosen multiple of $\arg(f_t'(z))$, then this similarly produces a martingale for SLE$_\kappa$ with $\kappa<4$. Again this leads to a coupling of such an SLE$_\kappa$ with a Dirichlet GFF, in which the curve drawn up to any time is locally coupled to the field. The difference is that now there is an extra ``twist'' term appearing in the harmonic function. 

Note that the argument function is the imaginary part of the complex logarithm. It should therefore not be too surprising that the real part of the complex logarithm gives rise to another martingale associated with a simple SLE$_\kappa$ curve. In fact, to observe this martingale one must ``grow'' the SLE$_\kappa$ in a slightly different way (using the so-called reverse Loewner flow under which new pieces of curve are iteratively ``added at the root'' rather than at the tip -- a little bit like when a plant grows from the bottom and ``pushes up'' the already existing part further up). Although we will not go into any more detail on this here, the consequence is another relationship, now between SLE$_\kappa$ and the Neumann GFF. To prove this one can use a very similar argument to that in the proof of Proposition \ref{sle4mgale}, modulo some minor tweaks concerning the use of the reverse Loewner flow and the Neumann GFF.

\begin{proposition}\label{reversecoupling}[Quantum Zipper, I]
	Let $\Gamma$ be a Neumann GFF in $\HH$ (viewed as a generalised function modulo additive constant) and let $\eta$ be an independent SLE$_\kappa$ curve from $0$ to $\infty$, for some $\kappa\in (0,4)$. Fix $t>0$ and write $f_t$ for the conformal map from $\HH\setminus \eta([0,t])\to \HH$ that sends $\eta(t)$ to $0$ and has $f_t(z) \sim z$ as $z\to \infty$. Finally define $\tilde{\Gamma}$ to be equal to $\Gamma$ plus the function $(2/\gamma)\log|z|$, where $\gamma=\sqrt{\kappa}>0$. Then the random generalised function $\tilde \Gamma_t$ modulo additive constant described by 
	$\tilde \Gamma_t := \tilde{\Gamma}\circ f_t^{-1} + Q \log |(f_t^{-1})|$ has the same 
	law as $\tilde{\Gamma}$, when $Q=2/\gamma+\gamma/2$. 
\end{proposition}

\begin {remark}
There exists a version of this result for $\kappa =4$, and actually also for $\kappa > 4$ as well -- but the LQG measures need to be properly defined (in the $\kappa > 4$, one uses the LQG measure with $\gamma^2 = 16 /\kappa$, and some new features have to be understood due to the fact that the SLE is no longer a simple curve.)
\end {remark}

The appearance of the term $Q \log |(f_t^{-1})|$ and the addition of the log singularity to the GFF at the origin may seem slightly odd here, but this should hopefully make a bit more sense in light the next section. We can already note that the way in which 
$\tilde \Gamma$ is defined out of $\Gamma$ is very much similar to the change of variables property of LQG measures (Lemma \ref{lem:coc}).

 Note that in the coupling of Proposition \ref{sle4mgale} between an SLE$_4$ curve and a Dirichlet GFF, the curve is very much \emph{not} independent of the field (indeed, we know that the SLE is determined by the field)! In the above proposition, the  SLE$_\kappa$ curve is by definition independent of the underlying Neumann GFF $\Gamma$ that one starts with. However, as we will briefly explain in the next section, it is very much not independent of the obtained fields $\tilde{\Gamma}\circ f_t^{-1} + Q \log |(f_t^{-1})|$.

\begin {remark} 
One way to rephrase Proposition \ref {reversecoupling} goes as follows: Consider $P$ to be the joint law of $\tilde \Gamma_0$ and an independent SLE$_\kappa$. Define $F_t ( \Gamma , \eta) := (\tilde \Gamma_t, f_t (\eta[t, \infty)))$. Then, the law $P$ is invariant under the flow $(F_t)_{t \ge 0}$. In the ``quantum zipper'' terminology coined in \cite {SheffieldZipper}, this corresponds to ``zipping down'' 
the Neumann GFF along an SLE curve.
\end {remark}

\subsection{LQG boundary length and conformal welding}

We now briefly discuss how to revert the previously described procedure and to understand the reverse of the flow $F_t$. 
A key-role will be played here by the so-called LQG boundary length measure for Neumann-GFFs. 

Consider first the case of a Neumann GFF $\Gamma$ in the upper half-plane. We want to associate to it a measure on the real line 
which could be interpreted as having a density proportional to  $\exp((\gamma/2)\Gamma(x)) $ with respect to the Lebesgue measure
(the reason for the use of $\gamma/2$ rather than $\gamma$ here will become clear in a moment -- it is related to the factor $2$ mentioned in Remark \ref {NeumannExample}), and it is an easy exercise to check that the ideas used to define the area measure can be adapted quite directly. One can, for instance,  approximate the field by semi-circle averages centred at boundary points and then take a limit of the natural ``approximate'' boundary length measures. There is also an analogue of the CLE$_4$ construction of the bulk LQG measures, but we will not get into this here. One noticeable difference is that because the Neumann GFF is intrinsically really only defined up to an additive constant, LQG boundary length measures are in turn only really defined up to a multiplicative constant (although this can of course be fixed in some way). We can also notice that it is also possible to then define the boundary length measure associated to fields such as $\tilde \Gamma_0$ in the previous section obtained by adding some (rather) nice deterministic functions to $\Gamma$.

Again it is the case that the LQG boundary measures can be nicely defined when $\gamma<2$. Moreover, for fixed $\gamma<2$ they will satisfy the same very same conformal covariance property as the area measures. To start with, if $\Gamma$ is a Neumann GFF defined in the upper half-plane and $f$ is a Moebius transformation of the half-plane onto itself, then if $\hat{\Gamma}=\Gamma\circ f + Q\log|f'(z)|$, the $\gamma$-boundary length measure for $\hat \Gamma$ is exactly the image under $f^{-1}$ of the $\gamma$-boundary LQG measure for ${\Gamma}$.

This behaviour under conformal maps allows for the definition of some LQG boundary length measures on much less well behaved domain boundaries. Suppose for instance that $f$ is a conformal map from some simply connected domain $D$ onto $\HH$. Suppose that $\Gamma$ is a Neumann GFF in $\HH$, and {\em define} the field $\hat \Gamma$ using the same formula as above (note that if the boundary of $D$ is not smooth, then it can happen that $\log |f_t'|$ is unbounded). One can the just {\em define} the $\gamma$-LQG boundary measure on $\partial D$ to be the image of the LQG boundary measure of $\Gamma$ under $f^{-1}$.

In the context of Proposition \ref{reversecoupling} above, this raises the following question. Suppose we have a GFF (Neumann or Dirichlet) in the upper half plane and draw an independent SLE$_\kappa$ curve $\eta$ on top of it. Then it is possible to define, via the discussion above and using the conformal map $f_t$, the $\gamma$-LQG boundary length according to $\tilde \Gamma_0$ on the boundary of $\HH \setminus \eta[0,t]$ (where the left-hand side of $\eta([0,t])$ and the right-hand side of $\eta([0,t])$ are treated as two different boundary parts), or equivalently of the field $\tilde \Gamma_t$ on the real line. 

\begin{proposition}
In this set-up, the boundary length measures of the images of the left-hand side and of the right-hand side of $\eta[0,t]$ under $f_t$ for the field $\tilde \Gamma_t$ do coincide. 
\end{proposition}

\begin {remark} 
The field $\tilde \Gamma$ is defined up to an additive constant, so the boundary length is defined up to a multiplicative constant -- and saying that two boundary lengths ``are equal'' is indeed something that does not depend on this multiplicative constant, 
\end {remark}

\begin {remark}
In other words, for each $t \ge 0$, the boundary length measure of the left-hand side of $\eta[0,t]$ and of the right-hand side of $\eta [0,t]$ for the field $\tilde \Gamma_0$ viewed in the domain $D \setminus \eta[0,t]$ (and defined via conformal invariance as described above) do coincide. It is then easy to see that actually, this holds also when viewed in $D \setminus \eta [0,t']$ for $t' > 0$. 
\end {remark}

This proposition is then the key to explain how to describe the ``zipping up'' flow $F_t^{-1}$: 
One observes the field $\tilde \Gamma_t$ in the upper half-plane, and the boundary length measure it defines on $\R$. In the way in which $\tilde \Gamma_t$ was defined, we see that $\eta$ has the property that if $f_t^+(\eta_s)$ and $f_t^- (\eta_s)$ denote the images of the left-side and right-side boundary point $\eta_s$ of $\HH \setminus \eta [0,t]$ for $s < t$, then the boundary length of 
$[f_t^- (\eta_s), 0 ]$ and of $[0, f_t^+ (\eta_s)]$ are the same (for $\tilde \Gamma_t$). 
The key property is that this feature in fact almost surely determines the curve $\eta$. 
In other words, an SLE$_{\kappa}$ curve for $\kappa < 4$ is obtained by ``welding together''  positive and negative real line segments according to their LQG boundary length measure. This is the starting point of a vast and far-reaching theory.

\section {The scaling limit of Wilson's algorithm in 2D} 
\label {SLE2}

Since we have described Wilson's algorithm and its relation to loop-soups and to the square of the GFF in the discrete setting, it is natural to say a few words about their continuum counterparts, even though this story is not so directly related to the GFF.
Most results in this section will be given without proofs. 

\subsection {A first remark in the discrete case} 

Suppose that we are dealing with the UST with wired boundary conditions in some subset of $D$, as described in the first section of Chapter \ref {Ch2}. 
Let us just focus on the first {branch constructed in} Wilson's algorithm. One starts a random walk $Z$ from a given point $x_1$ that is stopped at its exit time $\tau$ of $D$, and then defines its loop-erasure $L(Z)$ and the new domain $D \setminus L(Z)$.

- The exit point $Z_\tau$ of $D$ by $Z$ is chosen according to the so-called discrete {\em harmonic measure} of $\partial D$ seen from $x_1$. 
If one finishes off Wilson's algorithm and defines a UST ${\mathcal T}$, then this point $Z_\tau$ will be the one where the branch of the UST ${\mathcal T}$ starting from $x_1$ will join $\partial D$. Hence, we see that conditionally on $Z_\tau = z$, the conditional law of ${\mathcal T}$ is the uniform measure among all spanning trees (with wired boundary conditions) such that the branch that starts at $x_1$ joins $\partial D$ at $z$ -- we call this law ${\mathcal U}(D, x_1, z)$.

- If we now condition on $Z_{\tau}=z$ and $Z_{\tau -1}=z'$, then it is clear (just because conditioning a uniform measure to be on a smaller set gives the uniform measure on the smaller set), that the conditional law of the remaining part of ${\mathcal T}$ is described by ${\mathcal U}( D \setminus \{ z' \}, x_1, z')$. 

- However,  under this conditional law, the branch from $x_1$ to $z'$ can be obtained via Wilson's algorithm, as the loop-erasure of a random walk $Z'$ from $x_1$ conditioned to exit $D \setminus \{ z' \}$ through $z'$. 

\medbreak

All this shows that the time-reversal $W$ of $L(Z_0, \ldots, Z_\tau)= (L_0=x_1, \ldots, L_{\sigma} = z )$ has a nice ``Markovian-type'' property. More precisely, if we call ${\mathcal W} (D, z, x_1)$ the law of $W=(W_0 = z, W_1 = L_{\sigma -1}, \ldots , W_{\sigma} = x_1)$ when $Z$ is chosen according to the law of the random walk conditioned on $\{ Z_\tau =y \}$, then we have the following.

\begin {lemma}[Domain Markov property of time-reversed LERW] 
\label {DMP}
Suppose that $W=(W_0, \ldots, W_{\sigma})$ is chosen according to ${\mathcal W} (D, z, x_1)$. Then, the conditional law of $(W_1, \ldots, W_{\sigma})$ given $W_1 = z'$ is  given by ${\mathcal W} (D \setminus \{ z' \}, z', x_1 )$. 
\end {lemma}

This makes it natural to consider the conditional laws ${\mathcal U}(D, x_1, z)$ and ${\mathcal W} (D, x_1, z)$ and to progressively grow the LERW $L$ backwards, from its endpoint back to $x_1$.

\subsection {Radial Loewner chains} 

When trying to understand the scaling limit two-dimensional LERW, it appears very natural in view of Lemma \ref {DMP} to use the framework of radial Loewner chains. We now very briefly recall the definition and a few features of these objects.

- Suppose that $\gamma$ is a simple continuous curve that joins $1$ to the origin in the unit disc $\U$ (so apart from its starting point, the entire curve lies in the open unit disc). At each given time $t$, one can define the unique conformal transformation $g_t$ from $U_t := \U \setminus \gamma [0,t]$ onto $\U$, that is normalised at the origin by specifying that
$g_t (0) = 0 $ and $g_t' (0) \in \R_+$. 

It is easy to check that $|g_t'(0)|$ increases continuously from $1$ to infinity, so that it is possible to choose (and this choice is unique) the parametrisation of the path $\gamma : [0, \infty ) \to \U \cup \{ 1 \} $ such that $|g_t' (0)| 
= \exp (t)$ for all $t$.

We then define $\xi_t := g_t (\gamma_t)$ on the unit circle to be the {\em driving function} of the curve $\gamma$. The key point is that two different curves will necessarily have different driving functions. It is indeed possible to show that the functions $g_t$ satisfy the radial Loewner equation 
$$ \partial_t g_t (z) =  g_t (z) \frac { g_t (z) + \xi_t }{g_t (z) - \xi_t}$$ 
for all $z \in \U \setminus \gamma [0,t]$, which makes it possible to recover $\gamma$ from $\xi$. 

In view of Lemma \ref {DMP}, and the above considerations, it is natural to consider the case where the driving function is a Brownian motion on the unit circle. One then has the radial analogue of the first part of Proposition \ref {prop:slecurve}. 

\begin {proposition}[Definition of radial SLE$_\kappa$] 
When $\kappa \in [0,4]$ is fixed and $B$ denotes a one-dimensional Brownian motion, then the function $\xi_t = \exp (i B_{\kappa t})$ almost surely corresponds via the radial Loewner equation to a continuous simple curve $\gamma$ from $1$ to $0$ in the unit disc. This random curve is called radial SLE$_\kappa$.  
\end {proposition}

When $D$ is a simply connected domain and $\Phi$ is a conformal map from $D$ onto $\U$ with $x = \Phi^{-1} (0)$, we say that the image of a radial SLE$_\kappa$ (as defined above) under $\Phi^{-1}$ is an SLE$_\kappa$ from $\Phi^{-1} (1)$ to $x$ 
(here $\Phi^{-1} (1)$ is defined as a ``prime-end'' in case $\Phi^{-1}$ is not one-to-one -- we hope it will be clear in the next subsection what this will mean in the relevant context).  

\subsection {Scaling limit of LERW} 

In view of Lemma \ref {DMP}, and the fact that in two dimensions the LERW is obtained from the random walk which has a conformally invariant scaling limit (Brownian motion), it is natural to expect that any scaling limit of LERW should satisfy some continuum version of the domain Markov property and some conformal invariance features. This leads one directly to guess that the scaling limit in distribution of LERW should be one of the radial SLE curves. This turns out to be correct. 

Suppose that $D$ is a bounded simply connected domain. We assume that $\partial D$ is a continuous curve, in the sense that any conformal mapping $\Phi$ from $\U$ onto $D$ extends continuously into a mapping from $\U$ onto $\partial D$ (mind that this mapping is not necessarily one-to-one from $\partial \U$ onto $\partial D$ -- we for instance authorise the domain $D= \U \setminus [0,1]$ where the point $1$ would have two preimages. 

We fix an interior point $x$ of $D$ and a boundary point $z$ (together with a choice of $\Phi^{-1}(z)\in \partial U$ if necessary). Then, for each $\delta$, we choose (in some way) a lattice approximation $D_\delta$ of $D$ on $\delta \Z^2$, so that $D_\delta \subset D$. For each $\delta$, we choose $x_\delta\in D_\delta$ and $z_\delta$ to be points that are very close to $x$ and $z$: the point $z_\delta$ being a boundary point of $D_\delta$, and the image under $\Phi$ of $z_\delta$ begin close to $\Phi (z)$). 

We can now define for each $\delta$, the law of the loop-erasure of a random walk from {$x_\delta$} to the boundary of $D_\delta$, conditioned to hit this boundary at {$z_\delta$}.
Then, by interpolating it with linear segments of length $\delta$, we obtain a continuous function on the cable-system of $D_\delta$, and finally, can consider its time-reversal. This time-reversal $\gamma_{\delta}$ is a continuous curve from $z_\delta$ to $x_\delta$ for every $\delta$.

\begin {theorem} 
\label {lerwscalinglimit}
When $\delta \to 0$, the law of $\gamma_\delta$, viewed as a random compact subset of $\overline D$, does converge weakly towards the law of a radial SLE$_2$ from $z$ to $x$ in $D$. 
\end {theorem} 

It is possible to upgrade this convergence to stronger topologies. Actually, the proof in \cite {LSWLERW} gives convergence in the sup-norm when both are parametrised by log-conformal radius seen from $x$.

\subsection {UST scaling limits} 

We now list some further results about the scaling limit of 
the entire UST and of Wilson's algorithm: 

\begin {enumerate} 
 \item {[Scaling limit of LERW and erased loops]} 
In the previous framework, we have focused solely on the scaling limit of the LERW. In Wilson's algorithm, the collection of erased loops when performing the LERW correspond exactly to the loops in a random walk loop-soup that the LERW hits. Using the fact that the discrete random walk loop-soup converges to the Brownian loop-soup in the scaling limit (together with some a priori estimates on the time-lengths of the small random walk loops encountered by the LERW), it is actually possible to derive the following result. This can be viewed as the scaling limit of the first step in Wilson's algorithm (for a UST with wired boundary conditions) and can be stated as follows.
Suppose that $D$ is a bounded simply connected domain in the plane with $x \in D$. Choose $z \in \partial D$ according to the harmonic measure viewed from $x$ (i.e. distributed like the exit point from $D$ by a Brownian motion started from $x$). Then define a radial SLE$_2$ $\gamma$ from $z$ to $x$, that we ``time-reverse'' i.e., we look at it as a curve from $x$ to $z$. Consider further an independent Brownian loop-soup (of oriented Brownian loops) with intensity $\alpha = 1$ in $D$.
Then, it is not difficult to check that almost surely for every loop that $\gamma$ hits on its way from $x$ to $z$, the point on that loop that $\gamma$ hits for the first time is a simple point of this loop.
Furthermore, the sum of the time-lengths of the Brownian loops encountered by $\gamma$ is almost surely finite. 
This makes it possible to define the function $B$ obtained by concatenating the Brownian loops encountered by $\gamma$ in the order in which they are met by $\gamma$.
Then:  
\begin {proposition}[SLE$_2$+Brownian loops=BM]
\label {sle2prop}
The obtained path is a continuous path from $x$ to $\partial D$ that is distributed exactly like a Brownian motion started from $x$ until its first hitting time of $\partial D$.  
\end {proposition} 
In this way, one can indeed interpret $\gamma$ as a loop-erasure of this Brownian motion. However, it is worth noticing that this does not answer the following (still) open question: {\em In this coupling, is the SLE$_2$ a deterministic function of the Brownian motion?}.\vspace{0.1cm} 

\item {[Scaling limit of the finite-dimensional subtrees]} 
There are several ways to describe the scaling limit of the entire UST. In view of Wilson's algorithm, one natural way is via the law of its ``finite subtrees''. More precisely, for each finite collection of points $x_1, \ldots, x_n$, we can look at the ``subtree'' of the UST that connects these $n$ points and the boundary. In the discrete case, this would correspond to the tree obtained by performing Wilson's algorithm to successively discover the branches of the tree that connect these $n$ points to the boundary. 

By iteratively using Theorem \ref {lerwscalinglimit}, one can describe the scaling limit of these finite subtrees in terms of $n$ successive radial SLE$_2$ curves.\vspace{0.1cm}  

\item {[Scaling limit of the UST Peano curve]} 
Another natural way to describe the entire UST scaling limit at once is via the UST contour curve. Indeed, one can define the discrete contour curve of the UST, as the space-filling loop that draws the ``inside contour'' of the tree. It turns out that the scaling limit of this curve can be also described in terms of SLE curves: it is a variant of SLE$_8$. The information encapsulated by this space-filling loop turns out to be the same as the one provided by the collection of all finite-trees (say starting from points with rational coordinates) as described in the previous item. \vspace{0.1cm} 

\item {[Scaling limit of UST with other boundary conditions]} 
It is natural to ask what happens (for instance) in the scaling limit for a UST with free rather than wired boundary conditions. Again, there are two ways to go about this. One can either note that in the discrete case, the UST with free dual boundary conditions is the dual of a UST with wired boundary conditions (on the dual graph). In particular, the inside contour curve of the latter is the same as the outside contour of the former, so that one can just use the same SLE$_8$ loop to describe its scaling limit. 
If one would want to describe the law of the finite subtrees of the UST with free boundary conditions, one possibility is to control the 
Radon-Nikodym derivative of the law of the subtrees with respect to the ones with wired boundary conditions. This can be done using considerations on Brownian loop-soups (but it is not totally straightforward). 
\end {enumerate}

\begin {remark} 
In Proposition \ref {sle2prop}, we see that the SLE$_2$ is naturally coupled with an independent Brownian loop-soup with intensity $\alpha =1$ (or equivalently with an independent unoriented Brownian loop-soup 
with intensity $c=2$) in the sense that the union of these two independent objects can be used to define Brownian paths. 
This is a particular instance of the so-called restriction property: when $\kappa \in [0, 8/3)$, then SLE$_\kappa$ is naturally coupled with an independent loop-soup with intensity 
$$ c( \kappa) = \frac {(6 - \kappa) (8 - 3 \kappa)} {2 \kappa},$$
and the union of these two allow to construct ``restriction measures''.
This can be viewed as the $\kappa < 8/3$ counterpart of the construction of CLE$_\kappa$ from Brownian loop-soups for $\kappa \in (8/3, 4]$, as mentioned in (\ref {ckappa}). In particular, one can notice that the formula relating $|c|$ and $\kappa$ are the same.  
\end {remark}

\section*{Bibliographical comments}

The LQG area measures are a particular case of Gaussian multiplicative chaos constructed from ``log-correlated fields'', as pioneered in the work of Hoegh-Krohn and Kahane \cite {HK,K} or \cite {RhodesVargasPS} for a survey), and beginning with ideas of Mandelbrot (which gave rise to the name ``Mandelbrot multiplicative cascades'' as then later studied by Kahane and Peyri\`ere) -- see also \cite {AruReview}. The motivation from physics to investigate such measures is sometimes encapsulated by the term ``quantum gravity'', and the idea of using the exponential of the GFF in this setting is closely associated to the work of Polyakov. 
The relation to SLE and weldings (which is the perspective that we presented in Section \ref {S.zipper}) was initiated by Sheffield \cite {SheffieldZipper} and then considerably developed by Duplantier, Sheffield and Miller \cite {DS,DMS}. For the construction via nested CLE$_4$, see \cite {APS} and the references therein (in particular the unpublished preprint by Aidekon). For the branching Brownian motion, see \cite {ShiStFlour} and the references therein (in particular the seminal papers by Biggins). 
The fact that the LQG area measure corresponds to a metric is the outcome of a recent series of papers (see \cite {Metric2,Metric4} and the references therein). 

There are two closely related and important lines of research that we did not discuss here.  One is to do with discrete approximations of the LQG area measures and metrics via discrete planar maps (starting with the so-called Brownian map). The other is the construction of measures from a conformal field theory approach, based on a fine analysis of LQG correlation functions by David, Kupiainen, Rhodes and Vargas. Both topics could be the focus of entire books. We just point here to \cite {LeGallBG,MillerSheffieldLQGBM3} and \cite {VargasLN} and the references therein.    

The convergence of LERW to radial SLE$_2$ is the main result of \cite {LSWLERW}, that builds on some a priori estimates about LERW by Schramm \cite {Schramm}. Note also that a number of asymptotic results (such as the precise asymptotics of the determinant of the discrete Laplacian) had been derived by Kenyon \cite {K1,K2}. For results related to the scaling limits of USTs with free boundary conditions, one can look at the appendix of \cite {BDW}. Some references for the relation between Brownian motion, the erased loops and SLE$_2$ are 
\cite {LSWrest,LawlerWerner,SapoShi,Ambrosio}. For the final remark on restriction measures, see \cite {LSWrest}, or the survey \cite {Wrestriction}. 

\backmatter

\end{document}